\input amstex
\documentstyle{amsppt}
\nopagenumbers
\nologo
%
%
%
\catcode`@=11
\redefine\output@{%
  \def\break{\penalty-\@M}\let\par\endgraf
  \ifnum\pageno=1\global\voffset=90pt\else\global\voffset=25pt\fi
  \ifodd\pageno\global\hoffset=105pt\else\global\hoffset=8pt\fi  
  \shipout\vbox{%
    \ifplain@
      \let\makeheadline\relax \let\makefootline\relax
    \else
      \iffirstpage@ \global\firstpage@false
        \let\rightheadline\frheadline
        \let\leftheadline\flheadline
      \else
        \ifrunheads@ 
        \else \let\makeheadline\relax
        \fi
      \fi
    \fi
    \makeheadline \pagebody \makefootline}%
  \advancepageno \ifnum\outputpenalty>-\@MM\else\dosupereject\fi
}
\font\cpr=cmr7
\newcount\xnumber
\footline={\xnumber=\pageno
\divide\xnumber by 7
\multiply\xnumber by -7
\advance\xnumber by\pageno
\ifnum\xnumber>0\hfil\else\vtop{\vskip 0.5cm
\noindent\cpr CopyRight \copyright\ Sharipov R.A.,
2004.}\hfil\fi}
\def\setfirstpage{\global\firstpage@true}
\catcode`\@=\active

\fontdimen3\tenrm=3pt
\fontdimen4\tenrm=0.7pt

\def\leaderfill{\leaders\hbox to 0.3em{\hss.\hss}\hfill}
\font\tvbf=cmbx12
\font\tvrm=cmr12
\font\etbf=cmbx8
\font\tencyr=wncyr10


\Monograph
\def\negskp{\hskip -2pt}
\def\compos{\,\raise 1pt\hbox{$\sssize\circ$} \,}

\def\Cl{\operatorname{Cl}}
\def\tr{\operatorname{tr}}
\def\id{\operatorname{id}}

\def\Ker{\operatorname{Ker}}
\def\Map{\operatorname{Map}}
\def\Hom{\operatorname{Hom}}
\def\End{\operatorname{End}}
\def\Aut{\operatorname{Aut}}
\def\rank{\operatorname{rank}}
\def\Img{\operatorname{Im}}
\def\MatGrSO{\operatorname{SO}}
\def\MatGrO{\operatorname{O}}
\def\blue#1{#1}
\def\myref#1#2{\line{\vtop{\hsize 10pt\noindent #1}
\vtop{\advance\hsize -10pt\noindent #2}\hss}\vskip 4pt}
\pagewidth{360pt}
\pageheight{606pt}
\loadbold
\TagsOnRight
\document
\vbox to\vsize{

\centerline{\etbf RUSSIAN FEDERAL COMMITTEE}
\centerline{\etbf FOR HIGHER EDUCATION}
\bigskip
\centerline{\etbf BASHKIR STATE UNIVERSITY}
\vskip 3cm
\centerline{SHARIPOV\ R.\,A.}
\vskip 1.5cm
\centerline{\tvbf COURSE \ OF \ LINEAR \ ALGEBRA}
\vskip 0.2cm
\centerline{\tvbf AND \ MULTIDIMENSIONAL \ GEOMETRY}
\vskip 1.3cm
\centerline{\tvrm The Textbook}
\vfill
\centerline{Ufa 1996}}
\vbox to\vsize{
MSC 97U20\par
\par
UDC 517.9\par
\medskip
Sharipov R. A. {\bf Course of Linear Algebra and Multidimensional
Geometry}: the textbook / Publ\. of Bashkir State University
--- Ufa, 1996. --- pp\.~143. ISBN 5-7477-0099-5.
\bigskip
\bigskip
This book is written as a textbook for the course of multidimensional
geometry and linear algebra. At Mathematical Department of Bashkir
State University this course is taught to the first year students in the
Spring semester. It is a part of the basic mathematical education.
Therefore, this course is taught at Physical and Mathematical
Departments in all Universities of Russia.\par
     In preparing Russian edition of this book I used the computer
typesetting on the base of the \AmSTeX\ package and I used the Cyrillic
fonts of Lh-family distributed by the CyrTUG association of Cyrillic
\TeX\ users. English edition of this book is also typeset by means of
the \AmSTeX\ package.\par
\medskip
Referees:\ \ \ \
\vtop{\hsize 9.5cm\noindent Computational Mathematics and Cybernetics
group of Ufa State University for Aircraft and Technology ({\tencyr UGATU});
\vskip 0.1cm
\noindent Prof\.~S.~I.~Pinchuk, Chelyabinsk State University for
Technology ({\tencyr ChGTU}) and Indiana University.}
\medskip
\noindent {\bf Contacts to author}.
\medskip
\line{\vtop to 150pt{\hsize=300pt\settabs\+\indent Office:\ &\cr
\+ Office:\hss &Mathematics Department, Bashkir State University,\cr
\+\hss &32 Frunze street, 450074 Ufa, Russia\cr
\+ Phone:\hss &7-(3472)-23-67-18\cr
\+ Fax:\hss   &7-(3472)-23-67-74\cr
\medskip
\+ Home:\hss &5 Rabochaya street, 450003 Ufa, Russia\cr
\+ Phone:\hss &7-(917)-75-55-786\cr
\+ E-mails:\hss &{\catcode`_=11
   \catcode`\_=\active}%
   \blue{R\_\hskip 1pt Sharipov\@ic.bashedu.ru}%
   \cr
\+\hss &%
   \blue{r-sharipov\@mail.ru}\cr
\+\hss &{\catcode`_=11
   \catcode`\_=\active}%
   \blue{ra\_\hskip 1pt sharipov\@lycos.com}%
   \cr
\+\hss &{\catcode`_=11
   \catcode`\_=\active}%
   \blue{ra\_\hskip 1pt sharipov\@hotmail.com}%
   \cr
\+ URL:\hss &%
   \blue{http:/\negskp/www.geocities.com/r-sharipov}\cr
\vfil}\hss}
\vskip -100pt
\vfil
\line{\vtop{\hsize=300pt\settabs\+\indent\kern 180pt &\cr
\+ISBN 5-7477-0099-5\hss
  &\copyright\ Sharipov R.A., 1996\cr
\+\hss &\copyright\ Bashkir State University, 1996\cr
\+English translation\hss &\copyright\ Sharipov R.A., 2004\cr}\hss}
\vskip 1pt plus 1pt minus 1pt}
\topmatter
\title
CONTENTS.
\endtitle
\endtopmatter
\document
\vskip 30pt
\line{CONTENTS.\ \leaderfill\ 3.}
\medskip
\line{PREFACE.\ \leaderfill\ 5.}
\medskip
\line{CHAPTER~\uppercase\expandafter{\romannumeral 1}.
LINEAR VECTOR SPACES AND LINEAR MAPPINGS.\ \leaderfill\ 6.}
\medskip
\line{\S~1. The sets and mappings.\ \leaderfill\ 6.}
\line{\S~2. Linear vector spaces.\ \leaderfill\ 10.}
\line{\S~3. Linear dependence and linear independence.\ \leaderfill\ 14.}
\line{\S~4. Spanning systems and bases.\ \leaderfill\ 18.}
\line{\S~5. Coordinates. Transformation of the coordinates of a vector\hss}
\line{\qquad under a change of basis.\ \leaderfill\ 22.}
\line{\S~6. Intersections and sums of subspaces.\ \leaderfill\ 27.}
\line{\S~7. Cosets of a subspace. The concept
of factorspace.\ \leaderfill\ 31.}
\setfirstpage
\line{\S~8. Linear mappings.\ \leaderfill\ 36.}
\line{\S~9. The matrix of a linear mapping.\ \leaderfill\ 39.}
\line{\S~10. Algebraic operations with mappings.\hss}
\line{\qquad \ The space of homomorphisms $\Hom(V,W)$.\ \leaderfill\ 45.}
\medskip
\line{CHAPTER~\uppercase\expandafter{\romannumeral 2}.
LINEAR OPERATORS.\ \leaderfill\ 50.}
\medskip
\line{\S~1. Linear operators. The algebra of endomorphisms $\End(V)$\hss}
\line{\qquad and the group of automorphisms $\Aut(V)$.\ \leaderfill\ 50.}
\line{\S~2. Projection operators.\ \leaderfill\ 56.}
\line{\S~3. Invariant subspaces. Restriction and factorization of
operators.\ \leaderfill\ 61.}
\line{\S~4. Eigenvalues and eigenvectors.\ \leaderfill\ 66.}
\line{\S~5. Nilpotent operators.\ \leaderfill\ 72.}
\line{\S~6. Root subspaces. Two theorems on the sum of root
subspaces.\ \leaderfill\ 79.}
\line{\S~7. Jordan basis of a linear operator.
Hamilton-Cayley theorem.\ \leaderfill\ 83.}
\medskip
\line{CHAPTER~\uppercase\expandafter{\romannumeral 3}.
DUAL SPACE.\ \leaderfill\ 87.}
\medskip
\line{\S~1. Linear functionals. Vectors and covectors. Dual
space.\ \leaderfill\ 87.}
\line{\S~2. Transformation of the coordinates of a covector\hss}
\line{\qquad under a change of basis.\ \leaderfill\ 92.}
\line{\S~3. Orthogonal complements in a dual spaces.\ \leaderfill\ 94.}
\line{\S~4. Conjugate mapping.\ \leaderfill\ 97.}
\medskip
\line{CHAPTER~\uppercase\expandafter{\romannumeral 4}.
BILINEAR AND QUADRATIC FORMS.\ \leaderfill\ 100.}
\medskip
\line{\S~1. Symmetric bilinear forms and quadratic forms.
Recovery formula.\ \leaderfill\ 100.}
\line{\S~2. Orthogonal complements with respect to a quadratic
form.\ \leaderfill\ 103.}
\line{\S~3. Transformation of a quadratic form to its canonic form.\hss}
\line{\qquad Inertia indices and signature.\ \leaderfill\ 107.}
\line{\S~4. Positive quadratic forms. Silvester's
criterion.\ \leaderfill\ 114.}
\medskip\newpage
\line{CHAPTER~\uppercase\expandafter{\romannumeral 5}.
EUCLIDEAN SPACES.\ \leaderfill\ 119.}
\medskip
\line{\S~1. The norm and the scalar product. The angle between vectors.\hss}
\line{\qquad Orthonormal bases.\ \leaderfill\ 119.}
\line{\S~2. Quadratic forms in a Euclidean space. Diagonalization of a pair\hss}
\line{\qquad of quadratic forms.\ \leaderfill\ 123.}
\line{\S~3. Selfadjoint operators. Theorem on the spectrum and the basis\hss}
\line{\qquad of eigenvectors for a selfadjoint operator.\ \leaderfill\ 127.}
\line{\S~4. Isometries and orthogonal operators.\ \leaderfill\ 132.}
\medskip
\line{CHAPTER~\uppercase\expandafter{\romannumeral 6}.
AFFINE SPACES.\ \leaderfill\ 136.}
\medskip
\line{\S~1. Points and parallel translations. Affine
spaces.\ \leaderfill\ 136.}
\line{\S~2. Euclidean point spaces. Quadrics in a Euclidean space.
\ \leaderfill\ 139.}
\medskip
\line{REFERENCES.\ \leaderfill\ 143.}
\medskip
\newpage
\topmatter
\title
PREFACE.
\endtitle
\endtopmatter
\document
    There are two approaches to stating the linear algebra and the
multidimensional geometry. The first approach can be characterized
as the {\tencyr\char '074}coordinates and matrices approach{\tencyr
\char '076}. The second one is the {\tencyr\char '074}invariant 
geometric approach{\tencyr\char '076}.\par
\setfirstpage
     In most of textbooks the coordinates and matrices approach is
used. It starts with considering the systems of linear algebraic
equations. Then the theory of determinants is developed, the matrix
algebra and the geometry of the space $\Bbb R^n$ are considered.
This approach is convenient for initial introduction to the subject
since it is based on very simple concepts: the numbers, the sets of
numbers, the numeric matrices, linear functions, and linear equations.
The proofs within this approach are conceptually simple and mostly
are based on calculations. However, in further statement of the
subject the coordinates and matrices approach is not so advantageous.
Computational proofs become huge, while the intension to consider only
numeric objects prevents us from introducing and using new concepts.\par
    The invariant geometric approach, which is used in this book,
starts with the definition of abstract linear vector space. Thereby
the coordinate representation of vectors is not of crucial importance;
the set-theoretic methods commonly used in modern algebra become more
important. Linear vector space is the very object to which these methods
apply in a most simple and effective way: proofs of many facts can be
shortened and made more elegant.\par
    The invariant geometric approach lets the reader to get prepared to
the study of more advanced branches of mathematics such as differential
geometry, commutative algebra, algebraic geometry, and algebraic topology.
I prefer a self-sufficient way of explanation. The reader is assumed to
have only minimal preliminary knowledge in matrix algebra and in theory
of determinants. This material is usually given in courses of general
algebra and analytic geometry.\par
     Under the term {\tencyr\char '074}numeric field{\tencyr\char '076} 
in this book we assume one of the following three fields: the field of
rational numbers $\Bbb Q$, the field of real numbers $\Bbb R$, or the 
field of complex numbers $\Bbb C$. Therefore the reader should not know 
the general theory of numeric fields.\par
     I am grateful to E.~B.~Rudenko for reading and correcting the
manuscript of Russian edition of this book. 
\bigskip\bigskip
\line{\vbox{\hsize 7.5cm\noindent May, 1996;\newline May,
2004.}\hss R.~A.~Sharipov.}
\newpage
\setfirstpage
\topmatter
\title\chapter{1}
LINEAR VECTOR SPACES AND LINEAR MAPPINGS.
\endtitle
\endtopmatter
\document
\head
\S~1. The sets and mappings.
\endhead
\leftheadtext{CHAPTER~\uppercase\expandafter{\romannumeral 1}.
LINEAR VECTOR SPACES AND LINEAR MAPPINGS.}
    The concept of {\it a set} is a basic concept of modern mathematics.
It denotes any group of objects for some reasons distinguished from other
objects and grouped together. Objects constituting a given set are called
{\it the elements} of this set. We usually assign some literal names
(identificators) to the sets and to their elements. Suppose the set $A$
consists of three objects $m$, $n$, and $q$. Then we write
$$
A=\{m,\ n,\ q\}.
$$
The fact that $m$ is an element of the set $A$ is denoted by the
membership sign: $m\in A$. The writing $p\notin A$ means that the
object $p$ is not an element of the set $A$.\par
    If we have several sets, we can gather all of their elements
into one set which is called the {\it union\/} of initial sets. In
order to denote this gathering operation we use the union sign
$\cup$. If we gather the elements each of which belongs to all of 
our sets, they constitute a new set which is called the {\it
intersection\/} of initial sets. In order to denote this operation 
we use the intersection sign $\cap$.\par
     If a set $A$ is a part of another set $B$, we denote this fact
as $A\subset B$ or $A\subseteq B$ and say that the set $A$ is
{\it a subset} of the set $B$. Two signs $\subset$ and $\subseteq$
are equivalent. However, using the sign $\subseteq$, we emphasize
that the condition $A\subset B$ does not exclude the coincidence
of sets $A=B$. If $A\varsubsetneq B$, then we say that the set $A$
is {\it a strict subset} in the set $B$.\par
     The term {\it empty set} is used to denote the set $\varnothing$
that comprises no elements at all. The empty set is assumed to be a
part of any set: $\varnothing\subset A$.\par
\definition{Definition 1.1} The mapping $f\!:\,X\to Y$ from the set $X$
to the set $Y$ is a rule $f$ applicable to any element $x$ of the
set $X$ and such that, being applied to a particular element $x\in X$,
uniquely defines some element $y=f(x)$ in the set $Y$.
\enddefinition
     The set $X$ in the definition~1.1 is called {\it the domain} of the
mapping $f$. The set $Y$ in the definition~1.1 is called {\it the domain of
values} of the mapping $f$. The writing $f(x)$ means that the rule $f$
is applied to the element $x$ of the set $X$. The element $y=f(x)$
obtained as a result of applying $f$ to $x$ is called {\it the image}
of $x$ under the mapping $f$.\par
    Let $A$ be a subset of the set $X$. The set $f(A)$ composed by the
images of all elements $x\in A$ is called {\it the image} of the subset
$A$ under the mapping $f$:
$$
f(A)=\{y\in Y\!:\,\exists\,x\ ((x\in A)\and (f(x)=y))\}.
$$
If $A=X$, then the image $f(X)$ is called {\it the image of the mapping}
$f$. There is special notation for this image: $f(X)=\Img f$. {\it The
set of values} is another term used for denoting $\Img f=f(X)$; don't
confuse it with {\it the domain of values}.\par
     Let $y$ be an element of the set $Y$. Let's consider the set
$f^{-1}(y)$ consisting of all elements $x\in X$ that are mapped to
the element $y$. This set $f^{-1}(y)$ is called {\it the total preimage}
of the element $y$:
$$
f^{-1}(y)=\{x\in X\!:\,f(x)=y\}.
$$
Suppose that $B$ is a subset in $Y$. Taking the union of total preimages
for all elements of the set $B$, we get {\it the total preimage} of the
set $B$ itself:
$$
f^{-1}(B)=\{x\in X\!:\,f(x)\in B\}.
$$
It is clear that for the case $B=Y$ the total preimage $f^{-1}(Y)$
coincides with $X$. Therefore there is no special sign for denoting
$f^{-1}(Y)$.
\definition{Definition 1.2} The mapping $f\!:\,X\to Y$ is called
{\it injective} if images of any two distinct elements $x_1\neq
x_2$ are different, i\.\,e\. $x_1\neq x_2$ implies $f(x_1)\neq f(x_2)$.
\enddefinition
\definition{Definition 1.3} The mapping $f\!:\,X\to Y$ is called
{\it surjective} if total preimage $f^{-1}(y)$ of any element
$y\in Y$ is not empty.
\enddefinition
\definition{Definition 1.4} The mapping $f\!:\,X\to Y$ is called
a {\it bijective} mapping or a {\it one-to-one} mapping if total
preimage $f^{-1}(y)$ of any element $y\in Y$ is a set consisting
of exactly one element.
\enddefinition
\proclaim{Theorem 1.1} The mapping $f\!:\,X\to Y$ is bijective if
and only if it is injective and surjective simultaneously.
\endproclaim
\demo{Proof} According to the statement of theorem~1.1, simultaneous
injectivity and surjectivity is necessary and sufficient condition
for bijectivity of the mapping $f\!:\,X\to Y$. Let's prove the
necessity of this condition for the beginning.\par
     Suppose that the mapping $f\!:\,X\to Y$ is bijective. Then for
any $y\in Y$ the total preimage $f^{-1}(y)$ consists of exactly one
element. This means that it is not empty. This fact proves the
surjectivity of the mapping $f\!:\,X\to Y$. \par
     However, we need to prove that $f$ is not only surjective, but
bijective as well. Let's prove the bijectivity of $f$ by contradiction.
If the mapping $f$ is not bijective, then there are two distinct
elements $x_1\neq x_2$ in $X$ such that $f(x_1)=f(x_2)$. Let's denote
$y=f(x_1)=f(x_2)$ and consider the total preimage $f^{-1}(y)$. From
the equality $f(x_1)=y$ we derive $x_1\in f^{-1}(y)$. Similarly from
$f(x_2)=y$ we derive $x_2\in f^{-1}(y)$. Hence, the total preimage
$f^{-1}(y)$ is a set containing at least two distinct elements $x_1$
and $x_2$. This fact contradicts the bijectivity of the mapping $f\!:\,X
\to Y$. Due to this contradiction we conclude that $f$ is surjective
and injective simultaneously. Thus, we have proved the necessity of
the condition stated in theorem~1.1.\par
    Let's proceed to the proof of sufficiency. Suppose that the
mapping $f\!:\,X\to Y$ is injective and surjective simultaneously.
Due to the surjectivity the sets $f^{-1}(y)$ are non-empty for all
$y\in Y$. Suppose that someone of them contains more than one
element. If $x_1\neq x_2$ are two distinct elements of the set
$f^{-1}(y)$, then  $f(x_1)=y=f(x_2)$. However, this equality
contradicts the injectivity of the mapping $f\!:\,X\to Y$. Hence,
each set $f^{-1}(y)$ is non-empty and contains exactly one element.
Thus, we have proved the bijectivity of the mapping $f$.
\qed\enddemo
\proclaim{Theorem 1.2} The mapping $f\!:\,X\to Y$ is surjective
if and only if \ $\Img f=Y$.
\endproclaim
\demo{Proof} If the mapping $f\!:\,X\to Y$ is surjective, then for
any element $y\in Y$ the total preimage $f^{-1}(y)$ is not empty.
Choosing some element $x\in f^{-1}(y)$, we get $y=f(x)$. Hence, each
element $y\in Y$ is an image of some element $x$ under the mapping
$f$. This proves the equality $\Img f=Y$.\par
    Conversely, if \ $\Img f=Y$, then any element $y\in Y$ is an image
of some element $x\in X$, i\.\,e\. $y=f(x)$. Hence, for any $y\in Y$ the
total preimage $f^{-1}(y)$ is not empty. This means that $f$ is a
surjective mapping.
\qed\enddemo
    Let's consider two mappings $f\!:\,X\to Y$ and $g\!:\,Y\to Z$.
Choosing an arbitrary element $x\in X$ we can apply $f$ to it. As
a result we get the element $f(x)\in Y$. Then we can apply $g$ to
$f(x)$. The successive application of two mappings $g(f(x))$ yields
a rule that associates each element $x\in X$ with some uniquely
determined element $z=g(f(x))\in Z$, i\.\,e\. we have a mapping
$\varphi\!:\,X\to Z$. This mapping is called {\it the composition}
of two mappings $f$ and $g$. It is denoted as $\varphi=g\compos f$.
\proclaim{Theorem 1.3} The composition $g\compos f$ of two injective
mappings $f\!:\,X\to Y$ and $g\!:\,Y\to Z$ is an injective mapping.
\endproclaim
\demo{Proof} Let's consider two elements $x_1$ and $x_2$ of the
set $X$. Denote $y_1=f(x_1)$ and $y_2=f(x_2)$. Therefore $g\compos
f(x_1)=g(y_1)$ and $g\compos f(x_2)=g(y_2)$. Due to the injectivity
of $f$ from $x_1\neq x_2$ we derive $y_1\neq y_2$. Then due to the
injectivity of $g$ from $y_1\neq y_2$ we derive $g(y_1)\neq g(y_2)$.
Hence, $g\compos f(x_1)\neq g\compos f(x_2)$. The injectivity of
the composition $g\compos f$ is proved.
\qed\enddemo
\proclaim{Theorem 1.4} The composition $g\compos f$ of two surjective
mappings $f\!:\,X\to Y$ and $g\!:\,Y\to Z$ is a surjective mapping.
\endproclaim
\demo{Proof} Let's take an arbitrary element $z\in Z$. Due to the
surjectivity of $g$ the total preimage $g^{-1}(z)$ is not empty.
Let's choose some arbitrary vector $y\in g^{-1}(z)$ and consider
its total preimage $f^{-1}(y)$. Due to the surjectivity of $f$ it
is not empty. Then choosing an arbitrary vector $x\in f^{-1}(y)$,
we get $g\compos f(x)=g(f(x))=g(y)=z$. This means that
$x\in (g\compos f)^{-1}(z)$. Hence, the total preimage $(g\compos
f)^{-1}(z)$ is not empty. The surjectivity of $g\compos f$ is proved.
\qed\enddemo
    As an immediate consequence of the above two theorems we obtain
the following theorem on composition of two bijections.
\proclaim{Theorem 1.5} The composition $g\compos f$ of two bijective
mappings $f\!:\,X\to Y$ and $g\!:\,Y\to Z$ is a bijective mapping.
\endproclaim
    Let's consider three mappings $f\!:\,X\to Y$, $g\!:\,Y\to Z$,
and $h\!:\,Z\to U$. Then we can form two different compositions of
these mappings:
$$
\xalignat 2
&\hskip -2em
\varphi=h\compos(g\compos f),
&&\psi=(h\compos g)\compos f.
\tag1.1
\endxalignat
$$
The fact of coincidence of these two mappings is formulated as
the following theorem on associativity.
\proclaim{Theorem 1.6} The operation of composition for the mappings
is an associative operation, i\.\,e\. $h\compos(g\compos f)=(h\compos
g)\compos f$.
\endproclaim
\demo{Proof} According to the definition~1.1, the coincidence of two
mappings\linebreak $\varphi\!:\,X\to U$ and $\psi\!:\,X\to U$ is verified
by verifying the equality $\varphi(x)=\psi(x)$ for an arbitrary element
$x\in X$. Let's denote $\alpha=h\compos g$ and $\beta=g\compos f$.
Then
$$
\hskip -2em
\aligned
&\varphi(x)=h\compos\beta(x)=h(\beta(x))=h(g(f(x))),\\
&\psi(x)=\alpha\compos f(x)=\alpha(f(x))=h(g(f(x))).
\endaligned
\tag1.2
$$
Comparing right hand sides of the equalities \thetag{1.2}, we derive
the required equality $\varphi(x)=\psi(x)$ for the mappings \thetag{1.1}.
Hence, $h\compos(g\compos f)=(h\compos g)\compos f$.
\qed\enddemo
     Let's consider a mapping $f\!:\,X\to Y$ and the pair of identical
mappings $\id_X\!:\,X\to X$ and $\id_Y\!:\,Y\to Y$.
The last two mappings are defined as follows:
$$
\xalignat 2
&\id_X(x)=x,&&\id_Y(y)=y.
\endxalignat
$$
\definition{Definition 1.5} A mapping $l\!:\,Y\to X$ is called
{\it left inverse} to the mapping $f\!:\,X\to Y$ if \ $l\compos
f=\id_X$.
\enddefinition
\definition{Definition 1.6} A mapping $r\!:\,Y\to X$ is called
{\it right inverse} to the mapping $f\!:\,X\to Y$ if \ $f\compos
r=\id_Y$.
\enddefinition
The problem of existence of the left and right inverse mappings
is solved by the following two theorems.
\proclaim{Theorem 1.7} A mapping $f\!:\,X\to Y$ possesses the left inverse
mapping $l$ if and only if it is injective.
\endproclaim
\proclaim{Theorem 1.8} A mapping $f\!:\,X\to Y$ possesses the right inverse
mapping $r$ if and only if it is surjective.
\endproclaim
\demo{Proof of the theorem~1.7} Suppose that the mapping $f$ possesses
the left inverse mapping $l$. Let's choose two vectors $x_1$ and $x_2$
in the space $X$ and let's denote $y_1=f(x_1)$ and $y_2=f(x_2)$. The
equality $l\compos f=\id_X$ yields $x_1=l(y_1)$ and $x_2=l(y_2)$.
Hence, the equality $y_1=y_2$ implies $x_1=x_2$ and $x_1\neq x_2$ implies
$y_1\neq y_2$. Thus, assuming the existence of left inverse mapping $l$,
we defive that the direct mapping $f$ is injective.\par
     Conversely, suppose that $f$ is an injective mapping. First of all 
let's choose and fix some element $x_0\in X$. Then let's consider an
arbitrary element $y\in\Img f$. Its total preimage $f^{-1}(y)$ is not 
empty. For any $y\in\Img f$ we can choose and fix some element $x_y\in 
f^{-1}(y)$ in non-empty set $f^{-1}(y)$. Then we define the mapping 
$l\!:\,Y\to X$ by the following equality:
$$
l(y)=
\cases
x_y &\text{for \ }y\in\Img f,\\
x_0 &\text{for \ }y\not\in\Img f.
\endcases
$$
Let's study the composition $l{\ssize\circ}f$. It is easy to see
that for any $x\in X$ and for $y=f(x)$ the equality $l{\ssize\circ}f(x)
=x_y$ is fulfilled. Then $f(x_y)=y=f(x)$. Taking into account the 
injectivity of $f$, we get $x_y=x$. Hence, $l{\ssize\circ}f(x)=x$ for
any $x\in X$. The equality $l\compos f=\id_X$ for the mapping
$l$ is proved. Therefore, this mapping is a required left inverse mapping  
for $f$. Theorem is proved.
\qed\enddemo
\demo{Proof of the theorem~1.8} Suppose that the mapping $f$ possesses
the right inverse mapping $r$. For an arbitrary element $y\in Y$, from
the equality $f\compos r=\id_Y$ we derive $y=f(r(y))$. This means that
$r(y)\in f^{-1}(y)$, therefore, the total preimage $f^{-1}(y)$ is not
empty. Thus, the surjectivity of $f$ is proved.\par
     Now, conversely, let's assume that $f$ is surjective. Then for any
$y\in Y$ the total preimage $f^{-1}(y)$ is not empty. In each non-empty
set $f^{-1}(y)$ we choose and mark exactly one element $x_y\in f^{-1}(y)$.
Then we can define a mapping by setting $r(y)=x_y$. Since $f(x_y)=y$, we
get $f(r(y))=y$ and $f\compos r=\id_Y$. The existence of the right inverse
mapping $r$ for $f$ is established.
\qed\enddemo
     Note that the mappings $l\!:\,Y\to X$ and $r\!:\,Y\to X$ constructed
when proving theorems~1.7 and 1.8 in general are not unique. Even the method 
of constructing them contains definite extent of arbitrariness.\par
\definition{Definition 1.7} A mapping $f^{-1}\!:\,Y\to X$ is called 
{\it bilateral inverse} mapping or simply {\it inverse} mapping for the 
mapping $f\!:\,X\to Y$ if
$$
\xalignat 2
&\hskip -2em
f^{-1}\compos f=\id_X,
&&f\compos f^{-1}=\id_Y.
\tag1.3
\endxalignat
$$
\enddefinition
\proclaim{Theorem 1.9} A mapping $f\!:\,X\to Y$ possesses both left and right inverse mappings $l$ and $r$ if and only if it is bijective. In 
this case the mappings $l$ and $r$ are uniquely determined. They coincide
with each other thus determining the unique bilateral inverse mapping
$l=r=f^{-1}$.
\endproclaim
\demo{Proof} The first proposition of the theorem~1.9 follows from
theorems~1.7, 1.8, and 1.1. Let's prove the remaining propositions of
this theorem~1.9. The coincidence $l=r$ is derived from the following 
chain of equalities:
$$
l=l\,\compos\id_Y=l\,\compos(f\compos r)=
(l\,\compos f)\compos\,r=\id_X\!\compos\,r=r.
$$
The uniqueness of left inverse mapping also follows from the same chain of
equalities. Indeed, if we assume that there is another left inverse mapping 
$l'$, then from $l=r$ and $l'=r$ it follows that $l=l'$.\par
     In a similar way, assuming the existence of another right inverse mapping
$r'$, we get $l=r$ and $l=r'$. Hence, $r=r'$. Coinciding with each other,
the left and right inverse mappings determine the unique bilateral inverse mapping $f^{-1}=l=r$ satisfying the equalities \thetag{1.3}.
\qed\enddemo
\head
\S\,2. Linear vector spaces.
\endhead
     Let $M$ be a set. {\it Binary algebraic operation in $M$} is a rule 
that maps each ordered pair of elements $x,y$ of the set $M$ to some uniquely 
determined element $z\in M$. This rule can be denoted as a function $z=f(x,y)$. This notation is called a {\it prefix} notation for 
an algebraic operation: the operation sign $f$ in it precedes the elements $x$ 
and $y$ to which it is applied. There is another {\it infix} notation
for algebraic operations, where the operation sign is placed between the
elements $x$ and $y$. Examples are the binary operations of addition and
multiplication of numbers: $z=x+y$, $z=x\cdot y$. Sometimes special
brackets play the role of the operation sign, while operands are separated
by comma. The vector product of three-dimensional vectors yields an example
of such notation: $\bold z=[\bold x,\,\bold y]$.\par
    Let $\Bbb K$ be a numeric field. Under the numeric field in this book
we shall understand one of three such fields: the field of rational numbers
$\Bbb K=\Bbb Q$, the field of real numbers $\Bbb K=\Bbb R$, or the field of
complex numbers $\Bbb K=\Bbb C$. The operation of {\it multiplication by
numbers from the field $\Bbb K$} in a set $M$ is a rule that maps each pair
$(\alpha, x)$ consisting of a number $\alpha\in\Bbb K$ and of an element
$x\in M$ to some element $y\in M$. The operation of multiplication by
numbers is written in infix form: $y=\alpha\cdot x$. The multiplication
sign in this notation is often omitted: $y=\alpha\,x$.
\definition{Definition 2.1} A set $V$ equipped with binary operation of
addition and with the operation of multiplication by numbers from the field
$\Bbb K$, is called a {linear vector space over the field $\Bbb K$}, if the
following conditions are fulfilled:
\roster
\item $\bold u+\bold v=\bold v+\bold u$ for all $\bold u,\bold v\in V$;
\item $(\bold u+\bold v)+\bold w=\bold u+(\bold v+\bold w)$ for all 
      $\bold u,\bold v,\bold w\in V$;
\item there is an element $\bold 0\in V$ such that $\bold v+\bold 0=\bold v$ 
      for all $\bold v\in V$; any such element is called a {\it zero element};
\item for any $\bold v\in V$ and for any zero element $\bold 0$ there is an 
      element $\bold v'\in V$ such that $\bold v+\bold v'=0$; it is called an 
      {\it opposite element} for $\bold v$;
\item $\alpha\cdot(\bold u+\bold v)=\alpha\cdot\bold u+\alpha\cdot\bold v$ 
      for any number $\alpha\in\Bbb K$ and for any two elements $\bold u,
      \bold v\in V$;
\item $(\alpha+\beta)\cdot\bold v=\alpha\cdot\bold v+\beta\cdot\bold v$ for 
      any two numbers $\alpha,\beta\in\Bbb K$ and for any element $\bold v
      \in V$;
\item $\alpha\cdot(\beta\cdot\bold v)=(\alpha\beta)\cdot\bold v$ for any two 
      numbers $\alpha,\beta\in\Bbb K$ and for any element $\bold v\in V$;
\item $1\cdot\bold v=\bold v$ for the number $1\in\Bbb K$ and for any element 
      $\bold v\in V$.
\endroster
\enddefinition
     The elements of a linear vector space are usually called the 
{\it vectors}, while the conditions \therosteritem{1}-\therosteritem{8} 
are called the {\it axioms of a linear vector space}. We shall distinguish 
{\it rational}, {\it real}, and {\it complex} linear vector spaces depending 
on which numeric field $\Bbb K=\Bbb Q$, $\Bbb K=\Bbb R$, or $\Bbb K=\Bbb C$ 
they are defined over. Most of the results in this book are valid for any 
numeric field $\Bbb K$. Formulating such results, we shall not specify 
the type of linear vector space.\par
     Axioms \therosteritem{1} and \therosteritem{2} are the axiom of 
{\it commutativity\footnotemark} and the axiom of {\it associativity}
respectively. Axioms \therosteritem{5} and \therosteritem{6} express
the {\it distributivity}.\footnotetext{\ The system of axioms 
\therosteritem{1}-\therosteritem{8} is excessive: the axiom 
\therosteritem{1} can be derived from other axioms. I~am grateful to 
A.~B.~Muftakhov who communicated me this curious fact.}\par
\adjustfootnotemark{-1}
\proclaim{Theorem 2.1} Algebraic operations in an arbitrary linear
vector space $V$ possess the following properties:
\roster
\item[9] zero vector $\bold 0\in V$ is unique;
\item for any vector $\bold v\in V$ the vector $\bold v'$ opposite to 
      $\bold v$ is unique;
\item the product of the number $0\in\Bbb K$ and any vector $\bold v\in V$ 
      is equal to zero vector: $0\cdot v=\bold 0$;
\item the product of an arbitrary number $\alpha\in K$ and zero vector
      is equal to zero vector: $\alpha\cdot\bold 0=\bold 0$;
\item the product of the number $-1\in\Bbb K$ and the vector $\bold v\in V$
      is equal to the opposite vector: $(-1)\cdot\bold v=\bold v'$.
\endroster
\endproclaim
\demo{Proof} The properties \therosteritem{9}-\therosteritem{13} are
immediate consequences of 
the axioms\linebreak\therosteritem{1}-\therosteritem{8}.
Therefore, they are enumerated so that their numbers form successive 
series with the numbers of the axioms of a linear vector space.\par
     Suppose that in a linear vector space there are two elements 
$\bold 0$ and $\bold 0'$ with the properties of zero vectors. Then 
for any vector $\bold v\in V$ due to the axiom~\therosteritem{3} we 
have $\bold v=\bold v+\bold 0$ and $\bold v+\bold 0'=\bold v$. Let's
substitute $\bold v=\bold 0'$ into the first equality and substitute 
$\bold v=\bold 0$ into the second one. Taking into account the
axiom~\therosteritem{1}, we get 
$$
\bold 0'=\bold 0'+\bold 0=\bold 0+\bold 0'=\bold 0.
$$
This means that the vectors $\bold 0$ and $\bold 0'$ do actually 
coincide. The uniqueness of zero vector is proved.\par
     Let $\bold v$ be some arbitrary vector in a vector space $V$. 
Suppose that there are two vectors $\bold v'$ and $\bold v''$ 
opposite to $\bold v$. Then
$$
\xalignat 2
&\bold v+\bold v'=\bold 0, & &\bold v+\bold v''=\bold 0.
\endxalignat
$$
The following calculations prove the uniqueness of opposite vector:
$$
\align
\bold v''=\bold v''+\bold 0&=\bold v''+(\bold v+\bold v')=
(\bold v''+\bold v)+\bold v'=\\
&=(\bold v+\bold v'')+\bold v'=\bold 0+\bold v'=\bold v'+
\bold 0=\bold v'.
\endalign
$$
In deriving $\bold v''=\bold v'$ above we used the 
axiom~\therosteritem{4}, the associativity axiom~\therosteritem{2} 
and we used twice the commutativity axiom~\therosteritem{1}.
\par
     Again, let $\bold v$ be some arbitrary vector in a vector space 
$V$. Let's take $\bold x=0\cdot\bold v$, then let's add $\bold x$ with 
$\bold x$ and apply the distributivity axiom~\therosteritem{6}. 
As a result we get
$$
\bold x+\bold x=0\cdot\bold v+0\cdot\bold v=(0
+0)\cdot\bold v=0\cdot\bold v=\bold x.
$$
Thus we have proved that $\bold x+\bold x=\bold x$. Then we easily 
derive that $\bold x=\bold 0$:
$$
\bold x=\bold x+\bold 0=\bold x+(\bold x+\bold x')=(\bold x+
\bold x)+\bold x'=\bold x+\bold x'=\bold 0.
$$
Here we used the associativity axiom~\therosteritem{2}. The 
property~\therosteritem{11} is proved.\par
     Let $\alpha$ be some arbitrary number of a numeric field 
$\Bbb K$. Let's take $\bold x=\alpha\cdot\bold 0$, where $\bold 0$ 
is zero vector of a vector space $V$. Then
$$
\bold x+\bold x=\alpha\cdot\bold 0+\alpha\cdot\bold 0=
\alpha\cdot(\bold 0+\bold 0)=\alpha\cdot\bold 0=\bold x.
$$
Here we used the axiom~\therosteritem{5} and the property of zero 
vector from the axiom~\therosteritem{3}. From the equality $\bold x+
\bold x=\bold x$ it follows that $\bold x=\bold 0$ (see above). 
Thus, the property~\therosteritem{12} is proved.\par
     Let $\bold v$ be some arbitrary vector of a vector space $V$. 
Let $\bold x=(-1)\cdot\bold v$. Applying axioms \therosteritem{8} 
and \therosteritem{6}, for the vector $\bold x$ we derive 
$$
\bold v+\bold x=1\cdot\bold v+\bold x=1\cdot\bold v+(-1)\cdot
\bold v=(1+(-1))\cdot\bold v=0\cdot\bold v=0.
$$
The equality $\bold v+\bold x=0$ just derived means that 
$\bold x$ is an opposite vector for the vector $\bold v$ in the 
sense of the axiom \therosteritem{4}. Due to the uniqueness 
property \therosteritem{10} of the opposite vector we conclude that 
$\bold x=\bold v'$. Therefore, $(-1)\cdot\bold v=\bold v'$. The 
theorem is completely proved.
\qed\enddemo
     Due to the commutativity and associativity axioms we need not 
worry about setting brackets and about the order of the summands 
when writing the sums of vectors. The property \therosteritem{13} 
and the axioms\therosteritem{7} and \therosteritem{8} yield
$$
(-1)\cdot\bold v'=(-1)\cdot((-1)\cdot\bold v)=((-1)(-1))\cdot
\bold v=1\cdot\bold v=\bold v.
$$
This equality shows that the notation $\bold v'=-\bold v$ for an 
opposite vector is quite natural. In addition, we can write
$$
-\alpha\cdot\bold v=-(\alpha\cdot\bold v)=(-1)\cdot(\alpha
\cdot\bold v)=(-\alpha)\cdot\bold v.
$$
The operation of {\it subtraction} is an opposite operation for 
the vector addition. It is determined as the addition with the 
opposite vector: $\bold x-\bold y=\bold x+(-\bold y)$. The 
following properties of the operation of vector subtraction 
$$
\align
&(\bold a+\bold b)-\bold c=\bold a+(\bold b-\bold c),\\
&(\bold a-\bold b)+\bold c=\bold a-(\bold b-\bold c),\\
&(\bold a-\bold b)-\bold c=\bold a-(\bold b+\bold c),\\
&\alpha\cdot(\bold x-\bold y)=\alpha\cdot\bold x-\alpha
\cdot\,\bold y
\endalign
$$
make the calculations with vectors very simple and quite similar to 
the calculations with numbers. Proof of the above properties is left 
to the reader.\par
     Let's consider some examples of linear vector spaces. Real 
{\it arithmetic vector space} $\Bbb R^n$ is determined as a set of 
ordered $n$-tuples of real numbers $x^1,\ldots,x^n$. Such $n$-tuples 
are represented in the form of {\it column vectors}. Algebraic operations 
with column vectors are determined as the operations with their components:
$$
\xalignat 2
&\hskip -2em
 \Vmatrix x^1\\ x^2\\ \vdots\\ x^n\endVmatrix+
 \Vmatrix y^1\\ y^2\\ \vdots\\ y^n\endVmatrix=
 \Vmatrix x^1+y^1\\ x^2+y^2\\ \vdots\\ x^n+y^n\endVmatrix
&&\alpha\cdot
 \Vmatrix x^1\\ x^2\\ \vdots\\ x^n\endVmatrix=
 \Vmatrix \alpha\cdot x^1\\ \alpha\cdot x^2\\ \vdots\\
 \alpha\cdot x^n\endVmatrix
\tag2.1
\endxalignat
$$
We leave to the reader to check the fact that the set $\Bbb R^n$ of 
all ordered $n$-tuples with algebraic operations \thetag{2.1} is a
linear vector space over the field $\Bbb R$ of real numbers. Rational
arithmetic vector space $\Bbb Q^n$ over the field $\Bbb Q$ of rational 
numbers and complex arithmetic vector space $\Bbb C^n$ over the field
$\Bbb C$ of complex numbers are defined in a similar way.\par
     Let's consider the set of $m$-times continuously differentiable
real-valued functions on the segment $[-1,1]$ of real axis. This set
is usually denoted as $C^m([-1,1])$. The operations of addition and
multiplication by numbers in $C^m([-1,1])$ are defined as pointwise
operations. This means that the value of the function $f+g$ at a point
$a$ is the sum of the values of $f$ and $g$ at that point. In a similar
way, the value of the function $\alpha\cdot f$ at the point $a$ is the
product of two numbers $\alpha$ and $f(a)$. It is easy to verify that 
the set of functions $C^m([-1,1])$ with pointwise algebraic operations 
of addition and multiplication by numbers is a linear vector space over 
the field of real numbers $\Bbb R$. The reader can easily verify this 
fact.\par
\definition{Definition 2.2} A non-empty subset $U\subset V$ in a linear
vector space $V$ over a numeric field $\Bbb K$ is called a subspace of
the space $V$ if:
\roster
\item from $\bold u_1,\bold u_2\in U$ it follows that $\bold u_1
      +\bold u_2\in U$;
\item from $\bold u\in U$ it follows that $\alpha\cdot\bold u\in U$ for 
      any number $\alpha\in\Bbb K$.
\endroster
\enddefinition
     Let $U$ be a subspace of a linear vector space $V$. Let's regard
$U$ as an isolated set. Due to the above conditions \therosteritem{1}
and \therosteritem{2} this set is closed with respect to operations of
addition and multiplication by numbers. It is easy to show that zero 
vector is an element of $U$ and for any $\bold u\in U$ the opposite 
vector $\bold u'$ also is an element of $U$. These facts follow from
$\bold 0=0\cdot\bold u$ and $\bold u'=(-1)\cdot\bold u$. Relying upon 
these facts one can easily prove that any subspace $U\subset V$, when
considered as an isolated set, is a linear vector space over the field 
$\Bbb K$. Indeed, we have already shown that axioms \therosteritem{3} 
and \therosteritem{4} are valid for it. Verifying axioms \therosteritem{1}, 
\therosteritem{2} and remaining axioms \therosteritem{5}-\therosteritem{8}
consists in checking equalities written in terms of the operations of 
addition and multiplication by numbers. Being fulfilled for arbitrary
vectors of $V$, these equalities are obviously fulfilled for vectors
of subset $U\subset V$. Since $U$ is closed with respect to algebraic
operations, it makes sure that all calculations in these equalities are
performed within the subset $U$.\par
     As the examples of the concept of subspace we can mention the
following subspaces in the functional space $C^m([-1,1])$:
\roster
\item"--" the subspace of even functions ($f(-x)=f(x)$);
\item"--" the subspace of odd functions ($f(-x)=-f(x)$);
\item"--" the subspace of polynomials ($f(x)=a_n\,x^n+\ldots+a_1\,x+a_0$).
\endroster\par
\head
\S\,3. Linear dependence and linear independence. 
\endhead
     Let $\bold v_1,\,\ldots,\,\bold v_n$ be a system of vectors some from 
some linear vector space $V$. Applying the operations of multiplication 
by numbers and addition to them we can produce the following expressions
with these vectors:
$$
\hskip -2em
v=\alpha_1\cdot\bold v_1+\ldots+\alpha_n\cdot\bold v_n.
\tag3.1
$$
An expression of the form \thetag{3.1} is called a {\it linear 
combination} of the vectors $\bold v_1,\,\ldots,\,\bold v_n$. The numbers
$\alpha_1,\ldots,\alpha_n$ are taken from the field $\Bbb K$; they 
are called the {\it coefficients} of the linear combination 
\thetag{3.1}, while vector $\bold v$ is called the {\it value} of 
this linear combination. Linear combination is said to be {\it zero}
or {\it equal to zero} if its value is zero.\par
    A linear combination is called {\it trivial} if all its coefficients 
are equal to zero: $\alpha_1=\ldots=\alpha_n=0$. Otherwise it is called 
{\it nontrivial}.\par
\definition{Definition 3.1} A system of vectors $\bold v_1,\ldots,\bold
v_n$ in linear vector space $V$ is called {\it linearly dependent} if 
there exists some nontrivial linear combination of these vectors equal 
to zero.
\enddefinition
\definition{Definition 3.2} A system of vectors $\bold v_1,\,\ldots,\,
\bold v_n$ in linear vector space $V$ is called {\it linearly independent}
if any linear combination of these vectors being equal to zero is
necessarily trivial.
\enddefinition
     The concept of linear independence is obtained by direct logical
negation of the concept of linear dependence. The reader can give several
equivalent statements defining this concept. Here we give only one of such
statements which, to our knowledge, is most convenient in what follows.
\par
     Let's introduce one more concept related to linear combinations. We
say that vector $\bold v$ is {\it linearly expressed} through the vectors
$\bold v_1,\ldots,\bold v_n$ if $\bold v$ is the value of some linear
combination composed of \pagebreak $\bold v_1,\,\ldots,\,\bold v_n$.\par 
\proclaim{Theorem 3.1} The relation of linear dependence of vectors in a linear vector space has the following basic properties:
\roster
\item any system of vectors comprising zero vector is linearly dependent;
\item any system of vectors comprising linearly dependent subsystem is
      linearly dependent in whole;
\item if a system of vectors is linearly dependent, then at least one of 
      these vectors is linearly expressed through others;
\item if a system of vectors $\bold v_1,\,\ldots,\,\bold v_n$ is linearly
      independent and if adding the next vector $\bold v_{n+1}$ to it 
      we make it linearly dependent, then the vector $\bold v_{n+1}$ 
      is linearly expressed through previous vectors $\bold v_1,\,
      \ldots,\,\bold v_n$;
\item if a vector $\bold x$ is linearly expressed through the vectors 
      $\bold y_1,\,\ldots,\,\bold y_m$ and if each one of the vectors
      $\bold y_1,\,\ldots,\,\bold y_m$ is linearly expressed through 
      $\bold z_1,\,\ldots,\,\bold z_n$, then $\bold x$ is linearly
      expressed through $\bold z_1,\,\ldots,\,\bold z_n$.
\endroster
\endproclaim
\demo{Proof} Suppose that a system of vectors $\bold v_1,\,\ldots,\,
\bold v_n$ comprises zero vector. For the sake of certainty we can 
assume that $\bold v_k=0$. Let's compose the following linear combination
of the vectors $\bold v_1,\,\ldots,\,\bold v_n$:
$$
0\cdot\bold v_1+\ldots+0\cdot\bold v_{k-1}+1\cdot\bold v_k+0\cdot\bold
v_{k+1}+\ldots+0\cdot\bold v_n=\bold 0.
$$
This linear combination is nontrivial since the coefficient of vector 
$\bold v_k$ is nonzero. And its value is equal to zero. Hence, the vectors
$\bold v_1,\,\ldots,\,\bold v_n$ are linearly dependent. The property
\therosteritem{1} is proved.
     Suppose that a system of vectors $\bold v_1,\,\ldots,\,\bold v_n$
comprises a linear dependent subsystem. Since linear dependence is not
sensible to the order in which the vectors in a system are enumerated, 
we can assume that first $k$ vectors form linear dependent subsystem 
in it. Then there exists some nontrivial liner combination of these $k$
vectors being equal to zero:
$$
\alpha_1\cdot\bold v_1+\ldots+\alpha_k\cdot\bold v_k=\bold 0.
$$
Let's expand this linear combination by adding other vectors with
zero coefficients:
$$
\alpha_1\cdot\bold v_1+\ldots+\alpha_k\cdot\bold v_k+0\cdot\bold 
v_{k+1}+\ldots+0\cdot\bold v_n=\bold 0.
$$
It is obvious that the resulting linear combination is nontrivial and 
its value is equal to zero. Hence, the vectors $\bold v_1,\,\ldots,\,
\bold v_n$ are linearly dependent. The property \therosteritem{2} is
proved.\par
     Let assume that the vectors $\bold v_1,\,\ldots,\,\bold v_n$ are
linearly dependent. Then there exists a nontrivial linear combination 
of them being equal to zero:
$$
\hskip -2em
\alpha_1\cdot\bold v_1+\ldots+\alpha_n\cdot\bold v_n=\bold 0.
\tag3.2
$$
Non-triviality of the linear combination \thetag{3.2} means that at least
one of its coefficients is nonzero. Suppose that $\alpha_k\neq 0$. Let's
write \thetag{3.2} in more details:
$$
\pagebreak
\alpha_1\cdot\bold v_1+\ldots+\alpha_k\cdot\bold v_k+\ldots
+\alpha_n\cdot\bold v_n=\bold 0.
$$
Let's move the term $\alpha_k\cdot\bold v_k$ to the right hand side 
of the above equality, and then let's divide the equality by 
$-\alpha_k$:
$$
\bold v_k=-\frac{\alpha_1}{\alpha_k}\cdot\bold v_1
-\ldots-\frac{\alpha_{k-1}}{\alpha_k}\cdot\bold v_{k-1}
-\frac{\alpha_{k+1}}{\alpha_k}\cdot\bold v_{k+1}
-\ldots-\frac{\alpha_{n}}{\alpha_k}\cdot\bold v_n.
$$
Now we see that the vector $\bold v_k$ is linearly expressed through 
other vectors of the system. The property \therosteritem{3} is proved.
\par
     Let's consider a linearly independent system of vectors $\bold v_1,
\,\ldots,\,\bold v_n$ such that adding the next vector $\bold v_{n+1}$
to it we make it linearly dependent. Then there is some nontrivial linear
combination of vectors $\bold v_1,\,\ldots,\,\bold v_{n+1}$ being equal 
to zero:
$$
\alpha_1\cdot\bold v_1+\ldots+\alpha_n\cdot\bold v_n+\alpha_{n+1}\cdot
\bold v_{n+1}=\bold 0.
$$
Let's prove that $\alpha_{n+1}\neq 0$. If, conversely, we assume that
$\alpha_{n+1}=0$, we would get the nontrivial linear combination of 
$n$ vectors being equal to zero:
$$
\alpha_1\cdot\bold v_1+\ldots+\alpha_n\cdot\bold v_n=\bold 0.
$$
This contradicts to the linear independence of the first $n$ vectors
$\bold v_1,\,\ldots,\,\bold v_n$. Hence, $\alpha_{n+1}\neq 0$, and 
we can apply the trick already used above:
$$
\bold v_{n+1}=-\frac{\alpha_1}{\alpha_{n+1}}\cdot\bold v_1
-\ldots-\frac{\alpha_n}{\alpha_{n+1}}\cdot\bold v_n.
$$
This expression for the vector $\bold v_{n+1}$ completes the proof of 
the property \therosteritem{4}.\par
     Suppose that the vector $\bold x$ is linearly expressed through
$\bold y_1,\,\ldots,\,\bold y_m$, and each one of the vectors
$\bold y_1,\,\ldots,\,\bold y_m$ is linearly expressed through
$\bold z_1,\,\ldots,\,\bold z_n$. This fact is expressed by the 
following formulas:
$$
\xalignat 2
&x=\sum^m_{i=1}\alpha_i\cdot y_i,
&&y_i=\sum^n_{j=1}\beta_{ij}\cdot z_j.
\endxalignat
$$
Substituting second formula into the first one, for the vector $\bold x$
we get
$$
x=\sum^m_{i=1}\alpha_i\cdot\left(\,\shave{\sum^n_{j=1}}
\beta_{ij}\cdot z_j\right)=\sum^n_{j=1}
\left(\,\shave{\sum^m_{i=1}}\alpha_i\,\beta_{ij}\right)
\cdot z_j
$$
The above expression for the vector $\bold x$ shows that it is linearly
expressed through vectors $\bold z_1,\,\ldots,\,\bold z_n$. The property
\therosteritem{5} is proved. This completes the proof of theorem~3.1 in
whole.
\qed\enddemo
     Note the following important consequence that follows from the
property~\therosteritem{2} in the theorem~3.1.
\proclaim{Corollary} Any subsystem in a linearly independent system of
vectors is linearly\pagebreak independent.
\endproclaim
     The next property of linear dependence of vectors is known as 
Steinitz theorem. It describes some quantitative feature of this concept.
\proclaim{Theorem 3.2 (Steinitz)} If the vectors $\bold x_1,\,\ldots,\,
\bold x_n$ are linear independent and if each of them is expressed 
through the vectors $\bold y_1,\,\ldots,\,\bold y_m$, then 
$m\geqslant n$.
\endproclaim
\demo{Proof} We shall prove this theorem by induction on the number of
vectors in the system $\bold x_1,\,\ldots,\,\bold x_n$. Let's begin with 
the case $n=1$. Linear independence of a system with a single vector
$\bold x_1$ means that $\bold x_1\neq\bold 0$. In order to express the
nonzero vector $\bold x_1$ through the vectors of a system $\bold y_1,\,
\ldots,\,\bold y_m$ this system should contain at least one vector. Hence,
$m\geqslant 1$. The base step of induction is proved.\par
     Suppose that the theorem holds for the case $n=k$. Under this
assumption let's prove that it is valid for $n=k+1$. If $n=k+1$ we have
a system of linearly independent vectors $\bold x_1,\,\ldots,\,\bold
x_{k+1}$, each vector being expressed through the vectors of another
system $\bold y_1,\,\ldots,\,\bold y_m$. We express this fact by
formulas 
$$
\hskip -2em
\aligned
&\bold x_1=\alpha_{11}\cdot\,\bold y_1+\ldots+\alpha_{1m}\cdot\,\bold
y_m,\\
&\text{\hbox to 4.8cm{\leaders\hbox to 1em{\hss.\hss}\hfill}}\\
&\bold x_k=\alpha_{k1}\cdot\,\bold y_1+\ldots+\alpha_{km}\cdot\,\bold y_m.
\endaligned
\tag3.3
$$
We shall write the analogous formula expressing $\bold x_{k+1}$ through
$\bold y_1,\,\ldots,\,\bold y_m$ in a slightly different way:
$$
\bold x_{k+1}=\beta_1\cdot\,\bold y_1+\ldots+\beta_m\cdot\,\bold y_m.
$$
Due to the linear independence of vectors $\bold x_1,\,\ldots,\,\bold
x_{k+1}$ the last vector $x_{k+1}$ of this system is nonzero (as well 
as other ones). Therefore at least one of the numbers $\beta_1,\ldots,
\beta_m$ is nonzero. Upon renumerating the vectors $\bold y_1,\,\ldots,
\,\bold y_m$, if necessary, we can assume that $\beta_m\neq 0$. Then
$$
\hskip -2em
\bold y_m=\frac{1}{\beta_m}\cdot\,\bold x_{k+1}-\frac{\beta_1}{\beta_m}
\cdot\,\bold y_1-\ldots-\frac{\beta_{m-1}}{\beta_m}\cdot\,\bold y_{m-1}.
\tag3.4
$$
Let's substitute \thetag{3.4} into the relationships \thetag{3.3} and
collect similar terms in them. As a result the relationships \thetag{3.4}
are written as 
$$
\hskip -2em
\bold x_i-\frac{\alpha_{im}}{\beta_m}\cdot\,\bold x_{k+1}=
\sum^{m-1}_{j=1}\left(\alpha_{ij}-\beta_j\frac{\alpha_{im}}
{\beta_m}\right)\!\cdot\,\bold y_j,
\tag3.5
$$
where $i=1,\,\ldots,\,k$. In order to simplify \thetag{3.5} we introduce
the following notations:
$$
\xalignat 2
&\hskip -2em
x^*_i=x_i-\frac{\alpha_{im}}{\beta_m}\cdot x_{k+1},
&&\alpha^*_{ij}=\alpha_{ij}-\beta_j\frac{\alpha_{im}}{\beta_m}.
\tag3.6
\endxalignat
$$
In these notations the formulas \thetag{3.5} are written as
$$
\hskip -2em
\aligned
&\bold x^*_1=\alpha^*_{11}\cdot\,\bold y_1+\ldots+\alpha^*_{1\,m-1}
\cdot\,\bold y_{m-1},\\
&\text{\hbox to 5.4cm{\leaders\hbox to 1em{\hss.\hss}\hfill}}\\
&\bold x^*_k=\alpha^*_{k\,1}\cdot\,\bold y_1+\ldots+\alpha^*_{k\,m-1}
\cdot\,\bold y_{m-1}.
\endaligned
\tag3.7
$$
According to the above formulas, $k$ vectors $\bold x^*_1,\,\ldots,
\,\bold x^*_k$ are linearly expressed through $\bold y_1,\,\ldots,
\,y_{m-1}$. In order to apply the inductive hypothesis we need to
show that the vectors $\bold x^*_1,\,\ldots,\,\bold x^*_k$ are
linearly independent. Let's consider a linear combination of these
vectors being equal to zero:
$$
\hskip -2em
\gamma_1\cdot\bold x^*_1+\ldots+\gamma_k\cdot\bold x^*_k=\bold 0.
\tag3.8
$$
Substituting \thetag{3.6} for $x^*_i$ in \thetag{3.8}, upon collecting
similar terms, we get
$$
\gamma_1\cdot\bold x_1+\ldots+\gamma_k\cdot\bold x_k-
\left(\shave{\sum^k_{i=1}}\gamma_i\frac{\alpha_{im}}{\beta_m}
\right)\cdot\bold x_{k+1}=0.
$$
Due to the linear independence of the initial system of vectors 
$\bold x_1,\,\ldots,\,\bold x_{k+1}$ we derive $\gamma_1=\ldots=
\gamma_k=0$. Hence, the linear combination \thetag{3.8} is trivial,
which proves the linear independence of vectors $\bold x^*_1,\,
\ldots,\,\bold x^*_k$. Now, applying the inductive hypothesis to
the relationships \thetag{3.7}, we get $m-1\geqslant k$. The required
inequality $m\geqslant k+1$ proving the theorem for the case $n=k+1$ 
is an immediate consequence of $m\geqslant k+1$. So, the inductive step 
is completed and the theorem is proved.\qed\enddemo
\head
\S\,4. Spanning systems and bases.
\endhead
     Let $S\subset V$ be some non-empty subset in a linear vector space
$V$. The set $S$ can consist of either finite number of vectors, or of
infinite number of vectors. We denote by $\langle S\rangle$ the set of
all vectors, each of which is linearly expressed through some finite
number of vectors taken from $S$:
$$
\langle S\rangle=\{\bold v\in V\!:\ \exists\,n\ (\bold v=\alpha_1
\cdot\bold s_1+\ldots+\alpha_n\cdot\bold s_n\text{, where\ }
\bold s_i\in S)\}.
$$
This set $\langle S\rangle$ is called the {\it linear span} of a subset 
$S\subset V$.
\proclaim{Theorem 4.1}The linear span of any subset $S\subset V$ 
is a subspace in a linear vector space $V$.
\endproclaim
\demo{Proof} In order to prove this theorem it is sufficient to check
two conditions from the definition~2.2 for $\langle S\rangle$. Suppose
that $\bold u_1,\bold u_2\in\langle S\rangle$. Then
$$
\align
&\bold u_1=\alpha_1\cdot\bold s_1+\ldots+\alpha_n\cdot\bold s_n,\\
&\bold u_2=\beta_1\cdot\bold s^*_1+\ldots+\beta_m\cdot\bold s^*_m.
\endalign
$$
Adding these two equalities, we see that the vector $\bold u_1+\bold u_2$
also is expressed as a linear combination of some finite number of vectors
taken from $S$. Therefore, we have $\bold u_1+\bold u_2\in\langle S\rangle$.\par
     Now suppose that $\bold u\in\langle S\rangle$. Then $\bold u=
\alpha_1\cdot\bold s_1+\ldots+\alpha_n\cdot\bold s_n$. For the 
vector\linebreak $\alpha\cdot\bold u$, from this equality we derive
$$
\alpha\cdot\bold u=(\alpha\,\alpha_1)\cdot\bold s_1+\ldots+
(\alpha\,\alpha_n)\cdot\bold s_n.
$$
Hence, $\alpha\cdot u\in\langle S\rangle$. Both conditions \therosteritem{1} and \therosteritem{2} from the definition~2.2 
for $\langle S\rangle$ are fulfilled. Thus, the theorem is proved.
\qed\enddemo
\proclaim{Theorem 4.2} The operation of passing to the linear span in
a linear vector space $V$ possesses the following properties:
\roster
\item if $S\subset U$ and if $U$ is a subspace in $V$, then
      $\langle S\rangle\subset U$;
\item the linear span of a subset $S\subset V$ is the intersection of
      all subspaces comprising this subset $S$.
\endroster
\endproclaim
\demo{Proof} Let $\bold u\in\langle S\rangle$ and $S\subset U$, where 
$U$ is a subspace. Then for the vector $\bold u$ we have $\bold u=
\alpha_1\cdot\bold s_1+\ldots+\alpha_n\cdot\bold s_n$, where $\bold s_i
\in S$. But $\bold s_i\in S$ and $S\subset U$ implies $\bold s_i\in U$.
Since $U$ is a subspace, the value of any linear combination of its
elements again is an element of $U$. Hence, $\bold u\in U$. This proves
the inclusion $\langle S\rangle\subset U$.\par
     Let's denote by $W$ the intersection of all subspaces of $V$
comprising the subset $S$. Due to the property \therosteritem{1}, which
is already proved, the subset $\langle S\rangle$ is included into each
of such subspaces. Therefore, $\langle S\rangle\subset W$. On the other
hand, $\langle S\rangle$ is a subspace of $V$ comprising the subset $S$
(see theorem~4.1). Hence, $\langle S\rangle$ is among those subspaces
forming $W$. Then $W\subset\langle S\rangle$. From the two inclusions
$\langle S\rangle\subset W$ and $W\subset\langle S\rangle$ it follows
that $\langle S\rangle=W$. The theorem is proved.
\qed\enddemo
     Let $\langle S\rangle=U$. Then we say that the subset $S\subset V$
{\it spans} the subspace $U$, i\.\,e\. $S$ {\it generates} $U$ by means
of the linear combinations. This terminology is supported by the following 
definition.\par
\definition{Definition 4.1} A subset $S\subset V$ is called a 
{\it generating subset} or a {\it spanning system of vectors} in a linear
vector space $V$ if $\langle S\rangle=V$.
\enddefinition
     A linear vector space $V$ can have multiple spanning systems.
Therefore the problem of choosing of a minimal (is some sense)
spanning system is reasonable.\par
\definition{Definition 4.2} A spanning system of vectors $S\subset V$
in a linear vector space $V$ is called a {\it minimal spanning system} 
if none of smaller subsystems $S'\varsubsetneq S$ is a spanning system 
in $V$, i\.\,e\. if $\langle S'\rangle\neq V$ for all $S'\varsubsetneq S$.
\enddefinition
\definition{Definition 4.3} A system of vectors $S\subset V$ is called
{\it linearly independent} if any finite subsystem of vectors $\bold s_1,
\,\ldots,\,\bold s_n$ taken from $S$ is linearly independent.
\enddefinition
    This definition extends the definition~3.2 for the case of infinite
systems of vectors. As for the spanning systems, the relation of the
properties of minimality and linear independence for them is determined
by the following theorem.     
\proclaim{Theorem 4.3} A spanning system of vectors $S\subset V$ is 
minimal if and only if it is linearly independent.
\endproclaim
\demo{Proof} If a spanning system of vectors $S\subset V$ is linearly
dependent, then it contains some finite linearly dependent set of
vectors $\bold s_1, \,\ldots,\,\bold s_n$. Due to the item
\therosteritem{3} in the statement of theorem~3.1 one of these vectors
$\bold s_k$ is linearly expressed through others. Then the subsystem
$S'=S\,\diagdown\,\{\bold s_k\}$ obtained by omitting this vector 
$\bold s_k$ from $S$ is a spanning system in $V$. This fact obviously
contradicts the minimality of $S$ (see definition~4.2 above). Therefore
any minimal spanning system of vectors in $V$ is linearly independent.
\par
     If a spanning system of vectors $S\subset V$ is not minimal, then
there is some smaller spanning subsystem $S'\varsubsetneq S$, i\.\,e\.
subsystem $S'$ such that 
$$
\hskip -2em
\langle S'\rangle=\langle S\rangle=V.
\tag4.1
$$
In this case we can choose some vector $\bold s_0\in S$ such that
$\bold s_0\notin S'$. Due to \thetag{4.1} this vector is an element
of $\langle S'\rangle$. Hence, $\bold s_0$ is linearly expressed 
through some finite number of vectors taken from the subsystem $S'$:
$$
\hskip -2em
\bold s_0=\alpha_1\cdot\bold s_1+\ldots+\alpha_n\cdot\bold s_n.
\tag4.2
$$
One can easily transform \thetag{4.2} to the form of a linear 
combination equal to zero:
$$
\hskip -2em
(-1)\cdot\bold s_0+\alpha_1\,\cdot\bold s_1+\ldots+\alpha_n
\cdot\bold s_n=\bold 0.
\tag4.3
$$
This linear combination is obviously nontrivial. Thus, we have
found that the vectors $\bold s_0,\,\ldots,\,\bold s_n$ form a 
finite linearly dependent subset of $S$. Hence, $S$ is linearly
dependent (see the item~\therosteritem{2} in theorem~3.1 and the
definition~4.2). This fact means that any linearly independent
spanning system of vector in $V$ is minimal. 
\qed\enddemo
\definition{Definition 4.4} A linear vector space $V$ is called
{\it finite dimensional} if there is some finite spanning system
of vectors $S=\{\bold x_1,\ldots,\bold x_n\}$ in it. 
\enddefinition
     In an arbitrary linear vector space $V$ there is at lease one
spanning system, e\.\,g\. $S=V$. However, the problem of existence 
of minimal spanning systems in general case is nontrivial. The 
solution of this problem is positive, but it is not elementary and 
it is not constructive. This problem is solved with the use of the 
{\it axiom of choice} (see \cite{1}). Finite dimensional vector 
spaces are distinguished due to the fact that the proof of existence 
of minimal spanning systems for them is elementary.\par
\proclaim{Theorem 4.4} In a finite dimensional linear vector space
$V$ there is at least one minimal spanning system of vectors. Any
two of such systems $\{\bold x_1,\,\ldots,\,\bold x_n\}$ and 
$\{\bold y_1,\,\ldots,\,\bold y_n\}$ have the same number of elements 
$n$. This number $n$ is called the {\it dimension} of $V$, it is 
denoted as $n=\dim V$.
\endproclaim
\demo{Proof} Let $S=\{\bold x_1,\,\ldots,\,\bold x_k\}$ be some finite
spanning system of vectors in a finite-dimensional linear vector space
$V$. If this system is not minimal, then it is linear dependent. Hence,
one of its vectors is linearly expressed through others. This vector can
be omitted and we get the smaller spanning system $S'$ consisting of 
$k-1$ vectors. If $S'$ is not minimal again, then we can iterate the
process getting one less vectors in each step. Ultimately, we shall get
a minimal spanning system $S_{\text{min}}$ in $V$ with finite number of
vectors $n$ in it:
$$
\hskip -2em
S_{\text{min}}=\{\bold y_1,\,\ldots,\,\bold y_n\}.
\tag4.4
$$\par
     Usually, the minimal spanning system of vectors \thetag{4.4} is not
unique. Suppose that $\{\bold x_1,\,\ldots,\,\bold x_m\}$ is some other
minimal spanning system in $V$. Both systems $\{\bold x_1,\,\ldots,\,\bold
x_m\}$ and $\{\bold y_1,\,\ldots,\,\bold y_n\}$ are linearly independent
and 
$$
\hskip -2em
\aligned
&x^i\in\langle \bold y_1,\,\ldots,\,\bold y_n\rangle\text{\ \ for \ }i=1,\,\ldots,\,m,\\
&y^i\in\langle \bold x_1,\,\ldots,\,\bold x_m\rangle\text{\ \ for \ }i=1,\,\ldots,\,n.
\endaligned
\tag4.5
$$
Due to \thetag{4.5} we can apply Steinitz theorem~3.2 to the systems
of vectors $\{\bold x_1,\,\ldots,\,\bold x_m\}$ and $\{\bold y_1,\,
\ldots,\,\bold y_n\}$. As a result we get two inequalities
$n\geqslant m$ and $m\geqslant n$. Therefore, $m=n=\dim V$. The 
\pagebreak theorem is proved.
\qed\enddemo
     The dimension $\dim V$ is an integer invariant of a 
finite-dimensional linear vector space. If $\dim V=n$, then such a space
is called an $n$-dimensional space. Returning to the examples of linear
vector spaces considered in \S~2, note that $\dim\Bbb R^n=n$, while the
functional space $C^m([-1,1])$ is not finite-dimensional at all.\par
\proclaim{Theorem 4.5} Let $V$ be a finite dimensional linear
vector space. Then the following propositions are valid:
\roster
\rosteritemwd=5pt
\item the number of vectors in any linearly independent system of vectors
      $\bold x_1,\,\ldots,\,\bold x_k$ in $V$ is not greater than the
      dimension of $V$;
\item any subspace $U$ of the space $V$ is finite-dimensional and \ 
      $\dim U\leqslant\dim V$;
\item for any subspace $U$ in $V$ if \ $\dim U=\dim V$, then \ $U=V$;
\item any linearly independent system of $n$ vectors $\bold x_1,\,
      \ldots,\,\bold x_n$, where $n=\dim V$, is a spanning system in $V$.
\endroster
\endproclaim
\demo{Proof} Suppose that $\dim V=n$. Let's fix some minimal spanning
system of vectors $\bold y_1,\,\ldots,\,\bold y_n$ in $V$. Then each
vector of the linear independent system of vectors $\bold x_1,\,
\ldots,\,\bold x_k$ in proposition \therosteritem{1} is linearly
expressed through $\bold y_1,\,\ldots,\,\bold y_n$. Applying Steinitz
theorem~3.2, we get the inequality $k\leqslant n$. The first proposition
of theorem is proved.\par
     Let's consider all possible linear independent systems $\bold u_1,
\,\ldots,\,\bold u_k$ composed by the vectors of a subspace $U$. Due to
the proposition \therosteritem{1}, which is already proved, the number 
of vectors in such systems is restricted. It is not greater than 
$n=\dim V$. Therefore we can assume that $\bold u_1,\,\ldots,\,\bold 
u_k$ is a linearly independent system with maximal number of vectors:
$k=k_{\text{max}}\leqslant n=\dim V$. If $\bold u$ is an arbitrary
vector of the subspace $U$ and if we add it to the system $\bold u_1,\,
\ldots,\,\bold u_k$, we get a linearly dependent system; this is 
because $k=k_{\text{max}}$. Now, applying the property \therosteritem{4}
from the theorem~3.1, we conclude that the vector $\bold u$ is linearly
expressed through the vectors $\bold u_1,\,\ldots,\,\bold u_k$. Hence,
the vectors $\bold u_1,\,\ldots,\,\bold u_k$ form a finite spanning 
system in $U$. It is minimal since it is linearly independent (see 
theorem~4.3). Finite dimensionality of $U$ is proved. The estimate for
its dimension follows from the above inequality: $\dim U=k\leqslant n=
\dim V$.\par
     Let $U$ again be a subspace in $V$. Assume that $\dim U=\dim V=n$.
Let's choose some minimal spanning system of vectors $\bold u_1,\,
\ldots,\,\bold u_n$ in $U$. It is linearly independent. Adding an
arbitrary vector $\bold v\in V$ to this system, we make it linearly
dependent since in $V$ there is no linearly independent system with 
$(n+1)$ vectors (see proposition \therosteritem{1}, which is already 
proved). Furthermore, applying the property \therosteritem{3} from 
the theorem~3.1 to the system $\bold u_1,\,\ldots,\,\bold u_n,\,\bold
v$, we find that
$$
\bold v=\alpha_1\cdot\bold u_1+\ldots+\alpha_m\cdot\bold u_m.
$$
This formula means that $\bold v\in U$, where $\bold v$ is an arbitrary
vector of the space $V$. Therefore, $U=V$. The third proposition of the
theorem is proved.\par
     Let $\bold x_1,\,\ldots,\,\bold x_n$ be a linearly independent 
system of $n$ vectors in $V$, where $n$ is equal to the dimension of
the space $V$. Denote by $U$ the linear span of this system of
vectors: $U=\langle\bold x_1,\,\ldots,\,\bold x_n\rangle$. 
Since $\bold x_1,\,\ldots,\,\bold x_n$ are linearly independent, they 
form a minimal spanning system in $U$. Therefore, $\dim U=n=\dim V$.
Now, applying proposition \therosteritem{3} of the theorem, we get
$$
\pagebreak 
\langle\bold x_1,\,\ldots,\,\bold x_n\rangle=U=V.
$$
This equality proves the fourth proposition of theorem~4.5 and completes
the proof of the theorem in whole.
\qed\enddemo
\definition{Definition 4.5} A minimal spanning system $\bold e_1,\,
\ldots,\,\bold e_n$ with some fixed order of vectors in it is called 
a {\it basis} of a finite-dimensional vector space $V$.
\enddefinition
\proclaim{Theorem 4.6 (basis criterion)} An ordered system of vectors 
$\bold e_1,\,\ldots,\,\bold e_n$ is a basis in a finite-dimensional
vector space $V$ if and only if 
\roster
\item the vectors $\bold e_1,\,\ldots,\,\bold e_n$ are linearly 
      independent;
\item an arbitrary vector of the space $V$ is linearly expressed
      through $\bold e_1,\,\ldots,\,\bold e_n$.
\endroster
\endproclaim
     Proof is obvious. The second condition of theorem means that
the vectors $\bold e_1,\,\ldots,\,\bold e_n$ form a spanning system
in $V$, while the first condition is equivalent to its minimality.
\par
     In essential, theorem~4.6 simply reformulates the 
definition~4.5. We give it here in order to simplify the terminology.
The terms {\tencyr\char '074}spanning system{\tencyr\char '076} and
{\tencyr\char '074}minimal spanning system{\tencyr\char '076} are
huge and inconvenient for often usage.\par
\proclaim{Theorem 4.7} Let $\bold e_1,\,\ldots,\,\bold e_s$ be a
basis in a subspace $U\subset V$ and let $\bold v\in V$ be some 
vector outside this subspace: $\bold v\notin U$. Then the system 
of vectors $\bold e_1,\,\ldots,\,\bold e_s,\,\bold v$ is a
linearly independent system.
\endproclaim
\demo{Proof} Indeed, if the system of vectors $\bold e_1,\,\ldots,
\,\bold e_s,\,\bold v$ is linearly dependent, while $\bold e_1,\,
\ldots,\,\bold e_s$ is a linearly independent system, then 
$\bold v$ is linearly expressed through the vectors $\bold e_1,\,
\ldots,\,\bold e_s$, thus contradicting the condition $\bold v
\notin U$. This contradiction proves the theorem~4.7.
\qed\enddemo
\proclaim{Theorem 4.8 (on completing the basis)} Let $U$ be a
subspace in a finite-dimensional linear vector space $V$. Then
any basis $\bold e_1,\,\ldots,\,\bold e_s$ of $U$ can be 
completed up to a basis $\bold e_1,\,\ldots,\,\bold e_s,\,
\bold e_{s+1},\,\ldots,\,\bold e_n$ in $V$.
\endproclaim
\demo{Proof} Let's denote $U=U_0$. If $U_0=V$, then there is no need
to complete the basis since $\bold e_1,\,\ldots,\,\bold e_s$ is a basis
in $V$. Otherwise, if $U_0\neq V$, then let's denote by $\bold e_{s+1}$
some arbitrary vector of $V$ taken outside the subspace $U_0$. According
to the above theorem~4.7, the vectors $\bold e_1,\,\ldots,\,\bold e_s,\,
\bold e_{s+1}$ are linearly independent.\par
     Let's denote by $U_1$ the linear span of vectors $\bold e_1,\,\ldots,
\,\bold e_s,\,\bold e_{s+1}$. For the subspace $U_1$ we have the same
two mutually exclusive options $U_1=V$ or $U_1\neq V$, as we previously
had for the subspace $U_0$. If $U_1=V$, then the process of completing
the basis $\bold e_1,\,\ldots,\,\bold e_s$ is over. Otherwise, we can
iterate the process and get a chain of subspaces enclosed into each 
other:
$$
U_0\varsubsetneq U_1\varsubsetneq U_2\varsubsetneq\ldots\,.
$$
This chain of subspaces cannot be infinite since the dimension of every
next subspace is one as greater than the dimension of previous subspace,
and the dimensions of all subspaces are not greater than the dimension 
of $V$. The process of completing the basis will be finished in 
$(n-s)$-th step, where $U_{n-s}=V$.
\qed\enddemo
\head
\S\,5. Coordinates. Transformation of the coordinates of a 
vector under a change of basis.
\endhead
\rightheadtext{\S\,5. Transformation of the coordinates of
vectors\dots}
     Let $V$ be some finite-dimensional linear vector space over
the field $\Bbb K$ and let $\dim V=n$. In this section we shall
consider only finite-dimensional spaces. \pagebreak Let's choose 
a basis $\bold e_1,\,\ldots,\,\bold e_n$ in $V$. Then an arbitrary 
vector $x\in V$ can be expressed as linear combination of the 
basis vectors:
$$
\hskip -2em
\bold x=x^1\cdot\bold e_1+\ldots+x^n\cdot\bold e_n.
\tag5.1
$$
The linear combination \thetag{5.1} is called the {\it expansion} 
of the vector $\bold x$ in the basis $\bold e_1,\,\ldots,\,\bold e_n$.
Its coefficients $x^1,\,\ldots,\,x^n$ are the elements of the
numeric field $\Bbb K$. They are called the {\it components} or
the {\it coordinates} of the vector $\bold x$ in this basis.\par
     We use upper indices for the literal notations of the 
coordinates of a vector $\bold x$ in \thetag{5.1}. The usage of 
upper indices for the coordinates of vectors is determined by 
special convention, which is known as {\it tensorial notation}.
It was introduced for to simplify huge calculations in differential
geometry and in theory of relativity (see \cite{2} and \cite{3}).
Other rules of tensorial notation are discussed in coordinate
theory of tensors (see \cite{7}\footnotemark).\par
\footnotetext{\ The reference \cite{7} is added in 2004 to English 
translation of this book.}
\adjustfootnotemark{-1}
\proclaim{Theorem 5.1} For any vector $\bold x\in V$ its expansion
in a basis of a linear vector space $V$ is unique.
\endproclaim
\demo{Proof} The existence of an expansion \thetag{5.1} for a vector
$\bold x$ follows from the item~\therosteritem{2} of theorem~4.7.
Assume that there is another expansion
$$
\hskip -2em
\bold x=x'\vphantom{x}^1\cdot\bold e_1+\ldots+x'\vphantom{x}^n\cdot
\bold e_n.
\tag5.2
$$
Subtracting \thetag{5.1} from this equality, we get
$$
\hskip -2em
\bold x=(x'\vphantom{x}^1-x^1)\cdot\bold e_1+\ldots
+(x'\vphantom{x}^n-x^n)\cdot\bold e_n.
\tag5.3
$$
Since basis vectors $\bold e_1,\,\ldots,\,\bold e_n$ are linearly
independent, from the equality \thetag{5.3} it follows that the
linear combination \thetag{5.3} is trivial: $x'\vphantom{x}^i-x^i=0$. 
Then
$$
x'\vphantom{x}^1=x^1,\ \ldots,\ x'\vphantom{x}^n=x^n.
$$
Hence the expansions \thetag{5.1} and \thetag{5.2} do coincide. The
uniqueness of the expansion \thetag{5.1} is proved.
\qed\enddemo
     Having chosen some basis $\bold e_1,\,\ldots,\,\bold e_n$ in a
space $V$ and expanding a vector $\bold x$ in this base we can 
write its coordinates in the form of column vectors. Due to the
theorem~5.1 this determines a bijective map $\psi: V\to\Bbb K^n$.
It is easy to verify that 
$$
\xalignat 2
&\hskip -2em
\psi(x+y)=
 \Vmatrix x^1+y^1\\ \vdots\\ x^n+y^n\endVmatrix,
&&\psi(\alpha\cdot x)=
 \Vmatrix \alpha\cdot x^1\\ \vdots\\ \alpha\cdot x^n\endVmatrix.
\tag5.4
\endxalignat
$$
The above formulas \thetag{5.4} show that a basis is a very
convenient tool when dealing with vectors. In a basis algebraic
operations with vectors are replaced by algebraic operations 
with their coordinates, i\.\,e\. with numeric quantities.
However, coordinate approach has one disadvantage. The mapping
$\psi$ essentially depends on the basis we choose. And there is
no canonic choice of basis. In general, none of basis is preferable
with respect to another. Therefore we should be ready to consider
various bases and should be able to recalculate the coordinates of
vectors when passing from a basis to another basis.\par
     Let $\bold e_1,\,\ldots,\,\bold e_n$ and $\tilde\bold e_1,\,
\ldots,\,\tilde\bold e_n$ be two arbitrary bases in a linear vector
space $V$. We shall call them {\tencyr\char '074}wavy{\tencyr\char '076}
basis and {\tencyr\char '074}non-wavy{\tencyr\char '076} basis
(because of tilde sign we use for denoting the vectors of one of them).
The non-wavy basis will also be called the {\it initial basis} or the 
{\it old basis}, and the wavy one will be called the {\it new basis}.
Taking $i$-th vector of new (wavy) basis, we expand it in the old basis:
$$
\hskip -2em
\tilde\bold e_i=S^1_i\cdot\bold e_1+\ldots+S^n_i\cdot\bold e_n.
\tag5.5
$$
According to the tensorial notation, the coordinates of the vector
$\tilde\bold e_i$ in the expansion \thetag{5.5} are specified by 
upper index. The lower index $i$ specifies the number of the vector 
$\tilde\bold e_i$ being expanded. Totally in the expansion \thetag{5.5}
we determine $n^2$ numbers; they are usually arranged into a matrix:
$$
\hskip -2em
S=
\Vmatrix S^1_1 & \hdots & S^1_n \\
        \vdots & \ddots & \vdots\\
         S^n_1 & \hdots & S^n_n
\endVmatrix.
\tag5.6
$$
Upper index $j$ of the matrix element $S^j_i$ specifies the row number;
lower index $i$ specifies the column number. The matrix $S$ in \thetag{5.6}
the {\it direct transition matrix } for passing from the old basis
$\bold e_1,\,\ldots,\,\bold e_n$ to the new basis
$\tilde\bold e_1,\,\ldots,\,\tilde\bold e_n$.\par
     Swapping the bases $\bold e_1,\,\ldots,\,\bold e_n$ and $\tilde\bold
e_1,\,\ldots,\,\tilde\bold e_n$ we can write the expansion of the vector
$\bold e_j$ in wavy basis:
$$
\hskip -2em
\bold e_j=T^1_j\cdot\tilde\bold e_1+\ldots+T^n_j\cdot\tilde\bold e_n.
\tag5.7
$$
The coefficients of the expansion \thetag{5.7} determine the matrix 
$T$, which is called the {\it inverse transition matrix}. Certainly, 
the usage of terms {\tencyr\char '074}direct{\tencyr\char '076} and
{\tencyr\char '074}inverse{\tencyr\char '076} here is relative; it
depends on which basis is considered as an old basis and which one
is taken for a new one.\par
\proclaim{Theorem 5.2} The direct transition matrix $S$ and the 
inverse transition matrix $T$ determined by the expansions \thetag{5.5}
and \thetag{5.7} are inverse to each other.
\endproclaim
     Remember that two square matrices are inverse to each other if
their product is equal to unit matrix: $S\,T=1$. Here we do not define
the matrix multiplication assuming that it is known from the course of
general algebra.\par
\demo{Proof} Let's begin the proof of the theorem~5.2 by writing 
the relationships \thetag{5.5} and \thetag{5.7} in a brief symbolic
form:
$$
\xalignat 2
&\hskip -2em
\tilde\bold e_i=\sum^n_{k=1} S^k_i\cdot\bold e_k,
&&e_j=\sum^n_{i=1} T^i_j\cdot\tilde\bold e_i.
\tag5.8
\endxalignat
$$
Then we substitute the first relationship \thetag{5.8} into the second
one. This yields:
$$
\hskip -2em
\bold e_j=\sum^n_{i=1} T^i_j\cdot\!\left(\,\shave{\sum^n_{k=1}} S^k_i
\cdot\bold e_k\right)=\sum^n_{k=1}\left(\,\shave{\sum^n_{i=1}}
S^k_i\,T^i_j\right)\cdot\bold e_k.
\tag5.9
$$\par
    The symbol $\delta^k_j$, which is called the {\it Kronecker symbol},
is used for denoting the following numeric array:    
$$
\hskip -2em
\delta^k_j=
\cases
1 &\text{ for $k=j$,}\\
0 &\text{ for $k\neq j$.}
\endcases
\tag5.10
$$
We apply the Kronecker symbol determined in \thetag{5.10} in order
to transform left hand side of the equality \thetag{5.9}:
$$
\bold e_j=\sum^n_{k=1} \delta^k_j\cdot\bold e_k.
\tag5.11
$$
Both equalities \thetag{5.11} and \thetag{5.9} represent the same
vector $\bold e_j$ expanded in the same basis $\bold e_1,\,\ldots,
\,\bold e_n$. Due to the theorem~5.1 on the uniqueness of the
expansion of a vector in a basis we have the equality
$$
\sum^n_{i=1} S^k_i\,T^i_j=\delta^k_j.
$$
It is easy to note that this equality is equivalent to the matrix
equality $S\,T=1$. The theorem is proved.\qed\enddemo
\proclaim{Corollary} The direct transition matrix $S$ and the 
inverse transition matrix $T$ both are non-degenerate matrices
and \ $\det S\,\det T=1$.
\endproclaim
\demo{Proof} The relationship $\det S\,\det T=1$ follows from
the matrix equality $S\,T=1$, which was proved just above. This
fact is well known from the course of general algebra. If the 
product of two numbers is equal to unity, then none of these two
numbers can be equal to zero:
$$
\xalignat 2
&\det S\neq 0,
&&\det T\neq 0.
\endxalignat
$$
This proves the non-degeneracy of transition matrices $S$ and $T$. 
The corollary is proved.
\qed\enddemo
\proclaim{Theorem 5.3} Every non-degenerate $n\times n$ matrix $S$
can be obtained as a transition matrix for passing from some basis
$\bold e_1,\,\ldots,\,\bold e_n$ to some other basis $\tilde\bold
e_1,\,\ldots,\,\tilde\bold e_n$ in a linear vector space $V$ of the
dimension $n$.
\endproclaim
\demo{Proof} Let's choose an arbitrary $\bold e_1,\,\ldots,\,\bold 
e_n$ basis in $V$ and fix it. Then let's determine the other $n$ 
vectors $\tilde\bold e_1,\,\ldots,\,\tilde\bold e_n$ by means of 
the relationships \thetag{5.5} and prove that they are linearly
independent. For this purpose we consider a linear combination of
these vectors that is equal to zero:
$$
\hskip -2em
\alpha^1\cdot\,\tilde\bold e_1+\ldots+\alpha^n\cdot\,\tilde\bold e_n
=\bold 0.
\tag5.12
$$
Substituting \thetag{5.5} into this equality, one can transform it 
to the following one:
$$
\left(\,\shave{\sum^n_{i=1}} S^1_i\,\alpha^i\right)\cdot\bold e_1
+\ldots+
\left(\,\shave{\sum^n_{i=1}} S^n_i\,\alpha^i\right)\cdot\bold e_n=0.
$$
Since the basis vectors $\bold e_1,\,\ldots,\,\bold e_n$ are
linearly independent, it follows that all sums enclosed within 
the brackets in the above equality are equal to zero. Writing 
these sums in expanded form, we get a homogeneous system of
linear algebraic equations with respect to the variables $\alpha^1,
\,\ldots,\,\alpha^n$:
$$
\matrix
S^1_1\,\alpha^1 &+ &\hdots &+ &S^1_n\,\alpha^n &= &0,\\
\hdotsfor 7\\ \vspace{1\jot}
S^n_1\,\alpha^1 &+ &\hdots &+ &S^n_n\,\alpha^n &= &0.
\endmatrix
$$
The matrix of coefficients of this system coincides with $S$.
From the course of algebra we know that each homogeneous system
of linear equations with non-degenerate square matrix has unique
solution, which is purely zero: 
$$
\alpha^1=\ldots=\alpha^n=0.
$$
This means that an arbitrary linear combination \thetag{5.12},
which is equal to zero, is necessarily trivial. Hence, 
$\tilde\bold e_1,\,\ldots,\,\tilde\bold e_n$ is a linear 
independent system of vectors. Applying the proposition~\therosteritem{4}
from the theorem~4.5 to these vectors, we find that they form
a basis in $V$, while the matrix $S$ appears to be a direct transition 
matrix for passing from $\bold e_1,\,\ldots,\,\bold e_n$ to
$\tilde\bold e_1,\,\ldots,\,\tilde\bold e_n$. The theorem is proved.
\qed\enddemo
     Let's consider two bases $\bold e_1,\,\ldots,\,\bold e_n$ and 
$\tilde\bold e_1,\ldots,\tilde\bold e_n$ in a linear vector space 
$V$ related by the transition matrix $S$. Let $\bold x$ be some 
arbitrary vector of the space $V$. It can be expanded in each of 
these two bases:
$$
\xalignat 2
&\hskip -2em
\bold x=\sum^n_{k=1} x^k\cdot\,\bold e_k,
&&\bold x=\sum^n_{i=1}\tilde x^i\cdot\,\tilde\bold e_i.
\tag5.13
\endxalignat
$$
Once the coordinates of $\bold x$ in one of these two bases are
fixed, this fixes the vector $\bold x$ itself, and, hence, this
fixes its coordinates in another basis.
\proclaim{Theorem 5.4} The coordinates of a vector $\bold x$ in
two bases $\bold e_1,\,\ldots,\,\bold e_n$ and $\tilde\bold e_1,
\ldots,\tilde\bold e_n$ of a linear vector space $V$ are related
by formulas
$$
\xalignat 2
&\hskip -2em
x^k=\sum^n_{i=1} S^k_i\,\tilde x^i,
&&\tilde x^i=\sum^n_{k=1} T^i_k\,x^k,
\tag5.14
\endxalignat
$$
where $S$ and $T$ are direct and inverse transition matrices for the
passage from $\bold e_1,\,\ldots,\,\bold e_n$ to $\tilde\bold e_1,
\ldots,\tilde\bold e_n$, i\.\,e\. when $\bold e_1,\,\ldots,\,\bold 
e_n$ is treated as an old basis and $\tilde\bold e_1,\ldots,\tilde
\bold e_n$ is treated as a new one.
\endproclaim
The relationships \thetag{5.14} are known as {\it transformation 
formulas} for the coordinates of a vector under a change of basis.
\par
\demo{Proof} In order to prove the first relationship \thetag{5.14}
we substitute the expansion of the vector $\tilde\bold e_i$ taken
from \thetag{5.8} into the second relationship \thetag{5.13}: 
$$
\bold x=\sum^n_{i=1}\tilde x^i\cdot\left(\,\shave{\sum^n_{k=1}}
S^k_i\cdot\bold e_k\right)=\sum^n_{k=1}\left(\,\shave{\sum^n_{i=1}}
S^k_i\,\tilde x^i\right)\cdot\bold e_k.
$$
Comparing this expansion $\bold x$ with the first expansion 
\thetag{5.13} and applying the theorem on uniqueness of the 
expansion of a vector in a basis, we derive
$$
\pagebreak
x^k=\sum^n_{i=1} S^k_i\,\tilde x^i.
$$
This is exactly the first transformation formula \thetag{5.14}.
The second formula \thetag{5.14} is proved similarly. 
\qed\enddemo
\head
\S\,6. Intersections and sums of subspaces.
\endhead
     Suppose that we have a certain number of subspaces in a linear 
vector space $V$. In order to designate this fact we write $U_i\subset 
V$, where $i\in I$. The number of subspaces can be finite or infinite
enumerable, then they can be enumerated by the positive integers.
However, in general case we should enumerate the subspaces by the
elements of some indexing set $I$, which can be finite, infinite
enumerable, or even non-enumerable. Let's denote by $U$ and by $S$ the
intersection and the union of all subspaces that we consider:
$$
\xalignat 2
&\hskip -2em
U=\bigcap_{i\in I} U_i,
&&S=\bigcup_{i\in I} U_i.
\tag6.1
\endxalignat
$$
\proclaim{Theorem 6.1} The intersection of an arbitrary number of
subspaces in a linear vector space $V$ is a subspace in $V$.
\endproclaim
\demo{Proof} The set $U$ in \thetag{6.1} is not empty since zero
vector is an element of each subspace $U_i$. Let's verify the
conditions \therosteritem{1} and \therosteritem{2} from the
definition~2.2 for $U$.\par 
    Suppose that $\bold u_1$, $\bold u_2$, and $\bold u$ are the
vectors from the subset $U$. Then they belong to $U_i$ for each
$i\in I$. However, $U_i$ is a subspace, hence, $u_1+u_2\in U_i$
and $\alpha\cdot u\in U_i$ for any $i\in I$ and for any $\alpha\in
\Bbb K$. Therefore, $u_1+u_2\in U$ and $\alpha\cdot u\in U$.
The theorem is proved.
\qed\enddemo
     In general, the subset $S$ in \thetag{6.1} is not a subspace.
Therefore we need to introduce the following concept.
\definition{Definition 6.1} The linear span of the union of subspaces
$U_i$, $i\in I$, is called the {\it sum} of these subspaces.
\enddefinition
\noindent
For to denote the sum of subspaces $W=\langle S\rangle$ we use the
standard summation sign:
$$
W=\left<\bigcup U_i\right>=\sum_{i\in I} U_i.
$$
\proclaim{Theorem 6.2} A vector $\bold w$ of a linear vector space $V$
belongs to the sum of subspaces $U_i$, $i\in I$, if and  only if
it is represented as a sum of finite number of vectors each of 
which is taken from some subspace $U_i$:
$$
\hskip -2em
\bold w=\bold u_{i_{\ssize 1}}+\ldots+\bold u_{i_{\ssize k}}
\text{, where \ }\bold u_i\in U_i.
\tag6.2
$$
\endproclaim
\demo{Proof} Let $S$ be the union of subspaces $U_i\subset V$, $i\in I$.
Suppose that $\bold w\in W$. Then $\bold w$ is a linear combination of
finite number of vectors taken from $S$:
$$
\bold w=\alpha_1\cdot\bold s_1+\ldots+\alpha_k\cdot\bold s_k.
$$
But $S$ is the union of subspaces $U_i$. Therefore, $\bold s_m\in
U_{i_{\ssize m}}$ and $\alpha_m\cdot\,\bold s_m=\bold u_{i_{\ssize m}}
\in U_{i_{\ssize m}}$, where $m=1,\,\ldots,\,k$. This leads to the
equality \thetag{6.2} \pagebreak for the vector $\bold w$.\par
     Conversely, suppose that $\bold w$ is a vector given by formula
\thetag{6.2}. Then $u_{i_{\ssize m}}\in U_{i_{\ssize m}}$ and
$U_{i_{\ssize m}}\subset S$, i\.\,e\. $u_{i_{\ssize m}}\in S$.
Therefore, the vector $\bold w$ belongs to the linear span of $S$.
The theorem is proved.
\qed\enddemo
\definition{Definition 6.2} The sum $W$ of subspaces $U_i$, $i\in I$,
is called the {\it direct sum}, if for any vector $\bold w\in W$ the
expansion \thetag{6.2} is unique. In this case for the direct sum
of subspaces we use the special notation:
$$
W=\bigoplus_{i\in I} U_i.
$$
\enddefinition
\proclaim{Theorem 6.3} Let $W=U_1+\,\ldots\,+U_k$ be the sum a finite
number of finite-dimensional subspaces. The dimension of $W$ is equal
to the sum of dimensions of the subspaces $U_i$ if and only if $W$ is
the direct sum: $W=U_1\oplus\ldots\oplus U_k$.
\endproclaim
\demo{Proof} Let's choose a basis in each subspace $U_i$. Suppose 
that $\dim U_i=s_i$ and let $\bold e_{i\,1},\,\ldots,\,\bold e_{i\,
s_{\ssize i}}$ be a basis in $U_i$. Let's join the vectors of all bases
into one system ordering them alphabetically:
$$
\hskip -2em
\bold e_{1\,1},\,\ldots,\,\bold e_{1\,s_{\ssize 1}},\ \ \ldots,\ \
\bold e_{k\,1},\,\ldots,\,\bold e_{k\,s_{\ssize k}}.
\tag6.3
$$
Due to the equality $W=U_1+\ldots+U_k$ for an arbitrary vector $\bold w$
of the subspace $W$ we have the expansion \thetag{6.2}:
$$
\bold w=\bold u_1+\ldots+\bold u_k\text{, where \ }\bold u_i\in U_i.
\tag6.4
$$
Expanding each vector $\bold u_i$ of \thetag{6.4} in the basis of
corresponding subspace $U_i$, we get the expansion of $\bold w$ in
vectors of the system \thetag{6.3}. Hence, \thetag{6.3} is a spanning
system of vectors in $W$ (though, in general case it is not a minimal
spanning system).\par
     If $\dim W=\dim U_1+\ldots+\dim U_k$, then the number of vectors
in \thetag{6.3} cannot be reduced. Therefore \thetag{6.3} is a basis
in $W$. From any expansion \thetag{6.4} we can derive the following
expansion of the vector $\bold w$ in the basis \thetag{6.3}:
$$
\hskip -2em
\bold w=\left(\,\shave{\sum^{s_{\ssize 1}}_{j=1}}
\alpha_{1\,j}\cdot\bold e_{1\,j}\right)+\ldots+
\left(\,\shave{\sum^{s_{\ssize k}}_{j=1}}
\alpha_{k\,j}\cdot\bold e_{k\,j}\right).
\tag6.5
$$
The sums enclosed into the round brackets in \thetag{6.5} are determined 
by the expansions of the vectors $\bold u_1,\,\ldots,\,\bold u_k$ in
the bases of corresponding subspaces $U_1,\ldots,U_k$: 
$$
u_i=\sum^{s_{\ssize i}}_{j=1}\alpha_{i\,j}\cdot e_{i\,j}.
\tag6.6
$$
Due to \thetag{6.6} the existence of two different expansions
\thetag{6.4} for some vector $\bold w$ would mean the existence
of two different expansions \thetag{6.5} of this vector in the
basis \thetag{6.3}. Hence, the expansion \thetag{6.4} is unique
and the sum of subspaces $W=U_1+\ldots+U_k$ is the direct sum.
\par
     Conversely, suppose that $W=U_1\oplus\ldots\oplus U_k$.
We know that the vectors \thetag{6.3} span the subspace $W$.
Let's prove that they are linearly independent. For this 
purpose we consider a linear combination of these vectors
being equal to zero:
$$
\bold 0=\left(\,\shave{\sum^{s_{\ssize 1}}_{j=1}}
\alpha_{1\,j}\cdot\bold e_{1\,j}\right)+\ldots+
\left(\,\shave{\sum^{s_{\ssize k}}_{j=1}}
\alpha_{k\,j}\cdot\bold e_{k\,j}\right).
\tag6.7
$$
Let's denote by $\tilde\bold u_1,\,\ldots,\,\tilde\bold u_k$ the
values of sums enclosed into the round brackets in \thetag{6.7}.
It is easy to see that $\tilde\bold u_i\in U_i$, therefore, 
\thetag{6.7} is an expansion of the form \thetag{6.4} for the
vector $\bold w=\bold 0$. But $\bold 0=\bold 0+\ldots+\bold 0$
and $\bold 0\in U_i$. This is another expansion for the vector
$\bold w=\bold 0$. However, $W=U_1\oplus\ldots\oplus U_k$, therefore,
the expansion $\bold 0=\bold 0+\ldots+\bold 0$ is unique expansion
of the form \thetag{6.4} for zero vector $\bold w=\bold 0$. Then
we have the equalities 
$$
\bold 0=\sum^{s_{\ssize i}}_{j=1}\alpha_{i\,j}\cdot e_{i\,j}
\text{\ for all \ }i=1,\ldots,k.
$$
It's clear that these equalities are the expansions of zero vector
in the bases of the subspaces $U_i$. Hence, $\alpha_{i\,j}=0$. This
means that the linear combination \thetag{6.7} is trivial, and
\thetag{6.3} is a linearly independent system of vectors. Thus, being
a spanning system and being linearly independent, the system of
vectors \thetag{6.3} is a basis of $W$. Now we can find the dimension 
of the subspace $W$ by counting the number of vectors in \thetag{6.3}:
$\dim W=s_1+\ldots+s_k=\dim U_1+\ldots+\dim U_k$. The theorem is 
proved.
\qed\enddemo
     {\bf Note.} If the sum of subspaces $W=U_1+\ldots+U_k$ is not 
necessarily the direct sum, the vectors \thetag{6.3}, nevertheless, 
form a spanning system in $W$. But they do not necessarily form a
linearly independent system in this case. Therefore, we have
$$
\hskip -2em
\dim W\leqslant\dim U_1+\ldots+\dim U_k.
\tag6.8
$$
Sharpening this inequality in general case is sufficiently complicated.
We shall do it for the case of two subspaces.\par
\proclaim{Theorem 6.4} The dimension of the sum of two arbitrary 
finite-dimensional subspaces $U_1$ and $U_2$ in a linear vector
space $V$ is equal to the sum of their dimensions minus the dimension
of their intersection:
$$
\hskip -2em
\dim(U_1+U_2)=\dim U_1+\dim U_2-\dim(U_1\cap U_2).
\tag6.9
$$
\endproclaim
\demo{Proof} From the inclusion $U_1\cap U_2\subset U_1$ and from
the inequality \thetag{6.8} we conclude that all subspaces 
considered in the theorem are finite-dimensional. Let's denote
$\dim(U_1\cap U_2)=s$ and choose a basis $\bold e_1,\,\ldots,\,
\bold e_s$ in the intersection $U_1\cap U_2$.\par
     Due to the inclusion $U_1\cap U_2\subset U_1$ we can apply the
theorem~4.8 on completing the basis. This theorem says that we can
complete the basis $\bold e_1,\,\ldots,\,\bold e_s$ of the intersection
$U_1\cap U_2$ up to a basis $\bold e_1,\,\ldots,\,\bold e_s,\,
\bold e_{s+1},\,\ldots,\,\bold e_{s+p}$ in $U_1$. For the dimension
of $U_1$, we have $\dim U_1=s+p$. In a similar way, due to the inclusion
$U_1\cap U_2\subset U_2$ we can construct a basis $\bold e_1,\,\ldots,
\,\bold e_s,\,\bold e_{s+p+1},\,\ldots,\,\bold e_{s+p+q}$ in $U_2$.
For the dimension of $U_2$ this yields $\dim U_2=s+q$.\par
     Now let's join together the two bases constructed above with the
use of theorem~4.8 and consider the total set of vectors in them:
$$
\hskip -2em
\bold e_1,\,\ldots,\,\bold e_s,\,\bold e_{s+1},\,\ldots,\,
\bold e_{s+p},\,\bold e_{s+p+1},\,\ldots,\,\bold e_{s+p+q}.
\tag6.10
$$
Let's prove that these vectors \thetag{6.10} form a basis in the
sum of subspaces $U_1+U_2$. Let $\bold w$ be some arbitrary vector
in $U_1+U_2$. The relationship \thetag{6.2} for this vector is
written as $\bold w=\bold u_1+\bold u_2$. Let's expand the vectors
$\bold u_1$ and $\bold u_2$ in the above two bases of the subspaces 
$U_1$ and $U_2$ respectively:
$$
\align
&\bold u_1=\sum^s_{i=1}\alpha_i\cdot\bold  e_i+
     \sum^p_{j=1}\beta_{s+j}\cdot\bold e_{s+j},\\
&\bold u_2=\sum^s_{i=1}\tilde\alpha_i\cdot\bold e_i+
     \sum^q_{j=1}\gamma_{s+p+j}\cdot\bold e_{s+p+j}.
\endalign
$$
Adding these two equalities, we find that the vector $\bold w$ is
linearly expressed through the vectors \thetag{6.10}. Hence, 
\thetag{6.10} is a spanning system of vectors in $U_1+U_2$.\par
     In order to prove that \thetag{6.10} is a linearly independent
system of vectors we consider a linear combination of these vectors
being equal to zero:
$$
\hskip -2em
\sum^{s+p}_{i=1}\alpha_i\cdot\bold e_i+
\sum^{q}_{i=1}\alpha_{s+p+i}\cdot\bold e_{s+p+i}=\bold 0.
\tag6.11
$$
Then we transform this equality by moving the second sum to the
right hand side:
$$
\sum^{s+p}_{i=1}\alpha_i\cdot\bold e_i=-\sum^{q}_{i=1}\alpha_{s+p+i}
\cdot\bold e_{s+p+i}.
$$
Let's denote by $\bold u$ the value of left and right sides of this
equality. Then for the vector $\bold u$ we get the following two
expressions:
$$
\xalignat 2
&\hskip -2em
\bold u=\sum^{s+p}_{i=1}\alpha_i\cdot\bold e_i,
&&\bold u=-\sum^{q}_{i=1}\alpha_{s+p+i}\cdot\bold e_{s+p+i}.
\tag6.12
\endxalignat
$$
Because of the first expression \thetag{6.12} we have $\bold u\in U_1$,
while the second expression \thetag{6.12} yields $\bold u\in U_2$. Hence, $\bold u\in U_1\cap U_2$. This means that we can expand the vector $\bold
u$ in the basis $\bold e_1,\,\ldots,\,\bold e_s$:
$$
\hskip -2em
\bold u=\sum^s_{i=1}\beta_i\cdot\bold e_i.
\tag6.13
$$
Comparing this expansion with the second expression \thetag{6.12}, we
find that 
$$
\sum^s_{i=1}\beta_i\cdot\bold e_i+
\sum^{q}_{i=1}\alpha_{s+p+i}\cdot\bold e_{s+p+i}=\bold 0.
\tag6.14
$$
Note that the vectors $\bold e_1,\,\ldots,\,\bold e_s,\,\bold e_{s+p+1},
\,\ldots,\,\bold e_{s+p+q}$ form a basis in $U_2$. They are linearly
independent. Therefore, all coefficients in \thetag{6.14} are equal to
zero. In particular, we have the following equalities:
$$
\hskip -2em
\alpha_{s+p+1}=\ldots=\alpha_{s+p+q}=0.
\tag6.15
$$
Moreover, $\beta_1=\ldots=\beta_s=0$. Due to \thetag{6.13} this means that
$\bold u=0$. Now from the first expansion \thetag{6.12} we get the equality
$$
\sum^{s+p}_{i=1}\alpha_i\cdot\bold e_i=0.
$$
Since $\bold e_1,\,\ldots,\,\bold e_s,\,\bold e_{s+1},\,\ldots,\,
\bold e_{s+p}$ are linearly independent vectors, all coefficients
$\alpha_i$ in the above equality should be zero:
$$
\hskip -2em
\alpha_1=\ldots=\alpha_s=\alpha_{s+1}=\ldots=\alpha_{s+p}=0.
\tag6.16
$$
Combining \thetag{6.15} and \thetag{6.16}, we see that the linear
combination \thetag{6.11} is trivial. This means that the vectors
\thetag{6.10} are linearly independent. Hence, they form a basis in
$U_1+U_2$. For the dimension of the subspace $U_1+U_2$ this yields
$$
\align
\dim(U_1+U_2)&=s+p+q=(s+p)+(s+q)-s=\\
&=\dim U_1+\dim U_2-\dim(U_1\cap U_2).
\endalign
$$
Thus, the relationship \thetag{6.9} and the theorem~6.4 in whole
is proved.
\qed\enddemo 
\head
\S~7. Cosets of a subspace. The concept of factorspace.
\endhead
     Let $V$ be a linear vector space and let $U$ be a subspace in it.
A {\it coset} of the subspace $\bold U$ determined by a vector 
$\bold v\in V$ is the following set of vectors\footnote{\ We used the 
sign $\Cl$ for cosets since in Russia they are called {\it adjacency
classes.}}:\adjustfootnotemark{-1}
$$
\hskip -2em
\Cl_U(\bold v)=\{\bold w\in V\!:\ \bold w-\bold v\in U\}.
\tag7.1
$$
The vector $\bold v$ in \thetag{7.1} is called a {\it representative}
of the coset \thetag{7.1}. The coset $\Cl_U(\bold v)$ is a very simple
thing, it is obtained by adding the vector $\bold v$ with all vectors 
of the subspace $U$. The coset represented by zero vector is the 
especially simple thing since $\Cl_U(\bold 0)=U$. It is called a 
{\it zero coset}.
\proclaim{Theorem 7.1} The cosets of a subspace $U$ in a linear
vector space $V$ possess the following properties:
\roster
\item $\bold a\in\Cl_U(\bold a)$ for any $\bold a\in V$;
\item if $\bold a\in\Cl_U(\bold b)$, then $\bold b\in\Cl_U(\bold a)$;
\item if $\bold a\in\Cl_U(\bold b)$ and $\bold b\in\Cl_U(\bold c)$, 
then $\bold a\in\Cl_U(\bold c)$.
\endroster
\endproclaim
\demo{Proof} The first proposition is obvious. Indeed, the difference
$\bold a-\bold a$ is equal to zero vector, which is an element of any subspace: $\bold a-\bold a=\bold 0\in U$. Hence, due to the formula
\thetag{7.1}, which is the formal definition of cosets, we have
$\bold a\in\Cl_U(\bold a)$.\par
     Let $\bold a\in\Cl_U(\bold b)$. Then $\bold a-\bold b\in U$.
For $\bold b-\bold a$, we have $\bold b-\bold a=(-1)\cdot(\bold a
-\bold b)$. Therefore, $\bold b-\bold a\in U$ and $\bold b\in
\Cl_U(\bold a)$ (see formula \thetag{7.1} and the definition~2.2).
The second proposition is proved.\par
     Let $\bold a\in\Cl_U(\bold b)$ and $\bold b\in\Cl_U(\bold c)$. 
Then $\bold a-\bold b\in U$ and $\bold b-\bold c\in U$. Note that
$\bold a-\bold c=(\bold a-\bold b)+(\bold b-\bold a)$. Hence,
$\bold a-\bold c\in U$ and $a\in\Cl_U(c)$ (see formula \thetag{7.1} 
and the definition~2.2 again). The third proposition is proved. 
This completes the proof of the theorem in whole.
\qed\enddemo
     Let $\bold a\in\Cl_U(\bold b)$. This condition establishes some
kind of dependence between two vectors $\bold a$ and $\bold b$. This
dependence is not strict: the condition $\bold a\in\Cl_U(\bold b)$ 
does not exclude the possibility that $\bold a'\in\Cl_U(\bold b)$ for
some other vector $\bold a'$. Such non-strict dependences in mathematics
are described by the concept of {\it binary relation} (see details in
\cite{1} and \cite{4}). Let's write $\bold a\sim\bold b$ as an 
abbreviation for $\bold a\in\Cl_U(\bold b)$. Then the theorem~7.1
reveals the following properties of the binary relation $\bold a\sim
\bold b$, which is introduced just above:
\roster
\item {\it reflexivity}: $\bold a\sim\bold a$;
\item {\it symmetry}: $\bold a\sim\bold b$ implies $\bold b
\sim\bold a$;
\item {\it transitivity}: $\bold a\sim\bold b$ and $\bold b\sim
\bold c$ implies $\bold a\sim\bold c$.
\endroster
A binary relation possessing the properties of reflexivity, symmetry,
and transiti\-vity is called an {\it equivalence relation}. Each
equivalence relation determined in a set $V$ partitions this set into 
a union of mutually non-intersecting subsets, which are called the 
{\it equivalence classes}:
$$
\hskip -2em
\Cl(\bold v)=\{\bold w\in V\!:\ \bold w\sim \bold v\}.
\tag7.2
$$
In our particular case the formal definition \thetag{7.2} coincides
with the formal definition \thetag{7.1}. In order to keep the 
completeness of presentation we shall not use the notation $\bold
a\sim\bold b$ in place of $\bold a\in\Cl_U(\bold b)$ anymore, and 
we shall not refer to the theory of binary relations (though it is
simple and well-known). Instead of this we shall derive the result 
on partitioning $V$ into the mutually non-intersecting cosets from
the following theorem.
\proclaim{Theorem 7.2} If two cosets $\Cl_U(\bold a)$ and $\Cl_U(
\bold b)$ of a subspace $U\subset V$ are intersecting, then they 
do coincide.
\endproclaim
\demo{Proof} Assume that the intersection of two cosets $\Cl_U(
\bold a)$ and $\Cl_U(\bold b)$ is not empty. Then there is an
element $\bold c$ belonging to both of them: $\bold c\in \Cl_U(
\bold a)$ and $\bold c\in \Cl_U(\bold b)$. Due to the proposition
\therosteritem{2} of the above theorem~7.1 we derive
$\bold b\in \Cl_U(\bold c)$. Combining $\bold b\in \Cl_U(\bold c)$
and $\bold c\in \Cl_U(\bold a)$ and applying the proposition
\therosteritem{3} of the theorem~7.1, we get $\bold b\in \Cl_U(\bold 
a)$. The opposite inclusion $\bold a\in \Cl_U(\bold b)$ then is 
obtained by applying the proposition~\therosteritem{2} of the
theorem~7.1.\par
     Let's prove that two cosets $\Cl_U(\bold a)$ and $\Cl_U(\bold b)$
do coincide. For this purpose let's consider an arbitrary vector
$\bold x\in\Cl_U(\bold a)$. From $\bold x\in\Cl_U(\bold a)$ and
$\bold a\in\Cl_U(\bold b)$ we derive $\bold x\in\Cl_U(\bold b)$.
Hence, $\Cl_U(\bold a)\subset\Cl_U(\bold b)$. The opposite inclusion
$\Cl_U(\bold b)\subset\Cl_U(\bold a)$ is proved similarly. From these
two inclusions we derive $\Cl_U(\bold a)=\Cl_U(\bold b)$. The theorem
is proved.
\qed\enddemo
     The set of all cosets of a subspace $U$ in a linear vector space
$V$ is called the {\it factorset} or  {\it quotient set} $V/U$. Due to
the theorem proved just \pagebreak above any two different cosets $Q_1$
and $Q_2$ from the factorset $V/U$ have the empty intersection
$Q_1\cap Q_2=\varnothing$, while the union of all cosets coincides with 
$V$:
$$
V=\bigcup_{Q\in V/U}Q.
$$
This equality is a consequence of the fact that any vector $\bold v\in V$
is an element of some coset: $\bold v\in Q$. This coset $Q$ is determined 
by $\bold v$ according to the formula $Q=\Cl_U(\bold v)$. For this reason
the following theorem is a simple reformulation of the definition of
cosets.\par
\proclaim{Theorem 7.3} Two vectors $\bold v$ and $\bold w$ belong 
to the same coset of a subspace $U$ if and only if their difference 
$\bold v-\bold w$ is a vector of $U$.
\endproclaim
\definition{Definition 7.1} Let $Q_1$ and $Q_2$ be two cosets of
a subspace $U$. The {\it sum} of cosets $Q_1$ and $Q_2$ is a coset
$Q$ of the subspace $U$ determined by the equality $Q=\Cl_U(\bold v_1
+\bold v_2)$, where $\bold v_1\in Q_1$ and $\bold v_2\in Q_2$.
\enddefinition
\definition{Definition 7.2} Let $Q$ be a coset of a subspace $U$. The product of $Q$ and a number $\alpha\in\Bbb K$ is a coset $P$ of the 
subspace $U$ determined by the relationship $P=\Cl_U(\alpha\cdot\bold
v)$, where $\bold v\in Q$.
\enddefinition
     For the addition of cosets and for the multiplication of them by
numbers we use the same signs of algebraic operations as in case of vectors,
i\.\,e\. $Q=Q_1+Q_2$ and $P=\alpha\cdot Q$. The definitions~7.1 and 7.2
can be expressed by formulas
$$
\aligned
&\Cl_U(\bold v_1)+\Cl_U(\bold v_2)=\Cl_U(\bold v_1+\bold v_2),\\
&\alpha\cdot\Cl_U(\bold v)=\Cl_U(\alpha\cdot\bold v).
\endaligned
\tag7.3
$$
These definitions require some comments. Indeed, the coset $Q=Q_1+Q_2$
in the definition~7.1 and the coset $P=\alpha\cdot Q$ in the definition~7.2
both are determined using some representative vectors $\bold v_1\in Q_1$,
$\bold v_2\in Q_2$, and $\bold v\in Q$. The choice of a representative
vector in a coset is not unique; therefore, we need especially to prove the
uniqueness of the results of algebraic operations determined in the
definitions~7.1 and 7.2. This proof is called the {\it proof of
correctness}.\par
\proclaim{Theorem 7.4} The definitions~7.1 and 7.2 are correct and the
results of the algebraic operations of coset addition and of coset
multiplication by numbers do not depend on the choice of representatives
in cosets.
\endproclaim
\demo{Proof} For the beginning we study the operation of coset addition.
Lat's take consider two different choices of representatives within cosets
$Q_1$ and $Q_2$. Let $\bold v_1,\tilde\bold v_1$ be two vectors of $Q_1$
and let $\bold v_1,\tilde\bold v_1$ be two vectors of $Q_2$. Then due to 
the theorem~7.3 we have the following two equalities:
$$
\xalignat 2
&\tilde\bold v_1-\bold v_1\in U,
&&\tilde\bold v_2-\bold v_2\in U.
\endxalignat
$$
Hence, $(\tilde\bold v_1+\tilde\bold v_2)-(\bold v_1+\bold v_2)=(\tilde
\bold v_1-\bold v_1)+(\tilde\bold v_2-\bold v_2)\in U$. This means that
the cosets determined by vectors $\tilde\bold v_1+\tilde\bold v_2$ and
$\bold v_1+\bold v_2$ do coincide with each other:
$$
\pagebreak
\Cl_U(\tilde v_1+\tilde v_2)=\Cl_U(v_1+v_2).
$$
This proves the correctness of the definition~7.1 for the operation of
coset addition.\par
     Now let's consider two different representatives $\bold v$ and
$\tilde\bold v$ within the coset $Q$. Then $\tilde\bold v-\bold v\in 
U$. Hence, $\alpha\cdot\tilde\bold v-\alpha\cdot\bold v=\alpha\cdot(
\tilde\bold v-\bold v)\in U$. This yields
$$
\Cl_U(\alpha\cdot\tilde\bold v)=\Cl_U(\alpha\cdot\bold v),
$$
which proves the correctness of the definition~7.2 for the operation of
multiplication of cosets by numbers.
\qed\enddemo
\proclaim{Theorem 7.5} The factorset $V/U$ of a linear vector space $V$
over a subspace $U$ equipped with algebraic operations \thetag{7.3} is
a linear vector space. This space is called the {\it factorspace} or the
{\it quotient space} of the space $V$ over its subspace $U$.
\endproclaim
\demo{Proof} The proof of this theorem consists in verifying the axioms
\therosteritem{1}-\therosteritem{8} of a linear vector space for $V/U$.
The commutativity and associativity axioms for the operation of coset addition follow from the following calculations:
$$
\align
&\aligned
\Cl_U(\bold v_1)+\Cl_U(\bold v_2)&=\Cl_U(\bold v_1+\bold v_2)=\\
  &=\Cl_U(\bold v_2+\bold v_1)=\Cl_U(\bold v_2)+\Cl_U(\bold v_1),
 \endaligned\\
\vspace{1ex}
&\aligned
  (\Cl_U(\bold v_1)+\Cl_U(\bold v_2))+\Cl_U(\bold v_3)&=
  \Cl_U(\bold v_1+\bold v_2)+\Cl_U(\bold v_3)=\\
  =\Cl_U((\bold v_1+\bold v_2)+v_3)&=\Cl_U(\bold v_1+(\bold v_2
  +\bold v_3))=\\
  \Cl_U(\bold v_1)+\Cl_U(\bold v_2+\bold v_3)&=\Cl_U(\bold v_1)+
  (\Cl_U(\bold v_2)+\Cl_U(\bold v_3)).
  \endaligned
\endalign
$$
In essential, they follow from the corresponding axioms for the
operation of vector addition (see definition~2.1).\par
     In order to verify the axiom \therosteritem{3} we should 
have a zero element in $V/U$. The zero coset $\bold 0=\Cl_U(\bold 0)$
is the best pretender for this role:
$$
\Cl_U(\bold v)+\Cl_U(\bold 0)=\Cl_U(\bold v+\bold 0)=\Cl_U(\bold v).
$$\par
     In verifying the axiom \therosteritem{4} we should indicate
the opposite coset $Q'$ for a coset $Q=\Cl_U(\bold v)$. We define
it as follows: $Q'=\Cl_U(\bold v')$. Then
$$
Q+Q'=\Cl_U(\bold v)+\Cl_U(\bold v')=\Cl_U(\bold v+\bold v')
=\Cl_U(\bold 0)=\bold 0.
$$\par
     The rest axioms \therosteritem{5}-\therosteritem{8} are verified 
by direct calculations on the base of formula \thetag{7.3} for coset operations. Here are these calculations:
$$
\align
 &\aligned
   &\alpha\cdot(\Cl_U(\bold v_1)+\Cl_U(\bold v_2))=
    \alpha\cdot\Cl_U(\bold v_1+\bold v_2)=\\
   &=\Cl_U(\alpha\cdot(\bold v_1+\bold v_2))=
     \Cl_U(\alpha\cdot\bold v_1+\alpha\cdot\bold v_2)=\\
   &=\Cl_U(\alpha\cdot\bold v_1)+\Cl_U(\alpha\cdot\bold v_2)
    =\alpha\cdot\Cl_U(\bold v_1)+\alpha\cdot\Cl_U(\bold v_2),
  \endaligned\\
\vspace{1.7ex}
 &\aligned
  &(\alpha+\beta)\cdot\Cl_U(\bold v)=\Cl_U((\alpha+\beta)
    \cdot\bold v)=\Cl_U(\alpha\cdot\bold v+\beta\cdot\bold v)=\\
  &=\Cl_U(\alpha\cdot\bold v)+\Cl_U(\beta\cdot\bold v)
   =\alpha\cdot\Cl_U(\bold v)+\beta\cdot\Cl_U(\bold v),
  \endaligned\\
\vspace{1.7ex}
 &\aligned
  &\alpha\cdot(\beta\cdot\Cl_U(\bold v))=\alpha\cdot\Cl_U(\beta
  \cdot\bold v)=
  \Cl_U(\alpha\cdot(\beta\cdot\bold v))=\\
  &=\Cl_U((\alpha\beta)\cdot\bold v)=(\alpha\beta)\cdot\Cl_U(\bold v),
  \endaligned\\
\vspace{1.7ex}
 &1\cdot\Cl_U(\bold v)=\Cl_U(1\cdot\bold v)=\Cl_U(\bold v).
\endalign
$$
The above equalities complete the verification of the fact that
the factorset $V/U$ possesses the structure of a linear vector space.
\qed\enddemo
     Note that verifying the axiom \therosteritem{4} we have defined
the opposite coset $Q'$ for a coset $Q=\Cl_U(\bold v)$ by means of
the relationship $Q'=\Cl_U(\bold v')$, where $\bold v'$ is the opposite
vector for $\bold v$. One could check the correctness of this definition. 
However, this is not necessary since due to the property~\therosteritem{10},
see theorem~2.1, the opposite coset $Q'$ for $Q$ is unique.\par
     The concept of factorspace is equally applicable to 
finite-dimensional and to infinite-dimensional spaces $V$. The finite
or infinite dimensionality of a subspace $U$ also makes no difference.
The only simplification in finite-dimensional case is that we can
calculate the dimension of the factorspace $V/U$.\par
\proclaim{Theorem 7.6} If a linear vector space $V$ is finite-dimensional,
then for any its subspace $U$ the factorspace $V/U$ also is 
finite-dimensional and its dimension is determined by the following
formula:
$$
\dim U+\dim(V/U)=\dim V.
\tag7.4
$$
\endproclaim
\demo{Proof} If $U=V$ then the factorspace $V/U$ consists of zero
coset only: $V/U=\{\bold 0\}$. The dimension of such zero space is
equal to zero. Hence, the equality \thetag{7.4} in this trivial case
is fulfilled.\par
     Let's consider a nontrivial case $U\varsubsetneq V$. Due to
the theorem~4.5 the subspace $U$ is finite-dimensional. Denote
$\dim V=n$ and $\dim U=s$, then $s<n$. Let's choose a basis 
$\bold e_1,\,\ldots,\,\bold e_s$ in $U$ and, according to the 
theorem~4.8, complete it with vectors $\bold e_{s+1},\,\ldots,\,
\bold e_n$ up to a basis in $V$. For each of complementary vectors
$\bold e_{s+1},\,\ldots,\,\bold e_n$ we consider the corresponding 
coset of a subspace $U$:
$$
\hskip -2em
\bold E_1=\Cl_U(\bold e_{s+1}),\ \ldots,\ \bold E_{n-s}
=\Cl_U(\bold e_n).
\tag7.5
$$\par
    Now let's show that the cosets \thetag{7.5} span the factorspace
$V/U$. Indeed, let $Q$ be an arbitrary coset in $V/U$ and let $\bold 
v\in Q$ be some representative vector of this coset. Let's expand the 
vector $\bold v$ in the above basis of $V$:
$$
\bold v=(\alpha_1\cdot\bold e_1+\ldots+\alpha_s\cdot\bold e_s)+
\beta_1\cdot\bold e_{s+1}+\ldots+\beta_{n-s}\cdot\bold e_n.
$$
Let's denote by $\bold u$ the initial part of this expansion:
$\bold u=\alpha_1\cdot\bold e_1+\ldots+\alpha_s\cdot\bold e_s$. 
It is clear that $\bold u\in U$. Then we can write
$$
\bold v=\bold u+\beta_1\cdot\bold e_{s+1}+\ldots+\beta_{n-s}
\cdot\bold e_n.
$$
Since $\bold u\in U$, we have $Cl_U(\bold u)=\bold 0$. For the 
coset $Q=\Cl_U(\bold v)$ this equality yields $Q=\beta_1\cdot
\Cl_U(\bold e_{s+1})+\ldots+\beta_{n-s}\cdot\Cl_U(\bold e_n)$. 
Hence, we have
$$
Q=\beta_1\cdot\bold E_1+\ldots+\beta_{n-s}\cdot\bold E_{n-s}.
$$
This means that $\bold E_1,\,\ldots,\,\bold E_{n-s}$ is a 
finite spanning system in $V/U$. Therefore, $V/U$ is a 
finite-dimensional linear vector space. \pagebreak To determine its
dimension we shall prove that the cosets \thetag{7.5} are linearly
independent. Indeed, let's consider a linear combination of these
cosets being equal to zero:
$$
\hskip -2em
\gamma_1\cdot\bold E_1+\ldots+\gamma_{n-s}\cdot\bold E_{n-s}
=\bold 0.
\tag7.6
$$
Passing from cosets to their representative vectors, from 
\thetag{7.6} we derive
$$
\align
\gamma_1&\cdot\Cl_U(\bold e_{s+1})+\ldots+\gamma_{n-s}\cdot
 \Cl_U(\bold e_n)=\\
&=\Cl_U(\gamma_1\cdot\bold e_{s+1}+\ldots+\gamma_{n-s}\cdot\bold e_n)=
\Cl_U(\bold 0).
\endalign
$$
Let's denote $\bold u=\gamma_1\cdot\bold e_{s+1}+\ldots+\gamma_{n-s}
\cdot\bold e_n$. From the above equality for this vector we get
$\Cl_U(\bold u)=\Cl_U(\bold 0)$, which means $\bold u\in U$. Let's
expand $\bold u$ in the basis of subspace $U$: $\bold u=\alpha_1\cdot
\bold e_1+\ldots+\alpha_s\cdot\bold e_s$. Then, equating two expression
for the vector $\bold u$, we get the following equality:
$$
-\alpha_1\cdot\bold e_1-\ldots-\alpha_s\cdot\bold e_s+
\gamma_1\cdot\bold e_{s+1}+\ldots+\gamma_{n-s}\cdot\bold e_n=\bold 0.
$$
This is the linear combination of basis vectors of $V$, which is equal
to zero. Basis vectors $\bold e_1,\,\ldots,\,\bold e_n$ are linearly
independent. Hence, this linear combination is trivial and
$\gamma_1=\ldots=\gamma_{n-s}=0$. This proves the triviality of
linear combination \thetag{7.6} and, therefore, the linear independence
of cosets \thetag{7.5}. Thus, for the dimension of factorspace this
yields $\dim(V/U)=n-s$, which proves the equality \thetag{7.4}. The
theorem is proved.
\qed\enddemo
\head
\S\,8. Linear mappings.
\endhead
\definition{Definition 8.1} Let $V$ and $W$ be two linear vector spaces
over a numeric field $\Bbb K$. A mapping $f\!:\ V\to W$ from the space
$V$ to the space $W$ is called a {\it linear mapping} if the following
two conditions are fulfilled:
\roster
\item $f(\bold v_1+\bold v_2)=f(\bold v_1)+f(\bold v_2)$ for any two
      vectors $\bold v_1,\bold v_2\in V$;
\item $f(\alpha\cdot\bold v)=\alpha\cdot f(\bold v)$ for any vector 
      $\bold v\in V$ and for any number $\alpha\in\Bbb K$.
\endroster
\enddefinition
     The relationship $f(\bold 0)=\bold 0$ is one of the simplest and
immediate consequences of the above two properties  \therosteritem{1}
and \therosteritem{2} of linear mappings. Indeed, we have
$$
\hskip -2em
f(\bold 0)=f(\bold 0+(-1)\cdot\bold 0)=f(\bold 0)+(-1)\cdot
f(\bold 0)=\bold 0.
\tag8.1
$$
\proclaim{Theorem 8.1} Linear mappings possess the following three
properties:
\roster
\item the identical mapping $\id_V\!:\,V\to V$ of a linear vector space
      $V$ onto itself is a linear mapping;
\item the composition of any two linear mappings $f\!:\,V\to W$ and 
      $g\!:\,W\to U$ is a linear mapping $g\compos f\!:\,V\to U$;
\item if a linear mapping $f\!:\,V\to W$ is bijective, then the inverse
      mapping $f^{-1}\!:\,W\to V$ also is a linear mapping.
\endroster
\endproclaim
\demo{Proof} The linearity of the identical mapping is obvious. Indeed,
here is the verification of the conditions \therosteritem{1} and
\therosteritem{2} from the definition~8.1 for $\id_V$:
$$
\align
&\id_V(\bold v_1+\bold v_2)=\bold v_1+\bold v_2=\id_V(\bold
v_1)+\id_V(\bold v_2),\\
&\id_V(\alpha\cdot\bold v)=\alpha\cdot\bold v=\alpha\cdot\id_V(\bold v).
\endalign
$$\par
     Let's prove the second proposition of the theorem~8.1. Consider
the composition $g\compos f$ of two linear mappings $f$ and $g$. For 
this composition the conditions \therosteritem{1} and \therosteritem{2} from the definition~8.1 are verified as follows:
$$
\align
&\aligned
  g\compos f(\bold v_1&+\bold v_2)=g(f(\bold v_1
  +\bold v_2)=g(f(\bold v_1)+f(\bold v_2))=\\
  &=g(f(\bold v_1))+g(f(\bold v_2))=g\compos f(\bold v_1)
  +g\compos f(\bold v_2),
 \endaligned\\
\vspace{1ex}
&\aligned
 g\compos f(\alpha\cdot\bold v)&=g(f(\alpha\cdot\bold v))=
 g(\alpha\cdot f(\bold v))=\alpha\cdot g(f(\bold v))\\
 &=\alpha\cdot g\compos f(\bold v).
 \endaligned
\endalign
$$\par
     Now let's prove the third proposition of the theorem~8.1.
Suppose that $f\!:\,V\to W$ is a bijective linear mapping. Then 
it possesses unique bilateral inverse mapping $f^{-1}\!:\,W\to V$
(see theorem~1.9). Let's denote
$$
\align
&\bold z_1=f^{-1}(\bold w_1+\bold w_2)-f^{-1}(\bold w_1)
-f^{-1}(\bold w_2),\\
&\bold z_2=f^{-1}(\alpha\cdot\bold w)
-\alpha\cdot f^{-1}(\bold w).
\endalign
$$
It is obvious that the linearity of the inverse mapping $f^{-1}$
is equivalent to vanishing $\bold z_1$ and $\bold z_2$. Let's
apply $f$ to these vectors:
$$
\align
&\aligned
  f(\bold z_1)&=f(f^{-1}(\bold w_1+\bold w_2)-f^{-1}(\bold w_1)
  -f^{-1}(\bold w_2))=\\
  &=f(f^{-1}(\bold w_1+\bold w_2))-f(f^{-1}(\bold w_1))
  -f(f^{-1}(\bold w_2))=\\
  &=(\bold w_1+\bold w_2)-\bold w_1-\bold w_2=\bold 0,
 \endaligned\\
\vspace{1ex}
&\aligned
 f(\bold z_2)&=f(f^{-1}(\alpha\cdot\bold w)
 -\alpha\cdot f^{-1}(\bold w))=
 f(f^{-1}(\alpha\cdot\bold w))-\\
 &-\alpha\cdot f(f^{-1}(\bold w))=\alpha\cdot\bold w-\alpha
 \cdot\bold w=\bold 0.
 \endaligned
\endalign
$$
A bijective mapping is injective. Therefore, from the equalities
$f(\bold z_1)=\bold 0$ and $f(\bold z_2)=\bold 0$ just derived 
and from the equality $f(\bold 0)=\bold 0$ derived in \thetag{8.1}
it follows that $\bold z_1=\bold z_2=\bold 0$. The theorem is 
proved.
\qed\enddemo
     Each linear mapping $f\!:\,V\to W$ is related with two 
subsets: the {\it kernel\/} $\Ker f\subset V$ and the {\it image\/}
$\Img f\subset W$. The image $\Img f=f(V)$ of a linear mapping is
defined in the same way as it was done for a general mapping in \S~1:
$$
\Img f=\{\bold \in W\!:\,\exists\,\bold v\ ((\bold v\in A)\and 
(f(\bold v)=\bold w))\}.
$$
The kernel of a linear mapping $f\!:\,V\to W$ is the set of vectors
in the space $V$ that map to zero under the action of $f$:
$$
\Ker f=\{\bold v\in V: f(\bold v)=\bold 0\}
$$
\proclaim{Theorem 8.2} The kernel and the image of a linear mapping
$f\!:\,V\to W$ both are subspaces in $V$ and $W$ respectively.
\endproclaim
\demo{Proof}In order to prove this theorem we should check the
conditions \therosteritem{1} and \therosteritem{2} from the
definition~2.2 as applied\pagebreak to the subsets $\Ker f\subset V$ 
and $\Img f\subset W$. Suppose that $\bold v_1,\bold v_2\in\Ker f$. 
Then $f(\bold v_1)=\bold 0$ and $f(\bold v_2)=\bold 0$. Suppose
also that $\bold v\in\Ker f$. Then $f(\bold v)=\bold 0$. As a
result we derive
$$
\align
&f(\bold v_1+\bold v_2)=f(\bold v_1)+f(\bold v_2)
=\bold 0+\bold 0=\bold 0,\\
&f(\alpha\cdot\bold v)=\alpha\cdot f(\bold v)
=\alpha\cdot\bold 0=\bold 0.
\endalign
$$
Hence, $\bold v_1+\bold v_2\in\Ker f$ and $\alpha\cdot\bold v\in\Ker f$.
This proves the proposition of the theorem concerning the kernel $\Ker
f$.\par
     Let $\bold w_1,\bold w_2,\bold w\in\Img f$. Then there are three
vectors $\bold v_1,\bold v_2,\bold v$ in $V$ such that $f(\bold v_1)=
\bold w_1$, $f(\bold v_2)=\bold w_2$, and $f(\bold v)=\bold w$. Hence,
we have
$$
\align
&\bold w_1+\bold w_2=f(\bold v_1)+f(\bold v_2)
=f(\bold v_1+\bold v_2),\\
&\alpha\cdot\bold w=\alpha\cdot f(\bold v)
=f(\alpha\cdot\bold v).
\endalign
$$
This meant that $\bold w_1+\bold w_2\in\Img f$ and $\alpha\cdot\bold
w\in\Img f$. The theorem is proved.
\qed\enddemo
     Remember that, according to the theorem~1.2, a linear mapping
$f\!:\,V\to W$ is surjective if and only if $\Img f=W$. There is a
similar proposition for $\Ker f$.
\proclaim{Theorem 8.3} A linear mapping $f\!:\,V\to W$ is injective
if and only if its kernel is zero, i\.\,e\. $\Ker f=\{\bold 0\}$.
\endproclaim
\demo{Proof} Let $f$ be injective and let $\bold v\in\Ker f$. Then
$f(\bold 0)=\bold 0$ and $f(\bold v)=\bold 0$. But if $\bold v\neq
\bold 0$, then due to injectivity of $f$ it would be $f(\bold v)
\neq f(\bold 0)$. Hence, $\bold v=\bold 0$. This means that the
kernel of $f$ consists of the only one element: $\Ker f=\{\bold0\}$.
\par
     Now conversely, suppose that $\Ker f=\{\bold 0\}$. Let's consider 
two different vectors $\bold v_1\neq\bold v_2$ in $V$. Then $\bold v_1
-\bold v_2\neq\bold 0$ and $\bold v_1-\bold v_2\not\in\Ker f$. Therefore,
$f(\bold v_1-\bold v_2)\neq\bold 0$. Applying the linearity of $f$, from
this inequality we derive $f(\bold v_1)-f(\bold v_2)\neq\bold 0$, i\.\,e\.
$f(\bold v_1)\neq f(\bold v_2)$. Hence, $f$ is an injective mapping. The
theorem is proved.
\qed\enddemo
     The following theorem is known as the theorem on the linear
independence of preimages. Here is its statement.
\proclaim{Theorem 8.4} Let $f\!:\,V\to W$ be a linear mapping and let
$\bold v_1,\,\ldots,\,\bold v_s$ be some vectors of a linear vector
space $V$ such that their images $f(\bold v_1),\,\ldots,\,f(\bold v_n)$
in $W$ are linearly independent. Then the vectors $\bold v_1,\,\ldots,
\,\bold v_s$ themselves are also linearly independent.
\endproclaim
\demo{Proof} In order to prove the theorem let's consider a linear
combination of the vectors $\bold v_1,\,\ldots,\,\bold v_s$ being
equal to zero:
$$
\alpha_1\cdot\bold v_1+\ldots+\alpha_s\cdot\bold v_s=\bold 0.
$$
Applying $f$ to both sides of this equality and using the fact that
$f$ is a linear mapping, we obtain quite similar equality for the
images
$$
\alpha_1\cdot f(\bold v_1)+\ldots+\alpha_s\cdot f(\bold v_s)=\bold 0.
$$
However, these images $f(\bold v_1),\,\ldots,\,f(\bold v_n)$ are
linearly independent. Hence, all coefficients in the above linear
combination are equal to zero: \pagebreak $\alpha_1=\ldots=\alpha_s=0$.
Then the initial linear combination is also necessarily trivial. This 
proves that the vectors $\bold v_1,\,\ldots,\,\bold v_s$ are linearly
independent.
\qed\enddemo
     A linear vector space is a set. But it is not simply a set ---
it is a structured set. It is equipped with algebraic operations 
satisfying the axioms \therosteritem{1}-\therosteritem{8}. Linear
mappings are those being concordant with the structures of a linear 
vector space in the spaces they are acting from and to. In algebra
such mappings concordant with algebraic structures are called 
{\it morphisms}. So, in algebraic terminology, linear mappings are
morphisms of linear vector spaces. 
\definition{Definition 8.2} Two linear vector spaces $V$ and $W$ 
are called {\it isomorphic} if there is a bijective linear mapping
$f\!:\,V\to W$ binding them.
\enddefinition
     The first example of an isomorphism of linear vector spaces
is the mapping $\psi\!:\,V\to\Bbb K^n$ in \thetag{5.4}. Because
of the existence of such mapping we can formulate the following 
theorem.\par
\proclaim{Theorem 8.5} Any $n$-dimensional linear vector space 
$V$ is isomorphic to the arithmetic linear vector space $\Bbb K^n$.
\endproclaim
     Isomorphic linear vector spaces have many common features. 
Often they can be treated as undistinguishable. In particular, we
have the following fact.\par
\proclaim{Theorem 8.6} If a linear vector space $V$ is isomorphic 
to a finite-dimensional vector space $W$, then $V$ is also finite-dimensional and the dimensions of these two spaces do coincide:
$\dim V=\dim W$.
\endproclaim
\demo{Proof} Let $f\!:\,V\to W$ be an isomorphism of spaces $V$ and
$W$. Assume for the sake of certainty that $\dim W=n$ and choose 
a basis $\bold h_1,\,\ldots,\,\bold h_n$ in $W$. By means of inverse
mapping $f^{-1}\!:\,W\to V$ we define the vectors $\bold e_i=f^{-1}
(\bold h_i)$, $i=1,\ldots,n$. Let $\bold v$ be an arbitrary vector 
of $V$. Let's map it with the use of $f$ into the space $W$ and then
expand in the basis:
$$
f(\bold v)=\alpha_1\cdot\bold h_1+\ldots+\alpha_n\cdot\bold  h_n.
$$
Applying the inverse mapping $f^{-1}$ to both sides of this equality,
due to the linearity of $f^{-1}$ we get the expansion
$$
\bold v=\alpha_1\cdot\bold e_1+\ldots+\alpha_n\cdot\bold e_n.
$$
From this expansion we derive that $\{\bold e_1,\,\ldots,\,\bold e_n\}$
is a finite spanning system in $V$. The finite dimensionality of $V$ 
is proved. The linear independence of $\bold e_1,\,\ldots,\,\bold e_n$
follows from the theorem~8.4 on the linear independence of preimages.
Hence, $\bold e_1,\,\ldots,\,\bold e_n$ is a basis in $V$ and 
$\dim V=n=\dim W$. The theorem is proved.
\qed\enddemo
\head
\S\,9. The matrix of a linear mapping.
\endhead
     Let $f\!:\,V\to W$ be a linear mapping from $n$-dimensional
vector space $V$ to $m$-dimensional vector space $W$. Let's choose
a basis $\bold e_1,\,\ldots,\,\bold e_n$ in $V$ and a 
basis $\bold h_1,\,\ldots,\,\bold h_m$ in $W$. \pagebreak Then 
consider the images of basis vectors $\bold e_1,\,\ldots,\,\bold e_n$ 
in $W$ and expand them in the basis $\bold h_1,\,\ldots,\,\bold h_m$:
$$
\hskip -2em
\matrix
f(e_1) &= &F^1_1\cdot h_1 &+ &\hdots &+ &F^m_1\cdot h_m,\\
\hdotsfor 7\\ \vspace{1\jot}
f(e_n) &= &F^1_n\cdot h_1 &+ &\hdots &+ &F^m_n\cdot h_m.
\endmatrix
\tag9.1
$$
Totally in \thetag{9.1} we have $n$ expansions that define $n\,m$ numbers
$F^i_j$. These numbers are arranged into a rectangular $m\times n$ matrix
which is called the {\it matrix of the linear mapping} $f$ in a pair of
bases $\bold e_1,\,\ldots,\,\bold e_n$ and $\bold h_1,\,\ldots,\,\bold h_m$:
$$
\hskip -2em
F=\Vmatrix F^1_1 & \hdots & F^1_n \\
        \vdots & \hdots & \vdots\\
         F^m_1 & \hdots & F^m_n
\endVmatrix.
\tag9.2
$$
When placing the element $F^i_j$ into the matrix \thetag{9.2}, the upper
index determines the row number, while the lower index determines the
column number. In other words, the matrix $F$ is composed by the 
column vectors formed by coordinates of the vectors $f(\bold e_1),\,
\ldots,\,f(\bold e_n)$ in the basis $\bold h_1,\,\ldots,\,\bold h_m$.
The expansions \thetag{9.1}, which determine the components of this 
matrix, are convenient to write as follows:
$$
\hskip -2em
f(\bold e_j)=\sum^m_{i=1} F^i_j\cdot\bold h_i.
\tag9.3
$$\par
     Let $\bold x$ be an arbitrary vector of $V$ and let $\bold y=
f(\bold x)$ be its image under the mapping $f$. If we expand the 
vector $\bold x$ in the basis: $\bold x=x^1\cdot\bold e_1+\ldots
+x^n\cdot\bold e_n$, then, taking into account \thetag{9.3}, for the
vector $\bold y$ we get 
$$
\bold y=f(\bold x)=\sum^n_{j=1} x^j\cdot f(\bold e_j)=\sum^n_{j=1}
x^j\cdot\left(\,\shave{\sum^m_{i=1}}F^i_j\cdot\,\bold h_i\right).
$$
Changing the order of summations in the above expression, we get the
expansion of the vector $\bold y$ in the basis $\bold h_1,\,\ldots,\,
\bold h_m$:
$$
\bold y=f(\bold x)=\sum^m_{i=1}\left(\,\shave{\sum^n_{j=1}}F^i_j\,x^j
\right)\!\cdot\bold h_i.
$$
Due to the uniqueness of such expansion for the coordinates of the
vector $\bold y$ in the basis $\bold h_1,\,\ldots,\,\bold h_m$ we get
the following formula:
$$
\hskip -2em
y^i=\sum^n_{j=1}F^i_j\,x^j.
\tag9.4
$$
This formula \thetag{9.4} is the basic application of the matrix of
a linear mapping. It is used for calculating the coordinates of the
vector $f(\bold x)$ through the coordinates of $\bold x$. In matrix
form this formula is written as
$$
\hskip -2em
\Vmatrix y^1\\ \vdots\\ y^m\endVmatrix=
\Vmatrix F^1_1 & \hdots & F^1_n \\
        \vdots & \hdots & \vdots\\
         F^m_1 & \hdots & F^m_n
\endVmatrix\cdot
\Vmatrix \vphantom{y^1}x^1\\ \vdots\\ x^n\endVmatrix.
\tag9.5
$$\par
     Remember that when composing a column vector of the coordinates
of a vector $\bold x$, we negotiated to understand this procedure as a
linear mapping $\psi\!:\,V\to\Bbb K^n$ (see formulas \thetag{5.4} and
the theorem~8.5). Denote by $\tilde\psi\!:\,W\to\Bbb K^m$ the analogous
mapping for a vector $\bold y$ in $W$. Then the matrix relationship
\thetag{9.5} can be treated as a mapping $F\!:\,\Bbb K^n\to\Bbb K^m$.
These three mappings $\psi,\,\tilde\psi, F$ and the initial mapping $f$ 
can be written in a diagram:
$$
\hskip -2em
\CD
V @>f>> W\\
@V\psi VV @VV\tilde\psi V\\
\Bbb K^n @>>F>\Bbb K^m
\endCD
\tag9.6
$$
Such diagrams are called {\it commutative diagrams} if the compositions
of mappings {\tencyr\char '074}when passing along arrows{\tencyr\char '076}
from any node to any other node do not depend on a particular path
connecting these two nodes. When applied to the diagram \thetag{9.6}, the
commutativity means $\tilde\psi\compos f=F\compos\psi$. Due to the 
bijectivity of linear mappings $\psi$ and $\tilde\psi$ the condition of
commutativity of the diagram \thetag{9.6} can be written as
$$
\xalignat 2
&\hskip -2em
F=\tilde\psi\,{\ssize\circ}f{\ssize\circ}\,\psi^{-1},
&&f=\tilde\psi^{-1}{\ssize\circ}\,F{\ssize\circ}\,\psi.
\tag9.7
\endxalignat
$$
The reader can easily check that the relationships \thetag{9.7} are
fulfilled due to the way the matrix $F$ is constructed. Hence, the
diagram  \thetag{9.6} is commutative.\par
     Now let's look at the relationships \thetag{9.7} from a little
bit other point of view. Let $V$ and $W$ be two spaces of the 
dimensions $n$ and $m$ respectively. Suppose that we have an arbitrary
$m\times n$ matrix $F$. Then the relationship \thetag{9.5} determines
a linear mapping $F:\Bbb K^n\to\Bbb K^m$. Choosing bases
$\bold e_1,\,\ldots,\,\bold e_n$ and $\bold h_1,\,\ldots,\,\bold h_m$ 
in $V$ and $W$ we can use the second relationship \thetag{9.7} in order 
to define the linear mapping $f\!:\,V\to W$. The matrix of this mapping 
in the bases $\bold e_1,\,\ldots,\,\bold e_n$ and $\bold h_1,\,\ldots,
\,\bold h_m$ coincides with $F$ exactly. Thus, we have proved the
following theorem.
\proclaim{Theorem 9.1} Any rectangular $m\times n$ matrix $F$ can be
constructed as a matrix of a linear mapping $f\!:\,V\to W$ from  $n$-dimensional vector space $V$ to $m$-dimensional vector space $W$
in some pair of bases in these spaces.
\endproclaim
     A more straightforward way of proving the theorem~9.1 than we 
considered above can be based on the following theorem.
\proclaim{Theorem 9.2} For any basis $\bold e_1,\,\ldots,\,\bold e_n$
in $n$-dimensional vector space $V$ and for any set of $n$ vectors 
$\bold w_1,\,\ldots,\,\bold w_n$ in another vector space $W$ there 
is a linear mapping $f\!:\,V\to W$ such that $f(\bold e_i)=\bold w_i$ 
for $i=1,\,\ldots,\,n$. 
\endproclaim
\demo{Proof} Once the basis $\bold e_1,\,\ldots,\,\bold e_n$ in $V$ 
is chosen, this defines the mapping $\psi\!:\,V\to\Bbb K^n$ (see \thetag{5.4}). In order to construct the required mapping $f$ we 
define a mapping $\varphi\!:\,\Bbb K^n\to W$ by the following
relationship:
$$
\varphi:\
\Vmatrix x^1\\ \vdots\\ x^n\endVmatrix\mapsto
x^1\cdot w_1+\ldots+x^n\cdot w_n.
$$
Now it is easy to verify \pagebreak that the required mapping is the composition $f=\varphi\compos\psi$.
\qed\enddemo
     Let's return to initial situation. Suppose that we have a mapping
$f\!:\,V\to W$ that determines a matrix $F$ upon choosing two bases 
$\bold e_1,\,\ldots,\,\bold e_n$ and $\bold h_1,\,\ldots,\,\bold h_m$
in $V$ and $W$ respectively. The matrix $F$ essentially depends on the
choice of bases. In order to describe this dependence we consider four 
bases --- two bases in $V$ and other two bases in $W$. Suppose that
$S$ and $P$ are direct transition matrices for that pairs of bases.
Their components are defined as follows:
$$
\xalignat 2
&\tilde\bold e_k=\sum^n_{j=1} S^j_k\cdot\bold e_j,
&&\tilde\bold h_r=\sum^m_{i=1} P^i_r\cdot\bold h_i.
\endxalignat
$$
The inverse transition matrices $T=S^{-1}$ and $Q=P^{-1}$ are defined
similarly:
$$
\xalignat 2
&\bold e_j=\sum^n_{k=1} T^k_j\cdot\tilde\bold e_k,
&&\bold h_i=\sum^m_{r=1} Q^r_i\cdot\tilde\bold h_r.
\endxalignat
$$
We use these relationships and the above relationships \thetag{9.3} 
in order to carry out the following calculations for the vector
$f(\tilde\bold e_k)$:
$$
\align
f(\tilde\bold e_k)&=\sum^n_{j=1}S^j_k\cdot f(\bold e_j)=
\sum^n_{j=1}S^j_k\cdot\!\left(\,\shave{\sum^m_{i=1}}
F^i_j\cdot\bold h_i\right)=\\
&=\sum^n_{j=1}S^j_k\cdot\!\left(\,\shave{\sum^m_{i=1}}
F^i_j\cdot\!\left(\,\shave{\sum^m_{r=1}}Q^r_i\cdot\tilde\bold h_r
\right)\right).
\endalign
$$
Upon changing the order of summations this result is written as
$$
f(\tilde\bold e_k)=\sum^m_{r=1}\left(\,\shave{\sum^m_{i=1}
\sum^n_{j=1}}Q^r_i\,F^i_j\,S^j_k\right)\!\cdot\tilde\bold h_r.
$$
The double sums in round brackets are the coefficients of the 
expansion of the vector $f(\tilde\bold e_k)$ in the basis 
$\tilde\bold h_1,\,\ldots,\,\tilde\bold h_m$. They determine
the matrix of the linear mapping $f$ in wavy bases 
$\tilde\bold e_1,\,\ldots,\,\tilde\bold e_n$ and $\tilde\bold
h_1,\,\ldots,\,\tilde\bold h_m$:
$$
\hskip -2em
\tilde F^r_k=\sum^m_{i=1}\sum^n_{j=1}Q^r_i\,F^i_j\,S^j_k.
\tag9.8
$$
In a similar way one can derive the converse relationship expressing
$F$ through $\tilde F$:
$$
\hskip -2em
F^i_j=\sum^m_{r=1}\sum^n_{k=1}P^i_r\,\tilde F^r_k\,T^k_j.
\tag9.9
$$
The relationships \thetag{9.8} and \thetag{9.9} are called the {\it
transformation formulas for the matrix of a linear mapping under a
change of bases}. They can be written as
$$
\xalignat 2
&\tilde F=P^{-1}\,F\,S,
&&F=P\,\tilde F\,S^{-1}.
\tag9.10
\endxalignat
$$
This is the matrix form of the relationships \thetag{9.8} and \thetag{9.9}.
\par
     The transformation formulas like \thetag{9.10} lead us to the broad
class of problems of {\tencyr\char '074}bringing to a canonic
form{\tencyr\char '076}. In our particular case a change \pagebreak 
of bases in the spaces $V$ and $W$ changes the matrix of the linear mapping
$f\!:\,V \to W$. The problem of bringing to a canonic form in this case
consists in finding the optimal choice of bases, where the matrix $F$ has
the most simple (canonic) form. The following theorem solving this
particular problem is known as the theorem on bringing to the {\it almost
diagonal} form.\par
\proclaim{Theorem 9.3} Let $f\!:\,V\to W$ be some nonzero linear mapping
from $n$-dimensional vector space $V$ to $m$-dimensional vector space $W$.
Then there is a choice of bases in $V$ and $W$ such that the matrix $F$ of
this mapping has the following almost diagonal form:
$$
\hskip -2em
F=\aligned
&\hphantom{x}\overbrace{\hphantom{xxxxxxxxxxx}}^s\\
&\Vmatrix
1 & 0 & \hdots & 0 & 0 & \hdots & 0 & 0\\
0 & 1 & \hdots & 0 & 0 & \hdots & 0 & 0\\
\vdots & \vdots & \ddots & \vdots &\vdots &\, & \vdots &\vdots\\
0 & 0 & \hdots & 1 & 0 & \hdots & 0 & 0\\
0 & 0 & \hdots & 0 & 0 & \hdots & 0 & 0\\
\vdots & \vdots &\,  & \vdots &\vdots &\, & \vdots & \vdots\\
0 & 0 & \hdots & 0 & 0 & \hdots & 0 & 0\\
\endVmatrix
\endaligned
\aligned
&\vphantom{\vrule height 1ex depth 1.7ex}\\
&\left.\vphantom{\vrule height 5ex depth 5ex}\right\} s\\
&\vphantom{\vrule height 4ex depth 4ex}
\endaligned
\tag9.11
$$
\endproclaim
\demo{Proof} The purely zero mapping $0\!:\,V\to W$ maps each
vector of the space $V$ to zero vector in $W$. The matrix of such
mapping consists of zeros only. There is no need to formulate the 
problem of bringing it to a canonic form.\par
     Let $f\!:\,V\to W$ be a nonzero linear mapping. The integer 
number $s=\dim(\Img f)$ is called the {\it rank\/} of the mapping
$f$. The rank of a nonzero mapping is not equal to zero. We begin
constructing a canonic base in $W$ by choosing a base $\bold h_1,\,
\ldots,\,\bold h_s$ in the image space $\Img f$. For each basis
vector $\bold h_i\in\Img f$ there is a vector $\bold e_i\in V$ such
that $f(\bold e_i)=\bold h_i$, $i=1,\,\ldots,\,s$. These vectors
$\bold e_1,\,\ldots,\,\bold e_s$ are linearly independent due to
the theorem~8.4. Let $r=\dim(\Ker f)$. We choose a basis in $\Ker f$
and denote the basis vectors by $\bold e_{s+1},\,\ldots,\,\bold 
e_{s+r}$. Then we consider the vectors
$$
\hskip -2em
\bold e_1,\,\ldots,\,\bold e_s,\,\bold e_{s+1},\,\ldots,\,
\bold e_{s+r}
\tag9.12
$$
and prove that they form a basis in $V$. For this purpose we use the
theorem~4.6.\par
     Let's begin with checking the condition \therosteritem{1} in the
theorem~4.6 for the vectors \thetag{9.12}. In order to prove the linear independence of these vectors we consider a linear combination of them
being equal to zero:
$$
\hskip -2em
\alpha_1\cdot\bold e_1+\ldots+\alpha_s\cdot\bold e_s
+\alpha_{s+1}\cdot\bold e_{s+1}+\ldots+\alpha_{s+r}\cdot\bold
e_{s+r}=\bold 0.
\tag9.13
$$
Let's apply the mapping $f$ to both sides of the equality \thetag{9.13}
and take into account that $f(\bold e_i)=\bold h_i$ for $i=1,\ldots,s$.
Other vectors belong to the kernel of the mapping $f$, therefore,
$f(\bold e_{s+i})=\bold 0$ for $i=1,\ldots,r$. Then from \thetag{9.13}
we derive
$$
\alpha_1\cdot\bold h_1+\ldots+\alpha_s\cdot\bold h_s=\bold 0.
$$
The vectors \pagebreak $\bold h_1,\,\ldots,\,\bold h_s$ form a basis
in $\Img f$. They are linearly independent. Hence, $\alpha_1=\ldots=
\alpha_s=0$. Taking into account this fact, we reduce \thetag{9.13} to
$$
\alpha_{s+1}\cdot\bold e_{s+1}+\ldots+\alpha_{s+r}\cdot\bold 
e_{s+r}=\bold 0.
$$
The vectors $\bold e_{s+1},\,\ldots,\,\bold e_{s+r}$ form a basis in
$\Ker f$. They are linearly independent, therefore, $\alpha_{s+1}=
\ldots=\alpha_{s+r}=0$. As a result we have proved that all coefficients
of the linear combination \thetag{9.13} are necessarily zero. Hence,
the vectors \thetag{9.12} are linearly independent.\par
     Now lets check the second condition of the theorem~4.6 for the
vectors \thetag{9.12}. Assume that $\bold v$ is an arbitrary vector
in $V$. Then $f(\bold v)$ belongs to $\Img f$. Let's expand $f(\bold v)$
in the basis $\bold h_1,\,\ldots,\,\bold h_s$:
$$
\hskip -2em
f(\bold v)=\beta_1\cdot\bold h_1+\ldots+\beta_s\cdot\bold h_s.
\tag9.14
$$
Remember that $f(\bold e_i)=\bold h_i$ for $i=1,\,\ldots,\,s$. Then
from \thetag{9.14} we derive 
$$
\aligned
\bold 0=f(\bold v)&-\beta_1\cdot f(\bold e_1)-\ldots-\beta_s\cdot
f(\bold e_s)=\\
  &=f(\bold v-\beta_1\cdot\bold e_1-\ldots-\beta_s\cdot\bold e_s).
\endaligned
\tag9.15
$$
Let's denote $\tilde\bold v=\bold v-\beta_1\cdot\bold e_1-\ldots
-\beta_s\cdot\bold e_s$. From \thetag{9.15} we derive $f(\tilde
\bold v)=\bold 0$ for this vector $\tilde\bold v$. Hence, $\tilde
\bold v\in\Ker f$. Let's expand $\tilde\bold v$ in the basis of
$\Ker f$:
$$
\tilde\bold v=\beta_{s+1}\cdot\bold e_{s+1}+\ldots+
\beta_{s+r}\cdot\bold e_{s+r}.
$$
From the formula $\tilde\bold v=\bold v-\beta_1\cdot\bold e_1
-\ldots-\beta_s\cdot\bold e_s$ and the above expansion we get
$$
\bold v=\beta_1\cdot\bold e_1+\ldots+\beta_s\cdot\bold e_s
+\beta_{s+1}\cdot\bold e_{s+1}+\ldots+\beta_{s+r}\cdot\bold 
e_{s+r}.
$$
This means that the vectors \thetag{9.12} form a spanning system
in $V$. The condition \therosteritem{2} of the theorem~4.6 for
them is also fulfilled. Thus, the vectors \thetag{9.12} form a
basis in $V$. This yields the equality
$$
\hskip -2em
\dim V=s+r.
\tag9.16
$$
In order to complete the proof of the theorem we need to complete
the basis $\bold h_1,\,\ldots,\,\bold h_s$ of $\Img f$ up to a basis
$\bold h_1,\,\ldots,\,\bold h_s,\,\bold h_{s+1},\,\ldots,\,\bold
h_m$ in the space $W$. For the vector $f(\bold e_j)$ with 
$j=1,\,\ldots,\,s$ we have the expansion
$$
f(\bold e_j)=\bold h_j=\sum^s_{i=1}\delta^i_j\cdot\bold h_i+
\sum^m_{i=s+1} 0\cdot\bold h_i.
$$
If $j=s+1,\,\ldots,\,s+r$, the expansion for $f(\bold e_j)$ is
purely zero:
$$
f(\bold e_j)=\bold 0=\sum^s_{i=1} 0\cdot\bold h_i+\sum^m_{i=s+1}0
\cdot\bold h_i.
$$
Due to these expansions the matrix of the mapping $f$ in the
bases that we have constructed above has the required almost 
diagonal form \thetag{9.10}.
\qed\enddemo
     In proving this theorem we have proved simultaneously the
next \pagebreak one.\par
\proclaim{Theorem 9.4} Let $f\!:\,V\to W$ be a linear mapping from
$n$-dimensional space $V$ to an arbitrary linear vector space $W$. 
Then
$$
\hskip -2em
\dim(\Ker f)+\dim(\Img f)=\dim V.
\tag9.17
$$
\endproclaim
    This theorem~9.4 is known as the theorem {\it on the sum of
dimensions of the kernel and the image} of a linear mapping. The
proposition of the theorem in the form of the relationship 
\thetag{9.17} immediately follows from \thetag{9.16}.\par
\head
\S\,10. Algebraic operations with mappings.
The space of homomorphisms $\Hom(V,W)$.
\endhead
\rightheadtext{\S\,10. Algebraic operations with mappings.}
\definition{Definition 10.1} Let $V$ and $W$ be two linear vector
spaces and let $f\!:\,V\to W$ and $g\!:\,V\to W$ be two linear
mappings from $V$ to $W$. The linear mapping $h\!:\,V\to W$ defined
by the relationship $h(\bold v)=f(\bold v)+g(\bold v)$, where
$\bold v$ is an arbitrary vector of $V$, is called the {\it sum\/}
of the mappings $f$ and $h$.
\enddefinition
\definition{Definition 10.2} Let $V$ and $W$ be two linear vector
spaces over a numeric field $\Bbb K$ and let $f\!:\,V\to W$ be a
linear mapping from $V$ to $W$. The linear mapping $h\!:\,V\to W$
defined by the relationship $h(\bold v)=\alpha\cdot f(\bold v)$,
where $\bold v$ is an arbitrary vector of $V$, is called the
{\it product\/} of the number $\alpha\in\Bbb K$ and the mapping $f$.
\enddefinition
     The algebraic operations introduced by the definitions~10.1 
and 10.2 are called {\it pointwise addition} and {\it pointwise
multiplication by a number}. Indeed, they are calculated 
{\tencyr\char '074}pointwise{\tencyr\char '076} by adding the
values of the initial mappings and by multiplying them by a number
for each specific argument $\bold v\in V$. These operations are
denoted by the same signs as the corresponding operations with
vectors: $h=f+g$ and $h=\alpha\cdot f$. The writing $(f+g)(\bold v)$
is understood as the sum of mappings applied to the vector $\bold v$.
Another writing $f(\bold v)+g(\bold v)$ denotes the sum of the 
results of applying $f$ and $g$ to $\bold v$ separately. Though the
results of these calculations do coincide, their meanings are different.
In a similar way one should distinguish the meanings of left and right 
sides of the following equality:
$$
(\alpha\cdot f)(\bold v)=\alpha\cdot f(\bold v).
$$\par
     Let's denote by $\Map(V,W)$ the set of all mappings from the space
$V$ to the space $W$. Sometimes this set is denoted by $W^V$.
\proclaim{Theorem 10.1} Let $V$ and $W$ be two linear spaces over a 
numeric field $\Bbb K$. Then the set of mappings $\Map(V,W)$ equipped 
with the operations of pointwise addition and pointwise multiplication 
by numbers fits the definition of a linear vector space over the
numeric field $\Bbb K$.
\endproclaim
\demo{Proof} Let's verify the axioms of a linear vector space for the
set of mappings $\Map(V,W)$. In the case of the first axiom we should
verify the coincidence of the mappings $f+g$ and $g+f$. Remember that
the coincidence of two mappings is equivalent to the coincidence of
their values when applied to an arbitrary vector $v\in V$. The 
following calculations establish the latter coincidence:
$$
(f+g)(\bold v)=f(\bold v)+g(\bold v)=g(\bold v)
+f(\bold v)=(g+f)(\bold v).
$$
As we see in the above calculations, the equality $f+g=g+f$ follows
from the commutativity axiom for the addition of vectors in $W$ due
to pointwise nature of the addition of mappings. The same arguments
are applicable when verifying the axioms \therosteritem{2}, \therosteritem{5}, and \therosteritem{6} for the algebraic operations 
with mappings:
$$
\align
&\aligned
  &((f+g)+h)(\bold v)=(f+g)(\bold v)+h(\bold v)=(f(\bold v)
  +g(\bold v))+h(\bold v)=\\
  &=f(\bold v)+(g(\bold v)+h(\bold v))=f(\bold v)+(g+h)(\bold
  v)=(f+(g+h))(\bold v)
 \endaligned\\
\vspace{1.7ex}
&\aligned
  &(\alpha\cdot(f+g))(\bold v)=\alpha\cdot(f+g)(\bold v)
  =\alpha\cdot(f(\bold v)+g(\bold v))=\\
  &=\alpha\cdot f(\bold v)+\alpha\cdot g(\bold v)
  =(\alpha\cdot f)(\bold v)+(\alpha\cdot g)(\bold v)
  =(\alpha\cdot f+\alpha\cdot g)(\bold v)
 \endaligned\\
\vspace{1.7ex}
&\aligned
  ((\alpha+\beta)\cdot f)(\bold v)=(&\alpha+\beta)\cdot f(\bold
  v)=\alpha\cdot f(\bold v)+\beta\cdot f(\bold v)=\\
  &=(\alpha\cdot f)(\bold v)+(\beta\cdot f)(\bold v)
  =(\alpha\cdot f+\beta\cdot f)(\bold v)
 \endaligned\\
\endalign
$$
For the axioms \therosteritem{7} these calculations look like
$$
(\alpha\cdot(\beta\cdot f))(\bold v)=\alpha\cdot(\beta\cdot 
f)(\bold v)=\alpha\cdot(\beta\cdot f(\bold v))=(\alpha\beta)
\cdot f(\bold v)=((\alpha\beta)\cdot f)(\bold v).
$$
In the case of the axiom \therosteritem{8} the calculations
are even more simple:
$$
(1\cdot f)(\bold v)=1\cdot f(\bold v)=f(\bold v).
$$\par
     Now let's consider the rest axioms \therosteritem{3} and
\therosteritem{4}. The zero mapping is the best pretender for 
the role of zero element in the space $\Map(V,W)$, it maps
each vector $\bold v\in V$ to zero vector of the space $W$.
For this mapping we have
$$
(f+0)(\bold v)=f(\bold v)+0(\bold v)=f(v)+\bold 0=f(\bold v).
$$
As we see, the axiom \therosteritem{3} in $\Map(V,W)$ is fulfilled.
\par
     Suppose that $f\in\Map(V,W)$. We define the opposite mapping
$f'$ for $f$ as follows: $f'=(-1)\cdot f$. Then we have
$$
\align
(f&+f')(\bold v)=(f+(-1)\cdot f)(\bold v)=f(\bold v)\,+\\
&+\,((-1)\cdot f)(\bold v)=f(\bold v)+(-1)\cdot f(\bold v)
=\bold 0=0(\bold v).
\endalign
$$
The axiom \therosteritem{4} in $\Map(V,W)$ is also fulfilled. This
completes the proof of the theorem~10.1. 
\qed\enddemo
     In typical situation the space $\Map(V,W)$ is very large. Even
for the finite-dimensional spaces $V$ and $W$ usually it is an
infinite-dimensional space. In linear algebra the much smaller subset 
of $\Map(V,W)$ is studied. This is the set of all linear mappings 
from $V$ to $W$. It is denoted $\Hom(V,W)$ and is called the {\it set 
of homomorphisms}. The following two theorems show that $\Hom(V,W)$
is closed with respect to algebraic operations in $\Map(V,W)$.
Therefore, we can say that $\Hom(V,W)$ is the {\it space of 
homomorphisms}.
\par
\proclaim{Theorem 10.2} The pointwise sum of two linear mappings
$f\!:\,V\to W$ and $g\!:\,V\to W$ is a linear mapping from the 
space $V$ to the space $W$.
\endproclaim
\proclaim{Theorem 10.3} The pointwise product of a linear mapping
$f\!:\,V\to W$ by a number $\alpha\in\Bbb K$ is a linear mapping
from the space $V$ to the space $W$.
\endproclaim
\demo{Proof} Let $h=f+g$ be the sum of two linear mappings $f$ and $g$.
The following calculations prove the linearity of the mapping $h$:
$$
\align
&\aligned
  h(\bold v_1+\bold v_2)=f&(\bold v_1+\bold v_2)+g(\bold v_1
  +\bold v_2)=(f(\bold v_1)+\\
  +f(\bold v_2))&+((g(\bold v_1)+g(\bold v_2)=(f(\bold v_1)+\\
  &+g(\bold v_1))+(f(\bold v_2)+g(\bold v_2))=h(\bold v_1)
  +h(\bold v_2),
 \endaligned\\
\vspace{1.7ex}
&\aligned
  h(\beta\cdot\bold v)=f(\beta&\cdot\bold v)+g(\beta\cdot\bold
  v)=\beta\cdot f(\bold v)+\\
  &+\beta\cdot g(\bold v)=\beta\cdot(f(\bold v)+g(\bold v)
  =\beta\cdot h(\bold v).
 \endaligned
\endalign
$$
Now let's consider the product of the mapping $f$ and the number
$\alpha$. Let's denote it by $h$, i\.\,e\. let's denote $h=\alpha
\cdot f$. Then the following calculations 
$$
\align
&\aligned
  h(\bold v_1+\bold v_2)&=\alpha\cdot f(\bold v_1+\bold v_2)
  =\alpha\cdot(f(\bold v_1)+\\
  &+f(\bold v_2))=\alpha\cdot f(\bold v_1)+\alpha\cdot f(\bold v_2)
  =h(\bold v_1)+h(\bold v_2),
 \endaligned\\
\vspace{1.7ex}
&\aligned
  h(\beta\cdot\bold v)&=\alpha\cdot f(\beta\cdot\bold v)=\alpha
  \cdot(\beta\cdot f(\bold v))=\\
  &=(\alpha\beta)\cdot f(\bold v)=(\beta\alpha)\cdot f(\bold v)=
    \beta\cdot(\alpha\cdot f(\bold v))=\beta\cdot\,h(\bold v).
 \endaligned
\endalign
$$
prove the linearity of the mapping $h$ and thus complete the proofs 
of both theorems~10.2 and 10.3.
\qed\enddemo
     The space of homomorphisms $\Hom(V,W)$ is a subspace in the
space of all mappings $\Map(V,W)$. It is much smaller and it consists
of objects which are in the scope of linear algebra. For 
finite-dimensional spaces $V$ and $W$ the space of homomorphisms 
$\Hom(V,W)$ is also finite-dimensional. This is the result of the
following theorem.\par
\proclaim{Theorem 10.4} For finite-dimensional spaces $V$ and $W$ 
the space of homomorphisms $\Hom(V,W)$ is also finite-dimensional.
Its dimension is given by formula
$$
\hskip -2em
\dim(\Hom(V,W))=\dim(V)\cdot\,\dim(W).
\tag10.1
$$
\endproclaim
\demo{Proof} Let $\dim V=n$ and $\dim W=m$. We choose a basis 
$\bold e_1,\,\ldots,\,\bold e_n$ in the space $V$ and another
basis $\bold h_1,\,\ldots,\,\bold h_m$ in the space $W$. Let
$1\leqslant i\leqslant n$ and $1\leqslant j\leqslant m$. For
each fixed pair of indices $i,j$ within the above ranges we
consider the following set of $n$ vectors in the space $W$:
$$
\bold w_1=\bold 0,\ \ldots,\ \bold w_{i-1}=\bold 0,\ 
\bold w_i=\bold h_j,\ \bold w_{i+1}=\bold 0,\ \ldots,
\ \bold w_n=\bold 0.
$$
All vectors in this set are equal to zero, except for the $i$-th
vector $\bold w_i$ which is equal to $j$-th basis vector $\bold 
h_j$. Now we apply the theorem~9.2 to the basis $\bold e_1,\,
\ldots,\,\bold e_n$ in $V$ and to the set of vector $\bold w_1,
\,\ldots,\,\bold w_n$. This defines the linear mapping $E^i_j\!:
\,V\to W$ such that $E^i_j(\bold e_s)=\bold w_s$ for all $s=1,\,
\ldots,\,n$. We write this fact as
$$
\hskip -2em
E^i_j(\bold e_s)=\delta^i_s\cdot\bold h_j,
\tag10.2
$$
where $\delta^i_s$ is the Kronecker symbol. As a result we have
constructed $n\cdot\,m$ mappings $E^i_j$ satisfying the relationships \thetag{10.2}:
$$
\hskip -2em
E^i_j\!:\,V\to W,\text{\ where \ } 1\leqslant i\leqslant n,\ 
1\leqslant j\leqslant m.
\tag10.3
$$\par
     Now we show that the mapping \thetag{10.3} span the space of
homomorphisms $\Hom(V,W)$. For this purpose we take a linear mapping
$f\in\Hom(V,W)$. Suppose that $F$ is its matrix in the pair of bases
$\bold e_1,\,\ldots,\,\bold e_n$ and $\bold h_1,\,\ldots,\,\bold h_m$.
Denote by $F^j_i$ the elements of this matrix. Then the result of
applying $f$ to an arbitrary vector $\bold v\in V$ is determined by
coordinates of this vector according to the formula
$$
\hskip -2em
f(\bold v)=\sum^n_{i=1}v^i\cdot f(\bold e_i)=\sum^n_{i=1}\sum^m_{j=1}
(F^j_i\,v^i)\cdot\bold h_j.
\tag10.4
$$
Applying $E^i_j$ to the same vector $\bold v$ and taking into account
\thetag{10.2}, we derive 
$$
\hskip -2em
E^i_j(\bold v)=\sum^n_{s=1}v^s\cdot\,E^i_j(\bold e_s)
=\sum^n_{s=1}(v^s\,\delta^i_s)\cdot\bold h_j=v^i\cdot\,\bold h_j.
\tag10.5
$$
Now, comparing the relationships \thetag{10.4} and \thetag{10.5},
we find
$$
f(\bold v)=\sum^n_{i=1}\sum^m_{j=1}F^j_i\cdot\,E^i_j(\bold v).
$$
Since $\bold v$ is an arbitrary vector of the space $V$, this formula
means that $f$ is a linear combination of the mappings \thetag{10.3}:
$$
f=\sum^n_{i=1}\sum^m_{j=1}F^j_i\cdot\,E^i_j.
$$
Hence, the mappings \thetag{10.3} span the space of homomorphisms 
$\Hom(V,W)$. This proves the finite-dimensionality of the space
$\Hom(V,W)$.\par
     In order to calculate the dimension of $\Hom(V,W)$ we shall
prove that the mappings \thetag{10.3} are linearly independent.
Let's consider a linear combination of these mappings, which is 
equal to zero:
$$
\hskip -2em
\sum^n_{i=1}\sum^m_{j=1}\gamma^j_i\cdot E^i_j=0.
\tag10.6
$$
Both left and right hand sides of the equality \thetag{10.6} 
represent the zero mapping $0\!:\,V\to W$. Let's apply this 
mapping to the basis vector $\bold e_s$. Then
$$
\sum^n_{i=1}\sum^m_{j=1}\gamma^j_i\cdot E^i_j(\bold e_s)=
\sum^n_{i=1}\sum^m_{j=1}(\gamma^j_i\,\delta^i_s)\cdot\bold h_j
=\bold 0.
$$
The sum in the index $i$ can be calculated explicitly. As a result
we get the linear combinations of basis vectors in $W$, which are
equal to zero:
$$
\sum^m_{j=1}\gamma^j_i\cdot\bold h_j=\bold 0.
$$
Due to the linear independence of the vectors $\bold h_1,\,\ldots,\,
\bold h_m$ we derive $\gamma^j_i=0$. This means that the linear 
combination \thetag{10.6} is necessarily trivial. Hence, the mappings
\thetag{10.3} are linearly independent. They form a basis in 
$\Hom(V,W)$. Now, by counting these mappings we find that the
required formula \thetag{10.2} is valid.
\qed\enddemo
     The meaning of the above theorem becomes transparent in terms of
the matrices of linear mappings. Indeed, upon choosing the bases in
$V$ and $W$ the linear mappings from $\Hom(V,W)$ are represented by
rectangular $m\times n$ matrices. The sum of mappings corresponds to
the sum of matrices, and the product of a mapping by a number 
corresponds to the product of the matrix by that number. Note that
rectangular $m\times n$ matrices form a linear vector space isomorphic
to the arithmetic linear vector space $\Bbb K^{mn}$. This space is
denoted as $\Bbb K^{m\times n}$. So, the choice of bases in $V$ and
$W$ defines an isomorphism of $\Hom(V,W)$ and $\Bbb K^{m\times n}$.
\par
\newpage
\topmatter
\title\chapter{2}
Linear operators.
\endtitle
\endtopmatter
\document
\head
\S\,1. Linear operators. The algebra of endomorphisms $\End(V)$
and the group of automorphisms $\Aut(V)$.
\endhead
\rightheadtext{\S\,1. Endomorphisms and automorphisms.}
\leftheadtext{CHAPTER~\uppercase\expandafter{\romannumeral 2}.
LINEAR OPERATORS.}
\setfirstpage
     A linear mapping $f\!:\,V\to V$ acting from a linear vector space
$V$ to the same vector space $V$ is called a {\it linear 
operator\footnote{\ This terminology is not common, however, in this
book we strictly follow this terminology.}}. 
\adjustfootnotemark{-1}Linear operators are
special forms of linear mappings. Therefore, we can apply to them all
results of previous chapter. However, the less generality the more
specific features. Therefore, the theory of linear operators appears
to be more rich and more complicated than the theory of linear 
mappings. It contains not only the strengthening of previous theorems
for this particular case, but a class of problems that cannot be
formulated for the case of general linear mappings.\par
     Let's consider the space of homomorphisms $\Hom(V,W)$. If $W=V$,
this space is called the {\it space of endomorphisms\/} $\End(V)=
\Hom(V,V)$. It consists of linear operators $f\!:\,V\to V$ which are
also called {\it endomorphisms\/} of the space $V$. Unlike the space
of homomorphisms $\Hom(V,W)$, the space of endomorphisms $\End(V)$ 
is equipped with the additional binary algebraic operation. Indeed,
if we have two linear operators $f,g\in\End(V)$, we can not only add
them and multiply them by numbers, but we can also construct two
compositions $f\compos g\in\End(V)$ and $g\compos f\in\End(V)$.
\par
\proclaim{Theorem 1.1} Let $\End(V)$ be the space of endomorphisms
of a linear vector space $V$. Here, apart from the axioms 
\therosteritem{1}-\therosteritem{8} of a linear vector space, the
following relationships are fulfilled:\newline
\line{\vbox{\hsize=179pt
\roster
\item[9] $(f+g)\compos h=f\compos h+
         g\,\compos h$;
\item    $(\alpha\cdot f)\,\compos h=\alpha\cdot(f\compos h)$;
\endroster}\hss
\vbox{\hsize=179pt
\roster
\item[11] $f\,\compos (g+h)=f\compos\,g+
          f\compos h$;
\item     $f\compos (\alpha\cdot g)=\alpha\cdot(f\compos g)$;
\endroster}}
\endproclaim
\demo{Proof} Each of the equalities \therosteritem{9}-\therosteritem{12}
is an operator equality. As we know, the equality of two operators means
that these operators yield the same result when applied to an arbitrary
vector $v\in V$:
$$
\allowdisplaybreaks
\align
&\aligned
  ((f+g)\compos h)(\bold v)=(f&+g)(h(\bold v))
  =f(h(\bold v))+g(h(\bold v))=\\
  &=(f\compos\,h)(\bold v)+(g\compos\,h)(\bold v)=
    (f\compos\,h+g\compos\,h)(\bold v)
  \endaligned\\
\vspace{1.7ex}
&\aligned
  ((\alpha\cdot f)\compos h)(\bold v)
  =(\alpha\cdot f)(h(\bold v)&)=
  \alpha\cdot f(h(\bold v))=\\
  &=\alpha\cdot(f\compos\,h)(\bold v)=(\alpha\cdot
  (f\compos\,h))(\bold v)
  \endaligned\\
\vspace{1.7ex}
&\aligned
  &(f\compos (g+h))(\bold v)=f((g+h)(\bold v))
  =f(g(\bold v)+h(\bold v))=\\
  &=f(g(\bold v))+f(h(\bold v))=(f\compos\,g)(\bold v)+
  (f\compos\,h)(\bold v)=(f\compos\,g
  +f\compos\,h)(\bold v)
  \endaligned\\
\vspace{1.7ex}
&\aligned
  (f\compos (\alpha\cdot g))(\bold v)&=f((\alpha
  \cdot g)(\bold v))=
  f(\alpha\cdot g(\bold v))=\\
  &=\alpha\cdot f(g(\bold v))=\alpha\cdot
  (f\compos\,g)(\bold v)=(\alpha\cdot(f\compos\,g))(\bold v)
  \endaligned
\endalign
$$
The above calculations prove the properties
\therosteritem{9}-\therosteritem{12} of the composition of linear
operators.
\qed\enddemo
     Let's fix the operator $h\in\End(V)$ and consider the composition
$f{\ssize\circ}\,h$ as a rule that maps each operator $f$ to the other
operator $g=f{\ssize\circ}\,h$. Then we get a mapping:
$$
R_h:\End(V)\to\End(V).
$$
The first two properties \therosteritem{9} and \therosteritem{10} from
the theorem~1.1 mean that $R_h$ is a linear mapping. This mapping is 
called the {\it right shift by\/} $h$ since it acts as a composition,
where $h$ is placed on the right side. In a similar way we can define
another mapping, which is called the {\it left shift by\/} $h$:
$$
L_h:\End(V)\to\End(V).
$$
It acts according to the rule $L_h(f)=h\,\compos\,f$. This mapping is linear
due to the properties \therosteritem{11} and \therosteritem{12} from the
theorem~1.1.\par
     The operation of composition is an additional binary operation 
in the space of endomorphisms $\End(V)$. The linearity of the mapping
$R_h$ is interpreted as the {\it linearity\/} of this binary operation 
in its {\it first argument}, while the {\it linearity\/} of $L_h$ is 
said to be the linearity of composition in its {\it second argument}.
A binary algebraic operation linear in both arguments is called a
{\it bilinear\/} operation. A situation, where a linear vector space
is equipped with an additional bilinear algebraic operation, is rather
typical.  
\definition{Definition 1.1} A linear vector space $A$ over a numeric
field $\Bbb K$ equipped with a bilinear binary operation of vector
multiplication is called an {\it algebra over the field\/} $\Bbb K$
or simply a $\Bbb K$-{\it algebra}.
\enddefinition
     The operation of multiplication in algebras is usually denoted
by some sign like a dot 
{\tencyr\char '074}$\ssize\bullet${\tencyr\char '076} or a circle
{\tencyr\char '074}$\ssize\circ${\tencyr\char '076}, but very often 
this sign is omitted at all. The algebra $A$ is called a {\it 
commutative algebra} if the multiplication in it is commutative:
$a\,b=b\,a$. Similarly, the algebra $A$ is called an {\it associative
algebra} if the operation of multiplication is associative:
$(a\,b)\,c=a\,(b\,c)$.\par
     From the definition~1.1 and from the theorem~1.1 we conclude
that the linear space $\End(V)$ with the operation of composition
taken for multiplication is an algebra over the same numeric field 
$\Bbb K$ as the initial vector space $V$. This algebra is called 
the {\it algebra of endomorphisms} of a linear vector space $V$.
It is associative due to the theorem~1.6 from Chapter~\uppercase\expandafter{\romannumeral 1}. However, this algebra
is not commutative in general case.\par
     The operation of composition is treated as a multiplication in the
algebra of endomorphisms $\End(V)$. Therefore, it is usually omitted when
written in this context. The multiplication of operator is higher 
priority operation as compared to addition. The priority of operator
multiplication as compared to the multiplication by numbers makes no
difference at all. This follows from the axiom \therosteritem{7} for the
space $\End(V)$ and from the properties \therosteritem{10} and
\therosteritem{12} of the multiplication in $\End(V)$. Now we can consider
positive integer powers of linear operators:
$$
\xalignat 3
&f^2=f\,f,
&&f^3=f^2\,f,
&&f^{n+1}=f^n\,f.
\endxalignat
$$
If an operator $f$ is bijective, then we have the inverse operator 
$f^{-1}$ and we can consider negative inverse powers of $f$ as well:
$$
\xalignat 3
&f^{-2}=f^{-1}\,f^{-1},
&&f^{n}\,f^{-n}=\id_V,
&&f^{n+m}=f^n\,f^m.
\endxalignat
$$
The latter equality is valid either for positive and negative values
of integer constants $n$ and $m$.\par
\definition{Definition 1.2} An algebra $A$ over the field $\Bbb K$
is called an {\it algebra with unit element} or an {\it algebra with
unity} if there is an element $1\in A$ such that $1\cdot a=a$
and $a\cdot 1=a$ for all $a\in A$.
\enddefinition
     The algebra of endomorphisms $\End(V)$ is an algebra with unity.
The identical operator plays the role of unit element in this algebra:
$1=\id_V$. Therefore, this operator is also called the {\it unit 
operator\/} or the {\it operator unity}.
\definition{Definition 1.3} A linear operator $f\!:\,V\to V$ is called
a {\it scalar operator} if it is obtained by multiplying the unit operator
$1$ by a number $\lambda\in\Bbb K$, i\.\,e\. if $f=\lambda\cdot 1$.
\enddefinition
     The basic purpose of operators from the space $\End(V)$ is to 
act upon vectors of the space $V$. Suppose that $a,b\in\End(V)$ and 
let $\bold x,\bold y\in V$. Then
\roster
\item $(a+b)(\bold x)=a(\bold x)+b(\bold x)$;
\vskip 1ex
\item $a(\bold x+\bold y)=a(\bold x)+a(\bold y)$.
\endroster
These two relationships are well known: the first one follows from
the definition of the sum of two operators, the second relationship
follows from the linearity of the operator $a$. The question is
why the vectors $\bold x$ and $\bold y$ in the above formulas are
surrounded by brackets. This is the consequence of 
{\tencyr\char '074}functional{\tencyr\char '076} form of writing the
action of an operator upon a vector: the operator sign is put on the
left and the vector sign is put on the right and is enclosed into
brackets like an argument of a function: $\bold w=f(\bold v)$.
Algebraists use the more {\tencyr\char '074}deliberate{\tencyr\char '076}
form of writing: $\bold w=f\,\bold v$. The operator sign is on the left
and the vector sign on the right, but no brackets are used. If we know
that $f\in\End(V)$ and $\bold v\in V$, then such a writing makes no
confusion. In more complicated case even if we know that $\alpha\in
\Bbb K$, $f,g\in\End(V)$, and $\bold v\in V$, the writing $\bold w
=\alpha\,f\,g\,\bold v$ admits several interpretations:
$$
\xalignat 3
&\bold w=\alpha\cdot f(g(\bold v)),
&&\bold w=(\alpha\cdot f)(g(\bold v)),\\
\vspace{1ex}
&\bold w=(\alpha\cdot(f\compos g))(\bold v),
&&\bold w=((\alpha\cdot f)\compos g)(\bold v).
\endxalignat
$$
However, for any one of these interpretations we get the same vector 
$\bold w$. Therefore, in what follows we shall use the algebraic form
of writing the action of an operator upon a vector, especially in huge
calculations.\par
     Let $f\!:\,V\to V$ be a linear operator in a finite-dimensional
vector space $V$. According to general scheme of constructing the 
matrix of a linear mapping we should choose two bases $\bold e_1,\,
\ldots,\,\bold e_n$ and $\bold h_1,\,\ldots,\,\bold h_n$ in $V$ and
consider the expansions similar to \thetag{9.1} in
Chapter~\uppercase\expandafter{\romannumeral 1}. No doubt that this 
approach is valid, it could be very fruitful in some cases. However,
to have two bases in one space --- it\linebreak is certainly excessive.
Therefore, when constructing the matrix of a linear operator the second
basis $\bold h_1,\,\ldots,\,\bold h_n$ \pagebreak is chosen to be
coinciding with the first one. The matrix $F$ of an operator $f$ is
determined from the expansions
$$
\hskip -2em
\matrix
f(\bold e_1) &= &F^1_1\cdot\bold e_1 &+ &\hdots &+ &F^n_1\cdot\bold e_n,\\
\hdotsfor 7\\ \vspace{1\jot}
f(\bold e_n) &= &F^1_n\cdot\bold e_1 &+ &\hdots &+ &F^n_n\cdot\bold e_n,
\endmatrix
\tag1.1
$$
which can be expressed in brief form by the formula
$$
\hskip -2em
f(\bold e_j)=\sum^n_{i=1}F^i_j\cdot\bold e_i.
\tag1.2
$$
The matrix $F$ determined by the expansions \thetag{1.1} or by the
expansions \thetag{1.2} is called the {\it matrix of a linear operator\/}
$f$ in the basis $\bold e_1,\,\ldots,\,\bold e_n$. This is a square
$n\times n$ matrix, where $n=\dim V$.\par
\proclaim{Theorem 1.2} Matrices related to operators $f\in\End(V)$
in some fixed basis $\bold e_1,\,\ldots,\,\bold e_n$ possess the
following properties:
\roster
\rosteritemwd=5pt
\item the sum of two operators is represented by the sun of their
      matrices;
\item the product of an operator by a number is represented by
      the product of its matrix by that number;
\item the composition of two matrices is represented by the product
      of their matrices.
\endroster
\endproclaim
\demo{Proof} Consider the operators $f$, $g$, and $h$ from $\End(V)$.
Let $F$, $G$, and $H$ be their matrices in the basis $\bold e_1,\,
\ldots,\,\bold e_n$. Proving the first proposition in the theorem~1.2,
let's denote $h=f+g$. Then
$$
\align
h(\bold e_j)&=(f+g)\,\bold e_j=f(\bold e_j)+h(\bold e_j)=\\
&=\sum^n_{i=1} F^i_j\cdot\bold e_i+\sum^n_{i=1} G^i_j\cdot\bold
e_i=\sum^n_{i=1}(F^i_j+G^i_j)\cdot\bold e_i=\sum^n_{i=1}H^i_j
\cdot\bold e_i.
\endalign
$$
Due to the uniqueness of the expansion of a vector in a basis
we have $H^i_j=F^i_j+G^i_j$ and $H=F+G$. The first proposition
of the theorem is proved.\par
     The proof of the second proposition is similar. Let's denote
$f=\alpha\cdot h$. Then
$$
\align
h(\bold e_j)&=(\alpha\cdot f)\,\bold e_j=\alpha\cdot f(\bold e_j)=\\
&=\alpha\cdot\left(\,\shave{\sum^n_{i=1}}F^i_j\cdot\bold e_i\right)=
\sum^n_{i=1}(\alpha\,F^i_j)\cdot\bold e_i=\sum^n_{i=1} H^i_j\cdot\bold e_i.
\endalign
$$
Therefore, $H^i_j=\alpha\,F^i_j$ and $H=\alpha\cdot F$. The proof of 
the third proposition requires a little bit more efforts. Denote $h=
f\compos g$. Then
$$
\gather
h(\bold e_j)=(f\compos\,g)\,\bold e_j=f(g(\bold e_g))=
f\left(\,\shave{\sum^n_{i=1}} G^i_j\,\cdot\bold e_i\right)
=\sum^n_{i=1} G^i_j\cdot f(\bold e_i)=\\
=\sum^n_{i=1} G^i_j\,\cdot\!
\left(\,\shave{\sum^n_{s=1}} F^s_i\cdot\,\bold e_s\right)
=\sum^n_{s=1}\left(\,\shave{\sum^n_{i=1}} F^s_i\,G^i_j
\right)\!\cdot\bold e_s=\sum^n_{s=1}H^s_i\cdot\bold e_s.
\endgather
$$
Due to the uniqueness of the expansion of a vector in a basis
we derive
$$
H^s_i=\sum^n_{i=1} F^s_i\,G^i_j.
$$
The right side of this equality is easily interpreted as the
product of two matrices written in terms of the components 
of these matrices. Therefore, $H=F\,G$. The theorem is proved.
\qed\enddemo
    From the theorem that was proved just above we conclude that 
when relating an operator $f\in\End(V)$ with its matrix we establish
the isomorphism of the algebra $\End(V)$ and the matrix algebra
$\Bbb K^{n\times n}$ with standard matrix multiplication.\par
     Now let's study how the matrix of a linear operator $f\!:\,V\to V$
changes under the change of the basis $\bold e_1,\,\ldots,\,\bold e_n$ 
for some other basis $\tilde\bold e_1,\,\ldots,\,\tilde\bold e_n$. Let
$S$ be the direct transition matrix and let $T$ be the inverse one.
Note that we need not derive the transformation formulas again. We can
adapt the formulas \thetag{9.10} from
Chapter~\uppercase\expandafter{\romannumeral 1} for our present purpose.
Since the basis $\bold h_1,\,\ldots,\,\bold h_n$ coincides with $\bold
e_1,\,\ldots,\,\bold e_n$ and the basis $\tilde\bold h_1,\,\ldots,\,\tilde
\bold h_n$ coincides with $\tilde\bold e_1,\,\ldots,\,\tilde\bold e_n$, 
we have $P=S$. Then transformation formulas are written as
$$
\xalignat 2
&\hskip -2em
\tilde F=S^{-1}\,F\,S,
&&F=S\,\tilde F\,S^{-1}.
\tag1.3
\endxalignat
$$
These are the required formulas for transforming the matrix of a
linear operator under a change of basis. Taking into account that
$T=S^{-1}$ we can write \thetag{1.3} as 
$$
\xalignat 2
&\hskip -2em
\quad\tilde F^q_p=\sum^n_{i=1}\sum^n_{j=1} T^q_i\,S^j_p\,F^i_j,
&&F^i_j=\sum^n_{q=1}\sum^n_{p=1} S^i_q\,T^p_j\,F^q_p.\quad
\tag1.4
\endxalignat
$$\par
     The relationships \thetag{1.3} yield very important formula relating
the determinants of the matrices $F$ and $\tilde F$. Indeed, we have
$$
\det\tilde F=\det(S^{-1})\,\det F\,\det S=
(\det S)^{-1}\,\det F\,\det S=\det F.
$$
The coincidence of determinants of the matrices of a linear operator $f$
in two arbitrary bases mean that they represent a number which does not
depend on a basis at all. 
\definition{Definition 1.4} The {\it determinant\/} $\det f$ of a linear
operator $f$ is the number equal to the determinant of the matrix $F$ of
this linear operator in some basis.
\enddefinition
     A {\it numeric invariant\/} of a geometric object in a linear 
vector space $V$ is a number determined by this geometric object such 
that it does not depend on anything else other than that geometric 
object itself. The determinant of a linear operator $\det f$ is an 
example of such numeric invariant. Coordinates of a vector or components
of the matrix of a linear operator are not numeric invariants. Another
example of a numeric invariant of a linear operator is its rank:
$$
\rank f=\dim(\Img f).
$$
Soon we shall define a lot of other numeric invariants of a linear
\pagebreak operator.\par
     From the third proposition of the theorem~1.2 we derive the 
following formula for the determinant of a linear operator:
$$
\hskip -2em
\det(f\compos g)=\det(f)\cdot\det(g).
\tag1.5
$$\par
\proclaim{Theorem 1.3} A linear operator $f:V\to V$ in a 
finite-dimensional linear vector space $V$ is injective if and only
if it is surjective.
\endproclaim
\demo{Proof} In order to prove this theorem we apply the 
theorem~1.2 and two theorems 8.3 and 9.4 from
Chapter~\uppercase\expandafter{\romannumeral 1}. The injectivity
of the linear operator $f$ is equivalent to the condition $\Ker f
=\{\bold 0\}$, the surjectivity of the operator $f$ is equivalent to 
$\Img f=V$, while the theorem~9.4 from
Chapter~\uppercase\expandafter{\romannumeral 1} relates the 
dimensions of these two subspaces $\Ker f$ and $\Img f$:
$$
\dim(\Ker f)+\dim(\Img f)=\dim(V).
$$\par
     If the operator $f$ is injective, then $\Ker f=\{\bold 0\}$ and
$\dim(\Ker f)=0$. Then $\dim(\Img f)=\dim(V)$. Applying the
third proposition of the theorem~4.5 from
Chapter~\uppercase\expandafter{\romannumeral 1}, we get $\Img f=V$, 
which proves the surjectivity of the operator $f$.\par
     Conversely, if the operator $f$ is surjective, then $\Img f=V$
and $\dim(\Img f)=\dim(V)$. Hence, $\dim(\Ker f)=0$ and $\Ker f=\{\bold
0\}$. This proves the injectivity of the operator $f$.
\qed\enddemo
\proclaim{Theorem 1.4} A linear operator $f\!:\,V\to V$ in a
finite-dimensional linear vector space $V$ is bijective if and
only if \ $\det f\neq 0$.
\endproclaim
\demo{Proof} Let $\bold x$ be a vector of $V$ and let $\bold y
=f(\bold x)$. Expanding $\bold x$ and $\bold y$ in some basis 
$\bold e_1,\,\ldots,\,\bold e_n$, we get the following formula 
relating their coordinates:
$$
\hskip -2em
\Vmatrix y^1\\ \vdots\\ y^n\endVmatrix=
\Vmatrix F^1_1 & \hdots & F^1_n \\
        \vdots & \ddots & \vdots\\
         F^n_1 & \hdots & F^n_n
\endVmatrix\cdot
\Vmatrix \vphantom{y^1}x^1\\ \vdots\\ x^n\endVmatrix.
\tag1.6
$$
The formula \thetag{1.6} can be derived independently or one 
can derive it from the formula \thetag{9.5} of
Chapter~\uppercase\expandafter{\romannumeral 1}. From this
formula we derive that $\bold x$ belong to the kernel of the 
operator $f$ if and only if its coordinates $x^1,\,\ldots,\,x^n$
satisfy the homogeneous system of linear equations
$$
\hskip -2em
\Vmatrix F^1_1 & \hdots & F^1_n \\
        \vdots & \ddots & \vdots\\
         F^n_1 & \hdots & F^n_n
\endVmatrix\cdot
\Vmatrix \vphantom{y^1}x^1\\ \vdots\\ x^n\endVmatrix=
\Vmatrix \vphantom{y^1}0\\ \vdots\\ 0\endVmatrix,
\tag{1.7}
$$
The matrix of this system of equations coincides with the matrix
of the operator $f$ in the basis $\bold e_1,\,\ldots,\,\bold e_n$.
Therefore, the kernel of the operator $f$ is nonzero if and only if
the system of equations \thetag{1.7} has nonzero solution. Here we
use the well-known result from the theory of determinants: a 
homogeneous system of linear algebraic equations with square matrix 
$F$ has nonzero solution if and only if $\det F=0$. The proof of
this fact can be found in \cite{5}. From this result immediately
get that the condition $\Ker f\neq\{\bold 0\}$ is equivalent to 
$\Ker f=\{\bold 0\}$. Due to the previous theorem and due to the
theorem~1.1 from Chapter~\uppercase\expandafter{\romannumeral 1}
the latter equality $\Ker f\neq\{\bold 0\}$ is equivalent to bijectivity 
of $f$. \pagebreak The theorem is proved.
\qed\enddemo
     An operator $f$ with zero determinant $det f=0$ is called a
{\it degenerate operator}. Using this terminology we can formulate
the following corollary of the theorem~1.4.
\proclaim{Corollary} A linear operator $f\!:\,V\to V$ in a
finite-dimensional space $V$ has a nontrivial kernel $\Ker f\neq\{\bold 
0\}$ if and only if it is degenerate. Otherwise this linear operator is
bijective.
\endproclaim
     Remember that a bijective linear mapping $f$ from $V$ to $W$ 
is called an isomorphism. If $W=V$ such a mapping establishes an
isomorphism of the space $V$ with itself. Therefore, it is called
an {\it automorphism\/} of the space $V$. The set of all
automorphisms of the space $V$ is denoted by $\Aut(V)$. It is
obvious that $\Aut(V)$ possesses the following properties:
\roster
\item if $f,g\in\Aut(V)$, then $f{\ssize\circ}\,g\in\Aut(V)$;
\item if $f\in\Aut(V)$, then $f^{-1}\in\Aut(V)$;
\item $1\in\Aut(V)$, where $1$ is the identical operator.
\endroster
It is easy to see that due to the above three properties the set 
of automorphisms $\Aut(V)$ is equipped with a structure of a group.
{\it The group of automorphisms} $\Aut(V)$ is a subset in the
algebra of endomorphisms $\End(V)$, however, it does not inherit
the structure of an algebra, nor even the structure of a linear 
vector space. It is clear because, for instance, the zero operator
does not belong to $\Aut(V)$. In the case of finite-dimensional
space $V$ the group of automorphisms consists of all non-degenerate
operators.
\head
\S\,2. Projection operators.
\endhead
     Let $V$ be a linear vector space expanded into a direct sum of 
two subspaces:
$$
\hskip -2em
V=U_1\oplus U_2.
\tag2.1
$$
Due to the expansion \thetag{2.1} each vector $\bold v\in V$ is expanded
into a sum
$$
\hskip -2em
\bold v=\bold u_1+\bold u_2\text{, \ where \ }\bold u_1\in U_1
\text{\ \ and \ }\bold u_2\in U_2,
\tag2.2
$$
the components $\bold u_1$ and $\bold u_2$ in \thetag{2.2} being
uniquely determined by the vector $\bold v$.
\definition{Definition 2.1} The operator $P\!:\,V\to V$ mapping each
vector $v\in V$ to its first component $\bold u_1$ in the expansion
\thetag{2.2} is called the {\it operator of projection} onto the
subspace $U_1$ parallel to the subspace $U_2$.
\enddefinition
\proclaim{Theorem 2.1} For any expansion of the form \thetag{2.1}
the operator of projection onto the subspace $U_1$ parallel to the 
subspace $U_2$ is a linear operator.
\endproclaim
\demo{Proof} Let's consider a pair of vectors $\bold v_1,\bold v_2$ 
from the space $V$, and for each of them consider the expansion
like \thetag{2.2}:
$$
\aligned
\bold v_1&=\bold u_1+\bold u_2,\\
\bold v_2&=\tilde\bold u_1+\tilde\bold u_2.
\endaligned
$$
Then $P(\bold v_1)=\bold u_1$ and $P(\bold v_2)=\tilde\bold u_1$. Let's
add the above two expansions and write
$$
\hskip -2em
\bold v_1+\bold v_2=(\bold u_1+\tilde\bold u_1)+(\bold u_2
+\tilde\bold u_2).
\tag2.3
$$
From $\bold u_1,\tilde\bold u_1\in U_1$ and from $\bold u_2,
\tilde\bold u_2\in U_2$ we derive $\bold u_1+\tilde\bold u_1\in 
U_1$ and $\bold u_2+\tilde\bold u_2\in U_2$. Therefore, \thetag{2.3} 
is an expansion of the form \thetag{2.2} for the vector 
$\bold v_1+\bold v_2$. Then
$$
\hskip -2em
P(\bold v_1+\bold v_2)=\bold u_1+\tilde\bold u_1=P(\bold v_1)
+P(\bold v_2).
\tag2.4
$$\par
     Now let's consider the expansion \thetag{2.2} for an arbitrary
vector $\bold v\in V$ and multiply it by a number $\alpha\in\Bbb K$:
$$
\alpha\cdot\bold v=(\alpha\cdot\bold u_1)+(\alpha\cdot\bold u_2).
$$
Then $\alpha\cdot\bold u_1\in U_1$ and $\alpha\cdot\bold u_2\in U_2$,
therefore, due to the definition of $P$ we get
$$
\hskip -2em
P(\alpha\cdot\bold v)=\alpha\cdot\bold u_1=\alpha\cdot P(\bold v).
\tag2.5
$$
The relationships \thetag{2.4} and \thetag{2.5} are just the very 
relationships that mean the linearity of the operator $P$.
\qed\enddemo
     Suppose that $\bold v$ in the expansion \thetag{2.2} is chosen
to be a vector of the subspace $U_1$. Then the expansion \thetag{2.2} 
for this vector is $\bold v=\bold v+\bold 0$, therefore, $P(\bold v)
=\bold v$. This means that all vectors of the subspace $U_1$ are
projected by $P$ onto themselves. This fact has an important
consequence $P^2=P$. Indeed, for any $\bold v\in V$ we have 
$P(\bold v)\in U_1$, therefore, $P(P(\bold v))=P(\bold v)$.\par
     Besides $P$, by means of \thetag{2.2} we can define the other
operator $Q$ such that $Q(\bold v)=\bold u_2$. It is also a projection
operator: it projects onto $U_2$ parallel to $U_1$. Therefore, $Q^2=Q$. 
For the sum of these two operators we get $P+Q=1$. Indeed, for
any vector $v\in V$ we have
$$
P(\bold v)+Q(\bold v)=\bold u_1+\bold u_2=\bold v
=\id_V(\bold v)=1(\bold v).
$$\par
     If $\bold v\in U_1$, then the expansion \thetag{2.2} for this
vector is $\bold v=\bold +\bold 0$, therefore, $Q(\bold v)=\bold 0$.
Similarly, $P(\bold v)=\bold 0$ for all $\bold v\in U_2$. 
Hence, we derive $Q(P(\bold v))=\bold 0$ and $P(Q(\bold v))=\bold 0$ 
for any $v\in V$. Summarizing these results, we write
$$
\xalignat 2
&\hskip -2em
P^2=P, && P+Q=\bold 1,\\
\vspace{-1.9ex}
&&&\tag2.6\\
\vspace{-1.9ex}
&\hskip -2em
Q^2=Q, && P\,Q=Q\,P=0.
\endxalignat
$$
A pair of projection operators satisfying the relationships
\thetag{2.6} is called a {\it con\-cordant pair of projectors}.
\par
     in order to get a concordant pair of projectors it is sufficient
to define only one of them, for instance, the operator $P$. The second
operator $Q$ then is given by formula $Q=1-P$. All of the relationships
\thetag{2.6} thereby will be automatically fulfilled.
Indeed, we have the relationships
$$
\align
P\,Q&=P\compos (1-P)=P-P^2=P-P=0,\\
Q\,P&=(1-P)\,\compos P=P-P^2=P-P=0.
\endalign
$$
The relationship \ $Q^2=Q$ \ for \ $Q$ \ is derived in a similar way:
$$
Q^2=(1-P)\compos(1-P)=1-2P+P=1-P=Q.
$$
\proclaim{Theorem 2.2} An operator $P\!:\,V\to V$ is a projector
onto a subspace parallel to another subspace if and only if \ $P^2=P$.
\endproclaim
\demo{Proof} We have already shown that any projector satisfies the
equality $P^2=P$. Let's prove the converse proposition. Suppose that
$P^2=P$. Let's denote $Q=1-P$. Then for operators $P$ and $Q$ all of
the relationships \thetag{2.6} are fulfilled. Let's consider two
subspaces
$$
\xalignat 2
&U_1=\Img P, &&U_2=\Ker P.
\endxalignat
$$
For an arbitrary vector $\bold v\in V$ we have the expansion 
$$
\hskip -2em
\bold v=1(\bold v)=(P+Q)\bold v=P(\bold v)+Q(\bold v),
\tag2.7
$$
where $\bold u_1=P(\bold v)\in\Img P$. From the relationship $P\,Q=0$
for the other vector $\bold u_2=Q(\bold v)$ in \thetag{2.7} we get
the equality
$$
P(\bold u_2)=P(Q(\bold v))=\bold 0.
$$
This means $u_2\in\Ker P$. Hence, $V=\Img P+\Ker P$. Let's prove that
this is a direct sum of subspaces. We should prove the uniqueness of the expansion
$$
\hskip -2em
\bold v=\bold u_1+\bold u_2,
\tag2.8
$$
where $\bold u_1\in\Img P$ and $\bold u_2\in\Ker P$. From $\bold u_1\in
\Img P$ we conclude that $\bold u_1=P(\bold v_1)$ for some 
vector $\bold v_1\in V$. From $\bold u_2\in\Ker P$ we derive $P(\bold
u_2)=\bold 0$. Then from \thetag{2.8} we derive the following formulas:
$$
\align
P(\bold v)&=P(\bold u_1)+P(\bold u_2)=P(P(\bold v_1))=P^2(\bold
v_1)=P(\bold v_1)=\bold u_1,\\
Q(\bold v)&=(1-P)\,\bold v=\bold v-P(\bold v)=\bold v-\bold u_1=\bold u_2.
\endalign
$$
The relationships derives just above mean that any expansion \thetag{2.8}
coincides with \thetag{2.7}. Hence, it is unique and we have
$$
V=\Img P\oplus\Ker P.
$$
The operator $P$ maps an arbitrary vector $v\in V$ into the first component
of the expansion \thetag{2.8}. Hence, $P$ is an operator of projection
onto the subspace $\Img P$ parallel to the subspace $\Ker P$.
\qed\enddemo
     Now suppose that a linear vector space $V$ is expanded into the 
direct sum of several its subspaces $U_1,\,\ldots,\,U_s$:
$$
\hskip -2em
V=U_1\oplus\ldots\oplus U_s.
\tag2.9
$$
This expansion of the space $V$ implies the unique expansion for each
vector $\bold v\in V$:
$$
\hskip -2em
\bold v=\bold u_1+\ldots+\bold u_s\text{, \ where \ }
\bold u_i\in U_i
\tag2.10
$$
\definition{Definition 2.2} The operator $P_i\!:\,V\to V$ that maps
each vector $\bold v\in V$ to its $i$-th component $\bold u_i$ in
the expansion \thetag{2.10} is called the {\it operator of projection}
onto $U_i$ parallel to other subspaces.
\enddefinition
     The proof of linearity of the operators $P_i$ is practically the same as in case of two subspaces considered in theorem~2.1. It is based on the
uniqueness of the expansion \thetag{2.10}.\par
     Let's choose a vector $\bold u\in U_i$. Then its expansion
\thetag{2.10} looks like:
$$
\bold u=\bold 0+\ldots+\bold 0+\bold u+\bold 0+\ldots+\bold 0.
$$
Therefore for any such vector $\bold u$ we have $P_i(\bold u)=\bold u$ and
$P_j(\bold u)=\bold 0$ for $j\neq i$. For the projection operators $P_i$
this yields 
$$
\xalignat 2
&\hskip -2em
(P_i)^2=P_i, 
&& P_i\,\compos P_j=0\text{\ \ for \ }i\neq j.
\tag2.11
\endxalignat
$$
Moreover, from the definition of $P_i$ we get 
$$
\hskip -2em
P_1+\ldots+P_s=1.
\tag2.12
$$\par
     Due to the first relationship \thetag{2.11} the theory of separate
operators $P_i$ does not differ from the theory of projectors defined by 
two component expansions of the space $V$. In the case of multicomponent
expansions the collective behavior of projectors is of particular
interest. A family of projection operators $P_1,\,\ldots,\,P_s$ is called 
a {\it concordant family of projectors} if the operators of this family
satisfy the relationships \thetag{2.11} and \thetag{2.12}.\par
\proclaim{Theorem 2.3} A family of projection operators $P_1,\,\ldots,\,
P_s$ is determined by an expansion of the form \thetag{2.9} if and only
if it is concordant, i\.\,e\. if these operators satisfy the relationships 
\thetag{2.11} and \thetag{2.12}.
\endproclaim
\demo{Proof} We already know that a family of projectors determined by
an expansion \thetag{2.9} satisfy the relationships \thetag{2.11} and \thetag{2.12}. Let's prove the converse proposition. Suppose that we have
a family of operators $P_1,\,\ldots,\,P_s$ satisfying the relationships \thetag{2.11} and \thetag{2.12}. Then we define the subspaces $U_i=\Img P_i$. Due to the relationship \thetag{2.12} for an arbitrary vector 
$\bold v\in V$ we get
$$
\hskip -2em
\bold v=P_1(\bold v)+\ldots+P_s(\bold v),
\tag2.13
$$
where $P_i(\bold v)\in\Img P_i$. Hence, we have the expansion of $V$ 
into a sum of subspaces 
$$
\hskip -2em
V=\Img P_1+\ldots+\Img P_s.
\tag2.14
$$
Let's prove that the sum \thetag{2.14} is a direct sum. For this
purpose we consider an expansion of some arbitrary vector $\bold v
\in V$ corresponding to the expansion \thetag{2.14}:
$$
\pagebreak
\hskip -2em
\bold v=\bold u_1+\ldots+\bold u_s\text{, \ where \ }
\bold u_i\in\Img P_i
\tag2.15
$$
From $\bold u_i\in\Img P_i$ we conclude that $\bold u_i=P(\bold v_i)$,
where $\bold v_i\in V$. Then from the expansion \thetag{2.15} we derive
the following equality:
$$
P_i(\bold v)=P_i(\bold u_1+\ldots+\bold u_s)=\sum^s_{j=1}
P_i(P_j(\bold v_j)).
$$
Due to \thetag{2.11} only one term in the above sum is nonzero.
Therefore, we have 
$$
P_i(\bold v)=(P_i)^2\,\bold v_i=P_i(\bold v_i)=\bold u_i.
$$
This equality show that an arbitrary expansion \thetag{2.15} 
should coincide with \thetag{2.13}. This means that \thetag{2.13}
is the unique expansion of that sort. Hence, the sum \thetag{2.14}
is a direct sum and $P_i$ is the projection operator onto the
$i$-th component of the sum \thetag{2.14} parallel to its other
components. The theorem is proved.
\qed\enddemo
     Now we consider a projection operator $P$ as an example for the
first approach to the problem of bringing the matrix of a linear
operator to a canonic form.
\proclaim{Theorem 2.4} For any nonzero projection operator in a
finite-dimensional\linebreak vector space $V$ there is a basis $\bold
e_1,\,\ldots,\,\bold e_n$ such that the matrix of the operator $P$ has 
the following form in that basis:
$$
\hskip -2em
\Cal P=
\aligned
&\hphantom{x}\overbrace{\hphantom{xxxxxxxxxxx}}^s\\
&\Vmatrix
1 & 0 & \hdots & 0 & 0 & \hdots & 0 \\
0 & 1 & \hdots & 0 & 0 & \hdots & 0 \\
\vdots & \vdots & \ddots & \vdots &\vdots &\, & \vdots \\
0 & 0 & \hdots & 1 & 0 & \hdots & 0 \\
0 & 0 & \hdots & 0 & 0 & \hdots & 0 \\
\vdots & \vdots &\,  & \vdots &\vdots &\ddots & \vdots \\
0 & 0 & \hdots & 0 & 0 & \hdots & 0 \\
\endVmatrix
\endaligned
\aligned
&\vphantom{\vrule height 1ex depth 1.7ex}\\
&\left.\vphantom{\vrule height 5ex depth 5ex}\right\} s\\
&\vphantom{\vrule height 4ex depth 4ex}
\endaligned\lower 10pt\hbox {\kern -5pt .}
\tag2.16
$$
\endproclaim
\demo{Proof} Let's consider the subspaces $\Img P$ and $\Ker P$.
From the condition $P\neq 0$ we conclude that $s=\dim(\Img P)\neq 0$. 
Then we choose a basis $\bold e_1,\,\ldots,\,\bold e_s$ in $U_1=\Img P$
and if $U_1\neq V$, we complete it by choosing a basis in $\bold e_{s+1},
\,\ldots,\,\bold e_n$ in $U_2=\Ker P$. The sum of these two subspaces is
a direct sum: $V=U_1\oplus U_2$, therefore, joining together two bases
in them, we get a basis of $V$ (see the proof of theorem~6.3 in
Chapter~\uppercase\expandafter{\romannumeral 1}).\par
     Now let's apply the operator $P$ to the vectors of the basis we 
have constructed just above. This operator projects onto $U_1$ parallel 
to $U_2$, therefore, we have
$$
P(\bold e_i)=
\cases
\bold e_i&\text{ for \ }i=1,\ldots,s,\\
\bold 0  &\text{ for \ }i=s+1,\ldots,n.
\endcases
$$
Due to this formula it's clear that \thetag{2.16} is the matrix of the
projection operator $P$ the basis $\bold e_1,\,\ldots,\,\bold e_s$.
\qed\enddemo
\head
\S\,3. Invariant subspaces.\\
Restriction and factorization of operators.
\endhead
\rightheadtext{\S\,3. Invariant subspaces.}
     Let $f\!:\,V\to V$ be a linear operator and let $U$ be a
subspace of $V$. Let's restrict the domain of $f$ to the subspace
$U$. Thereby the image of $f$ shrinks to $f(U)$. However, in general, 
the subspace $f(U)$ is not enclosed into the subspace $U$. For this
reason in general case we should treat the {\it restricted} operator
$f$ as a linear mapping $f{\!\lower4pt\hbox{$\ssize U$}}\!:U\,\to V$,
rather than a linear operator.
\definition{Definition 3.1} A subspace $U$ is called an {\it invariant
subspace} of a linear operator $f\,:\,V\to V$ if $f(U)\subseteq U$, 
i\.\,e\. if \ $\bold u\in U$ \ implies \ $f(\bold u)\in U$.
\enddefinition
     If $U$ is an invariant subspace of $f$, the restriction $f{\!\lower4pt\hbox{$\ssize U$}}$ can be treated as a linear operator
in $U$. Its action upon vectors $\bold u\in U$ coincides with the
action of $f$ upon $u$. As for the vectors outside the subspace $U$,
the operator $f{\!\lower4pt\hbox{$\ssize U$}}$ cannot be applied to
them at all.
\proclaim{Theorem 3.1} The kernel and the image of a linear operator
$f\!:\,V\to V$ are invariant subspaces of $f$.
\endproclaim
\demo{Proof} Let's consider the kernel of $f$ for the first. If
$\bold u\in\Ker f$, then $f(\bold u)=\bold 0$. Hence, $f(\bold u)\in U$,
since the zero vector $\bold 0$ is an element of any subspace of $V$. 
The invariance of the kernel $Ker f$ is proved.\par
     Now let $\bold u\in\Img f$. Denote $\bold w=f(\bold u)$. 
Then $\bold w$ is the image of the vector $\bold u$, hence, $\bold w
=f(\bold u)\in\Img f$. The invariance of the image $\Img f$ is proved.
\qed\enddemo
\proclaim{Theorem 3.2} The intersection and the sum of an arbitrary 
number of invariant subspaces of a linear operator $f\!:\,V\to V$ 
both are the invariant subspaces of $f$.
\endproclaim
\demo{Proof} Let $U_i$, $i\in I$ be a family of invariant subspaces
of a linear operator $f\!:\,V\to V$. Let's consider the intersection
and the sum of these subspaces:
$$
\xalignat 2
&U=\bigcap_{i\in I} U_i,
&
&W=\sum_{i\in I} U_i.
\endxalignat
$$
In \S\,6 of Chapter~\uppercase\expandafter{\romannumeral 1} we have
proved that $U$ and $W$ are the subspaces of $V$. Now we should prove 
that they are invariant subspaces. For the first, let's prove that
$U$ is an invariant subspace. Consider a vector $\bold u\in U$. This
vector belongs to all subspaces $U_i$, which are invariant subspaces
of $f$. Therefore, $f(\bold u)$ also belongs to all subspaces $U_i$.
This means that $f(\bold u)$ belongs to their intersection $U$. 
The invariance of $U$ is proved.\par
     Now let's consider a vector $\bold w\in W$. According to the
definition of the sum of subspaces, this vector admits the expansion
$$
\bold w=\bold u_{i_1}+\ldots+\bold u_{i_s},\text{, where\ }
\bold u_{i_r}\in U_{i_r}.
$$
Applying the operator $f$ to both sides of this equality, we get:
$$
f(\bold w)=f(\bold u_{i_1})+\ldots+f(\bold u_{i_s}).
$$
Due to the invariance of $U_i$ we have $f(u_{i_r})\in U_{i_r}$. 
Hence, $f(w)\in W$. This yields the invariance of the sum $W$ of 
the invariant subspaces $U_i$.
\qed\enddemo
     Let $U$ be an invariant subspace of a linear operator $f\!:\,V\to
V$. Let's consider the factorspace $V/U$ and define the operator
$f{\!\lower4pt\hbox{$\ssize V/U$}}$ in this factorspace by formula
$$
\hskip -2em
f{\!\lower4pt\hbox{$\ssize V/U$}}(Q)=\Cl_U(f(\bold v))\text{, where \ }
Q=\Cl_U(\bold v).
\tag3.1
$$
The operator $f{\!\lower4pt\hbox{$\ssize V/U$}}\!:\,V/U\to V/U$ acting 
according to the rule \thetag{3.1} is called the {\it factoroperator}
of the {\it quotient operator} of the operator $f$ by the subspace
$U$. We can rewrite the formula \thetag{3.1} in shorter form as follows:
$$
\hskip -2em
f{\!\lower4pt\hbox{$\ssize V/U$}}(\Cl_U(\bold v))=\Cl_U(f(\bold v)).
\tag3.2
$$
Like formulas \thetag{7.3} in 
Chapter~\uppercase\expandafter{\romannumeral 1}, the formulas
\thetag{3.1} and \thetag{3.2} comprise the definite amount of
uncertainty due to the uncertainty of the choice of a representative
$\bold v$ in a coset $Q=\Cl_U(\bold v)$. Therefore, we need to prove
their correctness.
\proclaim{Theorem 3.3} The formula \thetag{3.1} and the equivalent 
formula \thetag{3.2} both are correct. They define a linear operator 
$f{\!\lower4pt\hbox{$\ssize V/U$}}$ in factorspace $V/U$.
\endproclaim
\demo{Proof} Let's conside two different representative vectors in a
coset $Q$, i\.\,e\. let $\bold v,\tilde\bold v\in Q$. Then $\tilde\bold
v-\bold v\in U$. According to the formula \thetag{3.1}, we consider two
possible results of applying the operator $f{\!\lower4pt\hbox{$\ssize
V/U$}}$ to $Q$:
$$
\xalignat 2
&f{\lower4pt\hbox{$\ssize V/U$}}(Q)=\Cl_U(f(\bold v)),
&&f{\lower4pt\hbox{$\ssize V/U$}}(Q)=\Cl_U(f(\tilde\bold v)).
\endxalignat
$$
Let's calculate the difference of these two possible results:
$$
\Cl_U(f(\tilde\bold v))-\Cl_U(f(\bold v))
=\Cl_U(f(\tilde\bold v)-f(\bold v))
=\Cl_U(f(\tilde\bold v-\bold v)).
$$
Note that the vector $\bold u=\tilde\bold v-\bold v$ belongs to the
subspace $U$. Since $U$ is an invariant subspace, we have $\tilde\bold
u=f(\bold u)\in U$. Therefore, we get
$$
\Cl_U(f(\tilde\bold v))-\Cl_U(f(\bold v))
=\Cl_U(\tilde\bold u)=\bold 0.
$$
This coincidence $\Cl_U(f(\tilde\bold v))=\Cl_U(f(\bold v))$ that
we have proved just above proves the correctness of the formula 
\thetag{3.1} and the formula \thetag{3.2} as well.\par
     Now let's prove the linearity of the factoroperator $f{\!\lower4pt\hbox{$\ssize V/U$}}\!:\,V/U\to V/U$. We shall carry out
the appropriate calculations on the base of formula \thetag{3.1}:
$$
\pagebreak
\align
&\aligned
  f{\!\lower4pt\hbox{$\ssize V/U$}}(Q_1&+Q_2)=
  f{\!\lower4pt\hbox{$\ssize V/U$}}(\Cl_U(\bold v_1)
  +\Cl_U(\bold v_2))=\\
  =f&{\!\lower4pt\hbox{$\ssize V/U$}}(\Cl_U(\bold v_1
  +\bold v_2))=\Cl_U(f(\bold v_1+\bold v_2))=\\
  &=\Cl_U(f(\bold v_1))+\Cl_U(f(\bold v_2))
  =f{\!\lower4pt\hbox{$\ssize V/U$}}(Q_1)+
  f{\!\lower4pt\hbox{$\ssize V/U$}}(Q_2),
  \endaligned\\
\vspace{1.7ex}
&\aligned
  f{\!\lower4pt\hbox{$\ssize V/U$}}(\alpha\cdot Q)&=
  f{\!\lower4pt\hbox{$\ssize V/U$}}(\alpha\cdot\Cl_U(\bold v))=\\
  =f&{\!\lower4pt\hbox{$\ssize V/U$}}(\Cl_U(\alpha\cdot\bold v))=
  \Cl_U(f(\alpha\cdot\bold v))=\\
  &=\Cl_U(\alpha\cdot f(\bold v))=\alpha\cdot\Cl_U(f(\bold v))=
  \alpha\cdot f{\!\lower4pt\hbox{$\ssize V/U$}}(Q).
  \endaligned
\endalign
$$
These calculations show that $f{\!\lower4pt\hbox{$\ssize V/U$}}$
is a linear operator. The theorem is proved.
\qed\enddemo
\proclaim{Theorem 3.4} Suppose that $U$ is a common invariant
subspace of two linear operators $f,g\in\End(V)$. Then $U$ is
an invariant subspace of the operators $f+g$, $\alpha\cdot f$
and $f\compos g$ as well. For their restrictions to the subspace
$U$ and for the corresponding factoroperators we have the following
relationships:
\vskip -4ex
$$
\xalignat 2
&(f+g){\lower4pt\hbox{$\ssize U$}}=f{\!\lower4pt\hbox{$\ssize U$}}+
 g{\lower4pt\hbox{$\ssize U$}}\,;&
&(f+g){\lower4pt\hbox{$\ssize V/U$}}=f{\!\lower4pt\hbox{$\ssize V/U$}}+
 g{\lower4pt\hbox{$\ssize V/U$}}\,;\\
&(\alpha\,\cdot f){\lower4pt\hbox{$\ssize U$}}=\alpha\,\cdot
 f{\!\lower4pt\hbox{$\ssize U$}}\,;&
&(\alpha\,\cdot f){\lower4pt\hbox{$\ssize V/U$}}=\alpha\,\cdot
 f{\!\lower4pt\hbox{$\ssize V/U$}}\,;\\
&(f\compos g){\lower4pt\hbox{$\ssize U$}}=f{\!\lower4pt
 \hbox{$\ssize U$}}\compos g{\lower4pt\hbox{$\ssize U$}}\,;&
&(f\compos g){\lower4pt\hbox{$\ssize V/U$}}=f{\!\lower4pt
 \hbox{$\ssize V/U$}}\compos g{\lower4pt\hbox{$\ssize V/U$}}\,.
\endxalignat
$$
\endproclaim
\demo{Proof} Let's begin with the first case. Denote $h=f+g$ and
assume that $\bold u$ is an arbitrary vector of $U$. Then $f(\bold 
u)\in U$ and $g(\bold u)\in U$ since $U$ is an invariant subspace
of both operators $f$ and $g$. For this reason we obtain
$h(\bold u)=f(\bold u)+g(\bold u)\in U$. This proves that $U$ is
an invariant subspace of $h$. The relationship $h{\lower4pt
\hbox{$\ssize U$}}=f{\!\lower4pt\hbox{$\ssize U$}}+
g{\lower4pt\hbox{$\ssize U$}}$ follows from $h=f+g$ since 
the results of applying the restricted operators to $\bold u$ do 
not differ from the results of applying $f$, $g$, and $h$ to
$\bold u$. The corresponding relationship for the factoroperators 
is proved as follows:
$$
\align
h{\lower4pt\hbox{$\ssize V/U$}}&(\Cl_U(\bold v))=
\Cl_U(h(\bold v))=\Cl_U(f(\bold v)+h(\bold v))=\\
=&\Cl_U(f(\bold v))+\Cl_U(g(\bold v))=f{\!\lower4pt
\hbox{$\ssize V/U$}}(\Cl_U(\bold v))+
g{\lower4pt\hbox{$\ssize V/U$}}(\Cl_U(\bold v)).
\endalign
$$\par
     The second case, where we denote $h=\alpha\cdot f$, is not
quite different from the first one. From $\bold u\in U$ it follows
that $f(u)\in U$, hence, $h(\bold u)=\alpha\cdot f(\bold u)\in U$.
The relationship $h{\lower4pt\hbox{$\ssize U$}}=\alpha\cdot
f{\!\lower4pt\hbox{$\ssize U$}}$ now is obvious due to the same 
reasons as above. For the factoroperators we perform the following
calculations:
$$
\align
h{\lower4pt\hbox{$\ssize V/U$}}&(\Cl_U(\bold v))=
\Cl_U(h(\bold v))=\Cl_U(\alpha\cdot f(\bold v))=\\
=\alpha&\cdot\Cl_U(f(\bold v))=\alpha\cdot f{\!\lower4pt
\hbox{$\ssize V/U$}}(\Cl_U(\bold v))=(\alpha\cdot
f{\!\lower4pt\hbox{$\ssize V/U$}})(\Cl_U(\bold v)).
\endalign
$$\par
     Now we consider the third case. Here we denote $h=f\compos g$. 
From $\bold u\in U$ we derive $\bold w=g(\bold u)\in U$, then
from $\bold w\in U$ we derive $f(\bold w)\in U$, which means that
$U$ is an invariant subspace of $h$. Indeed, $h(\bold u)=f(g(\bold u))=f(\bold w)\in U$. For the restricted operators this yields the
equality
$$
h{\lower4pt\hbox{$\ssize U$}}(\bold u)=h(\bold u)=f(g(\bold u))=
f{\!\lower4pt\hbox{$\ssize U$}}(g{\lower4pt\hbox{$\ssize U$}}
(\bold u)).
$$
Hence, $h{\lower4pt\hbox{$\ssize U$}}=f{\!\lower4pt\hbox{$\ssize U$}}
\compos g{\lower4pt\hbox{$\ssize U$}}$\,. Passing to factoroperators,
we obtain
$$
\align
h{\lower4pt\hbox{$\ssize V/U$}}&(\Cl_U(v))=\Cl_U(h(\bold v))
=\Cl_U(f(g(\bold v))=f{\!\lower4pt\hbox{$\ssize V/U$}}
(\Cl(g(\bold v)))=\\
&=f{\!\lower4pt\hbox{$\ssize V/U$}}(g{\lower4pt\hbox{$\ssize
V/U$}}\,(\Cl_U(\bold v)))=f{\!\lower4pt\hbox{$\ssize V/U$}}
\compos g{\lower4pt\hbox{$\ssize V/U$}}(\Cl_U(\bold v)).
\endalign
$$
The above calculations prove the last relationship of the 
theorem~3.4.
\qed\enddemo
\proclaim{Theorem 3.5} Let $V=U_1\oplus\ldots\oplus U_s$ be an
expansion of a linear vector space $V$ into a direct sum of its
subspaces. \pagebreak The subspaces $U_1,\,\ldots,\,U_s$ are 
invariant subspaces of an operator $f\!:\,V\to V$ if and only 
if the projection operators $P_1,\,\ldots,\,P_s$ associated with 
the expansion $V=U_1\oplus\ldots\oplus U_s$ commute with the 
operator $f$, i\.\,e\. if \ $f\compos P_i=P_i\,\compos f$, where 
\ $i=1,\,\ldots,\,s$.
\endproclaim
\demo{Proof} Suppose that all subspaces $U_i$ are invariant under 
the action of the operator $f$. For an arbitrary vector $\bold v
\in V$ we consider the expansion determined by the direct sum
$V=U_1\oplus\ldots\oplus U_s$:
$$
\bold v=\bold u_1+\ldots+\bold u_s.
$$
Here $\bold u_i=P_i(\bold v)\in U_i$. From this expansion we derive
$$
P_i(f(\bold v))=P_i(f(\bold u_1)+\ldots+f(\bold u_s))
=f(\bold u_i)=f(P_i(\bold v))
$$
We used the inclusion $\bold w_j=f(\bold u_j)\in U_j$ that follows 
from the invariance of the subspace $U_j$ under the action of $f$. 
We also used the following properties of projection operators
(they follow from \thetag{2.11} and $U_i=\Img P_i$, \ see \S~2 above):
$$
P_i(\bold w_j)=
\cases\bold w_i &\text{ for \ }j=i,\\
      \bold 0   &\text{ for \ }j\neq i.
\endcases
$$
Since $\bold v$ is an arbitrary vector of the space $V$, from the
above equality $P_i(f(\bold v))=f(P_i(\bold v))$ we derive $f\compos
P_i=P_i\compos f$.\par
     Conversely, suppose that the operator $f$ commute with all 
projection operators $P_1,\,\ldots,\,P_s$ associated with the
expansion $V=U_1\oplus\ldots\oplus U_s$. Let $\bold u$ be an
arbitrary vector of the subspace $U_i$. Then we denote $\bold w
=f(\bold u)$ and for $\bold w$ we derive
$$
P_i(\bold w)=P_i(f(\bold v))=f(P_i(\bold u))=f(\bold u)=\bold w.
$$
Remember that $P_i$ projects onto the subspace $U_i$. Hence, 
$P_i(\bold w)\in U_i$. But due to the above equality we find that
$P_i(\bold w)=\bold w=f(\bold u)\in U_i$. Thus we have shown
that the space $U_i$ is invariant under the action of the operator 
$f$. The theorem is completely proved.
\qed\enddemo
     Let's consider a linear operator $f$ in a finite-dimensional
linear vector space $V$ and possessing an invariant subspace $U$.
Suppose that $\dim V=n$ and $\dim U=s$. Let's choose a basis 
$\bold e_1,\,\ldots,\,\bold e_s$ in $U$ and then, if $s<n$, 
complete this basis up to a basis in $V$. Denote by $\bold e_{s+1},\,
\ldots,\,\bold e_n$ the complementary vectors. For $j\leqslant s$ 
due to the invariance of the subspace $U$ under the action of $f$ 
we have $f(\bold e_j)\in U$. Therefore, in the expansions
of these vectors 
$$
f(\bold e_j)=\sum^s_{i=1}F^i_j\cdot\bold e_i\text{, \ where \ }
j\leqslant s,
$$
the summation index $i$ runs  from $1$ to $s$, but not from $1$ to 
$n$ as it should in general case, where we expand an arbitrary vector
of $V$.  This means that if we construct the matrix of the operator 
$f$ in the basis $\bold e_1,\,\ldots,\,\bold e_n$, this matrix would 
be mounted of blocks with the lower left block in it being zero:
$$
\hskip -2em
F=
\aligned
&\hphantom{x}\overbrace{\hphantom{xxxxxxxxxxxxxx}}^{\dsize s}\\
&\Vmatrix
F^1_1 & F^1_2 & \hdots & F^1_s & F^1_{s+1} & \hdots & F^1_n \\
\vspace{1ex}
F^2_1 & F^2_2 & \hdots & F^2_s & F^2_{s+1} & \hdots & F^2_n \\
\vdots & \vdots & \ddots & \vdots &\vdots &\, & \vdots \\
F^s_1 & F^s_2 & \hdots & F^4_s & F^s_{s+1} & \hdots & F^s_n \\
\vspace{1ex}
0 & 0 & \hdots & 0 & F^{s+1}_{s+1} & \hdots & F^{s+1}_n \\
\vdots & \vdots &\,  & \vdots &\vdots &\ddots & \vdots \\
0 & 0 & \hdots & 0 & F^n_{s+1} & \hdots & F^n_n \\
\endVmatrix
\endaligned
\aligned
&\vphantom{\vrule height 1ex depth 1ex}\\
&\left.\vphantom{\vrule height 6ex depth 6ex}\right\} s\\
&\vphantom{\vrule height 5ex depth 4ex}
\endaligned
\tag3.3
$$
Matrices of this form are called {\it blockwise-triangular} matrices. 
The upper left diagonal block in the matrix \thetag{3.3} coincides
with the matrix of restricted operator $f{\!\lower4pt\hbox{$\ssize
U$}}\!:\,U\to U$ in the invariant subspace $U$.\par
     The lower right diagonal block of the matrix \thetag{3.3} can also 
be interpreted in a special way. In order to find this interpretation
let's consider the cosets of complementary vectors in the basis $\bold e_1,\,\ldots,\,\bold e_n$:
$$
\hskip -2em
\bold E_1=\Cl_U(\bold e_{s+1}),\ \ldots,\ \bold E_{n-s}
=\Cl_U(\bold e_n).
\tag3.4
$$
When proving the theorem~7.6 in 
Chapter~\uppercase\expandafter{\romannumeral 1}, we have found that
these cosets form a basis in the factorspace $V/U$. Applying the
factoroperator $f{\!\lower4pt\hbox{$\ssize V/U$}}$ to 
\thetag{3.4}, we get
$$
\align
f{\!\lower4pt\hbox{$\ssize V/U$}}\,&(\bold E_j)=
f{\!\lower4pt\hbox{$\ssize V/U$}}\,\Cl_U(\bold e_{s+j})=
\Cl_U(f(\bold e_{s+j}))=\\
=&\sum^s_{i=1}F^i_{s+j}\cdot\Cl_U(\bold e_i)+
\sum^n_{i=s+1}F^i_{s+j}\cdot\Cl_U(\bold e_i).
\endalign
$$
The first sum in the above expression is equal to zero since the
vectors $\bold e_1,\,\ldots,\,\bold e_s$ belong to $U$. Then, shifting 
the index $i+s\to i$, we find
$$
f{\!\lower4pt\hbox{$\ssize V/U$}}\,(\bold E_j)=
\sum^{n-s}_{i=1}F^{s+i}_{s+j}\cdot\bold E_i.
$$
Looking at this formula, we see that the matrix of the factoroperator
$f{\!\lower4pt\hbox{$\ssize V/U$}}$ in the basis \thetag{3.4} coincides
with the lower right diagonal block in the matrix \thetag{3.3}.
\proclaim{Theorem 3.6} Let $f\!:\,V\to V$ be a linear operator in a
finite-dimensional space and let $U$ be an invariant subspace of
this operator. Then the determinant of $f$ is equal to the product
of two determinants --- the determinant of the restricted operator
$f{\!\lower4pt\hbox{$\ssize U$}}$ and that of the factoroperator $f{\!\lower4pt\hbox{$\ssize V/U$}}$\,:
$$
\det f=\det(f{\lower4pt\hbox{$\ssize U$}})\cdot
\det(f{\lower4pt\hbox{$\ssize V/U$}}).
$$
\endproclaim
    The proof of this theorem is immediate from the following fact
well-known in the theory of determinants: the determinant of
blockwise-triangular matrix is equal to the product of determinants
of all its diagonal blocks.
\head
\S\,4. Eigenvalues and eigenvectors.
\endhead
     Let $f\!:\,V\to V$ be a linear operator. A nonzero vector 
$\bold v\neq\bold 0$ of the space $V$ is called an {\it eigenvector}
of the operator $f$ if $f\,\bold v=\lambda\cdot\bold v$, where
$\lambda\in\Bbb K$. The number $\lambda$ is called the {\it 
eigenvalue} of the operator $f$ associated with the eigenvector
$\bold v$.\par
     One eigenvalue $\lambda$ of an operator $f$ can be associated
with several or even with infinite number of eigenvectors. But
conversely, if an eigenvector is given, the associated eigenvalue
$\lambda$ for this eigenvector is unique. Indeed, from the equality
$f\,\bold v=\lambda\cdot\bold v=\lambda'\cdot\bold v$ and from 
$\bold v\neq \bold 0$ it follows that $\lambda=\lambda'$.\par
     Let $\bold v$ be an eigenvector of the operator $f\!:\,V\to V$.
Let's consider the other operator $h_\lambda=f-\lambda\cdot 1$. Then the equation $f\,\bold v=\lambda\cdot\bold v$ can be rewritten as
$$
\hskip -2em
(f-\lambda\cdot 1)\,\bold v=\bold 0.
\tag4.1
$$
Hence, $\bold v\in\Ker(f-\lambda\cdot 1)$. The condition $\bold v\neq
\bold 0$ means that the kernel of this operator is nonzero: $\Ker(f-
\lambda\cdot 1)\neq\{\bold 0\}$.
\definition{Definition 4.1} A number $\lambda\in\Bbb K$ is called
an {\it eigenvalue } of a linear operator $f\!:\,V\to V$ if the
subspace $V_\lambda=\Ker(f-\lambda\cdot 1)$ is nonzero. This subspace
$V_\lambda=\Ker(f-\lambda\cdot 1)\neq\{\bold 0\}$ is called the {\it
eigenspace\/} associated with the eigenvalue $\lambda$, while any nonzero
vector of $V_\lambda$ is called an {\it eigenvector\/} of the operator
$f$ associated with the eigenvalue $\lambda$.
\enddefinition
     The collection of all eigenvalues of an operator $f$ is sometimes
called the {\it spectrum} of this operator, while the brunch of mathematics
studying the spectra of linear operators is known as the {\it spectral 
theory of operators}. The spectral theory of linear operators in 
finite-dimensional spaces is the most simple one. This is the very theory
that is usually studied in the course of linear algebra.\par
     Let $f\!:\,V\to V$ be a linear operator in a finite-dimensional
linear vector space $V$. In order to find the spectrum of this operator
we apply the corollary of theorem~1.4. Due to this corollary a number
$\lambda\in\Bbb K$ is an eigenvalue of the operator $f$ if and only if
it satisfies the equation
$$
\hskip -2em
\det(f-\lambda\cdot 1)=0.
\tag4.2
$$
The equation \thetag{4.2} is called the {\it characteristic equation}
of the operator $f$, its roots are called the {\it characteristic
numbers} of the operator $f$.\par
     Let $\dim V=n$. Then the determinant in formula \thetag{4.2} is
equal to the determinant of the square $n\times n$ matrix. The matrix
of the operator $h_\lambda=f-\lambda\cdot 1$ is derived from the matrix
of the operator $f$ by subtracting $\lambda$ from each element on the 
primary diagonal of this matrix:
$$
\pagebreak
\hskip -2em
H_\lambda=
\Vmatrix
F^1_1-\lambda & F^1_2 & \hdots & F^1_n\\
\vspace{2ex}
F^2_1 & F^2_2-\lambda & \hdots & F^2_n\\
\vspace{1.7ex}
\vdots & \vdots & \ddots & \vdots\\
\vspace{1.7ex}
F^n_1 & F^n_2 & \hdots & F^n_n-\lambda
\endVmatrix.
\tag4.3
$$
The determinant of the matrix \thetag{4.3} is a polynomial of
$\lambda$: 
$$
\hskip -2em
\det(f-\lambda\cdot 1)=(-\lambda)^n+F_1\,(-\lambda)^{n-1}+\ldots+F_n.
\tag4.4
$$
The polynomial in right hand side of \thetag{4.4} is called the
{\it characteristic polynomial} of the operator $f$. If $F$ is the
the matrix of the operator $f$ in some basis, then the coefficients
$F_1,\,\ldots,\,F_n$ of characteristic polynomial \thetag{4.4} are
expressed through the elements of the matrix $F$. However, note that
left hand side of of \thetag{4.4} is basis independent, therefore, 
the coefficients $F_1,\,\ldots,\,F_n$ do not actually depend on the
choice of basis. They are scalar invariants of the operator $f$.
The fires and the last invariants in \thetag{4.4} are the most 
popular ones:
$$
\xalignat 2
&F_1=\tr f, &&F_n=\det f.
\endxalignat
$$
The invariant $F_1$ is called the {\it trace} of the operator $f$.
It is calculated through the matrix of this operator according to
the following formula:
$$
\hskip -2em
\tr f=\sum^n_{i=1} F^i_i.
\tag4.5
$$\par
     We shall not derive this formula \thetag{4.5} since it is well-known
in the theory of determinants. We shall only derive the invariance of the
trace immediately on the base of formula \thetag{1.4} which describes the transformation of the matrix of a linear operator under a change of basis:
$$
\sum^n_{p=1}\tilde F^p_p=\sum^n_{i=1}\sum^n_{j=1}\left(\,\shave{
\sum^n_{p=1}} T^p_i\,S^j_p\right)F^i_j=\sum^n_{i=1}\sum^n_{j=1}
\delta^j_i\,F^i_j=\sum^n_{i=1}F^i_i.
$$\par
     Upon substituting \thetag{4.4} into \thetag{4.2} we see that
the characteristic equation \thetag{4.2} of the operator $f$ is a
polynomial equation of $n$-th order with respect to $\lambda$:
$$
(-\lambda)^n+F_1\,(-\lambda)^{n-1}+\ldots+F_n=0.
\tag4.6
$$
Therefore we can estimate the number of eigenvalues of the operator
$f$. Any eigenvalue $\lambda\in\Bbb K$ is a root of characteristic
equation \thetag{4.6}. However, not any root of the equation 
\thetag{4.6} is an eigenvalue of the operator $f$. The matter is that
a polynomial equation with coefficients in the numeric field
$\Bbb K$ can have roots in some larger field $\tilde\Bbb K$ (e\.\,g\.
$\Bbb Q\subset\Bbb R$ or $\Bbb R\subset\Bbb C$). For the characteristic
number $\lambda$ of the operator $f$ to be an eigenvalue of this
operator it should belong to $\Bbb K$. From the course of general algebra
we know that the total number of roots of the equation \thetag{4.6}
counted according to their multiplicity and including those belonging
to the extensions of the field $\Bbb K$ is equal to $n$ (see \cite{4}).
\proclaim{Theorem 4.1} The number of eigenvalues of a linear operator
$f\!:\,V\to V$ equals to the dimension of the space $V$ at most.
\endproclaim
     Consider the case $\Bbb K=\Bbb Q$. The roots of a polynomial equation with rational coefficients are not necessarily rational numbers: the
equation $\lambda^2-3=0$ is an example. In the case of real numbers 
$\Bbb K=\Bbb R$ a polynomial equation with real coefficients can also
have non-real roots, e\.\,g\. the equation $\lambda^2+\sqrt{3}=0$. However,
the field of complex numbers $\Bbb K=\Bbb C$ is an exception.
\proclaim{Theorem 4.2} An arbitrary polynomial equation of $n$-th 
order with complex coefficients has exactly $n$ complex roots 
counted according to their multiplicity.
\endproclaim
     We shall not prove here this theorem referring the reader to the
course of general algebra (see \cite{4}). The theorem~4.2 is known as
the {\tencyr\char '074}basic theorem of algebra{\tencyr\char '076},
while the property of complex numbers stated in this theorem 
is called the {\it algebraic closure} of $\Bbb C$, i\.\,e\.
$\Bbb C$ is an algebraically closed numeric field.
\definition{Definition 4.2} A numeric field $\Bbb K$ is called an
{\it algebraically closed} field if the roots of any polynomial 
equation with coefficients from $\Bbb K$ are again in $\Bbb K$.
\enddefinition
     Certainly, $\Bbb C$ is not the unique algebraically closed field.
However, in the list of numeric fields $\Bbb Q$, $\Bbb R$, $\Bbb C$
that we consider in this book, only the field of complex numbers
is algebraically closed.\par
     Let $\lambda$ be an eigenvalue of a linear operator $f$. Then
$\lambda$ is a root of the equation \thetag{4.6}. The multiplicity
of this root $\lambda$ in the equation \thetag{4.6} is called the
{\it multiplicity} of the eigenvalue $\lambda$.
\proclaim{Theorem 4.3} For a linear operator $f\!:\,V\to V$ in a complex
linear vector space $V$ the number of its eigenvalues counted according
to their multiplicities is exactly equal to the dimension of $V$.
\endproclaim
     This proposition strengthen the theorem~4.1. It is an immediate
consequence of the algebraic closure of the field of complex numbers
$\Bbb C$. In the case $\Bbb K=\Bbb C$ the characteristic polynomial
\thetag{4.4} is factorized into a product of terms linear in $\lambda$:
$$
\hskip -2em
\det(f-\lambda\cdot 1)=\prod^n_{i=1}(\lambda_i-\lambda).
\tag4.7
$$
For some operators such an expansion can occur in the case $\Bbb K=\Bbb Q$
or $\Bbb K=\Bbb R$, however, it is not a typical situation. If $\lambda_1,
\,\ldots,\,\lambda_n$ are understood as characteristic numbers of the
operator $f$, then the formula \thetag{4.7} is always valid.\par
     Due to the formula \thetag{4.7} we can present the numeric 
invariants $F_1,\,\ldots,\,F_n$ of the operator $f$ as elementary 
symmetric polynomials of its characteristic numbers:
$$
F_i=\sigma_i(\lambda_1,\ldots,\lambda_n).
$$
In particular, for the trace and for the determinant 
of the operator $f$ we have
$$
\xalignat 2
&\hskip -2em
\tr f=\sum^n_{i=1}\lambda_i,
&&\det f=\prod^n_{i=1}\lambda_i.
\tag4.8
\endxalignat
$$
The theory of symmetric polynomials is given in the course of general
algebra (see, for example, the book \cite{4}).\par
\proclaim{Theorem~4.4} For any eigenvalue $\lambda$ of a linear operator
$f\!:\,V\to V$ the asso\-ciated eigenspace $V_\lambda$ is invariant under
the action of $f$.
\endproclaim
\demo{Proof} The definition~4.1 of an eigenspace $V_\lambda$ of a linear
operator $f$ can be reformulated as $V_\lambda=\{\bold v\in V\!:\ f(\bold 
v)=\lambda\cdot\bold v\}$. Therefore, $\bold v\in V_\lambda$ implies 
$f(\bold v)=\lambda\cdot\bold v\in V_\lambda$, which proves the invariance
of $V_\lambda$.
\qed\enddemo
     We know that the set of linear operators in a space $V$ form the
algebra $\End(V)$ over the numeric field $\Bbb K$. However, this algebra
is too big. Let's consider some operator $f\in\End(V)$ and complement it
with the identical operator $1$. Within the algebra $\End(V)$ we can take
positive integer powers of the operator $f$, we can multiply them by 
numbers from $\Bbb K$, we can add such products, and we can add to them scalar operators obtained by multiplying the identical operator $1$ by various numbers from $\Bbb K$. As a result we obtain various operators
of the form
$$
\hskip -2em
P(f)=\alpha_p\cdot f^p+\ldots+\alpha_1\cdot f+\alpha_0\cdot 1.
\tag4.9
$$
The set of all operators of the form \thetag{4.9} is called the {\it
polynomial envelope} of the operator $f$; it is denoted $\Bbb K[f]$. 
This is a subset of $\End(V)$ closed with respect to all algebraic operations in $\End(V)$. Such subsets are used to be called {\it subalgebras}. It is important to say that the subalgebra $\Bbb K[f]$ 
is commutative, i\.\,e\. for any two polynomials $P$ and $Q$ the corresponding operators \thetag{4.9} commute:
$$
\hskip -2em
P(f)\compos Q(f)=Q(f)\compos P(f).
\tag4.10
$$
The equality \thetag{4.10} is verified by direct calculation. Indeed,
let $P(f)$ and $Q(f)$ be two operator polynomials of the form:
$$
\xalignat 2
&P(f)=\sum^p_{i=0}\alpha_i\cdot f^i,
&&Q(f)=\sum^q_{j=0}\beta_j\cdot f^j.
\endxalignat
$$
Here we denote: $f^{\,0}=1$. This relationship should be treated as
the definition of zeroth power of the operator $f$. Then 
$$
P(f)\compos Q(f)=\sum^p_{i=0}\sum^q_{j=0}(\alpha_i\beta_j)\cdot
f^{i+j}=Q(f)\compos P(f).
$$
These calculations prove the relationship \thetag{4.10}.\par
\proclaim{Theorem~4.5} Let $U$ be an invariant subspace of an operator 
$f$. Then it is invariant under the action of any operator from the polynomial envelope $\Bbb K[f]$.
\endproclaim
\demo{Proof} Let $\bold u$ be an arbitrary vector of $U$. Let's 
consider the following vectors $\bold u_0=\bold u$, $\bold u_1=
f(\bold u)$, $\bold u_2=f^2(\bold u)$, $\ldots$, $\bold u_p=
f^p(\bold u)$. Every next vector in this sequence is obtained
by applying the operator $f$ to the previous one: $\bold u_{i+1}
=f(\bold u_i)$. Therefore, from $\bold u_0\in U$ it follows that 
$\bold u_1\in U$ since $U$ is an invariant subspace of $f$. Then,
in turn, we successively obtain $\bold u_2\in U$, $\bold u_3\in U$,
and so on up to $\bold u_p\in U$. Applying the operator $P(f)$ of
the form \thetag{4.9} to the vector $\bold u$, we get
$$
P(f)\,\bold u=\alpha_p\cdot\bold u_p+\ldots+\alpha_1\cdot\bold
u_1+\alpha_0\cdot\bold u_0.
$$
Hence, due to $\bold u_i\in U$ we find that $P(f)\,\bold u\in U$, 
which proves the invariance of $U$ under the action of the operator
$P(f)$.
\qed\enddemo
     The following fact is curious: if $\lambda$ is an eigenvalue 
of the operator $f$ and if $\bold v$ is an associated eigenvector,
then $P(f)\,\bold v=P(\lambda)\cdot\bold v$. \pagebreak Therefore, 
any eigenvector $\bold v$ of the operator $f$ is an eigenvector of 
the operator $P(f)$. The converse proposition, however, is not 
true.\par
     Let $\lambda_1,\,\ldots,\,\lambda_s$ be a set of mutually distinct
eigenvalues of the operator $f$. Let's consider the operators $h_i=f-
\lambda_i\cdot\bold 1$, which certainly belong to the polynomial envelope of $f$. The permutability of any two such operators follows from
\thetag{4.10}. The eigenspace $V_{\lambda_{\ssize i}}$ of the operator
$f$ is determined as the kernel of the operator $h_i$. According to the
definition~4.1, it is nonzero. Moreover, the theorems~4.4 and 4.5
say that $V_{\lambda_{\ssize i}}$ is invariant under the action
of $f$ and of all other operators $h_j$.\par
\proclaim{Theorem 4.6} Let $\lambda_1,\,\ldots,\,\lambda_s$ be a set of
mutually distinct eigenvalues of the operator $f\!:\,V\to V$. Then the
sum of associated eigenspaces $V_{\lambda_1},\,\ldots,\,V_{\lambda_s}$
is a direct sum: $V_{\lambda_1}+\ldots+ V_{\lambda_{\ssize s}}=
V_{\lambda_1}\oplus\ldots\oplus V_{\lambda_{\ssize s}}$.
\endproclaim
     Note that the set of mutually distinct eigenvalues $\lambda_1,\,
\ldots,\,\lambda_s$ of the operator $f$ in this theorem could be the
complete set of such eigenvalues, or it could include only a part of
such eigenvalues. This makes no difference for the result of the
theorem~4.6, it remains valid in either case.\par
\demo{Proof} Let's denote by $W$ the sum of eigenspaces of the operator
$f$:
$$
\hskip -2em
W=V_{\lambda_1}\oplus\ldots\oplus V_{\lambda_{\ssize s}}.
\tag4.11
$$
In order to prove that the sum \thetag{4.11} is a direct sum we need
to prove that for an arbitrary vector $\bold w\in W$ the expansion
$$
\hskip -2em
\bold w=\bold v_1+\ldots+\bold v_s\text{, \ where \ }\bold v_i\in
V_{\lambda_{\ssize i}},
\tag4.12
$$
is unique. For this purpose we consider the operator $f_i$ defined
by formula
$$
f_i=\prod^s_{r\neq i}h_r.
$$
The operator $f_i$ belongs to the polynomial envelope
of the operator $f$ and 
$$
\hskip -2em
f_i(\bold v_j)=\left(\,\shave{\prod^s_{r\neq i}}(\lambda_j-\lambda_r)
\right)\cdot\bold v_j.
\tag4.13
$$
This follows from $\bold v_j\in V_{\lambda_{\ssize j}}$, which implies
$h_r(\bold v_j)=(\lambda_j-\lambda_r)\cdot\bold v_j$. The formula
\thetag{4.13} means that $f_i(\bold v_j)=0$ for all $j\neq i$. Applying the operator $f_i$ to both sides of the expansion \thetag{4.12}, we get
the equality 
$$
f_i(\bold w)=\left(\,\shave{\prod^s_{r\neq i}}(\lambda_i-\lambda_r)
\right)\cdot\bold v_i.
$$
Hence, for the vector $\bold v_i$ in the expansion \thetag{4.12} we
derive
$$
\pagebreak 
\hskip -2em
\bold v_i=\frac{f_i(\bold w)}{\dsize\prod^s_{r\neq i}(\lambda_i
-\lambda_r)}.
\tag4.14
$$
The formula \thetag{4.14} uniquely determines all summands 
in the expansion \thetag{4.12} if the vector $\bold w\in W$ is given.
This means that the expansion \thetag{4.12} is unique and the sum of subspaces \thetag{4.11} is a direct sum.
\qed\enddemo
\definition{Definition4.3} A linear operator $f\!:\,V\to V$ in a 
linear vector space $V$ is called a {\it diagonalizable operator} if 
there is a basis $\bold e_1,\,\ldots,\,\bold e_n$ in the space $V$ 
such that the matrix of the operator $f$ is diagonal in this basis.
\enddefinition
\proclaim{Theorem 4.7} An operator $f\!:\,V\to V$ is diagonalizable 
if and only if the sum of all its eigenspaces coincides with $V$.
\endproclaim
\demo{Proof} Let $f$ be a diagonalizable operator. Then we can choose a basis $\bold e_1,\,\ldots,\,\bold e_n$ such that its matrix $F$ in this
basis is diagonal, i\.\,e\. only diagonal elements $F^i_i$ of this matrix
can be nonzero. Then the relationship \thetag{1.2}, which determines the
matrix $F$, is written as $f(\bold e_i)=F^i_i\cdot\bold e_i$. Hence,
each basis vector $\bold e_i$ is an eigenvector of the operator $f$,
while $\lambda_i=F^i_i$ is its associated eigenvalue. The expansion of
an arbitrary vector $\bold v$ in this base is an expansion by eigenvectors
of the operator $f$. Therefore, having collected together the terms with coinciding eigenvalues in this expansion, we get the expansion 
$$
\bold v=\bold v_1+\ldots+\bold v_s\text{, \ where \ }\bold v_i
\in V_{\lambda_{\ssize i}}.
$$
Since $\bold v$ is an arbitrary vector of $V$, this means that
$V_{\lambda_1}+\ldots+V_{\lambda_{\ssize s}}=V$. The direct proposition
of the theorem is proved.\par
     Conversely, suppose that $\lambda_1,\,\ldots,\,\lambda_s$ is the
total set of mutually distinct eigenvalues of the operator $f$ and assume
that $V_{\lambda_1}+\ldots+V_{\lambda_{\ssize s}}=V$. The theorem~4.6
says that this is a direct sum: $V=V_{\lambda_1}\oplus\ldots\oplus
V_{\lambda_{\ssize s}}=V$. Therefore, choosing a basis in each eigenspace
and joining them together, we get a basis in $V$ (see theorem~6.3 in
Chapter~\uppercase\expandafter{\romannumeral 1}). This is a basis composed
by eigenvectors of the operator $f$, the application of $f$ to each basis
vector reduces to multiplying this vector by its associated eigenvalue.
Therefore, the matrix $F$ of the operator $f$ in this basis is  diagonal.
Its diagonal elements coincide with the eigenvalues of the operator $f$.
The theorem is proved.
\qed\enddemo
     Assume that an operator $f\!:\,V\to V$ is diagonalizable and assume that we have chosen a basis where its matrix is diagonal. Then the matrix
$H_f$ in formula \thetag{4.3} is also diagonal. Hence, we immediately 
derive the following formula:
$$
\det(f-\lambda\cdot\bold 1)=\prod^n_{i=1}(F^i_i-\lambda).
$$
Due to this equality we conclude that the characteristic polynomial of a
diagona\-lizable operator is factorized into the product of a linear terms
and all roots of characteristic equation belong to the field $\Bbb K$
(not to its extension). This means that characteristic numbers of a
diagonalizable operator coincide with its eigenvalues. This is a necessary
condition for the operator $f$ to be diagonalizable. However, it is not
a sufficient condition. Even in the case of algebraically closed field
of complex numbers $\Bbb K=\Bbb C$ there are non-diagonalizable operators
in vector spaces over the field $\Bbb C$.\par
\head
\S\,5. Nilpotent operators.
\endhead
\definition{Definition 5.1} A linear  operator $f\!:\,V\to V$ is called
a {\it nilpotent operator} if for any vector $\bold v\in V$ there is a 
positive integer number $k$ such that $f^k(\bold v)=\bold 0$.
\enddefinition
     According to the definition~5.1 for any vector $v$ there is an integer number $k$ (depending on $\bold v$) such that $f^k(\bold v)=\bold 0$. The choice of such number has no upper bound, indeed, if $m>k$ and $f^k(\bold
v)=\bold 0$ then $f^m(\bold v)=\bold 0$. This means that there is a
minimal positive number $k=k_{\min}$ (depending on $\bold v$) such that 
$f^k(\bold v)=\bold 0$. This minimal number $k_{\min}$ is called the
{\it height} of the vector $\bold v$ respective to the nilpotent operator
$f$. The height of zero vector is taken to be zero by definition; for
any nonzero vector $\bold v$ its height is greater or equal to the unity.
Let's denote the height of $\bold v$ by $\nu(\bold v)$ and define the
number 
$$
\hskip -2em
\nu(f)=\max_{\lower2pt\hbox{$\ssize\bold v\in V$}}\nu(\bold v).
\tag5.1
$$
For each vector $\bold v\in V$ its height is finite, but the maximum 
in \thetag{5.1} can be infinite since the number of vectors in a linear
vector space usually is infinite.\par
\definition{Definition 5.2} In that case, where the maximum in the
formula \thetag{5.1} is finite, a nilpotent operator $f$ is called 
an {\it operator of finite height\/} and the number $\nu(f)$ is called
the {\it height\/} of a nilpotent operator $f$.
\enddefinition
\proclaim{Theorem 5.1} In a finite-dimensional linear vector space 
$V$ the height $\nu(f)$ of any nilpotent operator $f\!:\,V\to V$ is
finite.
\endproclaim
\demo{Proof} Let's choose a basis $\bold e_1,\,\ldots,\,\bold e_n$ 
in $V$ and consider the heights of all basis vectors $\nu(\bold e_1),
\,\ldots,\,\nu(\bold e_n)$ with respect to $f$. Then denote
$$
m=\max\{\nu(e_1),\ldots,\nu(e_n)\}.
$$
For an arbitrary vector $\bold v\in V$ consider its expansion
$\bold v=v^1\cdot\bold e_1+\ldots+v^n\cdot\bold e_n$. Then, applying
the operator $f^m$ to $\bold v$, we find
$$
\hskip -2em
f^m(\bold v)=\sum^n_{i=1}v^i\cdot f^m(\bold e_i)=\bold 0.
\tag5.2
$$
Due to the formula \thetag{5.2} we see that the heights of all vectors
of the space $V$ are restricted by the number $m$. This means that the
height of a nilpotent operator $f$ is finite: $\nu(f)=m<\infty$.
\qed\enddemo
\proclaim{Theorem 5.2} If $f\!:\,V\to V$ is a nilpotent operator 
and if $U$ is an invariant subspace of the operator $f$, then the
restricted operator $f{\!\lower4pt\hbox{$\ssize U$}}$ and the
factoroperator $f{\!\lower4pt\hbox{$\ssize V/U$}}$ both are nilpotent.
\endproclaim
\demo{Proof} Any vector $\bold u$ of the subspace $U\subset V$ is a
vector of $V$. Therefore, there is an integer number $k>0$ such that
$f^k(\bold u)=\bold 0$. However, the result of applying the restricted operator $f{\!\lower4pt\hbox{$\ssize U$}}$ the a vector of $U$
coincides with the result of applying the initial operator $f$ to
this vector. Hence, we have
$$
(f{\!\lower4pt\hbox{$\ssize U$}})^k\,\bold u=f^k(\bold u)=\bold 0.
$$
This proves that $f{\lower4pt\hbox{$\ssize U$}}$ is a nilpotent 
operator. In the case of factoroperator we consider an arbitrary coset 
$Q$ in the factorspace $V/U$. Let $Q=\Cl_U(\bold v)$, where 
$v$ is some fixed vector in $Q$, and let $k=\nu(v)$ be the height
of this vector $\bold v$ respective to the operator $f$. Then
we can calculate 
$$
(f{\!\lower4pt\hbox{$\ssize V/U$}})^k\,Q=\Cl_U(f^k(\bold v))=\bold 0.
$$
Now it is clear that the factoroperator $f{\!\lower4pt\hbox{$\ssize V/U$}}$
is a nilpotent operator. The theorem is completely proved.
\qed\enddemo
\proclaim{Theorem 5.3} A nilpotent operator $f$ cannot have a nonzero eigenvalue.
\endproclaim
\demo{Proof} Let $\lambda$ be an eigenvalue of a nilpotent operator 
$f$ and let $\bold v\neq\bold 0$ be an associated eigenvector. Then we
have $f(\bold v)=\lambda\cdot\bold v$. On the other hand, since $f$ is
nilpotent, there is a number $k>0$ such that $f^k(\bold v)=\bold 0$.
Then we derive
$$
f^k(\bold v)=\lambda^k\cdot\bold v=\bold 0.
$$
But $\bold v\neq\bold 0$, therefore, $\lambda^k=0$. This is the equation 
for $\lambda$ and $\lambda=0$ its unique root. The theorem is 
proved.
\qed\enddemo
     It the finite-dimensional case this theorem can be strengthened
as follows.
\proclaim{Theorem 5.4} In a finite-dimensional space $V$ of the 
dimension $\dim V=n$ any nilpotent operator $f$ has exactly one 
eigenvalue $\lambda=0$ with the multiplicity $n$.
\endproclaim
\demo{Proof} We shall prove this theorem by induction on $n=\dim V$. 
In the case $n=1$ we fix some vector $\bold v\neq\bold 0$ in $V$ and 
denote by $k=\nu(\bold v)$ its height. Then $f^k(\bold v)=\bold 0$ and
$f^{k-1}(\bold v)\neq \bold 0$. This means that $\bold w=f^{k-1}(\bold
v)\neq\bold 0$ is an eigenvector of $f$ with the eigenvalue $\lambda=0$
since $f(\bold w)=f^k(\bold v)=\bold 0=0\cdot\bold w$. The base of
the induction is proved.\par
     Suppose that the theorem is proved for any finite-dimensional
space of the dimension less than $n$ and consider a space $V$ of the dimension $n=\dim V$. As above, let's fix some vector $\bold v\neq\bold 0$
in $V$ and denote by $k=\nu(\bold v)$ its height respective to the
operator $f$. Then $f^k(\bold v)=\bold 0$ and $\bold w=f^{k-1}(\bold v)
\neq\bold 0$. Hence, for the nonzero vector $\bold w$ we get the following
series of equalities:
$$
f(\bold w)=f(f^{k-1}(\bold v))=f^k(\bold v)=\bold 0=0\cdot\bold w.
$$
Hence, $\bold w$ is an eigenvector of the operator $f$ and $\lambda=0$
is its associated eigenvalue. Let's consider the eigenspace $U=V_0$
corresponding to the eigenvalue $\lambda=0$. Let's denote $m=\dim U
\neq 0$. The restricted operator $f{\lower4pt\hbox{$\ssize U$}}$ is
zero, hence, for characteristic polynomial of this operator 
$f{\!\lower4pt\hbox{$\ssize U$}}=0$ we derive 
$$
\det(f{\!\lower4pt\hbox{$\ssize U$}}-\lambda\cdot 1)=
(-\lambda)^m.
$$
Now, applying the theorem~3.6, we derive the characteristic polynomial
of $f$:
$$
\hskip -2em
\det(f-\lambda\cdot 1)=(-\lambda)^m\,
\det(f{\!\lower4pt\hbox{$\ssize V/U$}}-\lambda\cdot 1).
\tag5.3
$$
The factoroperator $f{\!\lower4pt\hbox{$\ssize V/U$}}$ is an operator in
factorspace $V/U$ whose dimension $n-m$ is less than $n$. Due to the
theorem~5.2 the factoroperator $f{\!\lower4pt\hbox{$\ssize V/U$}}$ is
nilpotent, therefore, we can apply the inductive hypothesis to it. Then
for its characteristic polynomial of the factoroperator 
$f{\!\lower4pt\hbox{$\ssize V/U$}}$ we get
$$
\hskip -2em
\det(f{\lower4pt\hbox{$\ssize V/U$}}-\lambda\cdot 1)=
(-\lambda)^{n-m}.
\tag5.4
$$
Comparing the above relationships \thetag{5.3} and \thetag{5.4}, we
find the characteristic polynomial of the initial operator $f$:
$$
\det(f-\lambda\cdot 1)=(-\lambda)^n.
$$
This means that $\lambda=0$ is the only eigenvalue of the operator $f$
and its multiplicity is $n=\dim V$. The theorem is proved.
\qed\enddemo
     Let $f\!:\,V\to V$ be a linear operator. Consider a vector $\bold v
\in V$ and denote by $k=\nu(\bold v)$ its height respective to the operator
$f$. This vector $\bold v$ produces the chain of $k$ vectors according to
the following formulas:
$$
\hskip -2em
\bold v_1=f^{k-1}(\bold v),\ \ \bold v_2=f^{k-2}(\bold v),\ 
\ldots,\ \
\bold v_k=f^0(\bold v)=\bold v.
\tag5.5
$$
The chain vectors \thetag{5.5} are related with each other as follows:
$\bold v_i=f(\bold v_{i-1})$. Let's apply the operator $f$ to each vector
in the chain \thetag{5.5}. Then the first vector $\bold v_1$ vanished. 
Applying $f$ to the rest $k-1$ vectors we get another chain:
$$
\hskip -2em
\bold w_1=f^{k-1}(\bold v),\ \ \bold w_2=f^{k-2}(\bold v),\
\ldots,\ \bold w_{k-1}=f(\bold v).
\tag5.6
$$
Comparing these two chains \thetag{5.5} and \thetag{5.6}, we see that 
they are almost the same, but the second chain is shorter. It is 
obtained from the first one by removing the last vector $\bold v_k
=\bold v$.\par
     The vector $\bold v_1$ is called the {\it side vector\/} or the
{\it eigenvector\/} of the chain \thetag{5.5}. The other vectors are 
called the {\it adjoint vectors\/} of the chain. If the side vectors 
of two chains are different, then in these two chains there are no
coinciding vectors at all. However, there is even stronger result. It 
is known as the theorem on {\tencyr\char '074}linear independence of 
chains{\tencyr\char '076}.
\proclaim{Theorem 5.5} If the side vectors in several chains of the 
form \thetag{5.5} are linearly independent, then the whole set of 
vectors in these chains is linearly independent.
\endproclaim
\demo{Proof} We consider $s$ chains of the form \thetag{5.5}. In order
to specify the chain vector we use two indices $\bold v_{i,j}$. The
first index $i$ is the number of chain to which this vector $\bold v_{i,j}$ belongs, the second index $j$ specifies the number of this vector 
within the $i$-th chain. Denote by $k_1,\,\ldots,\,k_s$ the lengths of
our chains. Without loss of generality we can assume that the chains are
arranged in the order of decreasing their lengths, i\.\,e\. we have the
following inequalities:
$$
\hskip -2em
k_1\geqslant k_2\geqslant\ldots\geqslant k_s\geqslant 1.
\tag5.7
$$\par
     Let $k=\max\{k_1,\ldots,k_s\}$. We shall prove the theorem by induction on $k$. If $k=1$ then the lengths of all chains are equal 
to $1$. Therefore, they contain only the side vectors and have no
adjoint vectors at all. The proposition of the theorem in this case 
is obviously true.\par
     Suppose that the theorem is valid for the chains whose lengths 
are not greater than $k-1$. For our $s$ chains, whose lengths are restricted by the number $k$, we consider a linear combination of all
their vectors being equal to zero:
$$
\hskip -2em
\sum^s_{i=1}\sum^{k_{\ssize s}}_{j=1} \alpha_{i,j}\cdot
\bold v_{i,j}=\bold 0.
\tag5.8
$$
From this equality we should derive the triviality of the linear
combination in its left hand side. Let's apply the operator $f$ to
both sides of \thetag{5.8} and use the following quite obvious
relationships:
$$
f(\bold v_{i,j})=
\cases
 \bold 0,        &\text{\ for \ }j=1,\\
 \bold v_{i,j-1},&\text{\ for \ }j>1.
\endcases
$$
If we take into account \thetag{5.7}, then the result of applying $f$ 
to \thetag{5.8} is written as
$$
\hskip -2em
\sum^s_{i=1}\sum^{k_{\ssize s}}_{j=1}\alpha_{i,j}\cdot
f(\bold v_{i,j})=\sum^r_{i=1}\sum^{k_{\ssize r}}_{j=2}\alpha_{i,j}
\cdot\bold v_{i,j-1}=\bold 0.
\tag5.9
$$
In typical situation $r=s$. However, sometimes certain chains of vectors 
can drop from the above sums at all. This happens if a part of chains 
were of the length $1$. In this case $r<s$ and $k_{r+1}=\ldots=k_s=1$.
The lengths of all chains in \thetag{5.7} cannot be equal to $1$ since
$k>1$.\par
     Shifting the index $j+1\to j$ in the last sum we can write 
\thetag{5.9} as follows:
$$
\hskip -2em
\sum^r_{i=1}\sum^{k_{\ssize r}-1}_{j=1} \alpha_{i,j+1}\cdot
\bold v_{i,j}=\bold 0.
\tag5.10
$$
The left side of the relationship \thetag{5.10} is again a linear combination of chain vectors. Here we have $r$ chins with the lengths
$1$ less as compared to original ones in \thetag{5.8}. Now we can apply 
the inductive hypothesis, which yields the linear independence of all vectors presented in \thetag{5.10}. Hence, all coefficients of the
linear combination in left hand side of \thetag{5.10} are equal to zero.
When applied to \thetag{5.8} this fact means that the most part of terms 
in left hand side of this equality do actually vanish. The remainder is
written as follows:
$$
\hskip -2em
\sum^s_{i=1}\alpha_{i,1}\cdot\bold v_{i,1}=0.
\tag5.11
$$
Now in the linear combination \thetag{5.11} we have only the side vectors
of initial chains. The are linearly independent by the assumption of the
theorem. Therefore, the linear combination \thetag{5.11} is also trivial.
From triviality of \thetag{5.10} and \thetag{5.11} it follows that the
initial linear combination \thetag{5.8} is trivial too. We have completed
the inductive step and thus have proved the theorem in whole.
\qed\enddemo
     Let $f\!:\,V\to V$ be a nilpotent operator in a linear vector space 
$V$ and let $\bold v$ be a vector of the height $k=\nu(\bold v)$ in $V$.
Consider the chain of vectors \thetag{5.5} generated by $\bold v$ and 
denote by $U(\bold v)$ the linear span of chain vectors \thetag{5.5}:
$$
\hskip -2em
U(\bold v)=\langle\bold v_1,\,\ldots,\,\bold v_k\rangle.
\tag{5.12}
$$
Due to the theorem~5.5 the subspace $U(\bold v)$ is a finite-dimensional
subspace and $\dim U(\bold v)=k$. The chain vectors \thetag{5.5} form a basis in this subspace \thetag{5.12}. The following relationships
are derived directly from the definition of the chain \thetag{5.5}:
$$
\hskip -2em
\matrix
\format \l\\
f(\bold v_1)=\bold 0,\\
\vspace{0.5ex}
f(\bold v_2)=\bold v_1,\\
\hdotsfor 1\\
\vspace{0.4ex}
f(\bold v_k)=\bold v_{k-1}.
\endmatrix
\tag5.13
$$
Due to \thetag{5.13} the subspace \thetag{5.12} is invariant under the
action of the operator $f$. Hence, we can consider the restricted operator
$f{\!\lower4pt\hbox{$\ssize U(\bold v)$}}$ and, using \thetag{5.13}, we
can find the matrix of this restricted operator in the chain basis
$\bold v_1\,\ldots,\,\bold v_k$:
$$
\hskip -2em
J_k(0)=\Vmatrix
0 & 1 & \      &0 \\
\ & 0 & \ddots &\ \\
\ & \ & \ddots &1 \\
\ & \ & \      &0 \\
\endVmatrix.
\tag5.14
$$
A matrix of the form \thetag{5.14} is called a {\it Jordan block} or a
{\it Jordan cage} of a nilpotent operator. Its primary diagonal is filled
with zeros. The upward next diagonal parallel to the primary one is filled
with unities. All other space in the matrix \thetag{5.14} is filled with zeros again. The matrix \thetag{5.14} is a square $k\times k$ matrix, if 
$k=1$, this matrix degenerates and becomes purely zero matrix with the only
element: $J_1(0)=\Vmatrix 0\endVmatrix$.\par
     Let $f\!:\,V\to V$ again be a nilpotent operator. We continue to study
vector chains of the form \thetag{5.5}. For this purpose let's consider the
following subspaces:
$$
\hskip -2em
U_k=\Ker f\cap\Img f^{k-1}.
\tag5.15
$$
If $\bold u\in U_k$, then $\bold u\in\Img f^{k-1}$. Therefore, $\bold u
=f^{k-1}(\bold v)$ for some vector $\bold v$. This means that $\bold u$ 
is a chain vector in a chain of the form \thetag{5.5}. From the condition
$\bold u\in\Ker f$ we derive $f(\bold u)=f^k(\bold v)=\bold 0$. Hence,
$\bold v$ is a vector of the height $k$ and $\bold u$ is a side vector 
in the chain \thetag{5.5} initiated by the vector $\bold v$. For the
subspaces \thetag{5.15} we have the sequence of inclusions
$$
\hskip -2em
V_0=U_1\supseteq U_2\supseteq\ldots\supseteq U_k\supseteq\ldots,
\tag5.16
$$
where $V_0=\Ker f$ is the eigenspace corresponding to the unique
eigenvalue $\lambda=0$ of nilpotent operator $f$. The inclusions
\thetag{5.16} follow from the fact that any chain \thetag{5.5} of 
the length $k$ with the side vector $\bold u=f^{k-1}(\bold v)$ can 
be treated as a chain of the length $k-1$ by dropping the $k$-th
vector $\bold v_k=\bold v$ (see \thetag{5.5} and \thetag{5.6}). 
Then for the vector $\bold v'=f(\bold v)$ we have $\bold u=f^{k-2}
(\bold v')$. This yields the inclusion of subspaces $U_k\subset 
U_{k-1}$ for $k>1$.\par
     In a finite-dimensional space $V$ the height of any vector
$\bold v\in V$ is restricted by the height of the nilpotent operator 
$f$ itself:
$$
\nu(\bold v)\leq\nu(f)=m<\infty
$$
(see theorem~5.1). Therefore $U_{m+1}=\{\bold 0\}$. \pagebreak Hence, the
sequence of inclusions \thetag{5.16} terminates on $m$-th step, i\.\,e\.
we have a finite sequence of inclusions:
$$
\hskip -2em
V_0=U_1\supseteq U_2\supseteq\ldots\supseteq U_m\supseteq\{\bold 0\}.
\tag5.17
$$
Sequences of mutually enclosed subspaces of the form \thetag{5.16} or \thetag{5.17} are called {\it flags\/}, while each particular subspace
in a flag is called a {\it flag subspace}.
\proclaim{Theorem 5.6} For any nilpotent operator $f$ in a
finite-dimensional space $V$ there is a basis in $V$ composed by chain 
vectors of the form \thetag{5.5}. Such a basis is called a {\it canonic
basis} or a {\it Jordan basis} of a nilpotent operator $f$.
\endproclaim
\demo{Proof} The proof of the theorem is based on the fact that the
flag \thetag{5.17} is finite. We choose a basis in the smallest 
subspace $U_m$. Then we complete it up to a basis in $U_{m-1}$, in
$U_{m-2}$, and so on backward along the sequence \thetag{5.17}.
As a result we construct a basis $\bold e_1,\,\ldots,\,\bold e_s$ in
$V_0=\Ker f$. Note that each vector in such a basis is a side vector
of some chain of the form \thetag{5.5}. For basis vectors of the 
subspace $U_m$ the lengths of such chains are equal to $m$. For the 
complementary vectors from $U_{m-1}$ their chins are of the length
$m-1$ and further the length of chains decreases step by step until 
the unity for the complementary vectors in largest subspace $U_1
=V_0$.\par
    Let's join together all vectors of the above chains and let's 
enumerate them by means of double indices: $\bold e_{i,j}$. Here 
$i$ is the number of the chain and $j$ is the individual number 
of the vector within $i$-th chain. Then 
$$
\bold e_1=\bold e_{1,1},\ \ldots,\ \bold e_s=\bold e_{s,1}.
$$
Now let's prove that the set of all vectors from the above chains
form a basis in $V$. The linear independence of this set of vectors
follows from the theorem~5.5. We only have to prove that an arbitrary
vector $\bold v\in V$ can be represented as a linear combination 
of chain vectors $\bold e_{i,j}$. We shall prove this fact by induction
on the height of the vector $\bold v$.\par
     If $k=\nu(\bold v)=1$, then $\bold v\in\Ker f=V_0$. In this case 
$\bold v$ is expanded in the basis $\bold e_1,\,\ldots,\,\bold e_s$ 
of the subspace $V_0$. This is the base of induction.\par
    Now suppose that any vector of the height less than $k$ can be
represented as a linear combination of chain vectors $\bold e_{i,j}$.
Let's take a vector $\bold v$ of the height $k$ and denote $\bold u
=f^{k-1}(\bold v)$. Then $f(\bold u)=\bold 0$. This means that 
$\bold u$ is a side vector in a chain of the length $k$ initiated
by the vector $\bold v$. Therefore, $\bold u$ is an element of the
subspace $U_k$ (see formula \thetag{5.15}); this vector can be 
expanded in the basis of the subspace $U_k$, which we have 
constructed above: 
$$
\hskip -2em
\bold u=\sum^r_{i=1}\alpha_i\cdot\bold e_i.
\tag5.18
$$
Note that in the expansion \thetag{5.18} we have only a part of 
vectors $\bold e_1,\,\ldots,\,\bold e_s$, namely, we have only those
of them that belongs to $U_k$ and, hence, are side vectors in the
chains of the length not less than $k$. Therefore, we can write
$\bold e_i=f^{k-1}(\bold e_{i,k})$ for $i=1,\,\ldots,\,r$. 
Substituting these expressions into \thetag{5.18}, we obtain
$$
\hskip -2em
f^{k-1}(\bold v)=\sum^r_{i=1}\alpha_i\cdot f^{k-1}(\bold e_{i,k}).
\tag5.19
$$
By means of the coefficients of the expansion \thetag{5.19} we determine
the vector $\bold v'$:
$$
\hskip -2em
\bold v'=\bold v-\sum^r_{i=1}\alpha_i\cdot\bold e_{i,k}.
\tag5.20
$$
Applying the operator $f^{k-1}$ to $\bold v'$ and taking into
account \thetag{5.19}, we find
$$
f^{k-1}(\bold v')=f^{k-1}(\bold v)-\sum^r_{i=1}\alpha_i\cdot\bold
f^{k-1}(\bold e_{i,k})=\bold 0.
$$
Hence, the height of the vector $\bold v'$ is less than $k$ and we can
apply the inductive hypothesis to it. This means that $\bold v'$ can
be represented as a linear combination of chain vectors $\bold e_{i,j}$. 
But $\bold v$ is expressed through $\bold v'$ as follows:
$$
\bold v=\bold v'+\sum^r_{i=1}\alpha_i\cdot\bold e_{i,k}.
$$
Then $\bold v$ can also be expressed as a linear combination of chain
vectors $\bold e_{i,j}$. The inductive step is completed and the theorem
in whole is proved.
\qed\enddemo
     In the basis composed by chain vectors, the existence of which 
was proved in theorem~5.6, the matrix of nilpotent operator $f$ has 
the following form:
$$
\hskip -2em
\spreadmatrixlines{1.7ex}
F=\Vmatrix
J_{k_1}(0) &\          &\      & \\
\          &J_{k_2}(0) &\      &\ \\
\          &\          &\ddots &\ \\
\          &\          &\      &J_{k_{\ssize s}}(0)
\endVmatrix.
\tag5.21
$$
The matrix \thetag{5.21} is blockwise-diagonal, its diagonal blocks
are Jordan cages of the form \thetag{5.14}, all other space in this
matrix is filled with zeros. It is easy to understand this fact. 
Indeed, each chain with the side vector $\bold e_i$ produces the 
invariant subspace $U(v)$ of the form \thetag{5.12}, where $\bold v
=\bold e_{i,k_{\ssize i}}$. Due to the theorem~5.6 the space $V$
is the direct sum of such invariant subspaces:
$$
V=U(\bold e_{1,k_1})\oplus\ldots\oplus U(\bold e_{s,k_{\ssize s}}).
$$\par
     The matrix \thetag{5.21} is called a {\it Jordan form}
of the matrix of a nilpotent operator. The theorem~5.6 is known
as {\it the theorem on bringing the matrix of a nilpotent operator to
a canonic Jordan form}. If the chain basis $\bold e_1,\,\ldots,\,\bold 
e_s$ is constructed strictly according to the proof of the theorem~5.6,
then Jordan cages are arranged in the order of decreasing sizes:
$$
k_1\geqslant k_2\geqslant\ldots\geqslant k_s.
$$
However, the permutation of vectors $\bold e_1,\,\ldots,\,\ldots e_s$ 
can change this order, and this usually happens in practice.\par
\proclaim{Theorem 5.7} The height of a nilpotent operator $f$ in a
finite-dimensional space $V$ is less or equal to the dimension
$n=\dim V$ of this space and $f^n=0$.
\endproclaim
\demo{Proof} Above in proving the theorem~5.1 we noted that
the height $\nu(f)$ of a nilpotent operator $f$ coincides with
the greatest height of basis vectors. Due to the theorem~5.6 now
we can choose the chain basis. The height of a chain vector
is not greater than the length of the chain \thetag{5.5} to which 
it belongs. Therefore, the height of basis vectors in a chain basis
is not greater than the number of vectors in such a basis. This
yields $\nu(f)\leqslant n=\dim V$. The height of an arbitrary vector
$\bold v$ of $V$ is not greater than the height of the operator $f$.
Therefore, $f^n(\bold v)=\bold 0$ for all $\bold v\in V$. This means
that $f^n=0$. The theorem is proved.
\qed\enddemo
\head
\S\,6. Root subspaces.
Two theorems\\ on the sum of root subspaces.
\endhead
\rightheadtext{\S\,6. Root subspaces.}
\definition{Definition 6.1} The {\it root subspace} of a linear operator
$f\!:\,V\to V$ correspon\-ding to its eigenvalue $\lambda$ is the set
$$
V(\lambda)=\{\bold v\in V:\ \exists k\,((k\in\Bbb N)\and ((f-\lambda
\cdot 1)^k\,\bold v=\bold 0))\}
$$
that consist of vectors vanishing under the action of some positive
integer power of the operator $h_\lambda=f-\lambda\cdot 1$.
\enddefinition
     For each positive integer $k$ we define the subspace   
$V(k,\lambda)=\Ker (h_\lambda)^k$. For $k=1$ the subspace $V(1,\lambda)$ coincides with the eigenspace $V_\lambda$. Note that $(h_\lambda)^k\,\bold
v=\bold 0$ implies $(h_\lambda)^{k+1}\,\bold v=\bold 0$.
Therefore we have the sequence of inclusions
$$
\hskip -2em
V(1,\lambda)\subseteq V(2,\lambda)\subseteq\ldots\subseteq
V(k,\lambda)\subseteq\ldots
\tag6.1
$$
It is easy to see that all subspaces in the sequence \thetag{6.1}
are enclosed into the root subspace $V(\lambda)$. Moreover, 
$V(\lambda)$ is the union of the subspaces \thetag{6.1}:
$$
\hskip -2em
V(\lambda)=\bigcup^\infty_{k=1}V(k,\lambda)=
\sum^\infty_{k=1}V(k,\lambda).
\tag6.2
$$
In this case the sum of subspaces the sum of subspaces $V(k,\lambda)$
coincides with their union. Indeed, let $\bold v$ be a vector of the
sum of subspaces $V(k,\lambda)$. Then
$$
\hskip -2em
\bold v=\bold v_{k_1}+\ldots+\bold v_{k_{\ssize s}}\text{, where \ }
\bold v_{k_{\ssize s}}\in V(k_s,\lambda).
\tag6.3
$$
Let $k=\max\{k_1,\ldots,k_s\}$, then from the sequence of inclusions
\thetag{6.1} we derive $\bold v_{k_{\ssize i}}\in V(k,\lambda)$.
Therefore the vector \thetag{6.3} belongs to $V(k,\lambda)$, hence,
it belongs to the union of all subspaces $V(k,\lambda)$.\par
    The proof of coincidence of the sum and the union in \thetag{6.2}
is based on the inclusions \thetag{6.1}. Therefore, we have proved
the more general theorem.
\proclaim{Theorem 6.1} The sum of a growing sequence of mutually
enclosed subspaces coincides with their union.
\endproclaim
     The theorem~6.1 shows that the set $V(\lambda)$ in definition~6.1 
is actually a subspace in $V$. This subspace is nonzero since it 
comprises the eigenspace $V_\lambda$ as a subset.\par
\proclaim{Theorem 6.2} A root subspace $V(\lambda)$ of an operator $f$ 
is invariant under the action of $f$ and of all operators from its
polynomial envelope $P(f)$.
\endproclaim
\demo{Proof} Let $\bold v\in V(\lambda)$. Then there exists a positive
integer number $k$ such that $(h_\lambda)^k\,\bold v=\bold 0$. Let's
consider the vector $\bold w=f(\bold v)$. For this vector we have 
$$
(h_\lambda)^k\,\bold w=(h_\lambda)^k\compos f\,\bold v
=f\compos (h_\lambda)^k\,\bold v=f((h_\lambda)^k(\bold v))
=\bold 0.
$$
Here we used the permutability of the operators $h_\lambda$ and 
$f$, it follows from the inclusion $h_\lambda\in P(f)$. Due to
the above equality we have $\bold w=f(\bold v)\in V(\lambda)$.
The invariance of $V(\lambda)$ under the action of $f$ is proved. 
Its invariance under the action of operators from $P(f)$ now
follows from the theorem~4.5.
\qed\enddemo
\proclaim{Theorem 6.3} Let $\lambda$ and $\mu$ be two eigenvalues of 
a linear operator $f\!:\,V\to V$. Then the restriction of the operator
$h_\lambda=f-\lambda\cdot 1$ to the root subspace $V(\mu)$ is 
\roster
\item a bijective operator if $\mu\neq\lambda$;
\item a nilpotent operator if $\mu=\lambda$.
\endroster
\endproclaim
\demo{Proof} Let's prove the first proposition of the theorem.
We already know that the subspace $V(\mu)$ is invariant under the
action of $h_\lambda$. For the sake of convenience we denote by
$h_{\lambda,\mu}$ the restriction of $h_\lambda$ to the subspace
$V(\mu)$. This is an operator acting from $V(\mu)$ to $V(\mu)$. 
Let's find its kernel:
$$
\Ker h_{\lambda,\mu}=\{\bold v\in V(\mu)\!:\ h_\lambda(\bold v)
=\bold 0\}=\Ker h_\lambda\cap V(\mu).
$$
The kernel of the operator $h_\lambda$ by definition coincides with
the eigenspace $V_\lambda$. Therefore, $\Ker h_{\lambda,\mu}=V_\lambda
\cap V(\mu)$.\par
      Let $\bold v$ be an arbitrary vector of the kernel $\Ker h_{\lambda,
\mu}$. Due to the above result $\bold v$ belongs to $V_\lambda$. Therefore,
we have the equality 
$$
\hskip -2em
f(\bold v)=\lambda\cdot\bold v.
\tag6.4
$$
Simultaneously, we have the other condition $\bold v\in V(\mu)$ which 
means that there exists some integer number $k>0$ such that
$$
\hskip -2em
(h_\mu)^k\,\bold v=(f-\mu\cdot 1)^k\,\bold v=\bold 0.
\tag6.5
$$
From \thetag{6.4} we get $h_\mu(\bold v)=f(\bold v)-\mu\cdot\bold v
=(\lambda-\mu)\cdot\bold v$. Combining this equality with \thetag{6.5},
we obtain the following equality for $\bold v$:
$$
(h_\mu)^k\,\bold v=(\lambda-\mu)^k\cdot\bold v=\bold 0.
$$
Therefore, if $\lambda\neq\mu$, we immediately get $\bold v=0$, which means
that $\Ker h_{\lambda,\mu}=\{0\}$. Hence, in the case $\lambda\neq\mu$ the
operator $h_{\lambda,\mu}\!:\,V(\mu)\to V(\mu)$ is injective. 
The surjectivity of this operator and, hence, its bijectivity follows from
its injectivity due to the theorem~1.3.\par
    Now let's prove the second proposition of the theorem. In this case
$\mu=\lambda$, therefore, we consider the operator $h_{\lambda,\lambda}$
being the restriction of $h_\lambda$ to the subspace $V(\lambda)$. Note
that $h_{\lambda,\lambda}\,\bold v=h_\lambda\,\bold v$ for all $\bold v
\in V(\lambda)$. Therefore, from the definition of a root subspace we
conclude that for any vector $\bold v\in V(\lambda)$ \pagebreak 
there is a positive integer number $k$ such that $(h_{\lambda,\lambda})^k\,
\bold v=(f-\lambda\cdot 1)^k\,\bold v=\bold 0$. This equality means that
$h_{\lambda,\lambda}$ is a nilpotent operator in $V(\lambda)$. The theorem
is proved.
\qed\enddemo
\proclaim{Theorem 6.4} Let $\lambda_1,\,\ldots,\,\lambda_s$ be a set of
mutually distinct eigenvalues of a linear operator $f\!:\,V\to V$. Then
the sum of corresponding root subspaces is a direct sum:
$V(\lambda_1)+\ldots+V(\lambda_s)=V(\lambda_1)\oplus\ldots\oplus
V(\lambda_s)$.
\endproclaim
\demo{Proof} The proof of this theorem is similar to that of theorem~4.6.
Denote by $W$ the sum of subspaces specified in the theorem:
$$
\hskip -2em
W=V(\lambda_1)+\ldots+ V(\lambda_s).
\tag6.6
$$
In order to prove that the sum \thetag{6.6} is a direct sum we should
prove the uniqueness of the following expansion for an arbitrary vector
$\bold w\in W$:
$$
\hskip -2em
\bold w=\bold v_1+\ldots+\bold v_s\text{, \ where \ }\bold v_i
\in V(\lambda_i).
\tag6.7
$$
Consider another expansion of the same sort for the same vector 
$\bold w$:
$$
\hskip -2em
\bold w=\tilde\bold v_1+\ldots+\tilde\bold v_s\text{, \ where
\ }\tilde\bold v_i\in V(\lambda_i).
\tag6.8
$$
Then let's subtract the second expansion from the first one and for
the sake of brevity denote $\bold w_i=(\bold v_i-\tilde\bold v_i)
\in V(\lambda_i)$. As a result we get
$$
\hskip -2em
\bold w_1+\ldots+\bold w_s=\bold 0.
\tag6.9
$$\par
     Denote $h_r=f-\lambda_r\cdot 1$. According to the definition of
the root subspace $V(\lambda_r)$, for any vector $\bold w_r$ in the
expansion \thetag{6.9} there is some positive integer number $k_r$
such that $(h_r)^{k_{\ssize r}}\,\bold w_r=\bold 0$. We use this fact
and define the operators
$$
\hskip -2em
f_i=\prod^s_{r\neq i}(h_r)^{k_{\ssize r}}.
\tag6.10
$$
Due to the permutability of the operators $\bold h_1,\,\ldots,\,\bold 
h_s$ belonging to the polynomial envelope of the operator $f$ and due
to the equality $(h_r)^{k_{\ssize r}}\,\bold w_r=\bold 0$ we get
$$
f_i(\bold w_j)=\bold 0\text{ \ for all \ }j\neq i.
$$
Let's apply the operator \thetag{6.10} to both sides of the equality
\thetag{6.9}. Then all terms in the sum in left hand side of this
equality do vanish, except for $i$-th term only. This yields
$f_i(\bold w_i)=\bold 0$. Let's write this equality in expanded form:
$$
\hskip -2em
\left(\,\shave{\prod^s_{r\neq i}}(h_r)^{k_{\ssize r}}
\right)\bold w_i=\bold 0
\tag6.11
$$
The vector $\bold w_i$ belongs to the root space $V(\lambda_i)$, which
is invariant under the action of all operators $h_r$ in \thetag{6.11}.
Therefore we can replace the operators $h_r$ in \thetag{6.11} by their
restrictions $h_{r,i}$ to the subspace $V(\lambda_i)$:
$$
\pagebreak 
\hskip -2em
\left(\,\shave{\prod^s_{r\neq i}}(h_{r,i})^{k_{\ssize r}}
\right)\bold w_i=\bold 0.
\tag6.12
$$
According to the theorem~6.3, the restricted operators $h_{r,i}$ are
bijective if $r\neq i$. The product (the composition) of bijective 
operators is bijective. 
We also know that applying a bijective operator to nonzero vector 
we would get a nonzero result. Therefore, \thetag{6.12}
implies $\bold w_i=\bold 0$. Then $\bold v_i=\tilde\bold v_i$ and the
expansions \thetag{6.7} and \thetag{6.8} do coincide. The uniqueness 
of the above expansion \thetag{6.7} and the theorem in whole are proved.
\qed\enddemo
\proclaim{Theorem 6.5} Let $f$ be a linear operator in a finite-dimensional
space $V$ over the field $\Bbb K$ and suppose that its characteristic
polynomial factorizes into a product of linear terms in $\Bbb K$. Then the
sum of all root subspaces of the operator $f$ is equal to $V$, i\.\,e\.
$V(\lambda_1)\oplus\ldots\oplus V(\lambda_s)=V$, where 
$\lambda_1,\,\ldots,\,\lambda_s$ is the set of all mutually distinct
eigenvalues of the operator $f$.
\endproclaim
\demo{Proof} Since $\lambda_1,\,\ldots,\,\lambda_s$ is the set of all
mutually distinct eigenvalues of the operator $f$, for its characteristic
polynomial we get
$$
\det(f-\lambda\cdot 1)=\prod^s_{i=1}(\lambda_i-\lambda)^{n_{\ssize i}}.
$$
According to the hypothesis of theorem, it is factorized into a product
of linear polynomials of the form $\lambda_i-\lambda$, where $\lambda_i$
is an eigenvalue of $f$ and $n_i$ is the multiplicity of this eigenvalue.
Let's denote by $W$ the total sum of all root subspaces of the operator 
$f$, we know that this is a direct sum (see theorem~6.4): 
$$
W=V(\lambda_1)\oplus\ldots\oplus V(\lambda_s).
$$
The root subspaces are nonzero, hence, $W\neq\{\bold 0\}$.\par
     Further proof is by contradiction. Assume that the proposition of 
the theorem is false and $W\neq V$. The subspace $W$ is invariant under 
the action of $f$ as a sum of invariant subspaces $V(\lambda_i)$ (see
theorem~3.2). Due to the theorem~4.5 it is invariant under the action
of the operator $h_\lambda=f-\lambda\cdot 1$ as well. Let's apply the
theorem~3.5 to the operator $h_\lambda$. This yields
$$
\hskip -2em
\det(f-\lambda\cdot 1)=
\det(f{\!\lower4pt\hbox{$\ssize W$}}-\lambda\cdot 1)\,
\det(f{\!\lower4pt\hbox{$\ssize V/W$}}-\lambda\cdot 1).
\tag6.13
$$
Here we took into account that $1{\lower4pt\hbox{$\ssize W$}}=1$
and $1{\lower4pt\hbox{$\ssize V/W$}}=1$, we also used the theorem~3.4. The characteristic polynomial of the operator $f$ is the product of 
characteristic polynomial of restricted operator $f{\lower4pt\hbox{$\ssize
W$}}$ and that of factoroperator $f{\lower4pt\hbox{$\ssize V/W$}}$. The 
left hand side of \thetag{6.13} factorizes into a product of linear
polynomials in $\Bbb K$, therefore, each of the polynomials in right hand 
side of \thetag{6.13} should do the same. Let $\lambda_q$ be one of the
eigenvalues of the factoroperator $f{\!\lower4pt\hbox{$\ssize V/W$}}$
and let $Q\in V/W$ be the corresponding eigenvector. Due to \thetag{6.13}
the number $\lambda_q$ is in the list $\lambda_1,\ldots,\lambda_s$ of
eigenvalues of the operator $f$. Due to our assumption $W\neq V$ we conclude that the factorspace $V/W$ is nontrivial: $V/W\neq\{\bold 0\}$, and the
coset $Q$ is not zero. Suppose that $\bold v\in Q$ is a representative 
of this coset $Q$. Since $Q\neq\bold 0$, we have $\bold v\not\in W$. 
The coset $Q$ is an eigenvector of the factoroperator
$f{\lower4pt\hbox{$\ssize V/W$}}$, therefore, it should satisfy the 
following equality:
$$
\hskip -2em
(f{\lower4pt\hbox{$\ssize V/W$}}-\lambda_q\cdot 1)\,Q=
\Cl_W((f-\lambda_q\cdot 1)\,\bold v)=\bold 0.
\tag6.14
$$
Let's denote $h_r=f-\lambda_r\cdot 1$ for all $r=1,\,\ldots,\,s$ 
(we have already used this notation in proving the previous theorem).
The relationship \thetag{6.14} means that
$$
\hskip -2em
(f-\lambda_q\cdot 1)\,\bold v=h_q(\bold v)=\bold w\in W.
\tag6.15
$$
From the expansion $W=V(\lambda_1)\oplus\ldots\oplus V(\lambda_s)$ 
for the vector $\bold w$, which arises in formula \thetag{6.15},
we get the expansion
$$
\hskip -2em
h_q(\bold v)=\bold w=\bold v_1+\ldots+\bold v_s\text{, \ where \ }
\bold v_i\in V(\lambda_i).
\tag6.16
$$
Let's consider the restriction of the operator $h_q$ to the root
subspace $V(\lambda_i)$, this restriction is denoted
$h_{q,i}$ (see the proof of theorem~6.4). Due to the theorem~6.3 we
know that the operators $h_{q,i}\!:\,V(\lambda_i)\to V(\lambda_i)$ 
are bijective for all $i\neq q$. Therefore, for all $\bold v_1,\,
\ldots,\,\bold v_s$ in \thetag{6.16} other than $\bold v_q$ we can
find $\tilde\bold v_i\in V(\lambda_i)$ such that $\bold v_i=
h_{q,i}(\tilde\bold v_i)$. Let's substitute these expressions
into \thetag{6.16}. Then we get
$$
\hskip -2em
\bold w=h_q(\bold v)=\bold v_q+\sum^s_{i\neq q}h_q(\tilde\bold v_i).
\tag6.17
$$
Relying upon this formula \thetag{6.17}, we define the new vector
$\tilde\bold v_q$:
$$
\hskip -2em
\tilde\bold v_q=\bold v-\sum^s_{i\neq q}\tilde\bold v_i.
\tag6.18
$$
For this vector from \thetag{6.17} we derive $h_q(\tilde\bold v_q)
=\bold v_q\in V(\lambda_q)$. Due to the definition of the root subspace
$V(\lambda_q)$ there exists a positive integer number $k$ such that
$(h_q)^k\,\bold v_q=\bold 0$. Hence, $(h_q)^{k+1}\,\tilde\bold v_q
=\bold 0$ and, therefore, $\tilde\bold v_q\in V(\lambda_q)$. Returning 
back to the formula \thetag{6.18}, we derive
$$
\hskip -2em
\bold v=\sum^s_{i=1}\tilde\bold v_i\text{, \ where \ }\bold v_i
\in V(\lambda_i).
\tag6.19
$$
From the formula \thetag{6.19} and from the expansion$W=V(\lambda_1)\oplus
\ldots\oplus V(\lambda_s)$ it follows that $\bold v\in W$, but this
contradicts to our initial choice $\bold v\not\in W$, which was possible 
due to the assumption $W\neq V$. Hence, $W=V$. The theorem is proved.
\qed\enddemo
\head
\S\,7. Jordan basis of a linear operator.\\
Hamilton-Cayley theorem.
\endhead
\rightheadtext{\S\,7. Jordan basis of a linear operator.}
     Let $f\!:\,V\to V$ be a linear operator in finite-dimensional
linear vector space $V$. Suppose that $V$ is expanded into the sum
of root subspaces of the operator $f$:
$$
\hskip -2em
V=V(\lambda_1)\oplus\ldots\oplus V(\lambda_s).
\tag7.1
$$
Let's denote $h_i=f-\lambda_i\cdot 1$. Then denote by $h_{i,j}$ the
restriction of $h_i$ to $V(\lambda_j)$. According to the theorem~6.3,
the restriction $h_{i,i}$ is a nilpotent operator in $i$-th root 
subspace $V(\lambda_i)$. Therefore, in $V(\lambda_i)$ we can choose 
a canonic Jordan basis for this operator (see theorem~5.6). \pagebreak
The matrix
of the operator $h_{i,i}$ in canonic Jordan basis is a matrix of the
form \thetag{5.21} composed by diagonal blocks, where each diagonal
block is a matrix of the form \thetag{5.14}.
\definition{Definition 7.1} A {\it Jordan normal basis} of an operator
$f\!:\,V\to V$ is a basis composed by canonic Jordan bases of nilpotent
operators $h_{i,i}$ in the root subspaces $V(\lambda_i)$ of the operator
$f$.
\enddefinition
    Note that an operator $f$ in a finite dimensional space $V$ possesses 
a Jordan normal basis if and only if there $V$ is expanded into the sum
of root subspaces of the operator $f$, i\.\,e\. if we have \thetag{7.1}.
The theorem~6.5 yields a sufficient condition for the existence of a 
Jordan normal basis of a linear operator.\par
     Suppose that an operator $f$ in a finite-dimensional linear vector
space $V$ possesses a Jordan normal basis. The subspaces $V(\lambda_i)$
in \thetag{7.1} are invariant with respect to $f$. Let's denote by $f_i$
the restriction of $f$ to $V(\lambda_i)$. The matrix of the operator 
$f$ in a Jordan normal basis is a blockwise-diagonal matrix: 
$$
\hskip -2em
F=\Vmatrix
F_1 &\   &\      & \\
\   &F_2 &\      &\ \\
\   &\   &\ddots &\ \\
\   &\   &\      &F_s
\endVmatrix,
\tag7.2
$$
The diagonal blocks $F_i$ in \thetag{7.2} are determined by operators $f_i$.
Note that the operators $f_i$ and $h_{i,i}$ are related to each other by the
equality $f_i=h_{i,i}+\lambda_i\cdot 1$. Therefore, $F_i$ is also a 
blockwise-diagonal matrix:
$$
\spreadmatrixlines{1.7ex}
\hskip -2em
F_i=\Vmatrix
J_{k_1}(\lambda_i) &\          &\      & \\
\          &J_{k_2}(\lambda_i) &\      &\ \\
\          &\          &\ddots &\ \\
\          &\          &\      &J_{k_{\ssize r}}(\lambda_i)
\endVmatrix.
\tag7.3
$$
The number of diagonal blocks in \thetag{7.3} is determined the number
of chains in a canonic Jordan basis of the nilpotent operator $h_{i,i}$,
while these diagonal blocks themselves are matrices of the following form:
$$
\hskip -2em
J_k(\lambda)=\Vmatrix
\lambda & 1 & \      &0 \\
\ & \lambda & \ddots &\ \\
\ & \ & \ddots &1 \\
\ & \ & \      &\lambda \\
\endVmatrix.
\tag7.4
$$
A matrix of the form \thetag{7.4} is called a {\it Jordan block}
or a {\it Jordan cage} with $\lambda$ on the diagonal. This is
square $k\times k$ matrix; if $k=1$ this matrix degenerates and 
becomes a matrix with the single element $J_1(\lambda)=\Vmatrix
\lambda\endVmatrix$.\par
     The matrix of an operator $f$ in a Jordan normal base 
presented by the relationships \thetag{7.2}, \thetag{7.3}, and
\thetag{7.4} is called a {\it Jordan normal form} of the matrix 
of this operator. The problem of constructing a Jordan
normal basis for a linear operator $f$ and thus finding the
Jordan normal form $F$ of its matrix is known as the problem of
{\it bringing the matrix of a linear operator to a Jordan normal
form}.\par
     If the matrix of a linear operator can be brought to a Jordan
normal form, this fact has several important consequences. Note that
a matrix of the form\thetag{7.4} is upper-triangular. Hence, 
\thetag{7.3} and \thetag{7.2} all are upper-triangular matrices.
The entries on the diagonal of \thetag{7.2} are the eigenvalues
of the operator $f$, the $i$-th eigenvalue $\lambda_i$ being 
presented $n_i$ times, where $n_i=\dim V(\lambda_i)$. From the course
of algebra we know that the determinant of an upper-triangular matrix
is equal to the product of all its diagonal elements. Therefore, the
characteristic polynomial of an operator possessing a Jordan normal 
basis is given by the formula 
$$
\hskip -2em
\det(f-\lambda\cdot 1)=\prod^s_{i=1}
(\lambda_i-\lambda)^{n_{\ssize i}}.
\tag7.5
$$
\proclaim{Theorem 7.1} The matrix of a linear operator $f$ in a 
finite-dimensional linear vector space $V$ over a numeric field 
$\Bbb K$ can be brought to a Jordan normal form if and only if its 
characteristic polynomial factorizes into the product of linear 
polynomials in the field $\Bbb K$. 
\endproclaim
\demo{Proof} The necessity of the condition formulated in the 
theorem~7.1 is immediate from \thetag{7.5}; the sufficiency is 
provided by the theorems~5.6 and 6.5.
\qed\enddemo
     In the case of the field of complex numbers $\Bbb C$ any polynomial
factorizes into a product of linear terms. Therefore, the matrix of any
linear operator in a complex linear vector space can be brought to a
Jordan normal form.
\proclaim{Theorem 7.2} The multiplicity of an eigenvalue $\lambda$ of
a linear operator $f$ in a finite-dimensional linear vector space $V$
is equal to the dimension of the correspon\-ding root subspace 
$V(\lambda)$.
\endproclaim
     For the operator $f$, the characteristic polynomial of which 
factorizes into the product of linear terms, the proposition of 
theorem~7.2 immediately follows from the formula \thetag{7.5}. 
However, this fact is valid also in the case of partial factorization.
Such a case can be reduced to the case of complete factorization by
means of the field extension technique. We do not consider the field 
extension technique in this book. But it is worth to note that the
complete proof of the following Hamilton-Cayley theorem is also based 
on that technique.
\proclaim{Theorem 7.3} Let $P(\lambda)$ be the characteristic polynomial
of a linear operator $f$ in a finite-dimensional space $V$. Then $P(f)=0$.
\endproclaim
\demo{Proof} We shall prove the Hamilton-Cayley theorem for the case
where the characteristic polynomial $P(\lambda)$ factorizes into the 
product of linear terms:
$$
\hskip -2em
P(\lambda)=\prod^s_{i=1}(\lambda_i-\lambda)^{n_{\ssize i}}.
\tag7.6
$$
Denote $h_i:=f-\lambda_i\cdot 1$ and denote by $h_{i,j}$ the restriction
of $h_i$ to the root subspace $V(\lambda_j)$. Then from the formula \thetag{7.6} we derive 
$$
P(f)=\prod^s_{i=1}(h_i)^{n_{\ssize i}}.
$$
Let's apply $P(f)$ to an arbitrary vector $\bold v\in V$. Due to the
theorem \thetag{6.5} we can expand $\bold v$ into a sum $\bold v
=\bold v_1+\ldots+\bold v_s$, where $\bold v_i\in V(\lambda_i)$. 
Therefore, we have
$$
\hskip -2em
P(f)\,\bold v=P(f)\,\bold v_1+\ldots+P(f)\,\bold v_s.
\tag7.7
$$
The root subspace $V(\lambda_j)$ is invariant under the action of
the operators $h_i$. Then
$$
P(f)\,\bold v_j=\prod^s_{i=1}(h_i)^{n_{\ssize i}}\,\bold v_j=
\prod^s_{i=1}(h_{i,j})^{n_{\ssize i}}\,\bold v_j.
$$
Using permutability of the operators $h_i$ and their restrictions
$h_{i,j}$, we can bring the above expression for $P(f)\,\bold 
v_j$ to the following form:
$$
\hskip -2em
P(f)\,\bold v_j=\prod^s_{i\neq j}(h_{i,j})^{n_{\ssize i}}\,
(h_{j,j})^{n_{\ssize j}}\,\bold v_j.
\tag7.8
$$
The operator $h_{j,j}$ is a nilpotent operator in the subspace
$V(\lambda_j)$ and $n_j=\dim V(\lambda_j)$. Therefore, we can 
apply the theorem~5.7. As a result we obtain
$(h_{j,j})^{n_{\ssize j}}\,\bold v_j=\bold 0$. 
Now from \thetag{7.7} and \thetag{7.8} for an arbitrary vector 
$\bold v\in V$ we derive $P(f)\,\bold v=\bold 0$. This proves
the theorem for the special case, where the characteristic 
polynomial of an operator $f$ factorizes into a product of
linear terms. The general case is reduced to this special case 
by means of the field extension technique, which we do not
consider in this book.
\qed\enddemo
\newpage
\topmatter
\title\chapter{3}
Dual space.
\endtitle
\endtopmatter
\document
\head
\S\,1. Linear functionals.\\
Vectors and covectors. Dual space.
\endhead
\rightheadtext{\S\,1. Linear functionals. Dual space.}
\leftheadtext{CHAPTER~\uppercase\expandafter{\romannumeral 3}.
DUAL SPACE.}
\setfirstpage
\definition{Definition 1.1} Let $V$ be a linear vector space
over a numeric field $\Bbb K$. A numeric function $y=f(\bold v)$ 
with vectorial argument $\bold v\in V$ and with values $y\in
\Bbb K$ is called a {\it linear functional\/} if
\roster
\item $f(\bold v_1+\bold v_2)=f(\bold v_1)+f(\bold v_2)$ for any
      two $\bold v_1,\bold v_2\in V$;
\item $f(\alpha\cdot\bold v)=\alpha\,f(\bold v)$ for any $\bold v
       \in V$ and for any $\alpha\in\Bbb K$.
\endroster
\enddefinition
     The definition of a linear functional is quite similar to the
definition of a linear mapping (see definition~8.1 in 
Chapter~\uppercase\expandafter{\romannumeral 1}). Comparing these
two definitions, we see that any linear functional $f$ is a linear
mapping $f\!:\,V\to\Bbb K$ and, conversely, any such linear mapping
is a linear functional. Thereby the numeric field $\Bbb K$ is treated 
as a linear space of the dimension $1$ over itself.\par
     Linear functionals, as linear mappings from $V$ to $\Bbb K$, 
constitute the space $\Hom(V,\Bbb K)$, which is called the {\it dual 
space\/} or the {\it conjugate space\/} for the space $V$. The dual
space $\Hom(V,\Bbb K)$ is denoted by $V^*$. The space of homomorphisms 
$\Hom(V,W)$ is usually determined by two spaces $V$ and $W$. However,
the dual space is an exception $V^*=\Hom(V,\Bbb K)$, it is determined
only by $V$ since $\Bbb K$ is known whenever $V$ is given (see 
definition~2.1 in Chapter~\uppercase\expandafter{\romannumeral 1}).\par
     Thus, $V^*=\Hom(V,\Bbb K)$ is a linear vector space over the same
numeric field $\Bbb K$ as $V$. If $V$ is finite-dimensional, then the 
dimension of the conjugate space is determined by the theorem~10.4 in
Chapter ~\uppercase\expandafter{\romannumeral 1}: \ $\dim V^*=\dim V$. 
The structure of a linear vector space in $V^*=\Hom(V,\Bbb K)$ is
determined by two algebraic operations: the operation of pointwise
addition and pointwise multiplication by numbers (see definitions~10.1 
and 10.2 in Chapter~\uppercase\expandafter{\romannumeral 1}). However,
it would be worth to formulate these two definitions especially for the
present case of linear functionals.
\definition{Definition 1.2} Let $f$ and $g$ be two linear functionals
of $V^*$. The {\it sum\/} of functionals $f$ and $g$ is a functional 
$h$ whose values are determined by formula $h(\bold v)=f(\bold v)
+g(\bold v)$ for all $\bold v\in V$.
\enddefinition
\definition{Definition 1.3} Let $f$ be a linear functional of $V^*$. 
The product of the functional $f$ by a number $\alpha\in\Bbb K$ 
is a functional $h$ whose values are determined by formula
$h(\bold v)=\alpha\cdot f(\bold v)$ for all $\bold v\in V$.
\enddefinition
     Let $V$ be a finite-dimensional vector space over a field $\Bbb K$ 
and let $\bold e_1,\,\ldots,\,\bold e_n$ be a basis in $V$. Then each 
vector $\bold v\in V$ can be expanded in this basis:
$$
\hskip -2em
\bold v=v^1\cdot\bold e_1+\ldots+v^n\cdot\bold e_n.
\tag1.1
$$
Let's consider $i$-th coordinate of the vector $\bold v$. Due to the
uniqueness of the expansion \thetag{1.1}, when the basis is fixed, 
$v^i$ is a number uniquely determined by the vector $\bold v$. Hence,
we can consider a map $h^i\!:\,V\to\Bbb K$, defining it by formula
$h^i(\bold v)=v^i$. When adding vectors, their coordinates are added;
when multiplying a vector by a number, its coordinates are multiplied
by that number (see the relationships \thetag{5.4} in
Chapter~\uppercase\expandafter{\romannumeral 1}). Therefore, 
$h^i\!:\,V\to\Bbb K$ is a linear mapping. This means that each 
basis $\bold e_1,\,\ldots,\,\bold e_n$ of a linear vector space $V$
determines $n$ linear functionals in $V^*$. The functionals 
$h^1,\,\ldots,\,h^n$ are called the {\it coordinate functionals} 
the basis $\bold e_1,\,\ldots,\,\bold e_n$. They satisfy the
relationships
$$
\hskip -2em
h^i(\bold e_j)=\delta^i_j,
\tag1.2
$$
where $\delta^i_j$ is the Kronecker symbol. These relationships 
\thetag{1.2} are called the {\it relationships of biorthogonality}.
\par
     The proof of the relationships of biorthogonality is very simple.
If we expand the vector $\bold e_j$ in the basis $\bold e_1,\,\ldots,\,
\bold e_n$, then its $j$-th component is equal to unity, while all
other components are equal to zero. Note that $h^i(\bold e_j)$ is a
number equal to $i$-th component of the vector $\bold e_j$. Therefore,
$h^i(\bold e_j)=1$ if $i=j$ and $h^i(\bold e_j)=0$ in all other cases.
\par
\proclaim{Theorem 1.1} Coordinate functionals $\bold h^1,\,\ldots,\,
\bold h^n$ are linearly independent; they form a basis in dual space
$V^*$.
\endproclaim
\demo{Proof} Let's consider a linear combination of the coordinate
functionals associated with a basis $\bold e_1,\,\ldots,\,\bold e_n$
in $V$ and assume that it is equal to zero:
$$
\hskip -2em
\alpha_1\cdot h^1+\ldots+\alpha_n\cdot h^n=0.
\tag1.3
$$
Right hand side of \thetag{1.3} is zero functional. Its value when
applied to the base vector $\bold e_j$ is equal to zero. Hence, we 
have
$$
\hskip -2em
\alpha_1\,h^1(\bold e_j)+\ldots+\alpha_n\,h^n(\bold e_j)=0.
\tag1.4
$$
Now we use the relationships of biorthogonality \thetag{1.2}. Due to
these relationships among $n$ terms $h^1(\bold e_j),\,\ldots,\,
h^n(\bold e_j)$ in left hand side of the equality \thetag{1.4} only
one term is nonzero: $h_j(\bold e_j)=1$. Therefore, \thetag{1.4}
reduces to $\alpha_j=0$. But $j$ is an index that runs from $1$
to $n$. Hence, all coefficients of the linear combination \thetag{1.3} 
are zero, i\.\,e\. it is trivial and coordinate functionals
$h^1,\,\ldots,\,h^n$ are linearly independent.\par
     In order to complete the proof of the theorem now we could use
the equality $\dim V^*=\dim V=n$ and refer to the item 
\therosteritem{4} of the theorem~4.5 in
Chapter~\uppercase\expandafter{\romannumeral 1}. However, we choose
more explicit way and directly prove that coordinate functio\-nals
$h^1,\,\ldots,\,h^n$ span the dual space $V^*$. Let $f\in V^*$ be
an arbitrary linear functional and let $\bold v$ be an arbitrary 
vector of $V$. Then from \thetag{1.1} we derive
$$
f(\bold v)=v^1\,f(\bold e_1)+\ldots+v^n\,f(\bold e_n)
=f(\,e_1)\,h^1(\bold v)+\ldots+f(e_n)\,h^n(\bold v).
$$
Here $f(\bold e_1),\,\ldots,\,f(\bold e_n)$ are numeric \pagebreak
coefficients from $\Bbb K$ and $\bold v$ is an arbitrary vector 
of $V$. Therefore, the above equality can be rewritten as an 
equality of linear functionals in the conjugate space $V^*$:
$$
\hskip -2em
f=f(\bold e_1)\cdot h^1+\ldots+f(\bold e_n)\cdot h^n.
\tag1.5
$$
The formula \thetag{1.5} shows that an arbitrary function $f\in V^*$
can be represented as a linear combination of coordinate functionals
$h^1,\,\ldots,\,h^n$. Hence, being linearly independent, they form
a basis in $V^*$. The theorem is proved.
\qed\enddemo
\definition{Definition 1.4} The basis $h^1,\ldots,h^n$ in $V^*$ formed
by coordinate functionals associated with a basis $\bold e_1,\,\ldots,
\,\bold e_n$ in $V$ is called the {\it dual basis\/} or the 
{\it conjugate basis\/} for $\bold e_1,\,\ldots,\,\bold e_n$.
\enddefinition
\definition{Definition 1.5} Let $f$ be a linear functional in a 
finite-dimensional space $V$ and let $\bold e_1,\,\ldots,\,\bold 
e_n$ be a basis in this space. The numbers $f_1,\,\ldots,\,
f_n$ determined by the linear functional $f$ according to the
formula
$$
\hskip -2em
f_i=f(\bold e_i)
\tag1.6
$$
are called the {\it coordinates\/} or the {\it components\/}  
of $f$ in the basis $\bold e_1,\,\ldots,\,\bold e_n$.
\enddefinition
     As we see in formula \thetag{1.5}, the numbers \thetag{1.6}
are the coefficients of the expansion of $f$ in the conjugate basis
$h^1,\,\ldots,\,h^n$. However, in the definition~1.5 they are mentioned 
as the components of $f$ in the basis $\bold e_1,\,\ldots,\,\bold e_n$.
This is purely terminological trick, it means that we consider 
$\bold e_1,\,\ldots,\,\bold e_n$ as a primary basis, while the
conjugate basis is treated as an auxiliary and complementary thing.
\par
     The algebraic operations of addition and multiplication by numbers
in the spaces $V$ and $V^*$ are related to each other by the following
equalities:
$$
\xalignat 2
&\hskip -2em
f(\bold v_1+\bold v_2)=f(\bold v_1)+f(\bold v_2), 
&&f(\alpha\cdot\bold v)=\alpha\,f(\bold v);\\
\vspace{-1.9ex}
&&&\tag1.7\\
\vspace{-1.9ex}
&\hskip -2em
(f_1+f_2)(\bold v)=f_1(\bold v)+f_2(\bold v), 
&&(\alpha\cdot f)(\bold v)=\alpha\,f(\bold v).
\endxalignat
$$
Vectors and linear functionals enter these equalities in a quite
similar way. The fact that in the writing $f(\bold v)$ the 
functional plays the role of a function, while the vector is 
written as an argument is not so important. Therefore, sometimes
the quantity $f(\bold v)$ is denoted differently:
$$
\hskip -2em
f(\bold v)=\langle f\,|\,\bold v\rangle.
\tag1.8
$$
The writing \thetag{1.8} is associated with the special terminology.
Functionals from the dual space $V^*$ are called {\it covectors},
while the expression $\langle f\,|\,\bold v\rangle$ itself is called
the {\it pairing}, or the {\it contraction}, or even the {\it scalar
product\/} of a vector and a covector.\par
     The scalar product \thetag{1.8} possesses the property of 
bilinearity: it is linear in its first argument $f$ and in its second
argument $\bold v$. This follows from the relationships \thetag{1.7}, 
which are now written as 
$$
\xalignat 2
&\hskip -2em
\langle f_1+f_2\,|\,\bold v\rangle=\langle f_1\,|\,\bold v\rangle+
\langle f_2\,|\,\bold v\rangle, && \langle \alpha\cdot f\,|\,\bold v
\rangle=\alpha\,\langle f\,|\,\bold v\rangle;\\
\vspace{-1.9ex}
&&&\tag1.9\\
\vspace{-1.9ex}
&\hskip -2em
\langle f\,|\,\bold v_1+\bold v_2\rangle=\langle f\,|\,
\bold v_1\rangle+\langle f\,|\,\bold v_2\rangle, 
&& \langle f\,|\,\alpha\cdot v\rangle=\alpha\,\langle f\,|\,
\bold v\rangle.
\endxalignat
$$
We have already dealt with the concept of bilinearity earlier in this 
book (see theorem~1.1 in Chapter~\uppercase\expandafter{\romannumeral 2}).
\par
     The properties \thetag{1.9} of the scalar properties \thetag{1.8}
are analogous to the properties of the scalar product of geometric
vectors --- it is usually studied in the course analytic geometry
(see \cite{5}). However, in contrast to that 
{\tencyr\char '074}geometric{\tencyr\char '076} scalar product, the
scalar product \thetag{1.8} is not symmetric: its arguments belong to different spaces --- they cannot be swapped. Covectors in the scalar
product \thetag{1.8} are always written on the left and vectors are
always on the right.\par
     The following definition is dictated by the intension to strengthen
the analogy of \thetag{1.8} and traditional 
{\tencyr\char '074}geometric{\tencyr\char '076} scalar product.
\definition{Definition 1.7} A vector $\bold v$ and a covector $f$ are
called {\it orthogonal} to each other if their scalar product is
zero: $\langle f\,|\,\bold v\rangle=0$.
\enddefinition
\proclaim{Theorem 1.2} Let $U\varsubsetneq V$ be a subspace in a 
finite-dimensional vector space $V$ and let $\bold v\not\in U$. 
Then there exists a linear functional $f$ in $V^*$ such that 
$f(\bold v)\neq 0$ and $f(\bold u)=0$ for all $\bold u\in U$.
\endproclaim
\demo{Proof} Let $\dim V=n$ and $\dim U=s$. Let's choose a basis
$\bold e_1,\,\ldots,\,\bold e_s$ in a subspace $U$. Let's add 
the vector $\bold v$ to basis vectors $\bold e_1,\,\ldots,\,\bold e_s$
and denote it $\bold v=\bold e_{s+1}$. The extended system of vectors 
is linearly independent since $\bold v\not\in U$, see the item 
\therosteritem{4} of the theorem~3.1 in
Chapter~\uppercase\expandafter{\romannumeral 1}. Denote by
$W=\langle\bold e_1,\,\ldots,\,\bold e_{s+1}\rangle$ the linear span
of this system of vectors. It is clear that $W$ is a subspace of $V$
comprising the initial subspace $U$; its dimension is one as greater
than the dimension of $U$. The vectors $\bold e_1,\,\ldots,\,\bold 
e_{s+1}$ form a basis in $W$. If $W\neq V$, then we complete the
basis $\bold e_1,\,\ldots,\,\bold e_{s+1}$ up to a basis 
$\bold e_1,\,\ldots,\,\bold e_n$ in the space $V$ and consider the
coordinate functionals $h^1,\ldots,h^n$ associated with this base.
Let's denote $f=h^{s+1}$. Then from the relationships of
biorthogonality \thetag{1.2} we derive 
$$
f(\bold v)=h^{s+1}(\bold e_{s+1})=1\text{ \ and \ }f(\bold e_i)=0
\text{ \ for \ }i=1,\ldots,s.
$$
Being zero on the basis vectors of the subspace $U$, the functional
$f=h^{s+1}$ vanishes on all vector $\bold u\in U$. Its value on
the vector $\bold v$ is equal to unity.
\qed\enddemo
     Let's consider the case $U=\{\bold 0\}$ in the above theorem.
Then for any nonzero vector $\bold v$ we have $\bold v\not\in U$,
and we can formulate the following corollary of the theorem~1.2.
\proclaim{Corollary} For any vector $\bold v\neq\bold 0$ in a 
finite-dimensional space $V$ there is a linear functional $f$ in
$V^*$ such that $f(\bold v)\neq 0$.
\endproclaim
     Let $V$ be a linear vector space over the field $\Bbb K$ and let
$W=V^*$ be the conjugate space of $V$. We know that $W$ is also a linear
vector space over the field $\Bbb K$. Therefore, it possesses its own
conjugate space $W^*$. With respect to $V$ this is the double conjugate
space $V^{**}$. We can also consider triple conjugate, fourth conjugate,
etc. Thus we would have the infinite sequence of conjugate spaces.
However, soon we shall see, that in the case of finite-dimensional spaces
there is no need to consider the multiple conjugate spaces.\par
     Let $\bold v\in V$. To any $f\in V^*$ we associate the number 
$f(\bold v)\in\Bbb K$. Thus we define a mapping $\varphi_{\bold
v}\!:\,V^*\to\Bbb K$, which is linear due to the following relationships:
$$
\align
&\varphi_{\bold v}(f_1+f_2)=(f_1+f_2)(\bold v)=f_1(\bold v)+f_2(\bold v)=\varphi_{\bold v}(f_1)
 +\varphi_{\bold v}(f_2),\\
\vspace{1.7ex}
&\varphi_{\bold v}(\alpha\cdot f)=(\alpha\cdot f)(\bold v)=\alpha\,
f(\bold v)=\alpha\,\varphi_{\bold v}(f).
\endalign
$$
Hence, $\varphi_{\bold v}$ is a linear functional in the space $V^*$ or, 
in other words, it is an element of double conjugate space. The functional
$\varphi_{\bold v}$ is determined by a vector $\bold v\in V$. Therefore,
when associating $\varphi_{\bold v}$ with a vector $\bold v$, we define
a mapping 
$$
\hskip -2em
h:V\to V^{**}\text{, \ where \ }h(\bold v)=\varphi_{\bold v}
\text{ \ for all \ }\bold v\in V.
\tag1.10
$$
The mapping \thetag{1.10} is a linear mapping. In order to prove this
fact we should verify the following identities for this mapping $h$:
$$
\xalignat 2
&\hskip -2em
h(\bold v_1+\bold v_2)=h(\bold v_1)+h(\bold v_2), 
&&h(\alpha\cdot\bold v)=\alpha\,h(\bold v).
\tag1.11
\endxalignat
$$
The result of applying $h$ to a vector of the space $V$ is an element
of double conjugate space $V^{**}$. Therefore, in order to verify the
equalities \thetag{1.11} we should apply both sides of these equalities 
to an arbitrary covector $f\in V^*$ and check the coincidence of the 
results that we obtain: 
$$
\align
 &\aligned
  h(\bold v_1&+\bold v_2)(f)=\varphi_{\bold v_1+\bold v_2}(f)
  =f(\bold v_1)+f(\bold v_2)=\varphi_{\bold v_1}(f)+\\
\vspace{1ex}
  &+\varphi_{\bold v_2}(f)=h(\bold v_1)(f)+h(\bold
  v_2)(f)=(h(\bold v_1)+h(\bold v_2))(f),
  \endaligned\\
\vspace{1.7ex}
&\aligned
  h(\alpha\cdot\bold v)(f)=\varphi_{\alpha\cdot\bold v}&(f)
  =f(\alpha\cdot\bold v)=\alpha\,f(\bold v)=\\
\vspace{1ex}
  &=\alpha\,\varphi_{\bold v}(f)=\alpha\,h(\bold v)(f)
  =(\alpha\cdot h(\bold v))(f).
  \endaligned
\endalign
$$
\proclaim{Theorem 1.3} For a finite-dimensional linear vector 
space $V$ the mapping \thetag{1.10} is bijective. It is an 
isomorphism of the spaces $V$ and $V^{**}$. This isomorphism 
is called {\it canonic isomorphism\/} of these spaces.
\endproclaim
\demo{Proof} First of all we shall prove the injectivity of the 
mapping \thetag{1.10}. For this purpose we consider its kernel
$\Ker h$. Let $\bold v$ be an arbitrary vector of $\Ker h$. 
Then $\varphi_{\bold v}=h(\bold v)=0$. But $\varphi_{\bold v}
\in V^{**}$, this means that $\varphi$ is a linear functional 
in the space $V^*$. Therefore, the equality $\varphi_{\bold v}=0$ 
means that $\varphi(f)=0$ for any covector $f\in V^*$. Using this
equality, from \thetag{1.10} we derive
$$
h(v)(f)=\varphi_v(f)=f(v)=0\text{ \ for all \ }f\in V^*.
\tag1.12
$$
Now let's apply the corollary of the theorem~1.2. If the vector 
$\bold v$ would be nonzero, then there would be a functional 
$f$ such that $f(\bold v)\neq 0$. This would contradict the above 
condition \thetag{1.12}. Hence, $\bold v=0$ by contradiction. 
This means that $\Ker h=\{\bold 0\}$ and $h$ is an injective 
mapping.\par
     In order to prove the surjectivity of the mapping \thetag{1.10}
we use the theorem~9.4 from
Chapter~\uppercase\expandafter{\romannumeral 1}. According to this
theorem
$$
\dim(\Ker h)+\dim(\Img h)=\dim V.
$$
We have already proved that $\dim(\Ker h)=0$. Hence, $\dim(\Img h)
=\dim V$. Since $\Img h$ is a subspace of $V^{**}$ and $\dim V^{**}
=\dim V^*=\dim V$, we have $\Img h=V^{**}$ (see item \therosteritem{3}
of theorem~4.5 in
Chapter~\uppercase\expandafter{\romannumeral 1}). This completes the
proof of surjectivity of the mapping $h$ and the proof of the theorem
\pagebreak in whole.
\qed\enddemo
     Canonic isomorphism \thetag{1.10} possesses the property that for
any vector $\bold v\in V$ and for any covector $f\in V^*$ the following
equality holds:
$$
\hskip -2em
\langle h(v)\,|\,f\rangle=\langle f\,|\,v\rangle.
\tag1.13
$$
The equality \thetag{1.13} is derived from the definition of $h$.
Indeed, $\langle h(\bold v)\,|\,f\rangle=h(\bold v)(f)=
\varphi_{\bold v}(f)=f(\bold v)=\langle f\,|\,\bold v\rangle$. 
The relationship \thetag{1.13} distinguishes canonic isomorphism
among all other isomorphisms relating the spaces $V$ and $V^{**}$.
\head
\S\,2. Transformation of the coordinates of a covector
under a change of basis.
\endhead
\rightheadtext{\S\,2. Transformation of the coordinates of 
a covector\dots}
     Let $V$ be a finite-dimensional linear vector space and let
$V^*$ be the associated dual space. If we treat $V^*$ separately
forgetting its relation to $V$, then a choice of basis and a change
of basis in $V^*$ are quite the same as in any other linear vector
space. However, the conjugate space $V^*$ is practically never 
considered separately. The theory of this space should be understood 
as an extension of the theory of initial space $V$.\par
     Let $\bold e_1,\,\ldots,\,\bold e_n$ be a basis in a linear
vector space $V$. Each such basis of $V$ has the associated basis
of coordinate functionals in $V^*$. Choosing another basis 
$\tilde\bold e_1,\,\ldots,\,\tilde\bold e_n$ in $V$ we immediately
get another conjugate basis $\tilde\bold h^1,\,\ldots,\,\tilde
\bold h^n$ in $V^*$. Let $S$ be the transition matrix for passing
from the old basis $\bold e_1,\,\ldots,\,\bold e_n$ to the new 
basis $\tilde\bold e_1,\,\ldots,\,\tilde\bold e_n$. Similarly,
denote by $P$ the transition matrix for passing from the old
dual basis $h^1,\,\ldots,\,h^n$ to the new dual basis $\tilde h^1,
\,\ldots,\,\tilde h^n$. The components of these two transition 
matrices $S$ and $P$ are used to expand the vectors of {\tencyr\char
'074}wavy{\tencyr\char '076} bases in corresponding {\tencyr\char
'074}non-wavy{\tencyr\char '076} bases:
$$
\xalignat 2
&\hskip -2em
\tilde e_j=\sum^n_{i=1} S^i_j\cdot e_i,
&&\tilde h^r=\sum^n_{s=1} P^r_s\cdot h^s.
\tag2.1
\endxalignat
$$
Note that the second formula \thetag{2.1} differs from the standard 
given by formula \thetag{5.5} in
Chapter~\uppercase\expandafter{\romannumeral 1}: the vectors of dual 
bases in \thetag{2.1} are specified by upper indices despite to the 
usual convention of enumerating the basis vectors. The reason is that
the dual space $V^*$ and the dual bases are treated as complementary 
objects with respect to the initial space $V$ and its bases. We have 
already seen such deviations from the standard notations in constructing
the basis vectors $E^i_j$ in $\Hom(V,W)$ (see proof of the theorem~10.4
in Chapter~\uppercase\expandafter{\romannumeral 1}).\par
     In spite of the breaking the standard rules in indexing the basis
vectors, the formula \thetag{2.1} does not break the rules of tensorial
notation: the free index $r$ is in the same upper position in both sides of
the equality, the summation index $s$ enters twice --- once as an upper
index and for the second time as a lower index.\par
\proclaim{Theorem 2.1} The transition matrix $P$ for passing from the old
conjugate\linebreak basis $h^1,\,\ldots,\,h^n$ to the new conjugate basis
$\tilde h^1,\,\ldots,\,\tilde h^n$ is inverse to the transition matrix $S$
that is used for passing from the old basis $\bold e_1,\,\ldots,\,\bold e_n$
to the new\linebreak basis $\tilde\bold e_1,\,\ldots,\,\tilde\bold e_n$.
\endproclaim
\demo{Proof} In order to prove this theorem we use the biorthogonality relationships \thetag{1.2}. Substituting \thetag{2.1} into these relationships, we get
$$
\delta^r_j=h^r(e_j)=\sum^n_{s=1}\sum^n_{i=1} P^r_s\,S^i_j\,h^s(e_i)
=\sum^n_{s=1}\sum^n_{i=1} P^r_s\,S^i_j\,\delta^s_i=
\sum^n_{i=1}P^r_i\,S^i_j.
$$
The above relationship can be written in matrix form $P\,S=1$. This means
that $P=S^{-1}$. The theorem is proved.
\qed\enddemo
     Remember that the inverse transition matrix $T$ is also the inverse
matrix for $S$. Therefore, in order to write the complete set of formulas
relating two pairs of bases in $V$ and $V^*$ it is sufficient to know two
matrices $S$ and $T=S^{-1}$:
$$
\xalignat 2
&\hskip -2em
\tilde\bold e_j=\sum^n_{i=1}S^i_j\cdot\bold e_i,
&&\tilde h^r=\sum^n_{s=1}T^r_s\cdot h^s,\\
\vspace{-1.9ex}
&&&\tag2.2\\
\vspace{-1.9ex}
&\hskip -2em
\bold e_i=\sum^n_{j=1} T^j_i\cdot\tilde\bold e_j,
&&h^s=\sum^n_{r=1} S^s_r\cdot\tilde h^r.\\
\endxalignat
$$\par
     Let $f$ be a covector from the conjugate space $V^*$. Let's consider
its expansions in two conjugate bases $h^1,\,\ldots,\,h^n$ and $\tilde
h^1,\,\ldots,\,\tilde h^n$:
$$
\xalignat 2
&\hskip -2em
f=\sum^n_{s=1}f_s\cdot h^s,
&&f=\sum^n_{r=1}\tilde f_r\cdot\tilde h^r.
\tag2.3
\endxalignat
$$
The expansions \thetag{2.3} also differ from the standard introduced
by formula \thetag{5.1} in
Chapter~\uppercase\expandafter{\romannumeral 1}. To the coordinates 
of covectors the other standard is applied: they are specified by
lower indices and are written in row vectors.\par
\par
\proclaim{Theorem 2.2} The coordinates of a covector $f$ in two 
dual bases $h^1,\,\ldots,\,h^n$ and $\tilde h^1,\,\ldots,\,
\tilde h^n$ associated with the bases $\bold e_1,\,\ldots,\,\bold 
e_n$ and $\tilde\bold e_1,\,\ldots,\,\tilde\bold e_n$ in $V$ are
related to each other by formulas 
$$
\xalignat 2
&\hskip -2em
\tilde f_r=\sum^n_{s=1}S^s_r\,f_s,
&&f_s=\sum^n_{j=1}T^r_s\,\tilde f_r,
\tag2.4
\endxalignat
$$
where $S$ is the direct transition matrix for passing from $\bold e_1,\,\ldots,\,\bold e_n$ to the
{\tencyr\char '074}wavy{\tencyr\char '076} basis $\tilde\bold e_1,
\,\ldots,\,\tilde \bold e_n$, while $T=S^{-1}$ is the inverse transition
matrix.
\endproclaim
\demo{Proof} In order to prove the first relationship \thetag{2.4}
we substitute the fourth expression \thetag{2.2} for $h^s$ into the 
first expansion \thetag{2.3}:
$$
f=\sum^n_{s=1}f_s\cdot\left(\,\shave{\sum^n_{r=1}}
S^s_r\cdot\tilde h^r\right)=\sum^n_{r=1}\left(\,\shave{
\sum^n_{s=1}}S^s_r\,f_s\right)\cdot\tilde h^r.
$$
Then we compare the resulting expansion of $f$ with the second
expansion \thetag{2.3} and derive the first formula \thetag{2.4}. 
The second formula \thetag{2.4} is derived similarly.
\qed\enddemo
     Note that the formulas \thetag{2.4} can be derived immediately
from the definition~1.5 and from formula \thetag{1.6} without using
the conjugate bases.
\proclaim{Theorem 2.3} The scalar product of a vector $\bold v$ and
a covector $f$ is determined by their coordinates according to the
formula 
$$
\hskip -2em
\langle f\,|\,\bold v\rangle=\sum^n_{i=1}f_i\,v^i=
f_1\,v^1+\ldots+f_n\,v^n.
\tag2.5
$$
\endproclaim
\demo{Proof} In order to prove \thetag{2.5} we use the
relationship \thetag{1.6}:
$$
\langle f\,|\,\bold v\rangle=f(\bold v)=\sum^n_{i=1}f(\bold e_i)
\,v^i=\sum^n_{i=1}f_i\,v^i.
$$
In \thetag{2.5} and in the above calculations $f$ is assumed to 
be expanded in the basis $h^1,\,\ldots,\,h^n$ conjugated to the 
basis $\bold e_1,\,\ldots,\,\bold e_n$, where $\bold v$ is expanded. 
\qed\enddemo
\head
\S\,3. Orthogonal complements in a dual space.
\endhead
\rightheadtext{\S\,3. Orthogonal complements in a dual space.}
\definition{Definition 3.1} Let $S$ be a subset in a linear vector 
space $V$. The {\it orthogonal complement} of the subset $S$ in
the conjugate space $V^*$ is the set $S^{\sssize\perp}\subset V^*$ 
composed by covectors each of which orthogonal to all vectors of $S$.
\enddefinition
    The above definition of the orthogonal complement $S^{\sssize\perp}$
can be expressed by the formula $S^{\sssize\perp}=\{f\in V^*\!:\,\forall
\bold v\,((v\in S)\Rightarrow (\langle f\,|\,\bold v\rangle=0))\}$.
\proclaim{Theorem 3.1} The operation of constructing orthogonal 
complements of sub\-sets $S\subset V$ in the conjugate space $V^*$ 
possesses the following properties:
\roster
\item $S^{\sssize\perp}$ is a subspace in $V^*$;
\item $S_1\subset S_2$ implies $(S_2)^{\sssize\perp}
      \subset (S_1)^{\sssize\perp}$;
\item $\langle S\rangle^{\sssize\perp}=S^{\sssize\perp}$, where
      $\langle S\rangle$ is the linear span of $S$;
\vskip 1ex
\item $\dsize\left(\,\,\shave{\bigcup_{i\in I}} S_i
       \right)^{\kern -2pt\lower 2pt\hbox{$\ssize\perp$}}=
       \bigcap_{i\in I}\,(S_i)^{\sssize\perp}$.
\endroster
\endproclaim
\demo{Proof} Let's prove the first item of the theorem for the beginning.
For this purpose we should verify two conditions from the definition of
a subspace. Let $f_1,f_2\in S^{\sssize\perp}$, then $\langle f_1\,|\,
\bold v\rangle=0$ and $\langle f_2\,|\,\bold v\rangle=0$ for all $\bold v
\in S$. Therefore, for all vectors $\bold v\in S$ we derive the equality
$\langle f_1+f_2\,|\,\bold v\rangle=\langle f_1\,|\,\bold v\rangle+
\langle f_2\,|\,\bold v\rangle=0$ which means that $f_1+f_2\in S^{\sssize
\perp}$.\par
    Now assume that $f\in S^{\sssize\perp}$. Then $\langle f\,|
\,\bold v\rangle=0$ for all vectors $\bold v\in S$. Hence, for the 
covector $\alpha\cdot f$ we defive $\langle \alpha\cdot f\,|\,\bold v\rangle=\alpha\,\langle f\,|\,\bold v\rangle=0$. This means that $\alpha\cdot f\in S^{\sssize\perp}$. Thus, the first item
in the theorem~3.1 is proved.\par
     In order to prove the inclusion $(S_2)^{\sssize\perp}\subset
(S_1)^{\sssize\perp}$ in the second item of the theorem~3.1 we consider
an arbitrary covector $f$ of $(S_2)^{\sssize\perp}$. From the condition
$f\in (S_2)^{\sssize\perp}$ we derive $\langle f\,|\,\bold v\rangle=0$
for any $\bold v\in S_2$. But $S_1\subset S_2$, therefore, the 
equality $\langle f\,|\,\bold v\rangle=0$ holds for any $\bold v\in S_1$.
Then $f\in (S_1)^{\sssize\perp}$. This means that $f\in (S_2)^{\sssize
\perp}$ implies $f\in (S_1)^{\sssize\perp}$. The required inclusion
is proved.\par
     In order to prove the third item of the theorem note that the
linear span of $S$ comprises this set: $S\subset\langle S\rangle$.
Applying the item \therosteritem{2} of the theorem, which we have 
already proved, we obtain the inclusion $\langle S\rangle^{\sssize
\perp}\subset S^{\sssize\perp}$. Now we need the opposite inclusion
$S^{\sssize\perp}\subset\langle S\rangle^{\sssize\perp}$. In order to
prove it let's remember that the linear span $\langle S\rangle$ 
consists of all possible linear combinations of the form
$$
\hskip -2em
\bold v=\alpha_1\cdot\bold v_1+\ldots+\alpha_r\cdot\bold v_r
\text{, \ where \ }\bold v_i\in S.
\tag3.1
$$
Let $f\in S^{\sssize\perp}$, then $\langle f\,|\,\bold v\rangle=0$ 
for all $\bold v\in S$. In particular, this applies to the vectors
$\bold v_i$ in the expansion \thetag{3.1}, i\.\,e\. $\langle f\,|
\,\bold v_i\rangle=0$. Then from \thetag{3.1} we derive
$$
\langle f\,|\,\bold v\rangle=\alpha_1\,\langle f\,|\,\bold v_1
\rangle+\ldots+\alpha_r\,\langle f\,|\,\bold v_r\rangle=0.
$$
This means that $\langle f\,|\,\bold v\rangle=0$ for all $\bold v
\in\langle S\rangle$. This proves the opposite inclusion 
$S^{\sssize\perp}\subset\langle S\rangle^{\sssize\perp}$ and, thus, completes the proof of the equality $\langle S\rangle^{\sssize
\perp}=S^{\sssize\perp}$.\par
     Now let's proceed to the proof of the fourth item of the
theorem~3.1. For this purpose we introduce the following notations:
$$
\xalignat 2
&S=\bigcup_{i\in I}\,S_i,
&&\tilde S=\bigcap_{i\in I}\,(S_i)^{\sssize\perp}.
\endxalignat
$$
Let $f\in S^{\sssize\perp}$. Then $\langle f\,|\,\bold v\rangle
=0$ for all $\bold v\in S$. But $S_i\subset S$ for any $i\in I$.
Therefore, $\langle f\,|\,\bold v\rangle=0$ for all $\bold v\in 
S_i$ and for all $i\in I$. This means that $f$ belongs to each of
the orthogonal complement $(S_i)^{\sssize\perp}$, therefore, it 
belongs to their intersection. Thus, we have proved the inclusion
$S^{\sssize\perp}\subset\tilde S$.\par
     Conversely, from the inclusion $f\in (S_i)^{\sssize\perp}$ 
for all $i\in I$ we derive $\langle f\,|\,\bold v\rangle=0$
for all $\bold v\in S_i$ and for all $i\in I$. This means that the
equality $\langle f\,|\,\bold v\rangle=0$ holds for all vectors
$\bold v$ in the union of all sets $S_i$. This proves the converse
inclusion $\tilde S\subset S^{\sssize\perp}$. Thus, we have proved
that $S^{\sssize\perp}=\tilde S$. The theorem is proved.
\qed\enddemo
\definition{Definition 3.2} Let $S$ be a subset of conjugate space
$V^*$. The {\it orthogonal complement} of $S$ in $V$ is the set
$S^{\sssize\perp}\in V$ formed by vectors each of which is orthogonal 
to all covectors of the set $S$.
\enddefinition
     The above definition of the orthogonal complement $S^{\sssize
\perp}\subset V$ can be expressed by the formula $S^{\sssize\perp}
=\{\bold v\in V\!:\,\forall f\,((f\in S)\Rightarrow (\langle f\,|\,
\bold v\rangle=0))\}$. For this orthogonal complement one can
formulate a theorem quite similar to the theorem~3.1.
\proclaim{Theorem 3.2} The operation of constructing orthogonal 
complements of sub\-sets $S\subset V^*$ in $V$ possesses the 
following four properties:
\roster
\item $S^{\sssize\perp}$ is a subspace in $V$;
\item $S_1\subset S_2$ implies $(S_2)^{\sssize\perp}
      \subset (S_1)^{\sssize\perp}$;
\item $\langle S\rangle^{\sssize\perp}=S^{\sssize\perp}$, where
      $\langle S\rangle$ is a linear span of $S$;
\vskip 1ex
\item $\dsize\left(\,\,\shave{\bigcup_{i\in I}} S_i
       \right)^{\kern -2pt\lower 2pt\hbox{$\ssize\perp$}}=
       \bigcap_{i\in I}\,(S_i)^{\sssize\perp}$.
\endroster
\endproclaim
     The proof of this theorem almost literally coincides with the proof 
of the theorem~3.1. Therefore, here we omit this proof.\par
\proclaim{Theorem 3.3} Let $V$ be a finite-dimensional vector space and
suppose that we have a subspaces $U\subset V$ and a subspace $W\subset 
V^*$. The the condition $W=U^{\sssize\perp}$ in the sense of 
definition~3.1 is equivalent to the condition $U=W^{\sssize\perp}$ in the
sense of definition~3.2.
\endproclaim
\demo{Proof} Suppose that $W=U^{\sssize\perp}$ in the sense of 
definition~3.1. Then for any $w\in W$ and for any $\bold u\in U$ 
we have the orthogonality $\langle w\,|\,\bold u\rangle=0$. By definition
$W^{\sssize\perp}$ is the set of all vectors $\bold v\in V$ such that
$\langle w\,|\,\bold v\rangle=0$ for all covectors $w\in W$. Hence,
$\bold u\in U$ implies $\bold u\in W^{\sssize\perp}$ and we have the
inclusion $U\subset W^{\sssize\perp}$.\par
     However, we need to prove the coincidence $U=W^{\sssize\perp}$.
Let's do it by con\-tradiction. Suppose that $U\neq W^{\sssize\perp}$.
Then there is a vector $\bold v_0$ such that $\bold v_0\in
W^{\sssize\perp}$ and $\bold v_0\notin U$. In this case we can apply
the theorem~1.2 which says that there is a linear functional $f$ such 
that it vanishes on all vectors $\bold u\in U$ and is nonzero
on the vector $\bold v_0$. Then $f\in W$ and $\langle f\,|\,\bold v_0
\rangle\neq 0$, so we have the contradiction to the condition $\bold 
v_0\in W^{\sssize\perp}$. This contradiction proves that
$U=W^{\sssize\perp}$. As a result we have proved that 
$W=U^{\sssize\perp}$ implies $U=W^{\sssize\perp}$.\par
     Now, conversely, let $U=W^{\sssize\perp}$. Then for any $w\in W$ 
and for any $\bold u\in U$ we have the orthogonality $\langle w\,|\,\bold
u\rangle=0$. By definition $U^{\sssize\perp}$ is the set of all covectors
$f$ perpendicular to all vectors $\bold u\in U$. Hence, $w\in W$ implies
$w\in U^{\sssize\perp}$ and we have the inclusion $W\subset U^{\sssize
\perp}$.\par
     Next step is to prove the coincidence $W=U^{\sssize\perp}$. We shall
do it again by contradiction. Assume that $W\neq U^{\sssize\perp}$. Then
there is a covector $f_0\in U^{\sssize\perp}$ such that $f_0\notin W$.
Let's apply the theorem~1.2. It this case it says that there is a linear functional $\varphi$ in $V^{**}$ such that it vanishes on $W$ and is 
nonzero on the covector $f_0$. Remember that we have the canonic 
isomorphism $h\!:\,V\to V^{**}$. We apply $h^{-1}$ to $\varphi$ and get 
the vector $\bold v=h^{-1}(\varphi)$. Then we take into account \thetag{1.13}
which yields $\bold v\in U$ and $\langle f_0\,|\,\bold v\rangle\neq 0$.
This contradicts to the condition $f_0\in U^{\sssize\perp}$. Hence, by
contradiction, $U=W^{\sssize\perp}$ and $U=W^{\sssize\perp}$ implies
$W=U^{\sssize\perp}$. The theorem is\linebreak completely proved.
\qed\enddemo
     The proposition of the theorem~3.3 can be reformulated as follows:
in the case of a finite-dimensional space $V$ for any subspace $U\in V$
and for any subspace $W\in V^*$ the following relationships are valid:
$$
\xalignat 2
&\hskip -2em
(U^{\sssize\perp})^{\sssize\perp}=U,
&&(W^{\sssize\perp})^{\sssize\perp}=W.
\tag3.2
\endxalignat
$$
For arbitrary subsets $S\in V$ and $R\in V^*$ (not subspaces) in the
case of a finite-dimensional space $V$ we have the relationships
$$
\xalignat 2
&\hskip -2em
(S^{\sssize\perp})^{\sssize\perp}=\langle S\rangle,
&&(R^{\sssize\perp})^{\sssize\perp}=\langle R\rangle.
\tag3.3
\endxalignat
$$
These relationships \thetag{3.3} are derived from \thetag{3.2} by using
the item \therosteritem{3} in theorems~3.1 and 3.2.
\proclaim{Theorem 3.4} In the case of a finite-dimensional linear vector
space $V$ if $U$ is a subspace of $V$ or if $U$ is a subspace of $V^*$,
then \ $\dim U+\dim U^{\sssize\perp}=\dim V$.
\endproclaim
\demo{Proof} Due to the relationships \thetag{3.2} the second case
$U\subset V^*$ in the theorem~3.4 is reduced to the first case
$U\subset V$ if we replace $U$ by $U^{\sssize\perp}$. Therefore, 
we consider only the first case $U\subset V$.\par
     Let $\dim V=n$ and $\dim U=s$. We choose a basis $\bold e_1,\,
\ldots,\,\bold e_s$ in the subspace $U$ and complete it up to a basis
$\bold e_1,\,\ldots,\,\bold e_n$ in the subspace $V$. The basis $\bold
e_1,\,\ldots,\,\bold e_n$ determines the conjugate basis
$h^1,\,\ldots,\,h^n$ in $V^*$. If we specify vectors by their coordinates 
in the basis $\bold e_1,\,\ldots,\,\bold e_n$ and if we specify covectors
by their coordinates in dual basis $h^1,\,\ldots,\,h^n$, then we can apply
the formula \thetag{2.5}.\par
     By construction of the basis $\bold e_1,\,\ldots,\,\bold e_n$ the
subspace $U$ consists of vectors the initial $s$ coordinates of which
are deliberate, while the remaining $n-s$ coordinates are equal to zero.
Therefore the condition $f\in U^{\sssize\perp}$ means that the equality
$$
\langle f\,|\,v\rangle=\sum^s_{i=1} f_i\,v^i=0
$$
should be fulfilled identically for any numbers $v^1,\,\ldots,\,v^s$.
This is the case if and only if the first $s$ coordinates of the 
covector $f$ are zero. Other $n-s$ coordinates of $f$ are deliberate. 
This means that the subspace $U^{\sssize\perp}$ is the linear span of 
the last $n-s$ basis vectors of the conjugate basis:
$$
U^{\sssize\perp}=\langle h^{s+1},\ldots,h^n\rangle.
$$
For the dimension of the subspace $U^{\sssize\perp}$ this yields
$\dim U^{\sssize\perp}=n-s$, hence, we have the required identity
$\dim U+\dim U^{\sssize\perp}=\dim V$. The theorem is proved
\qed\enddemo
    The theorem~3.4 is known as the theorem {\it on the dimension of
orthogonal complements}. As an immediate consequence of this theorem
we get
$$
\xalignat 2
&\hskip -2em
\{0\}^{\sssize\perp}=V,   &&V^{\sssize\perp}=\{0\},\\
\vspace{-1.9ex}
&\hskip -2em&&\tag3.4\\
\vspace{-1.9ex}
&\hskip -2em
\{0\}^{\sssize\perp}=V^*, &&(V^*)^{\sssize\perp}=\{0\}.
\endxalignat
$$
All these equalities have the transparent interpretation. The first
three of the equalities \thetag{3.4} can be proved immediately without
using the finite-dimensionality of $V$. The proof of the last equality
\thetag{3.4} uses the corollary of the theorem~1.2, while this theorem
assumes $V$ to be a finite-dimensional space.
\proclaim{Theorem 3.5} In the case of a finite-dimensional space $V$
for any family of subspaces in $V$ or in $V^*$ the following relationships
are fulfilled
$$
\xalignat 2
&\hskip -2em
\left(\,\,\shave{\sum_{i\in I}} U_i\right)^{\kern -2pt
 \lower 2pt\hbox{$\ssize\perp$}}=
 \bigcap_{i\in I}\,(U_i)^{\sssize\perp},
&&\left(\,\,\shave{\bigcap_{i\in I}} U_i\right)^{\kern -2pt
 \lower 2pt\hbox{$\ssize\perp$}}=
 \sum_{i\in I}\,(U_i)^{\sssize\perp}.
\tag3.5
\endxalignat
$$
\endproclaim
\demo{Proof} The sum of subspaces is the span of their union.
Therefore, the first relationship \thetag{3.5} is an immediate
consequence of the items \therosteritem{3} and \therosteritem{4} 
in the theorems~3.1 and 3.2. The finite-dimensionality of $V$ here
is not used.\par
     The second relationship \thetag{3.5} follows from the first one
upon substituting $U_i$ by $(U_i)^{\sssize\perp}$. Indeed, applying
\thetag{3.2}, we derive the equality
$$
\qquad\left(\,\,\shave{\sum_{i\in I}} (U_i)^{\sssize\perp}
\right)^{\kern -2pt\lower 2pt\hbox{$\ssize\perp$}}=
\bigcap_{i\in I}\,((U_i)^{\sssize\perp})^{\sssize\perp}=
\bigcap_{i\in I}\,U_i.
$$
Now it is sufficient to pass to orthogonal complements in both
sides of this equality and apply \thetag{3.2} again. The theorem
is proved.
\qed\enddemo
\head
\S\,4. Conjugate mapping.
\endhead
\definition{Definition 4.1} Let $f\!:\,V\to W$ be a linear mapping
from $V$ to $W$. A linear mapping $\varphi\!:\,W^*\to V^*$ is called a
{\it conjugate mapping} for $f$ if for any $\bold v\in V$ and for 
any $w\in W^*$ the relationship $\langle\varphi(w)\,|\,\bold v\rangle
=\langle w\,|\,f(\bold v)\rangle$ is fulfilled.
\enddefinition
     The problem of the existence of a conjugate mapping is solved by
the definition~4.1 itself. Indeed, in order to define a mapping
$\varphi: W^*\to V^*$ for each functional $w\in W^*$ we should specify
the corresponding functional $h=\varphi(w)\in V^*$. But to specify a
functional in $V^*$ this means that we should specify its action upon
an arbitrary vector $\bold v\in V$. In the sense of this reasoning the
defining relationship for a conjugate mapping is written as follows:
$$
h(\bold v)=\langle h\,|\,\bold v\rangle=\langle \varphi(w)\,|\,\bold
v\rangle=\langle w\,|\,f(\bold v)\rangle.
$$
It is easy to verify that the above equality defines a linear functional
$h=h(\bold v)$:
$$
\align
&\aligned
h(\bold v_1+\bold v_2)=\langle w\,|\,&f(\bold v_1+\bold v_2)
\rangle=\langle w\,|\,f(\bold v_1)+f(\bold v_2)\rangle=\\
&=\langle w\,|\,f(\bold v_1)\rangle+\langle w\,|\,f(\bold v_2)
\rangle=h(\bold v_1)+h(\bold v_2),
\endaligned\\
\vspace{1.7ex}
&\aligned
h(\alpha\cdot\bold v)=\langle w\,|\,f(\alpha\cdot\bold v)\rangle
=\langle w\,|\,\alpha\cdot f(\bold v)\rangle=\alpha\,\langle
w\,|\,f(\bold v)\rangle=\alpha\,h(\bold v).
\endaligned
\endalign
$$
\proclaim{Theorem 4.1} For a linear mapping $f\!:\,V\to W$ from $V$ 
to $W$ the conjugate mapping $\varphi\!:\,W^*\to V^*$ is also linear.
\endproclaim
\demo{Proof} Due to the relationship~4.1 for the conjugate mapping 
$\varphi\!:\,W^*\to V^*$ we have the following relationships:
$$
\align
&\aligned
  \varphi(w_1&+w_2)(\bold v)=\langle w_1+w_2\,|\,f(\bold v)\rangle=
  \langle w_1\,|\,f(\bold v)\rangle+\\
  &+\langle w_2\,|\,f(\bold v)\rangle=\varphi(w_1)(\bold v)+
  \varphi(w_2)(\bold v)=
  (\varphi(w_1)+\varphi(w_2))(\bold v),
  \endaligned\\
\vspace{1.7ex}
&\aligned
  \varphi(\alpha\cdot w)(\bold v)=\langle\alpha\cdot w\,|\,f(\bold
  v)\rangle&=\alpha\,\langle w\,|\,f(\bold v)\rangle=\\
  &=\alpha\,\varphi(w)(\bold v)=(\alpha\cdot\varphi(w))(\bold v).
  \endaligned
\endalign
$$
Since $\bold v\in V$ is an arbitrary vector of $V$ from the above
calculations we obtain $\varphi(w_1+w_2)=\varphi(w_1)+\varphi(w_2)$
and $\varphi(\alpha\cdot w)=\alpha\cdot\varphi(w)$. This means that
the conjugate mapping $\varphi$ is a linear mapping.
\qed\enddemo
     As we have seen above, the conjugate mapping $\varphi\!:\,W^*\to 
V^*$ for a mapping $f\!:\,V\to W$ is unique. It is usually denoted
$\varphi=f^*$. The operation of passing from $f$ to its conjugate
mapping $f^*$ possesses the following properties:
$$
\xalignat 3
&(f+g)^*=f^*+g^*,
&&(\alpha\cdot f)^*=\alpha\cdot f^*,
&&(f{\ssize\circ}\,g)^*=g^*{\ssize\circ}\,f^*.
\endxalignat
$$
The first two properties are naturally called the {\it linearity\/}.
The last third property makes the operation $f\to f^*$ an analog of 
the matrix transposition. All three of the above properties are proved
by direct calculations on the base of the definition~4.1. We shall not
give these calculations here since in what follows we shall not use the
above properties at all.\par
\proclaim{Theorem 4.2} In the case of finite-dimensional spaces $V$ and
$W$ the kernels and images of the mappings $f:V\to W$ and $f^*:W^*\to V^*$
are related as follows:
$$
\xalignat 2
&\hskip -2em
\Ker f^*=(\Img f)^{\sssize\perp},
&&\Ker f=(\Img f^*)^{\sssize\perp},\\
\vspace{-1.9ex}
&&&\tag4.1\\
\vspace{-1.9ex}
&\hskip -2em
\Img f=(\Ker f^*)^{\sssize\perp},
&&\Img f^*=(\Ker f)^{\sssize\perp}.
\endxalignat
$$
\endproclaim
\demo{Proof} The kernel $\Ker f^*$ is the set of linear functionals
of $W^*$ that are mapped to the zero functional in $V^*$ under the
action of the mapping $f^*$. Therefore, $w\in\Ker f^*$ is equivalent 
to the equality $f^*(w)(\bold v)=0$ for all  $\bold v\in V$. As a result
of simple calculations we obtain 
$$
f^*(w)(\bold v)=\langle f^*(w)\,|\,\bold v\rangle=
\langle w\,|\,f(\bold v)\rangle=0.
$$
Hence, the kernel $\Ker f^*$ is the set of covectors orthogonal to the
vectors of the form $f(\bold v)$. But the vectors of the form $f(\bold v)
\in W$ constitute the image $\Img f$. Therefore, $\Ker f^*=(\Img
f)^{\sssize\perp}$. The first relationship \thetag{4.1} is proved. In
proving this relationship we did not use the finite-dimensionality of
$W$. It is valid for infinite dimensional spaces as well.\par
     In order to prove the second relationship we consider the orthogonal
complement $(\Img f^*)^{\sssize\perp}$. It is formed by the vectors 
orthogonal to all covectors of the form $f^*(w)$:
$$
0=\langle f^*(w)\,|\,\bold v\rangle=\langle w\,|\,f(\bold v)\rangle.
$$
Using the finite-dimensionality of $W$, we apply the corollary of the 
theorem~1.2. It says that if $\langle w\,|\,f(\bold v)\rangle=0$ for all
$w\in W^*$, then $f(\bold v)=\bold 0$. Therefore, we have
$(\Img f^*)^{\sssize\perp}=\Ker f$. The second relationship \thetag{4.1}
is proved. The third and the fourth relationships are derived from the
first and the second ones by means of the theorem~3.3. Thereby we use
the finite-dimensionality of the spaces $W$ and $V$.
\qed\enddemo
     Let the spaces $V$ and $W$ be finite-dimensional. Let's choose a
basis $\bold e_1,\,\ldots,\,\bold e_n$ in $V$ and a basis $\tilde\bold e_1,\,\ldots,\,\tilde\bold e_m$ in the space $W$. This choice uniquely
determines the conjugate bases $h^1,\,\ldots,\,h^n$ and $\tilde h^1,\,
\ldots,\,\tilde h^m$ in $V^*$ and $W^*$. Let's consider a mapping 
$f\!:\,V\to W$ and the conjugate mapping $f^*\!:\,W^*\to V^*$. The 
matrices of the mappings $f$ and $f^*$ are determined by the expansions:
$$
\xalignat 2
&\hskip -2em
\qquad f(\bold e_j)=\sum^m_{k=1}F^k_j\,\tilde\bold e_k,
&&f^*(\tilde h^i)=\sum^n_{q=1}\Phi^i_q\,h^q.
\tag4.2
\endxalignat
$$
The second relationship \thetag{4.1} is somewhat different by structure
from the first one. The matter is that the basis vectors of the dual 
basis are indexed differently (with upper indices). However, this 
relationship implement the same idea as the first one: the mapping is
applied to a basis vector of one space and the result is expanded in the
basis of another space.\par
\proclaim{Theorem 4.3} The matrices of the mappings $f$ and $f^*$
determined by the relationships \thetag{4.2} are the same, i\.\,e\.
$F^i_j=\Phi^i_j$.
\endproclaim
\demo{Proof} From the definition of the conjugate mapping we derive
$$
\hskip -2em
\langle\tilde h^i\,|\,f(\bold e_j)\rangle=
\langle f^*(\tilde h^i)\,|\,\bold e_j\rangle.
\tag4.3
$$
Let's calculate separately the left and the right hand sides of this
equality using the expansion \thetag{4.2} for this purpose:
$$
\align
&\langle\tilde h^i\,|\,f(\bold e_j)\rangle=
 \sum^m_{k=1} F^k_j\,\langle\tilde h^i\,|\,\tilde\bold e_k\rangle=
 \sum^m_{k=1} F^k_j\,\delta^i_k=F^i_j,\\
&\langle f^*(\tilde h^i)\,|\,\bold e_j\rangle=
 \sum^n_{q=1}\Phi^i_q\,\langle h^q\,|\,\bold e_j\rangle=
 \sum^n_{q=1}\Phi^i_q\,\delta^q_j=\Phi^i_j.
\endalign
$$
Substituting the above expressions back to the formula \thetag{4.3}, 
we get the required coincidence of the matrices: $F^i_j=\Phi^i_j$.
\qed\enddemo
\subhead
Remark
\endsubhead
In some theorems of this chapter the restrictions to the 
finite-dimensional case can be removed. However, the prove of 
such strengthened versions of these theorems is based on the 
axiom of choice (see \cite{1}).\par
\newpage
\topmatter
\title\chapter{4}
Bilinear and quadratic forms.
\endtitle
\endtopmatter
\document
\head
\S\,1. Symmetric bilinear forms\\
and quadratic forms. Recovery formula.
\endhead
\leftheadtext{CHAPTER~\uppercase\expandafter{\romannumeral 4}.
BILINEAR AND QUADRATIC FORMS.}
\rightheadtext{\S\,1. Symmetric bilinear forms \dots}
\setfirstpage
\definition{Definition 1.1} Let $V$ be a linear vector space
over a numeric field $\Bbb K$. A numeric function $y=f(\bold v,
\bold w)$ with two arguments $\bold v,\bold w\in V$ and with
the values in the field $\Bbb K$ is called a {\it bilinear form\/} 
if 
\roster
\item $f(\bold v_1+\bold v_2,\bold w)=f(\bold v_1,\bold w)
      +f(\bold v_2,\bold w)$ for any two $\bold v_1,\bold v_2\in V$;
\item $f(\alpha\cdot\bold v,\bold w)=\alpha\,f(\bold v,\bold w)$ for
      any $\bold v\in V$ and for any $\alpha\in\Bbb K$;
\item $f(\bold v,\bold w_1+\bold w_2)=f(\bold v,\bold w_1)
      +f(\bold v,\bold w_2)$ for any two $\bold v_1,\bold v_2\in V$;
\item $f(\bold v,\alpha\cdot\bold w)=\alpha\,f(\bold v,\bold w)$ for
      any $\bold v\in V$ and for any $\alpha\in\Bbb K$.
\endroster
\enddefinition
     The bilinear form $f(\bold v,\bold w)$ is linear in its first 
argument $\bold v$ when the second argument $\bold w$ is fixed; it 
is also linear in its second argument $\bold w$ when the first 
argument $\bold v$ is fixed.
\definition{Definition 1.2} A bilinear form $f(\bold v,\bold w)$ is called
a {\it symmetric bilinear form\/} if $f(\bold v,\bold w)=f(\bold w,\bold v)$.
\enddefinition
\definition{Definition 1.3} A bilinear form $f(\bold v,\bold w)$ is called
a {\it skew-symmetric bilinear form\/} or an {\it antisymmetric bilinear
form\/} if $f(\bold v,\bold w)=-f(\bold w,\bold v)$.
\enddefinition
     Having a bilinear form $f(\bold v,\bold w)$, one can produce a
symmetric bilinear form:
$$
\hskip -2em
f_{\sssize+}(\bold v,\bold w)=\frac{f(\bold v,\bold w)
+f(\bold w,\bold v)}{2}.
\tag1.1
$$
Similarly, one can produce a skew-symmetric bilinear form:
$$
\hskip -2em
f_{\sssize-}(\bold v,\bold w)=\frac{f(\bold v,\bold w)
-f(\bold w,\bold v)}{2}.
\tag1.2
$$
The operation \thetag{1.1} is called the {\it symmetrization\/} of 
the bilinear form $f$; the operation \thetag{1.2} is called the {\it
alternation\/} of this bilinear form. Thereby any bilinear form is
the sum of a symmetric bilinear form and a skew-symmetric one:
$$
\hskip -2em
f(\bold v,\bold w)=f_{\sssize+}(\bold v,\bold w)
+f_{\sssize-}(\bold v,\bold w).
\tag1.3
$$
\proclaim{Theorem 1.1} The expansion of a given bilinear form $f(\bold v,
\bold w)$ into the sum of a symmetric and a skew-symmetric bilinear forms
is unique.
\endproclaim
\demo{Proof} Let's consider an expansion of $f(\bold v,\bold w)$ into
the sum of a symmetric and a skew-symmetric bilinear forms
$$
\hskip -2em
f(\bold v,\bold w)=h_{\sssize+}(\bold v,\bold w)
+h_{\sssize-}(\bold v,\bold w).
\tag1.4
$$
By means of symmetrization and alternation from \thetag{1.4} we derive
$$
\align
&\aligned
  f(v,w)+f(w,v)&=(h_{\sssize+}(v,w)+h_{\sssize+}(w,v))+\\
  &+(h_{\sssize-}(v,w)+h_{\sssize-}(w,v))=2\,h_{\sssize+}(v,w),
  \endaligned\\
\vspace{1.7ex}
&\aligned
  f(v,w)-f(w,v)&=(h_{\sssize+}(v,w)-h_{\sssize+}(w,v))+\\
  &+(h_{\sssize-}(v,w)-h_{\sssize-}(w,v))=2\,h_{\sssize-}(v,w),
  \endaligned
\endalign
$$
Hence, $h_{\sssize+}=f_{\sssize+}$ and $h_{\sssize-}=f_{\sssize-}$.
Therefore, the expansion \thetag{1.4} coincides with the expansion
\thetag{1.3}. The theorem is proved.
\qed\enddemo
\definition{Definition 1.4} A numeric function $y=g(\bold v)$ with
one vectorial argument $\bold v\in V$ is called a {\it quadratic 
form} in a linear vector space $V$ if $g(\bold v)=f(\bold v,\bold v)$
for some bilinear form $f(\bold v,\bold w)$.
\enddefinition
    If $g(\bold v)=f(\bold v,\bold v)$, then the quadratic form $g$ 
is said to be {\it generated} by the bilinear form $f$. For a 
skew-symmetric bilinear form we have $f_{\sssize-}(\bold v,\bold v)
=-f_{\sssize-}(\bold v,\bold v)$. Hence, $f_{\sssize-}(\bold v,\bold v)
=0$. Then from the expansion \thetag{1.3} we derive
$$
\hskip -2em
g(\bold v)=f(\bold v,\bold v)=f_{\sssize+}(\bold v,\bold v).
\tag1.5
$$
The same quadratic form can be generated by several bilinear forms.
The relationship \thetag{1.5} shows that any quadratic form can be
generated by a symmetric bilinear form.
\proclaim{Theorem 1.2} For any quadratic form $g(\bold v)$ there is 
the unique bilinear form $f(\bold v,\bold w)$ that generates $g(\bold v)$.
\endproclaim
\demo{Proof} The existence of a symmetric bilinear form $f(\bold v,\bold
w)$ generating $g(\bold v)$ follows from \thetag{1.5}. Let's prove the
uniqueness of this form. From $g(\bold v)=f(\bold v,\bold v)$ and from
the symmetry of the form $f$ we derive
$$
\aligned
g(\bold v+\bold w)&=f(\bold v+\bold w,\bold v+\bold w)=f(\bold v,\bold v)
+f(\bold v,\bold w)+\\
&+f(\bold w,\bold v)+f(\bold w,\bold w)=f(\bold v,\bold v)+2\,f(\bold v,
\bold w)+f(\bold w,\bold w).
\endaligned
$$
Now $f(\bold v,\bold v)$ and $f(\bold w,\bold w)$ in right hand side of 
this formula can be replaced by $g(\bold v)$ and $g(\bold w)$ respectively. 
Hence, we get
$$
\hskip -2em
f(\bold v,\bold w)=\frac{g(\bold v+\bold w)-g(\bold v)-g(\bold w)}{2}.
\tag1.6
$$
Formula \thetag{1.6} shows that the values of the symmetric bilinear form
$f(\bold v,\bold w)$ are uniquely determined by the values of the quadratic
form $g(\bold v)$. This proves the uniqueness of the form $f$.
\qed\enddemo
     The formula \thetag{1.6} is called a {\it recovery formula}. Usually,
a quadratic form and an associated symmetric bilinear form for it both are denoted by the same symbol: $g(\bold v)=g(\bold v,\bold v)$. Moreover, when a quadratic form is given, we assume without stipulations that the associated symmetric bilinear form $g(\bold v,\bold w)$ is also \pagebreak
given.\par
     Let $f(\bold v,\bold w)$ be a bilinear form in a finite-dimensional
linear vector space $V$ and let $\bold e_1,\,\ldots,\,\bold e_n$ be a basis
in this space. The numbers $f_{ij}$ determined by formula
$$
\hskip -2em
f_{ij}=f(\bold e_i,\bold e_j)
\tag1.7
$$
are called the {\it coordinates\/} or the {\it components\/} of the form $f$
in the basis $\bold e_1,\,\ldots,\,\bold e_n$. The numbers \thetag{1.7} 
are written in form of a matrix 
$$
F=\Vmatrix f_{11} & \hdots & f_{1n}\\
\vspace{1.5ex}
\vdots & \ddots & \vdots\\
\vspace{1.5ex}
f_{n1} & \hdots & f_{nn}
\endVmatrix,
\tag1.8
$$
which is called the matrix of the bilinear form $f$ in the basis $\bold e_1,
\,\ldots,\,\bold e_n$. For the element $f_{ij}$ in the matrix \thetag{1.8} the first index $i$ specifies the row number, the second index $j$ specifies
the column number. The matrix of a symmetric bilinear form $g$ is also
symmetric: $g_{ij}=g_{ji}$. Further, saying the matrix of a quadratic
form $g(\bold v)$, we shall assume the matrix of an associated symmetric bilinear form $g(\bold v,\bold w)$.\par
     Let $v^1,\,\ldots,\,v^n$ and $w^1,\,\ldots,\,w^n$ be coordinates of
two vectors $\bold v$ and $\bold w$ in the basis $\bold e_1,\,\ldots,\,
\bold e_n$. Then the values $f(\bold v,\bold w)$ and $g(\bold v)$ of a
bilinear form and of a quadratic form respectively are calculated by the
following formulas:
$$
\xalignat 2
&\hskip -2em
f(v,w)=\sum^n_{i=1}\sum^n_{j=1} f_{ij}\,v^i\,w^j,
&&g(v)=\sum^n_{i=1}\sum^n_{j=1} g_{ij}\,v^i\,v^j.
\tag1.9
\endxalignat
$$
In the case when $g_{ij}$ is a diagonal matrix, the formula for $g(\bold v)$
contains only the squares of coordinates of a vector $\bold v$:
$$
\hskip -2em
g(\bold v)=g_{11}\,(v^1)^2+\ldots+g_{nn}\,(v^n)^2.
\tag1.10
$$
This supports the term {\tencyr\char '074}quadratic form{\tencyr\char '076}.
Bringing a quadratic form to the form \thetag{1.10} by means of choosing
proper basis $\bold e_1,\,\ldots,\,\bold e_n$ in a linear space $V$ is one
of the problems which are solved in the theory of quadratic form.\par
     Let $\bold e_1,\,\ldots,\,\bold e_n$ and $\tilde\bold e_1,\,\ldots,
\tilde\bold e_n$ be two bases in a linear vector space $V$. Let's denote
by $S$ the transition matrix for passing from the first basis to the second
one. Denote $T=S^{-1}$. From \thetag{1.7} we easily derive the formula relating the components of a bilinear form $f(\bold v,\bold w)$ these two
bases. For this purpose it is sufficient to substitute the relationship
\thetag{5.8} of Chapter~\uppercase\expandafter{\romannumeral 1} into
the formula \thetag{1.7} and use the bilinearity of the form $f(\bold
v,\bold w)$:
$$
f_{ij}=f(\bold e_i,\bold e_j)=\sum^n_{k=1}\sum^n_{q=1}T^k_i\,
T^q_j\,f(\tilde\bold e_k,\tilde\bold e_q)=\sum^n_{k=1}\sum^n_{q=1}
T^k_i\,T^q_j\,\tilde f_{kq}.
$$
The reverse formula expressing $\tilde f_{kq}$ through $f_{ij}$ is derived
similarly:
$$
\pagebreak
\xalignat 2
&\hskip -2em
f_{ij}=\sum^n_{k=1}\sum^n_{q=1}T^k_i\,T^q_j\,\tilde f_{kq},
&&\tilde f_{kq}=\sum^n_{i=1}\sum^n_{j=1}S^i_k\,S^j_q\,f_{ij}.
\tag1.11
\endxalignat
$$
In matrix form these relationships are written as follows:
$$
\xalignat 2
&\hskip -2em
F=T^{\tr}\,\tilde F\,T,
&&\tilde F=S^{\tr}\,F\,S.
\tag1.12
\endxalignat
$$
Here $S^{\tr}$ and $T^{\tr}$ are two matrices obtained from 
$S$ and $T$ by transposition.\par
\head
\S\,2. Orthogonal complements\\
with respect to a quadratic form.
\endhead
\rightheadtext{\S\,2. Orthogonal complements \dots}
\definition{Definition 2.1} Two vectors $\bold v$ and $\bold w$ 
in a linear vector space $V$ are called {\it orthogonal to each
other with respect to the quadratic form} $g$ if $g(\bold v,\bold w)
=0$.
\enddefinition
\definition{Definition 2.2} Let $S$ be a subset of a linear vector
space $V$. The {\it orthogonal complement\/} of the subset $S$ {\it 
with respect to a quadratic form\/} $g(\bold v)$ is the set of 
vectors each of which is orthogonal to all vectors of $S$ with 
respect that quadratic form $g$. The orthogonal complement of $S$
is denoted $S_{\sssize\perp}\subset V$.
\enddefinition
    The orthogonal complement of a subset $S$ with respect to a 
quadratic form $g$ can be defined formally: $S_{\sssize\perp}
=\{\bold v\in V\!:\,\forall \bold w\,((\bold w\in S)\Rightarrow 
(g(\bold v,\bold w)=0))\}$. For the orthogonal complements 
determined by a quadratic form $g(\bold v)$ there is a theorem
analogous to theorems~3.1 and 3.2 in 
Chapter~\uppercase\expandafter{\romannumeral 3}.
\proclaim{Theorem 2.1} The operation of constructing orthogonal 
complements of sub\-sets $S\subset V$ with respect to a quadratic
form $g$ possesses the following properties:
\roster
\item $S_{\sssize\perp}$ is a subspace in $V$;
\item $S_1\subset S_2$ implies $(S_2)_{\sssize\perp}
      \subset (S_1)_{\sssize\perp}$;
\item $\langle S\rangle_{\sssize\perp}=S_{\sssize\perp}$, where
      $\langle S\rangle$ is the linear span of $S$;
\vskip 1ex
\item $\dsize\left(\,\,\shave{\bigcup_{i\in I}} S_i
       \right)_{\kern -2pt\raise 3pt\hbox{$\ssize\perp$}}=
       \bigcap_{i\in I}\,(S_i)_{\sssize\perp}$.
\endroster
\endproclaim
\demo{Proof} Let's prove the first item in the theorem for the beginning.
For this purpose we should verify two conditions from the definition of a
subspace.\par
     Let $\bold v_1,\bold v_2\in S_{\sssize\perp}$. Then $g(\bold v_1,
\bold w)=0$ and $g(\bold v_2,\bold w)=0$ for all $\bold w\in S$. Hence,
for all $\bold w\in S$ we have $g(\bold v_1+\bold v_2,\bold w)=g(\bold
v_1,\bold w)+g(\bold v_2,\bold w)=0$. This means that $\bold v_1+\bold
v_2\in S_{\sssize\perp}$, so the first condition is verified.\par
    Now let $\bold v\in S_{\sssize\perp}$. Then $g(\bold v,\bold w)=0$ 
for all $\bold w\in S$. Hence, for the vector $\alpha\cdot\bold v$ we
derive $g(\alpha\cdot\bold v,\bold w)=\alpha\,g(\bold v,\bold w)=0$.
This means that $\alpha\cdot\bold v\in S_{\sssize\perp}$. Thus, the 
first item of the theorem~2.1 is proved.\par
     In order to prove the inclusion $(S_2)_{\sssize\perp}\subset
(S_1)_{\sssize\perp}$ in the second item of the theorem~2.1
we consider an arbitrary vector $\bold v$ in $(S_2)_{\sssize\perp}$. 
From the condition $\bold v\in (S_2)_{\sssize\perp}$ we get $g(\bold v,
\bold w)=0$ for any $\bold w\in S_2$. But $S_1\subset S_2$, therefore,
the equality $g(\bold v,\bold w)=0$ is fulfilled for any $\bold w\in S_1$.
Then $\bold v\in (S_1)_{\sssize\perp}$. Thus, $\bold v\in (S_2)_{\sssize
\perp}$ implies $\bold v\in(S_1)_{\sssize\perp}$. This proves the required
inclusion.\par
     Now let's proceed to the third item of the theorem. Note that the
linear span of $S$ comprises this set: $S\subset\langle S\rangle$. Applying
the second item of the theorem, which is already proved, we get the
inclusion $\langle S\rangle_{\sssize\perp}\subset S_{\sssize\perp}$.
In order to prove the coincidence $\langle S\rangle_{\sssize\perp}= S_{\sssize\perp}$ we have to prove the converse inclusion $S_{\sssize
\perp}\subset\langle S\rangle_{\sssize\perp}$. For this purpose let's 
remember that the linear span $\langle S\rangle$ is formed by the
linear combinations
$$
\hskip -2em
\bold w=\alpha_1\cdot\bold w_1+\ldots+\alpha_r\cdot\bold w_r\text{, \ 
where \ }\bold w_i\in S.
\tag2.1
$$
Let $\bold v\in S_{\sssize\perp}$, then $g(\bold v,\bold w)=0$ for all
$\bold w\in S$. In particular, this is true for the vectors 
$\bold w_i$ in the expansion \thetag{2.1}: $g(\bold v,\bold w_i)=0$. 
Then from \thetag{2.1} we derive
$$
g(\bold v,\bold w)=\alpha_1\,g(\bold v,\bold w_1)+\ldots+\alpha_r
\,g(\bold v,\bold w_r)=0.
$$
Hence, $g(\bold v,\bold w)=0$ for all $\bold w\in \langle S\rangle$. 
This proves the converse inclusion $S_{\sssize\perp}\subset\langle
S\rangle_{\sssize\perp}$ and thus completes the proof of the coincidence
$\langle S\rangle_{\sssize\perp}=S_{\sssize\perp}$.\par
      In proving the fourth item of the theorem we introduce the
following notations:
$$
\xalignat 2
&S=\bigcup_{i\in I}\,S_i,
&
&\tilde S=\bigcap_{i\in I}\,(S_i)_{\sssize\perp}.
\endxalignat
$$
Let $\bold v\in S_{\sssize\perp}$. Then $g(\bold v,\bold w)=0$ for
all $\bold w\in S$. But $S_i\subset S$ for any $i\in I$. Therefore,
$g(\bold v,\bold w)=0$ for all $\bold w\in S_i$ and for all $i\in I$. 
This means that $\bold v$ belongs to each of the orthogonal complements
$(S_i)_{\sssize\perp}$ and, hence, it belongs to their intersection. 
This proves the inclusion $S_{\sssize\perp}\subset\tilde S$.\par
     Conversely, if $\bold v\in (S_i)_{\sssize\perp}$ for all $i\in I$,
then $g(\bold v,\bold w)=0$ for all $\bold w\in S_i$ and for all
$i\in I$. Hence, $g(\bold v,\bold w)=0$ for any vector $\bold w$ in
the union of all sets $S_i$. This proves the converse inclusion
$\tilde S\subset S_{\sssize\perp}$.\par
     The above two inclusions $S_{\sssize\perp}\subset\tilde S$ and
$\tilde S\subset S_{\sssize\perp}$ prove the coincidence of two sets
$S_{\sssize\perp}=\tilde S$. The theorem~2.1 is proved.
\qed\enddemo
\definition{Definition 2.3} The {\it kernel\/} of a quadratic form
$g(\bold v)$ in a linear vector space $V$ is the set $\Ker g
=V_{\sssize\perp}$ formed by vectors orthogonal to each vector of
the space $V$ with respect to the form $g$.
\enddefinition
\definition{Definition 2.4} A quadratic form with nontrivial kernel
$\Ker g\neq\{\bold 0\}$ is called a {\it degenerate\/} quadratic 
form. Otherwise, if $\Ker g=\{\bold 0\}$, then the form $g$ is called 
a {\it non-degenerate\/} quadratic form.
\enddefinition
     Due to the theorem~2.1 the kernel of a form $g(\bold v)$ is a 
subspace of the space $V$, where it is defined. The term 
{\tencyr\char '074}kernel{\tencyr\char '076} is not an occasional 
choice for denoting the set $V_{\sssize\perp}$. Each quadratic form
is associated with some mapping, for which the subspace 
$V_{\sssize\perp}$ is the kernel.
\definition{Definition 2.5} An {\it associated mapping} of a quadratic
form $g$ is the mapping $a_g\!:\,V\to V^*$ that takes each vector
$\bold v$ of the space $V$ to the linear functional $f_{\bold v}$
in the conjugate space $V^*$ determined by the relationship
$$
f_{\bold v}(\bold w)=g(\bold v,\bold w)\text{\ \ for all \ }\bold w\in V.
\tag2.2
$$
\enddefinition
     The associated mapping $a_g\!:\,V\to V^*$ is linear, this fact is immediate from the bilinearity of the form $g$. Its kernel $\Ker a_g$
coincides with the kernel of the form $g$. Indeed, the condition
$\bold v\in \Ker a_g$ means that the functional $f_v$ determined by
\thetag{2.2} is identically zero. Hence, $\bold v$ is orthogonal to 
all vectors $\bold w\in V$ with respect to the quadratic form 
$g(\bold v)$.\par
     The associated mapping $a_g$ relates orthogonal complements 
$S_{\sssize\perp}$ determined by the quadratic form $g$ and and
orthogonal complements $S^{\sssize\perp}$ in a dual space, which
we considered earlier in Chapter~\uppercase\expandafter{\romannumeral 3}.
\proclaim{Theorem 2.2} For any subset $S\subset V$ and for any quadratic
form $g(\bold v)$ in a linear vector space $V$ the set $S_{\sssize\perp}$
is the total preimage of the set $S^{\sssize\perp}$ under the associated
mapping $a_g$, i\.\,e\. $S_{\sssize\perp}=a_g^{-1}(S^{\sssize\perp})$.
\endproclaim
\demo{Proof} The condition $\bold v\in S_{\sssize\perp}$ means that
$g(\bold v,\bold w)=0$ for all $\bold w\in S$. But this equality can
be rewritten in the following way:
$$
g(\bold v,\bold w)=f_{\bold v}(\bold w)=a_g(\bold v)(\bold w)
=\langle a_g(\bold v)\,|\,\bold w\rangle=0\text{\ \ for all \ }
\bold w\in S.
$$
Hence, the condition $\bold v\in S_{\sssize\perp}$ is equivalent to
$a_g(\bold v)\in S^{\sssize\perp}$. This proves the required 
equality $S_{\sssize\perp}=a_g^{-1}(S^{\sssize\perp})$.
\qed\enddemo
     According to the definition~2.3, vectors of the kernel $\Ker g$ 
are orthogonal to all vectors of the space $V$ with respect to the form
$g$. Therefore $(\Ker g)_{\sssize\perp}=V$. If we apply the result of the
theorem~2.2 to the kernel $S=\Ker g$, we get 
$$
a_g^{-1}((\Ker g)^{\sssize\perp})=(\Ker g)_{\sssize\perp}=V.
$$
This result becomes more clear if we write it in the following
equivalent form:
$$
\hskip -2em
\Img a_g=a_g(V)\subseteq (\Ker g)^{\sssize\perp}.
\tag2.3
$$
\proclaim{Corollary 1} The image of the associated mapping
$a_g$ is enclosed into the orthogonal complement to its kernel
$(\Ker a_g)^{\sssize\perp}$, i\.\,e\. $\Img a_g\subseteq
(\Ker a_g)^{\sssize\perp}$.
\endproclaim
     This corollary of the theorem~2.2 is derived from the formula
\thetag{2.3} if we take into account $\Ker g=\Ker a_g$. For a quadratic
form $g$ in a finite-dimensional space $V$ it can be strengthened.
\proclaim{Corollary 2} For a quadratic form $g(\bold v)$ in a
finite-dimensional linear vector space $V$ the image of the associated
mapping $a_g\!:\,V\to V^*$ coincides with the orthogonal complement of
its kernel $\Ker a_g$:
$$
\hskip -2em
\Img a_g=(\Ker a_g)^{\sssize\perp}.
\tag2.4
$$
\endproclaim
\demo{Proof} Using the theorem~9.4 from 
Chapter~\uppercase\expandafter{\romannumeral 1}, we calculate the 
dimension of the image \ $\Img a_g$ of the associated mapping:
$$
\dim(\Img a_g)=\dim V-\dim(\Ker a_g).
$$
The dimension of the orthogonal complement of $\Ker a_g$ in the dual 
space is determined by the theorem~3.4 in 
Chapter~\uppercase\expandafter{\romannumeral 3}:
$$
\dim(\Ker a_g)^{\sssize\perp}=\dim V-\dim(\Ker a_g).
$$
As we can see, the dimensions of these two subspaces are equal to each
other. Therefore. we can apply the above corollary~1 and the item
\therosteritem{3} of the theorem~4.5 from
Chapter~\uppercase\expandafter{\romannumeral 1}. As a result we get
the required equality \thetag{2.4}.
\qed\enddemo
\proclaim{Theorem~2.3} Let $U\varsubsetneq V$ be a subspace of a
finite-dimensional space $V$ comprising the kernel of a quadratic
form $g$. For any vector $\bold v\notin U$ there exists a vector 
$\bold w\in V$ such that $g(\bold v,\bold w)\neq 0$ and $g(\bold 
v,\bold u)=0$ for all $\bold u\in U$.
\endproclaim
\demo{Proof} This theorem is an analog of the theorem~1.2 from
Chapter~\uppercase\expandafter{\romannumeral 3}. It's proof
is  essentially based on that theorem. Applying the theorem~1.2 from
Chapter~\uppercase\expandafter{\romannumeral 3}, we get that there
exist a linear functional $f\in V^*$ such that $f(\bold v)\neq 0$ 
and $f(\bold u)=\langle f\,|\,\bold u\rangle=0$ for all $\bold u\in U$.
Due to the last condition this functional $f$ belongs to the
orthogonal complement $U^{\sssize\perp}$. From the inclusion $\Ker g
\subset U$, applying the item \therosteritem{2} of the theorem~3.1 
from Chapter~\uppercase\expandafter{\romannumeral 3}, we get
$U^{\sssize\perp}\subset (\Ker g)^{\sssize\perp}$. Hence, we conclude
that $f\in(\Ker g)^{\sssize\perp}$.\par
     Now we apply the corollary~2 from the theorem~2.2. From this corollary
we obtain that $(\Ker g)^{\sssize\perp}=\Img a_g$. Hence, $f\in\Img a_g$
and there is a vector $\bold w\in V$ that is taken to $f$ by the
associated mapping $a_g$, i\.\,e\. $f=a_g(\bold w)$. Then
$$
\align
&g(\bold v,\bold w)=a_g(\bold w)(\bold v)=f(\bold v)\neq 0,\\
&g(\bold v,\bold u)=a_g(\bold w)(\bold u)=f(\bold u)=0\text{ \ 
for all \ }\bold u\in U.
\endalign
$$
Due to these relationship we find that $\bold w$ is the very vector
that we need to complete the proof of the theorem.
\qed\enddemo
\proclaim{Theorem~2.4} Let $V$ be a finite-dimensional linear vector
space and let $U$ and $W$ be two subspaces of $V$ comprising the
kernel $\Ker g$ of a quadratic form $g$. Then the conditions 
$W=U_{\sssize\perp}$ and $U=W_{\sssize\perp}$ are equivalent to
each other.
\endproclaim
\demo{Proof} The theorem~2.4 is an analog of the theorem~3.3 from
Chapter~\uppercase\expandafter{\romannumeral 3}. The proofs of these 
two theorems are also very similar.\par
     Suppose that the condition $W=U_{\sssize\perp}$ is fulfilled. 
Then for any vector $\bold w\in W$ and for any vector $\bold u\in U$ 
we have the relationship $g(\bold w,\bold u)=0$. The set
$W_{\sssize\perp}$ is formed by vectors orthogonal to all vectors of
$W$ with respect to the quadratic form $g$. Therefore, we have the 
inclusion $U\subset W_{\sssize\perp}$.\par
     Further proof is by contradiction. Assume that $U\neq W_{\sssize
\perp}$. Then there is a vector $v_0$ such that $\bold v_0\in W_{\sssize
\perp}$ and $\bold v_0\not\in U$. In this situation we can apply the
theorem~2.3 which says that there exists a vector $\bold v$ such that
$g(\bold v,\bold v_0)\neq 0$ and $g(\bold v,\bold u)=0$ for all
$\bold u\in U$. The latter condition means that $\bold v\in U_{\sssize
\perp}=W$. Then the other condition $g(\bold v,\bold v_0)\neq 0$
contradicts to the initial choice $\bold v_0\in W_{\sssize\perp}$.
This contradiction shows that the assumption $U\neq W_{\sssize\perp}$
is not true and we have the coincidence $U=W_{\sssize\perp}$. Thus,
$W=U_{\sssize\perp}$ implies $U=W_{\sssize\perp}$. We can swap 
$U$ and $W$ and obtain that $U=W_{\sssize\perp}$ implies $W=U_{\sssize
\perp}$. Hence, these two conditions are equivalent.
\qed\enddemo
     The proposition of the theorem~2.3 can be reformulated as follows:
for a subspace $U\subset V$ in a finite-dimensional space $V$ the condition
$\Ker g\subset U$ means that double orthogonal complement of $U$ coincides 
with that space: $(U_{\sssize\perp})_{\sssize\perp}=U$. For an arbitrary 
subset $S\subset V$ of a finite-dimensional space $V$ one can derive 
$$
\hskip -2em
(S_{\sssize\perp})_{\sssize\perp}=\langle S\rangle+\Ker g.
\tag2.5
$$
Let's prove the relationship \thetag{2.5}. Note that vectors of the kernel
$\Ker g$ are orthogonal to all vectors of $V$. Therefore, joining the vectors
of the kernel \pagebreak $\Ker g$ to $S$, we do not change the orthogonal complement
of this subset:
$$
S_{\sssize\perp}=(S\cup\Ker g)_{\sssize\perp}.
$$
Now let's apply the item \therosteritem{3} of the theorem~2.1. This yields
$$
S_{\sssize\perp}=(S\cup\Ker g)_{\sssize\perp}
=\langle S\cup\Ker g\rangle_{\sssize\perp}
=(\langle S\rangle+\Ker g)_{\sssize\perp}.
$$
The subspace $U=\langle S\rangle+\Ker g$ comprises the kernel of the 
form $g$. Therefore, $(U_{\sssize\perp})_{\sssize\perp}=U$. This completes
the proof of the relationship \thetag{2.5}:
$$
(S_{\sssize\perp})_{\sssize\perp}=
((\langle S\rangle+\Ker g)_{\sssize\perp})_{\sssize\perp}=
\langle S\rangle+\Ker g.
$$
\proclaim{Theorem~2.5} In the case of finite-dimensional linear vector
space $V$ for any subspace $U$ of $V$ we have the equality
$$
\hskip -2em
\dim U+\dim U_{\sssize\perp}=\dim V+\dim(\Ker g\cap U),
\tag2.6
$$
where $U_{\sssize\perp}$ is the orthogonal complement of $U$ with respect
to the form $g$.
\endproclaim
\demo{Proof} The vectors of the kernel $\Ker g$ are orthogonal to all
vectors of the space $V$, therefore, joining them to $U$, we do not change 
the orthogonal complement $U_{\sssize\perp}$. Let's denote $W=U+\Ker g$. Then $U_{\sssize\perp}=W_{\sssize\perp}$. Applying the theorem~6.4 from
Chapter~\uppercase\expandafter{\romannumeral 1}, for the dimension
of $W$ we derive the formula 
$$
\hskip -2em
\dim W=\dim U+\dim(\Ker g)-\dim(\Ker g\cap U).
\tag2.7
$$
Now let's apply the theorem~2.2 to the subset $S=W$. This yields
$W_{\sssize\perp}=a_g^{-1}(W^{\sssize\perp})$. Note that $\Ker g
\subset W$, this differs $W$ from the initial subspace $U$. Let's
apply the item \therosteritem{2} of the theorem~3.1 to the inclusion
$\Ker g\subset W$ and take into account the corollary~2 of the 
theorem~2.2. This yields 
$$
W^{\sssize\perp}\subset (\Ker g)^{\sssize\perp}=\Img a_g.
$$
The inclusion $W^{\sssize\perp}\subset\Img a_g$ means that 
the preimage of each element $f\in W^{\sssize\perp}$ under the
mapping $a_g$ is not empty, while the equality $W_{\sssize\perp}=
a_g^{-1}(W^{\sssize\perp})$ shows that such preimage is enclosed
into $W_{\sssize\perp}$. Therefore, $W_{\sssize\perp}=
a_g^{-1}(W^{\sssize\perp})$ implies the equality $a_g(W_{\sssize
\perp})=W^{\sssize\perp}$.\par
    Now let's consider the restriction of the associated mapping 
$a_g$ to the subspace $W_{\sssize\perp}$. We denote this restriction
by $a$:
$$
\hskip -2em
a\!:\,W_{\sssize\perp}\to V^*.
\tag2.8
$$
The kernel of the mapping \thetag{2.8} coincides with the kernel of
the non-restricted mapping $a_g$ since $\Ker a_g=\Ker g\subset W_{\sssize\perp}$. For the image of this mapping we have 
$$
\Img a=a_g(W_{\sssize\perp})=W^{\sssize\perp}.
$$
Let's apply the theorem on the sum of dimensions of the kernel and the 
image (see theorem~9.4 in
Chapter~\uppercase\expandafter{\romannumeral 1}) to the mapping $a$:
$$
\pagebreak 
\hskip -2em
\dim(\Ker g)+\dim W^{\sssize\perp}=\dim W_{\sssize\perp}
\tag2.9
$$
In order to determine the dimension of $W^{\sssize\perp}$ we apply the
relationship
$$
\hskip -2em
\dim W+\dim W^{\sssize\perp}=\dim V
\tag2.10
$$
which follows from the theorem~3.4 of 
Chapter~\uppercase\expandafter{\romannumeral 3}. Now let's add the
relationships \thetag{2.7} and \thetag{2.9} and subtract the relationship
\thetag{2.10}. Taking into account the coincidence $W_{\sssize\perp}
=U_{\sssize\perp}$, we get the required equality \thetag{2.6}.
\qed\enddemo
     The analogs of the relationships \thetag{3.4} from
Chapter~\uppercase\expandafter{\romannumeral 3} in present case are
the relationships $\{\bold 0\}_{\sssize\perp}=V$ and 
$V_{\sssize\perp}=\Ker g$.\par
\proclaim{Theorem 2.6} In the case of finite-dimensional linear vector
space $V$ equipped with a quadratic form $g$ for any family of subspaces 
in $V$, each of which comprises the kernel $\Ker g$, the following
relationships are fulfilled:
$$
\xalignat 2
&\hskip -2em
 \left(\,\,\shave{\sum_{i\in I}} U_i\right)_{\kern -1pt
 \raise 3pt\hbox{$\ssize\perp$}}=
 \bigcap_{i\in I}\,(U_i)_{\sssize\perp},
&&\left(\,\,\shave{\bigcap_{i\in I}} U_i\right)_{\kern -1pt
 \raise 2pt\hbox{$\ssize\perp$}}=
 \sum_{i\in I}\,(U_i)_{\sssize\perp}.
\tag2.11
\endxalignat
$$
\endproclaim
\demo{Proof} In proving the first relationship \thetag{2.11} the
condition $\Ker g\subset U_i$ is inessential. This relationship
is derived from the items \therosteritem{3} and \therosteritem{4}
of the theo\-rem~2.1 if we take into account that the sum of 
subspaces is the linear span of the union of these subspaces.\par
     The second relationship \thetag{2.11} is derived from the
first one. From the condition $\Ker g\in U_i$ we derive that
$((U_i)_{\sssize\perp})_{\sssize\perp}=U_i$ (see theorem~2.4).
Let's denote $(U_i)_{\sssize\perp}=V_i$ and apply the first relationship
\thetag{2.11} to the family of subsets $V_i$:
$$
\left(\,\,\shave{\sum_{i\in I}} (U_i)_{\sssize\perp}
\right)_{\kern -1pt\raise 3pt\hbox{$\ssize\perp$}}=
\left(\,\,\shave{\sum_{i\in I}} V_i\right)_{\kern -1pt
\raise 3pt\hbox{$\ssize\perp$}}=
\bigcap_{i\in I}\,(V_i)_{\sssize\perp}=
\bigcap_{i\in I}\,U_i.
$$
Now it is sufficient to pass to orthogonal complements in left and
right hand sides of the above equality and apply the theorem~2.4 
again. This yields the required equality \thetag{2.11}. The theorem
is proved.
\qed\enddemo
\head
\S\,3. Transformation of a quadratic form\\to its canonic form.
Inertia indices and signature.
\endhead
\rightheadtext{\S\,3. Transformation to a canonic form.}
\definition{Definition 3.1} A subspace $U$ in a linear vector space
$V$ is called {\it regular with respect to a quadratic form} $g$ if
$U\cap U_{\sssize\perp}\subseteq\Ker g$.
\enddefinition
\proclaim{Theorem 3.1} Let $U$ be a subspace in a finite-dimensional
space $V$ regular with respect to a quadratic form $g$. Then 
$U+U_{\sssize\perp}=V$.
\endproclaim
\demo{Proof} Let's denote $W=U+U_{\sssize\perp}$ and then let's 
calculate the dimension of the subspace $W$ applying the theorem~6.4
from Chapter~\uppercase\expandafter{\romannumeral 1}:
$$
\dim W=\dim U+\dim U_{\sssize\perp}-
\dim(U\cap U_{\sssize\perp}).
$$
The vectors of the kernel $\Ker g$ are perpendicular to all vectors 
of the space $V$. Therefore, $\Ker g\subseteq U_{\sssize\perp}$.
Moreover, due to the regularity of $U$ with respect to the form $g$  
we have $U\cap U_{\sssize\perp}\subseteq\Ker g$. Therefore, we derive
$$
U\cap U_{\sssize\perp}=(U\cap U_{\sssize\perp})\cap\Ker g=
U\cap(U_{\sssize\perp}\cap\Ker g)=U\cap\Ker g.
$$
Because of the equality $U\cap U_{\sssize\perp}=U\cap\Ker g$ the above
formula for the dimension of the subspace $W$ can be written as follows:
$$
\hskip -2em
\dim W=\dim U+\dim U_{\sssize\perp}-
\dim(U\cap\Ker g).
\tag3.1
$$
Let's compare \thetag{3.1} with the formula \thetag{2.6} from the
theorem~2.5. This comparison yields $\dim W=\dim V$. Now, applying 
the item \therosteritem{3} of the theorem~4.5 from
Chapter~\uppercase\expandafter{\romannumeral 1}, we get $W=V$.
The theorem is proved.
\qed\enddemo
\proclaim{Theorem 3.2} Let $U$ be a subspace of a finite-dimensional
space $V$ regular with respect to a quadratic form $g$. If $U_{\sssize
\perp}\neq\Ker g$, then there exists a vector $\bold v\in U_{\sssize
\perp}$ such that $g(\bold v)\neq 0$.
\endproclaim
\demo{Proof} The proof is by contradiction. Assume that there is no
vector $\bold v\in U_{\sssize\perp}$ such that $g(\bold v)\neq 0$.
Then the numeric function $g(\bold v)$ is identically zero in the
subspace $U_{\sssize\perp}$. Due to the recovery formula \thetag{1.6}
the numeric function $g(\bold v,\bold w)$ is also identically zero 
for all $\bold v,\bold w\in U_{\ssize\perp}$.\par
     Now let's apply the theorem~3.1 and expand an arbitrary vector
$\bold x\in V$ into a sum of two vectors $\bold x=\bold u+\bold w$, 
where $\bold u\in U$ and $\bold w\in U_{\sssize\perp}$. Then for 
an arbitrary vector $\bold v$ of the subspace $U_{\ssize\perp}$ we 
derive
$$
g(\bold v,\bold x)=g(\bold v,\bold u+\bold w)=g(\bold v,\bold u)
+g(\bold v,\bold w)=0+0=0.
$$
The first summand $g(\bold v,\bold u)$ in right hand side of the above
equality is zero since the subspaces $U$ and $U_{\sssize\perp}$ are
orthogonal to each other. The second summand $g(\bold v,\bold w)$ is
zero due to our assumption in the beginning of the proof. Since 
$g(\bold v,\bold x)=0$ for an arbitrary vector $\bold x\in V$, we get 
$\bold v\in\Ker g$. But $\bold v$ is an arbitrary vector of the subspace
$U_{\ssize\perp}$. Therefore, $U_{\sssize\perp}\subseteq\Ker g$. 
The converse inclusion $\Ker g\subseteq U_{\sssize\perp}$ is always
valid. Hence, $U_{\sssize\perp}=\Ker g$, which contradicts the 
hypothesis of the theorem. This contradiction means that the assumption,
which we have made in the beginning of our proof, is not valid and, thus,
it proves the existence of a vector $\bold v\in U_{\sssize\perp}$ such 
that $g(\bold v)\neq 0$. The theorem is proved.
\qed\enddemo
\proclaim{Theorem 3.3} For any quadratic form $g$ in a finite-dimensional
vector space $V$ there exists a basis $\bold e_1,\,\ldots,\,\bold e_n$
such that the matrix of $g$ is diagonal in this basis.
\endproclaim
\demo{Proof} The case $g=0$ is trivial. The matrix of the zero quadratic
form $g$ is purely zero in any basis. The square $n\times n$ matrix, which
is purely zero, is obviously a diagonal matrix.\par
     Suppose that $g\not\equiv 0$. We shall prove the theorem by induction
on the dimension of the space $\dim V=n$. In the case $n=1$ the proposition
of the theorem is trivial: any $1\times 1$ matrix is a diagonal matrix.\par
     Suppose that the theorem is valid for any quadratic form in any space
of the dimension less than $n$. Let's consider the subspace $U=\Ker g$. 
It is regular with respect to the form $g$ and $U_{\sssize\perp}=V$.
Therefore, we can apply the theorem~3.2. According to this theorem, there
exists a vector $\bold v_0\not\in U$ such that $g(\bold v_0)\neq 0$.
Let's consider the subspace $W$ obtained by joining $\bold v_0$ to
$U=\Ker g$:
$$
\hskip -2em
W=\Ker g+\langle\bold v_0\rangle=U\oplus\langle\bold v_0\rangle.
\tag3.2
$$
This subspace $W$ determines the following two cases: $W=V$ or 
$W\neq V$.\par
     In the case $W=V$ we choose a basis $\bold e_1,\,\ldots,\,
\bold e_s$ in the kernel $\Ker g$ and complete it by one additional vector
$\bold e_{s+1}=\bold v_0$. As a result we get the basis in $V$. The matrix
of the quadratic form $g$ in this basis is a matrix almost completely filled with zeros, indeed, for $i=1,\,\ldots,\,s$ and $j=1,\,\ldots,\,s+1$ we have
$g_{ij}=g_{ji}=g(\bold e_i,\bold e_j)=0$ since $\bold e_i\in\Ker g$. 
The only nonzero element is $g_{s+1\,s+1}$, it is a diagonal element:
$g_{s+1\,s+1}=g(\bold e_{s+1},\bold e_{s+1})=g(\bold v_0)\neq 0$.\par
     In the case $W\neq V$ we consider the intersection $W\cap W_{\sssize
\perp}$. Let $\bold w\in W\cap W_{\sssize\perp}$. Then from \thetag{3.2} we derive $\bold w=\alpha\cdot\bold v_0+\bold u$, where $\bold u\in\Ker g$. 
Since $\bold w$ is a vector of $W$ ans simultaneously it is a vector of
$W_{\sssize\perp}$, it should be orthogonal to itself with respect to
the quadratic form $g$:
$$
\aligned
&g(\bold w,\bold w)=g(\alpha\cdot\bold v_0+\bold u,\alpha\cdot\bold v_0
+\bold u)=\\
&=\alpha^2\,g(\bold v_0,\bold v_0)+2\alpha\,g(\bold v_0,\bold u)
+g(\bold u,\bold u)=0.
\endaligned
\tag3.3
$$
But $\bold u\in\Ker g$, therefore, $g(\bold v_0,\bold u)=0$ and $g(\bold
u,\bold u)=0$, while $g(\bold v_0,\bold v_0)=g(\bold v_0)\neq 0$. 
Hence, from \thetag{3.3} we get $\alpha=0$. This means that $\bold w
=\bold u\in\Ker g$. Thus, we have proved the inclusion 
$W\cap W_{\sssize\perp}\subseteq\Ker g$, which means the regularity of
the subspace $W$ with respect to the quadratic form $g$.\par
     Now let's apply the theorem~3.1. It yields the expansion $V=W
+W_{\sssize\perp}$. Note that $\bold v_0\in W$, but $\bold v_0\notin W
\cap W_{\sssize\perp}$. This follows from $g(\bold v_0,\bold v_0)\neq 0$.
Hence, $\bold v_0\notin W_{\sssize\perp}$ and $W_{\sssize\perp}\neq V$.
This means that the dimension of the subspace $W_{\sssize\perp}$ is less
than $n$. The formula \thetag{2.6} yield the exact value of this dimension
$$
\dim W_{\ssize\perp}=\dim V+\dim(U\cap\Ker g)-\dim U=n-1.
$$
Let's consider the restriction of $g$ to the subspace $W_{\ssize\perp}$.
We can apply the inductive hypothesis to $g$ in $W_{\ssize\perp}$. Let
$\bold e_1,\,\ldots,\,\bold e_{n-1}$ be a basis of the subspace $W_{\ssize
\perp}$ in which the matrix of the restriction of $g$ to $W_{\ssize\perp}$ 
is diagonal:
$$
\hskip -2em
g_{ij}=g_{ji}=g(\bold e_i,\bold e_j)=0\text{ \ for \ }i<j\leqslant n-1.
\tag3.4
$$
We complete this basis by one vector $\bold e_n=\bold v_0$. Since
$\bold v_0\notin W_{\ssize\perp}$ the extended system of vectors
$\bold e_1,\,\ldots,\,\bold e_n$ is linearly independent and, hence,
is a basis of $V$. Let's find the matrix of the quadratic form $g$
in the extended basis. For the elements in the extension of this 
matrix we obtain the relationships
$$
\hskip -2em
g_{in}=g_{ni}=g(\bold e_i,\bold e_n)=0\text{ \ for \ }i<n.
\tag3.5
$$
They follow from the orthogonality of $\bold e_i$ and $\bold e_n$ in
\thetag{3.5}. Indeed, $\bold e_n\in W$ and $\bold e_i\in W_{\ssize
\perp}$. The relationships \thetag{3.4} and \thetag{3.5} taken together
mean that the matrix of the quadratic form $g$ is diagonal in the
basis $\bold e_1,\,\ldots,\,\bold e_n$. The inductive step is over
and the theorem is completely proved.
\qed\enddemo
     Let $g$ be a quadratic form in a finite dimensional space $V$ and
let $\bold e_1,\,\ldots,\,\bold e_n$ be a basis in which the matrix of
$g$ is diagonal. Then the value of $g(\bold v)$ can be calculated by 
formula \thetag{1.10}. A part of the diagonal elements $g_{11},\,\ldots,
\,g_{nn}$ can be equal to zero. Let's denote by $s$ the number of such
elements. We can renumerate the basis vectors $\bold e_1,\,\ldots,\,
\bold e_n$ so that 
$$
\hskip -2em
g_{11}=\ldots=g_{ss}=0.
\tag3.6
$$
The first $s$ vectors of the  basis, which correspond to the matrix 
elements \thetag{3.6}, belong to the kernel of the form $\Ker g$. 
Indeed, if $\bold w=\bold e_i$ for $i=1,\,\ldots,\,s$, then $g(\bold v,
\bold w)=0$ for all vectors $\bold v\in V$. This fact can be easily 
derived with the use of formulas \thetag{1.9}.\par
     Conversely, suppose that $\bold w\in\Ker g$. Then for an arbitrary
vector $\bold v\in V$ we have the following relationships:
$$
g(\bold v,\bold w)=\sum^n_{i=1}\sum^n_{j=1} g_{ij}\,v^i\,w^j=
\sum^n_{i=s+1} g_{ii}\,v^i\,w^i=0.
$$
Since $\bold v\in V$ is an arbitrary vector, the above equality should
be fulfilled identically in $v^{s+1},\,\ldots,\,v^n$. But $g_{ii}\neq 0$ 
for $i\geqslant s+1$, therefore, $w^{s+1}=\ldots=w^n=0$. From these equalities for the vector $\bold w$ we derive
$$
\bold w=w^1\cdot\bold e_1+\ldots+w^s\cdot\bold e_s.
$$
The conclusion is that any vector $\bold w$ of the kernel $\Ker g$ can
be expanded into a linear combination of the first $s$ basis vectors.
Hence, these basis vectors $\bold e_1,\,\ldots,\,\bold e_s$ form a
basis in $\Ker g$. The above considerations prove the following proposition
that we present in the form of a theorem.
\proclaim{Theorem 3.4} The number of zeros on the diagonal of the matrix
of a quadratic form $g$, brought to the diagonal form, is a {\it geometric
invariant} of the form $g$. It does not depend on the method used for bringing this matrix to a diagonal form and coincides with the dimension of the kernel of the quadratic form: $s=\dim(\Ker g)$.
\endproclaim
\definition{Definition 3.2} The number $s=\dim(\Ker g)$ is called the
{\it zero inertia index} of a quadratic form $g$.
\enddefinition
     Let $g$ be a quadratic form in a linear vector space over the field
of complex numbers $\Bbb C$ such that its matrix is diagonal a basis
$\bold e_1,\,\ldots,\,\bold e_n$. Suppose that $s$ is the zero inertia 
index of the quadratic form $g$. Without loss of generality we can assume
that the first $s$ basis vectors $\bold e_1,\,\ldots,\,\bold e_s$ form
a basis in the kernel $\Ker g$. We define the numbers 
$\gamma_1,\,\ldots,\,\gamma_n$ by means of formula
$$
\hskip -2em
\gamma_i=
\cases 1&\text{\ \ for \ }i\leqslant s,\\
\sqrt{g_{ii}}&\text{\ \ for \ }i> s.
\endcases
\tag3.7
$$
Remember that for any complex number one can take its square root which
is again a complex number. Complex numbers \thetag{3.7} are nonzero. We
use them in order to construct the new basis:
$$
\hskip -2em
\tilde\bold e_i=(\gamma_i)^{-1}\cdot\bold e_i,\ i=1,\ldots,n.
\tag3.8
$$
The matrix of the quadratic form $g$ in the new basis \thetag{3.8}
is again a diagonal matrix. Indeed, we can explicitly calculate the
matrix elements:
$$
\pagebreak
\tilde g_{ij}=g(\tilde\bold e_i,\tilde\bold e_j)=
(\gamma_i\,\gamma_j)^{-1}\,g_{ij}=0\text{\ \ for \ }i\neq j.
$$
For the diagonal elements of the matrix of $g$ we derive
$$
\tilde g_{ii}=g(\tilde\bold e_i,\tilde\bold e_i)=
(\gamma_i)^{-2}\,g_{ii}=
\cases 0&\text{\ \ for \ }i\leqslant s,\\
       1&\text{\ \ for \ }i> s.
\endcases
$$
The matrix of the quadratic form $g$ in the basis $\tilde\bold
e_1,\,\ldots,\,\tilde\bold e_n$ has the following form which is used
to be called the {\it canonic form} of the matrix of a quadratic form
over the field of complex numbers $\Bbb C$:
$$
\hskip -2em
\Cal G=
\aligned
&\Vmatrix
0 & \      & \ & \ & \      & \ \\
\ & \ddots & \ & \ & \      & \ \\
\ & \      & 0 & \ & \      & \ \\
\ & \      & \ & 1 & \      & \ \\
\ & \      & \ & \ & \ddots & \ \\
\ & \      & \ & \ & \      & 1 \\
\endVmatrix
\endaligned
\aligned
&\left.\vphantom{\vrule height 4ex depth 4ex}\right\} s\\
&\left.\vphantom{\vrule height 4ex depth 4ex}\right\} n-s
\endaligned
\tag3.9
$$
The matrix $\Cal G$ in \thetag{3.9} is a diagonal matrix, its diagonal
is filled with $s$ zeros and $n-s$ ones, where $s=\dim\Ker g$.\par
     In the case of a linear vector space over the field of real numbers
$\Bbb R$ the canonic form of the matrix of a quadratic form is 
different from \thetag{3.9}. Let $\bold e_1,\,\ldots,\,\bold e_n$ be a
basis in which the matrix of $g$ is diagonal. Diagonal elements of this
matrix now is subdivided into three groups: zero elements, positive
elements, and negative elements. If $s$ is the number of zero elements
and $r$ is the number of positive elements, then remaining $n-s-r$
elements on the diagonal are negative numbers.
Without loss of generality we can assume that the basis vectors 
$\bold e_1,\,\ldots,\,\bold e_n$ are enumerated so that $g_{ii}=0$ 
for $i=1,\,\ldots,\,s$ and $g_{ii}>0$ for $i=s+1,\,\ldots,\,s+r$. Then
$g_{ii}<0$ for $i=s+r+1,\,\ldots,\,n$. In the field of reals we can take
the square root only of non-negative numbers. Therefore, here we define
$\gamma_1,\,\ldots,\,\gamma_n$ a little bit differently than it was done
in \thetag{3.7} for complex numbers:
$$
\hskip -2em
\gamma_i=
\cases 1&\text{\ \ for \ }i\leqslant s,\\
\sqrt{|g_{ii}|}&\text{\ \ for \ }i> s.
\endcases
\tag3.10
$$\par
     By means of \thetag{3.10} we define new basis $\tilde\bold e_1,\,
\ldots,\,\tilde\bold e_n$ using the formulas \thetag{3.9}. Here is the
matrix of the quadratic form $g$ in this basis:
$$
\Cal G=
\aligned
&\Vmatrix
0\vphantom{\vrule height 3ex}
  & \      & \ & \ & \      & \ & \ & \      & \ \\
\ & \ddots & \ & \ & \      & \ & \ & \      & \ \\
\ & \      & 0 & \ & \      & \ & \ & \      & \ \\
\ & \      & \ & 1 & \      & \ & \ & \      & \ \\
\ & \      & \ & \ & \ddots & \ & \ & \      & \ \\
\ & \      & \ & \ & \      & 1 & \ & \      & \ \\
\ & \      & \ & \ & \      & \ &-1 & \      & \ \\
\ & \      & \ & \ & \      & \ & \ & \ddots & \ \\
\ & \      & \ & \ & \      & \ & \ & \      &-1 \\
\endVmatrix
\endaligned
\aligned
&\left.\vphantom{\vrule height 4ex depth 4ex}\right\} s\\
&\left.\vphantom{\vrule height 4ex depth 4ex}\right\} r_p\\
&\left.\vphantom{\vrule height 4ex depth 4ex}\right\} r_n
\endaligned
\tag3.11
$$
\definition{Definition 3.3} The formula \thetag{3.11} defines 
the {\it canonic form\/} of the matrix of a quadratic form $g$ 
in a space over the real numbers $\Bbb R$. The integers $r_p$ 
and $r_n$ that determine the number of plus ones and the number
of minus ones on the diagonal of the matrix \thetag{3.10} are 
called the {\it positive inertia index\/} and the {\it negative 
inertia index\/} of the quadratic form $g$ respectively.
\enddefinition
\proclaim{Theorem 3.5} The positive and the negative inertia indices
$r_p$ and $r_n$ of a quadratic form $g$ in a space over the field
of real numbers $\Bbb R$ are geometric invariants of $g$. They do not
depend on a particular way how the matrix of $g$ was brought to the
diagonal form. 
\endproclaim
\demo{Proof} Let $\bold e_1,\,\ldots,\,\bold e_n$ be a basis of a 
space $V$ in which the matrix of $g$ has the canonic form \thetag{3.11}.
Let's consider the following subspaces:
$$
\xalignat 2
&\hskip -2em
U_+=\langle e_1,\ldots,e_{s+r_p}\rangle,
&&U_-=\langle e_{s+r_p+1},\ldots,e_n\rangle.
\tag3.12
\endxalignat
$$
The intersection of $U_+$ and $U_-$ is trivial, $\dim U_+=s+r_p$, 
$\dim U_-=r_n$, and for their sum we have $U_+\oplus U_-=V$.
\par
     Let's take a vector $\bold v\in U_+$. The value of the quadratic form
$g$ for that vector is determined by the matrix \thetag{3.11} according to
the formula \thetag{1.10}:
$$
g(\bold v)=\sum^{s+r_p}_{i=s+1}(v^i)^2.
$$
The sum of squares in the right hand side of this equality is a 
non-negative quantity, i\.\,e\. $g(\bold v)\geqslant 0$ for all
$\bold v\in U_+$.\par
     Now let's take a vector $\bold v\in U_-$. For this vector 
the formula \thetag{1.10} is written as 
$$
g(\bold v)=\sum^n_{i=s+r_p+1}(-(v^i)^2).
$$
If $v\neq\bold 0$, then at least one summand in right hand side 
is nonzero. Hence,  $g(\bold v)<0$ for all nonzero vectors of 
the subspace $U_-$.\par
     Suppose that $\tilde\bold e_1,\,\ldots,\,\tilde\bold e_n$ is
some other basis in which the matrix of $g$ has the canonic form.
Denote by $\tilde s$, $\tilde r_p$, and $\tilde r_n$ the inertia 
indices of $g$ in this basis. The zero inertia indices in both
bases are the same $s=\tilde s$ since they are determined by
the kernel of $g$: $s=\dim(\Ker g)$ and $\tilde s=\dim(\Ker g)$.
\par
     Let's prove the coincidence of the positive and the negative
inertia indices in two bases. For this purpose we consider the
subspaces $\tilde U_+$ and $\tilde U_-$ determined by the relationships
of the form \thetag{3.12} but for the 
{\tencyr\char '074}wavy{\tencyr\char '076} basis $\tilde\bold e_1,\,
\ldots,\,\tilde\bold  e_n$. If we assume that $r_p\neq\tilde r_p$, then
$r_p>\tilde r_p$ or $r_p<\tilde r_p$. For the sake of certainty suppose 
that $r_p>\tilde r_p$. Then we calculate the dimensions of $U_+$ and
$\tilde U_-$:
$$
\xalignat 2
&\dim U_+=s+r_p,
&&\dim\tilde U_-=\tilde r_n=n-s-\tilde r_p.
\endxalignat
$$
For the sum of dimensions of these two subspaces $U_+$ and $\tilde U_-$
we get the equality $\dim U_++\dim \tilde U_-=n+(r_p-\tilde r_p)$. Due
to the above assumption $r_p>\tilde r_p$ we derive 
$$
\hskip -2em
\dim U_++\dim \tilde U_->\dim V
\tag3.13
$$
From the natural inclusion $U_++\tilde U_-\subseteq V$ we get 
$\dim(U_++\tilde U_-)\leqslant\dim V$. Using this estimate together 
with the inequality \thetag{3.13} and applying the theorem~6.4 of 
Chapter~\uppercase\expandafter{\romannumeral 1} to them, we derive
$\dim(U_+\cap\tilde U_-)>0$. Hence, the intersection $U_+\cap\tilde 
U_-$ is nonzero, it contains a nonzero vector $\bold v\in U_+\cap
\tilde U_-$. From the conditions $\bold v\in U_+$ and $\bold v\in 
U_-$ we obtain two inequalities
$$
\xalignat 2
&g(\bold v)\geqslant 0,
&&g(\bold v)<0
\endxalignat
$$
contradicting each other. This contradiction shows that our
assumption $r_p\neq\tilde r_p$ is not valid and the inertia indices
$r_p$ and $\tilde r_p$ do coincide. From $r_p=\tilde r_p$ and $s
=\tilde s$ then we derive $r_n=\tilde r_n$. The theorem is proved.
\qed\enddemo
\definition{Definition 3.4} The total set of inertia indices 
is called the {\it signature\/} of a quadratic form. In the case of
a quadratic form in complex space ($\Bbb K=\Bbb C$) the signature
is formed by two numbers $(s,n-s)$, in the case of real space
($\Bbb K=\Bbb R$) it is formed by three numbers $(s,r_p,r_n)$.
\enddefinition
     In the case of a linear space over the field of rational numbers
$\Bbb K=\Bbb Q$ we can also diagonalize the matrix of a quadratic form
and subdivide the diagonal elements into three parts: positive, negative,
and zero elements. This determines the numbers $s$, $r_p$, and $r_n$,
which are geometric invariants of $g$, and we can define its signature.
\par
     However, in the case $\Bbb K=\Bbb Q$ we cannot reduce the nonzero 
diagonal elements to plus ones and minus ones only. Therefore, the 
number of geometric invariants in this case is greater than $3$. We shall
not look for the complete set of geometric invariants of a quadratic form
in the case $\Bbb K=\Bbb Q$ and we shall not construct their theory since
this would lead us to the number theory toward the problems of divisibility,
primality, factorization of integers, etc.\par
\head
\S\,4. Positive quadratic forms.\\Silvester's criterion.
\endhead
\rightheadtext{\S\,4. Silvester's criterion.}
     In this section we consider quadratic forms in linear vector spaces
over the field of real numbers $\Bbb R$. However, almost all results of
this section remain valid for quadratic forms in rational vector spaces
as well.
\definition{Definition 4.1} A quadratic form $g$ in a space $V$ over the
field $\Bbb R$ is called a {\it positive\/} form if $g(\bold v)>0$ for
any nonzero vector $\bold v\in V$.
\enddefinition
\proclaim{Theorem 4.1} A quadratic form $g$ in a finite-dimensional
space $V$ is positive if and only if the numbers $s$ and $r_n$ in its
signature $(s,r_p,r_n)$ are equal to zero.
\endproclaim
\demo{Proof} Let $g$ be a positive quadratic form and let $\bold e_1,\,
\ldots,\,\bold e_n$ be a basis in which the matrix of $g$ has the
canonic form \thetag{3.11}. If $s\neq 0$ then for the basis vector 
$\bold e_1\neq 0$ we would get $g(\bold e_1)=g_{11}=0$, which would
contradict the positivity of $g$. If $r_n\neq 0$, then for the basis 
vector $\bold e_n\neq 0$ we would get $g(\bold e_n)=g_{nn}=-1$, which 
would also contradict the positivity of $g$. Hence, $s=r_n=0$.\par
     Now, conversely, let $s=r_n=0$. Then in the basis $\bold e_1,\,
\ldots,\,\bold e_n$, in which the matrix of $g$ has the form 
\thetag{3.11}, its value $g(\bold v)$ is the sum of squares
$$
g(v)=(v^1)^2+\ldots+(v^n)^2,
$$
where $v^1,\,\ldots,\,v^n$ are coordinates of a vector $\bold v$.
This formula follows from the formula \thetag{1.10}. For a nonzero
vector at least one of its coordinates is nonzero. Hence, $g(\bold v)
>0$. This proves the positivity of $g$ and thus completes the proof 
of the theorem in whole.
\qed\enddemo
     The condition $s=\Ker g$ obtained in the theorem~3.4 and the
condition $s=0$ mean that a positive form $g$ in a finite-dimensional 
space $V$ is non-degenerate: $\Ker g=\{\bold 0\}$. This fact is valid
for a form in an infinite-dimensional space as well.
\proclaim{Theorem 4.2} Any positive quadratic form $g$ is
non-degenerate.
\endproclaim
\demo{Proof} If $\Ker g\neq\{\bold 0\}$ then there is a nonzero vector
$\bold v\in\Ker g$. The vector $\bold v$ of the kernel is orthogonal
to all vectors of the space $V$. Hence, it is orthogonal to itself: 
$g(\bold v)=g(\bold v,\bold v)=0$. If so, this fact contradicts the 
positivity of the form $g$. Therefore, any positive form $g$ should
be non-degenerate. 
\qed\enddemo
\proclaim{Theorem 4.3} Any subspace $U\subset V$ is regular with respect
to a positive quadratic form $g$ in a linear vector space $V$.
\endproclaim
\demo{Proof} Since the kernel $\Ker g$ of a positive form $g$
is zero, the regularity of a subspace $U$ with respect to $g$ is
equivalent to the equality $U\cap U_{\sssize\perp}=\{\bold 0\}$ 
(see definition~3.1). Let's prove this equality. Let $\bold v$ be
an arbitrary vector of the intersection $U\cap U_{\sssize\perp}$. 
From $\bold v\in U_{\sssize\perp}$ we derive that it is orthogonal
to all vectors of $U$. Hence, it is also orthogonal to itself since
$\bold v\in U$. Therefore, $g(\bold v)=g(\bold v,\bold v)=0$. Due
to positivity of $g$ the equality $g(\bold v)=0$ holds only for the
zero vector $\bold v=\bold 0$. Thus, we get $U\cap U_{\sssize\perp}
=\{\bold 0\}$. The theorem is proved.
\qed\enddemo
\proclaim{Theorem 4.4} For any subspace $U\subset V$ and for any
positive quadratic form $g$ in a finite-dimensional space $V$ there
is an expansion $V=U\oplus U_{\sssize\perp}$.
\endproclaim
\demo{Proof} The expansion $V=U+U_{\sssize\perp}$ follows from the
theorem~3.1. We need only to prove that the sum in this expansion
is a direct sum. For the sum of the dimensions of $U$ and $U_{\sssize
\perp}$ from the theorem~2.5 due to the triviality of the kernel
$\Ker g=\{\bold 0\}$ of a positive quadratic form $g$ we derive
$$
\dim U+\dim U_{\sssize\perp}=\dim V.
$$
Due to this equality in order to complete the proof it is sufficient
to apply the theorem~6.3 of
Chapter~\uppercase\expandafter{\romannumeral 1}. 
\qed\enddemo
     Let $g$ be a quadratic form in a finite-dimensional space $V$
over the field of real numbers $\Bbb R$. Let's choose an arbitrary
basis $\bold e_1,\,\ldots,\,\bold e_n$ in $V$ and then let's construct
the matrix of the quadratic form $g$:
$$
\hskip -2em
\Cal G=
\Vmatrix g_{11} & \hdots & g_{1n}\\
           \vdots & \ddots & \vdots\\
           g_{n1} & \hdots & g_{nn}
\endVmatrix
\tag4.1
$$
Let's delete the last $n-k$ columns and the last $n-k$ raws in the above matrix \thetag{4.1}. \pagebreak The determinant of the matrix thus obtained is called the {\it $k$-th  principal minor} of the matrix $\Cal G$. We denote this determinant by $M_k$:
$$
\hskip -2em
M_k=\det
\vmatrix g_{11} & \hdots & g_{1k}\\
           \vdots & \ddots & \vdots\\
           g_{k1} & \hdots & g_{kk}
\endvmatrix
\tag4.2
$$
The $n$-th principal minor $M_n$ coincides with the determinant of the
matrix $\Cal G$.
\proclaim{Theorem 4.5} Let $g$ be a positive quadratic form in a
finite-dimensional space $V$. Then the determinant of the matrix 
of $g$ in an arbitrary basis of $V$ is positive.
\endproclaim
\demo{Proof} For the beginning we consider a canonic basis $\bold e_1,
\,\ldots,\,\bold e_n$ in which the matrix of $g$ has the canonic form
\thetag{3.11}. According to the theorem~4.1, the matrix of a positive quadratic form $g$ in a canonic basis is the unit matrix. Hence, its
determinant is equal to unity and thus it is positive: 
$\det\Cal G=1>0$.\par
     Now let $\tilde\bold e_1,\,\ldots,\,\tilde\bold e_n$ be an arbitrary
basis and let $S$ be the transition matrix for passing from $\bold e_1,\,
\ldots,\,\bold e_n$ to $\tilde\bold e_1,\,\ldots,\,\tilde\bold e_n$.
Applying the formula \thetag{1.12}, we get
$$
\det\tilde{\Cal G}=
\det S^{\tr}\,(\det\Cal G)\,\det S=(\det S)^2.
$$
In a linear vector space $V$ over the real numbers $\Bbb R$ the elements
of any transition matrix $S$ are real numbers. Its determinant is also a
nonzero real number. Therefore, $(\det S)^2$ is a positive number.
The theorem is proved.
\qed\enddemo
     Now again let $\bold e_1,\,\ldots,\,\bold e_n$ be an arbitrary 
basis of $V$ and let $g_{ij}$ be the matrix of a positive quadratic
form $g$ in this basis. Let's consider the subspace
$$
U_k=\langle\bold e_1,\,\ldots,\,\bold e_k\rangle.
$$
Let's denote by $h_k$ the restriction of $g$ to the subspace $U_k$.
The matrix of the form $h_k$ in the basis $\bold e_1,\,\ldots,\,
\bold e_k$ coincides with upper left diagonal block in the matrix of
the initial form $g$. This is the very block that determines the
$k$-th principal minor $M_k$ in the formula \thetag{4.2}. It is clear
that the restriction of a positive form $g$ to any subspace is again 
a positive quadratic form. Therefore, we can apply the theorem~4.5 to
the form $h_k$. This yields $M_k>0$.\par
     {\bf Conclusion}: the positivity of all principal minors \thetag{4.2}
is a necessary condition for the positivity of a quadratic form $g$ itself.
As appears, this condition is a sufficient condition as well. This fact
is known as {\it Silvester's criterion}.
\proclaim{Theorem 4.6 (Silvester)} A quadratic form $g$ in a 
finite-dimensional space $V$ is positive if and only if all principal
minors of its matrix are positive.
\endproclaim
\demo{Proof} The positivity of $g$ implies the positivity of all principal
minors in its matrix. This fact is already proved. Let's prove the converse
proposition. Suppose that all diagonal minors \thetag{4.2} in the matrix
of a quadratic form $g$ are positive. We should prove that $g$ is positive.
The proof is by induction on $n=\dim V$.\par
     The basis of the induction in the case $\dim V=1$ is obvious. Here the
matrix of $g$ consists of the only element $g_{11}$ that coincides with the
only principal minor: $g_{11}=M_1$. The value $g(\bold v)$ in one-dimensional
space is determined by the only coordinate of a vector $\bold v$ according
to the formula $g(\bold v)=g_{11}\,(v^1)^2$. Therefore $M_1>0$ implies the
positivity of the form $g$.\par
     Suppose that the proposition we are going to prove is valid for a
quadratic form in any space of the dimension less than $n=\dim V$. Let
$g_{ij}$ be the matrix of our quadratic form $g$ in some basis 
$\bold e_1,\,\ldots,\,\bold e_n$ of $V$. Let's denote
$$
U=\langle\bold e_1,\,\ldots,\,\bold e_{n-1}\rangle.
$$
Denote by $h$ the restriction of the form $g$ to the subspace $U$ of the
dimension $n-1$. The matrix elements $h_{ij}$ in the matrix of $h$ 
calculated in the basis $\bold e_1,\,\ldots,\,\bold e_{n-1}$ coincide with 
corresponding elements in the matrix of the initial form: $h_{ij}=g_{ij}$.
Therefore, the minors $M_1,\ldots,M_{n-1}$ can be calculated by means of
the matrix $h_{ij}$. Due to the positivity of these minors, applying the
inductive hypothesis, we find that $h$ is a positive quadratic form in 
$U$.\par
     Let $\tilde\bold e_1,\,\ldots,\,\tilde\bold e_{n-1}$ be a basis 
in which the matrix of the form $h$ has the canonic form \thetag{3.11}. 
Applying the theorem~4.1 to the form $h$, we conclude that 
the matrix $\tilde h_{ij}$ in the canonic basis $\tilde\bold e_1,\,\ldots,
\,\tilde\bold e_{n-1}$ is the unit matrix. Let's complete the
basis $\tilde\bold e_1,\,\ldots,\,\tilde\bold e_{n-1}$ of the subspace
$U$ by the vector $\bold e_n\notin U$. As a result we get the basis
$\tilde\bold e_1,\,\ldots,\,\tilde\bold e_{n-1},\,\bold e_n$ in which
the matrix of $g$ has the form
$$
\hskip -2em
\Cal G_1=
\Vmatrix
1      & \hdots & 0      & \tilde g_{1n}     \\
\vdots & \ddots & \vdots & \vdots            \\
0      & \hdots & 1      & \tilde g_{n-1\,n} \\
\vspace{1ex}
\tilde g_{n1} & \hdots & \tilde g_{n\,n-1} & g_{nn}
\endVmatrix.
\tag4.3
$$
The passage from the basis $\bold e_1,\,\ldots,\,\bold e_n$ to the basis
$\tilde\bold e_1,\,\ldots,\,\tilde\bold e_{n-1},\,\bold e_n$ is described
by a blockwise-diagonal matrix $S$ of the form
$$
\hskip -2em
S_1=
\Vmatrix
S^1_1     & \hdots & S^1_n     & 0      \\
\vdots    & \ddots & \vdots    & \vdots \\
\vspace{1ex}
S^{n-1}_1 & \hdots & S^{n-1}_n & 0      \\
\vspace{1ex}
0         & \hdots & 0         & 1
\endVmatrix.
\tag4.4
$$
The formula \thetag{1.12} relates the matrix \thetag{4.3} with the
matrix $\Cal G$ of the quadratic form $g$ in the initial basis:
$\Cal G_1=S^{\tr}\,\Cal G\,S$. From this formula we derive
$$
\hskip -2em
\det\Cal G_1=\det\Cal G\,(\det S)^2=M_n\,(\det S)^2.
\tag4.5
$$
Due to the above formula \thetag{4.5} the positivity of the principal 
minor $M_n=\det\Cal G$ in the initial matrix \thetag{4.1} implies the positivity of the determinant of the matrix \thetag{4.3}, i\.\,e\. \
$\det\Cal G_1>0$.\par
     Let's calculate the determinant of the matrix \thetag{4.3}
explicitly. For this purpose we multiply the first column of this
matrix by $\tilde g_{1n}$ and subtract it from the last column. 
Then we multiply the second column by $\tilde g_{2n}$ and subtract
it from the last one. We produce such an operation repeatedly for each
of the first $n-1$ columns of the matrix \thetag{4.3}. From the course
of algebra we know that such transformations do not change the 
determinant of a matrix. In present case they simplify the matrix 
\thetag{4.3} bringing it to a lower-triangular form. \pagebreak 
Therefore, we can calculate the determinant of the matrix \thetag{4.3} 
in explicit form:
$$
\hskip -2em
\det\Cal G_1=\det
\vmatrix
1      & \hdots & 0      & 0      \\
\vdots & \ddots & \vdots & \vdots \\
0      & \hdots & 1      & 0      \\
\vspace{1ex}
\tilde g_{n1} & \hdots & \tilde g_{n\,n-1} &\tilde g_{nn}
\endvmatrix
=\tilde g_{nn}.
\tag4.6
$$
The element $\tilde g_{nn}$ in the transformed matrix is given by the
formula
$$
\hskip -2em
\tilde g_{nn}=g_{nn}-\sum^{n-1}_{i=1} g_{ni}\,g_{in}=
g_{nn}-\sum^{n-1}_{i=1} (g_{in})^2.
\tag4.7
$$\par
     The matrix of the quadratic form $g$ in the basis $\tilde\bold e_1,
\,\ldots,\,\tilde\bold e_{n-1},\bold e_n$ is close to the diagonal matrix.
Let's complete the process of diagonalization replacing the vector
$\bold e_n$ by the vector $\tilde\bold e_n\not\in U$ such that
$$
\tilde\bold e_n=\bold e_n-\sum^{n-1}_{i=1} g_{in}\cdot\tilde\bold e_i.
$$
The passage from $\tilde\bold e_1,\,\ldots,\,\tilde\bold e_{n-1},\bold e_n$ 
to $\tilde\bold e_1,\,\ldots,\,\tilde\bold e_n$ changes only the last basis
vector. Therefore, the unit diagonal block in the matrix \thetag{4.3}
remains unchanged. For non-diagonal elements $g(\tilde\bold e_k,\tilde\bold
e_n)$ in the new basis we have 
$$
g(\tilde\bold e_k,\tilde\bold e_n)=\tilde g_{kn}-
\sum^{n-1}_{i=1}g_{in}\,g(\tilde\bold e_k,\tilde\bold e_i)=
\tilde g_{kn}-\sum^{n-1}_{i=1}g_{in}\,\tilde h_{ki}=0.
$$
The equality $g(\tilde\bold e_k,\tilde\bold e_n)=0$ in the above formula is due to the fact that the matrix of the restricted form $h$ in its canonic
basis $\tilde\bold e_1,\,\ldots,\,\tilde\bold e_{n-1}$ is the unit matrix.
For the diagonal element $g(\tilde\bold e_n,\tilde\bold e_n)$ from this
fact we derive 
$$
g(\tilde\bold e_n,\tilde\bold e_n)=g_{nn}-\sum^{n-1}_{i=1}
\sum^{n-1}_{k=1} g_{in}\,g_{kn}\,h_{ik}=
g_{nn}-\sum^{n-1}_{i=1} (g_{in})^2.
$$
Comparing this expression with \thetag{4.7}, we find that
$g(\tilde\bold e_n,\tilde\bold e_n)=\tilde g_{nn}$. Thus, the matrix 
of $g$ in the basis $\tilde\bold e_1,\,\ldots,\,\tilde\bold e_n$
is a diagonal matrix of the form
$$
\hskip -2em
\Cal G_2=
\Vmatrix
1      & \hdots & 0      & 0      \\
\vdots & \ddots & \vdots & \vdots \\
0      & \hdots & 1      & 0      \\
\vspace{1ex}
0      & \hdots & 0      &\tilde g_{nn}
\endVmatrix.
\tag4.8
$$
Combining \thetag{4.5} and \thetag{4.6}, for the element $\tilde 
g_{nn}$ in \thetag{4.8} we get $\tilde g_{nn}=M_n\,(\det S)^2$. Since 
the principal minor $M_n$ of the initial matrix \thetag{4.1} is positive,
we find that $\tilde g_{nn}$ in \thetag{4.8} is also positive. Hence,
$g$ is a positive quadratic form. Thus, we have completed the inductive 
step and have proved the theorem in whole.
\qed\enddemo
\newpage
\topmatter
\title\chapter{5}
Euclidean spaces.
\endtitle
\endtopmatter
\document
\head
\S\,1. The norm and the scalar product. The angle\\between vectors.
Orthonormal bases.
\endhead
\leftheadtext{CHAPTER~\uppercase\expandafter{\romannumeral 5}.
EUCLIDEAN SPACES.}
\rightheadtext{\S\,1. The norm and the scalar product.}
\setfirstpage
\definition{Definition 1.1} A {\it Euclidean vector space\/} is a
linear vector space $V$ over the field of reals $\Bbb R$ which is
equipped with some fixed positive quadratic form $g$.
\enddefinition
     Let $(V,g)$ be a Euclidean vector space. There many positive
quadratic forms in the linear vector space $V$, however, only one of 
them is associated with $V$ so that it defines the structure of 
Euclidean space in $V$. Two Euclidean vector spaces $(V,g_1)$ and
$(V,g_2)$ with $g_1\neq g_2$ coincide as linear vector spaces, but 
they are different when considered as Euclidean vector spaces.\par
     The structure of the Euclidean vector space $(V,g)$ is associated
with a special terminology and special notations. The value of the 
quadratic form $g(\bold v)$ is non-negative. The square root of 
$g(\bold v)$ is called the {\it norm\/} or the {\it length\/} of a 
vector $\bold v$. The norm of a vector $\bold v$ is denoted as follows:
$$
\hskip -2em
|\bold v|=\sqrt{g(\bold v)}.
\tag1.1
$$
The quadratic form $g(\bold v)$ produces the bilinear form $g(\bold v,
\bold w)$ determined by the recovery formula \thetag{1.6} of
Chapter~\uppercase\expandafter{\romannumeral 4}. The value of that
bilinear form is called the {\it scalar product\/} of two vectors
$\bold v$ and $\bold w$. The scalar product is denoted as follows:
$$
\hskip -2em
(\bold v\,|\,\bold w)=g(\bold v,\bold w).
\tag1.2
$$
Due to the notation \thetag{1.1} and \thetag{1.2}, when dealing with
some fixed Euclidean space $(V,g)$, we can omit the symbol $g$ 
at all.\par
     The scalar product \thetag{1.2} is defined for a pair of two vectors
$\bold v,\bold w\in V$. It is quite different from the scalar product
\thetag{1.8} of 
Chapter~\uppercase\expandafter{\romannumeral 3}, which is defined for
a pair of a vector and a covector. The scalar product \thetag{1.2}
of a Euclidean vector space possesses the following properties:
\roster
\item \quad$(\bold v_1+\bold v_2\,|\,\bold w)=(\bold v_1\,|\,
      \bold w)+(\bold v_2\,|\,\bold w)$ \ for all \ $\bold v_1,
      \bold v_2,\bold w\in V$;
\item \quad$(\alpha\cdot\bold v\,|\,\bold w)=\alpha\,(\bold v\,
      |\,\bold w)$ \ for all \ $\bold v,\bold w\in V$ \ and for
      all \ $\alpha\in\Bbb R$;
\item \quad$(\bold v\,|\,\bold w_1+\bold w_2)=(\bold v\,|\,
      \bold w_1)+(\bold v\,|\,\bold w_2)$ \ for all \ $\bold w_1,
      \bold w_2,\bold v\in V$;
\item \quad$(\bold v\,|\,\alpha\cdot\bold w)=\alpha\,(\bold v\,
      |\,\bold w)$ \ for all \ $\bold v,\bold w\in V$ \ and for
      all \ $\alpha\in\Bbb R$;
\item \quad$(\bold v\,|\,\bold w)=(\bold w\,|\,\bold v)$
       \ for all \ $\bold v,\bold w\in V$;
\item \quad$|\bold v|^2=(\bold v\,|\,\bold v)\geqslant 0$ \ for all
      $\bold v\in V$ \ and \ $|\bold v|=0$ implies $\bold v=\bold 0$.
\endroster
The properties \therosteritem{1}-\therosteritem{4} reflect the
bilinearity of the form $g$ in \thetag{1.2}. They are analogous
to that of the scalar product of a vector and a covector (see
formulas \thetag{1.9} in \pagebreak 
Chapter~\uppercase\expandafter{\romannumeral 3}).\par
     The properties \therosteritem{5} and \therosteritem{6} have no
such analogs. But they are the very properties that make the scalar
product \thetag{1.2} a generalization of the scalar product of
$3$-dimensional geometric vectors.
\proclaim{Theorem 1.1} The following two additional properties of the
scalar product \thetag{1.2} are derived from the properties
\therosteritem{1}-\therosteritem{6}:
\roster
\item [7] \quad$|(\bold v,\bold w)|\leqslant |\bold v|\,|\bold w|$
      \ for all \ $\bold v,\bold w\in V$;
\item     \quad$|\bold v+\bold w|\leqslant |\bold v|+|\bold w|$ 
      \ for all \ $\bold v,\bold w\in V$.
\endroster
The property \therosteritem{7} is known as the {\it 
Cauchy-Bunyakovsky-Schwarz inequality}, while the property 
\therosteritem{8} is called the {\it triangle inequality}.
\endproclaim
\demo{Proof} In order to prove the inequality \therosteritem{7}
we choose two arbitrary nonzero vectors $\bold v,\bold w\in V$ 
and consider the numeric function $f(\alpha)$ of a numeric 
argument $\alpha$ defined by the following explicit formula:
$$
\hskip -2em
f(\alpha)=|\bold v+\alpha\cdot\bold w|^2.
\tag1.3
$$
Using the properties \therosteritem{1}-\therosteritem{6} we find 
that $f(\alpha)$ is a polynomial of degree two:
$$
\align
f(\alpha)=|\bold v+\alpha\cdot\bold w|^2&=(\bold v+\alpha\cdot\bold w
\,|\,\bold v+\alpha\cdot\bold w)=\\
&=(\bold v\,|\,\bold v)+2\,\alpha\,(\bold v\,|\,\bold w)\,
+\alpha^2\,(\bold w\,|\,\bold w).
\endalign
$$
The function \thetag{1.3} has a lower bound: $f(\alpha)\geqslant 0$. 
This follows from the property \therosteritem{6}. Let's calculate the minimum point of the function $f(\alpha)$ by equating its derivative
$f'(\alpha)$ to zero. This yields the following equation:
$$
f'(\alpha)=2\,(\bold v\,|\,\bold w)\,+2\,\alpha\,(\bold w\,|\,\bold w)=0.
$$
Solving this equation, we find $\alpha_{\min}=-(\bold v\,|\,\bold
w)/(\bold w\,|\,\bold w)$. Now let's write the condition $f(\alpha)
\geqslant 0$ for the minimal value of the function $f(\alpha)$:
$$
\hskip -2em
f_{\min}=f(\alpha_{\min})=\frac{|\bold v|^2\,|\bold w|^2-(\bold v\,
|\,\bold w)^2}{|\bold w|^2}
\geqslant 0.
\tag1.4
$$
The denominator of the fraction \thetag{1.4} is positive, therefore,
from the inequality \thetag{1.4} we easily derive the property 
\therosteritem{7}.\par
     In order to prove the property \therosteritem{8} we consider
the square of the norm for the vector $\bold v+\bold w$. For this
quantity we derive
$$
\hskip -2em
|\bold v+\bold w|^2=(\bold v+\bold w\,|\,\bold v+\bold w)
=|\bold v|^2+2\,(\bold v\,|\,\bold w)+|\bold w|^2.
\tag1.5
$$
Applying the property \therosteritem{8}, which is already proved, for
the right hand side of the equality \thetag{1.5} we get the following
estimate:
$$
|\bold v|^2+2\,(\bold v\,|\,\bold w)+|\bold w|^2\leqslant
|\bold v|^2+2\,|\bold v|\,|\bold w|+|\bold w|^2=(|\bold v|
+|\bold w|)^2.
$$
From the relationship \thetag{1.5} and from the above inequality
we derive the other inequality $|\bold v+\bold w|^2\leq(|\bold v|
+|\bold w|)^2$. Now the property \therosteritem{8} is derived 
by taking the square root of both sides of this inequality. This
operation is correct since $y=\sqrt{x}$ is an increasing function 
of the real semiaxis $[0,\ +\infty)$. \pagebreak
The theorem is proved.
\qed\enddemo\par
     Due to the analogy of \thetag{1.2} and the scalar product of 
geometric vectors and due to the Cauchy-Bunyakovsky-Schwarz
inequality $|(\bold v,\bold w)|\leqslant |\bold v|\,|\bold w|$ 
we can introduce the concept of an {\it angle between vectors\/} 
in a Euclidean vector space.
\definition{Definition 1.2} The number $\varphi$ from the interval $0\leqslant\varphi\leqslant\pi$, which is determined by the following 
implicit formula
$$
\hskip -2em
\cos(\varphi)=\frac{(\bold v\,|\,\bold w)}{|\bold v|\,|\bold w|},
\tag1.6
$$
is called the {\it angle\/} between two nonzero vectors $\bold v$ and 
$\bold w$ in a Euclidean space $V$.
\enddefinition
     Due to the property \therosteritem{7} from the theorem~1.1 the
modulus of the fraction in left hand side of \thetag{1.6} is not greater
than $1$. Therefore, the formula \thetag{1.6} is correct. It determines 
the unique number $\varphi$ from the specified interval $0\leqslant\varphi
\leqslant \pi$.
\definition{Definition 1.3} Two vectors $\bold v$ and $\bold w$ in a
Euclidean space $V$ are called {\it orthogonal vectors\/} if they form a right
angle ($\varphi=\pi/2$).
\enddefinition
     The definition~1.3 applies only to nonzero vectors $\bold v$ and
$\bold w$. The definition ~2.1 of
Chapter~\uppercase\expandafter{\romannumeral 4} is more general. Let's
reformulate for the case of Euclidean spaces.
\definition{Definition1.4} Two vectors $\bold v$ and $\bold w$ in
a Euclidean space $V$ are called {\it orthogonal vectors\/} if their
scalar product is zero: $(\bold v\,|\,\bold w)=0$.
\enddefinition
     For nonzero vectors $\bold v$ and $\bold w$ these two definitions
1.3 and 1.4 are equivalent.\par
     Let $\bold v_1,\,\ldots,\,\bold v_m$ be a system of vectors 
in a Euclidean space $(V,g)$. The matrix $g_{ij}$ composed by the 
mutual scalar products of these vectors 
$$
\hskip -2em
g_{ij}=(\bold v_i\,|\,\bold v_j),
\tag1.7
$$
is called the {\it Gram matrix} of the system of vectors $\bold
v_1,\,\ldots,\,\bold v_m$.
\proclaim{Theorem 1.2} A system of vectors $\bold v_1,\,\ldots,\,\bold v_m$
in a Euclidean space is linearly dependent if and only if the determinant of
their Gram matrix is equal to zero.
\endproclaim
\demo{Proof} Suppose that the vectors $\bold v_1,\,\ldots,\,\bold v_m$
are linearly dependent. Then there is a nontrivial linear combination
of these vectors which is equal to zero:
$$
\hskip -2em
\alpha_1\cdot\bold v_1+\ldots+\alpha_m\cdot\bold v_m=\bold 0.
\tag1.8
$$
Using the coefficients of the linear combination \thetag{1.8}, we construct
the following expression with the components of Gram matrix \thetag{1.7}:
$$
\sum^m_{j=1}g_{ij}\,\alpha_j=\sum^m_{j=1}(\bold v_i\,|\,\bold v_j)\,
\alpha_j=(\bold v_i\,|\,\alpha_1\cdot\bold v_1+\ldots+\alpha_m\cdot
\bold v_m)=(\bold v_i\,|\,\bold 0)=0.
$$
Since $i$ is a free index running over the interval of integer numbers
from $1$ to $m$, this formula means that the columns of Gram matrix 
$g_{ij}$ are linearly dependent. Hence, its determinant is equal to zero
(this fact is known from the course of general algebra).\par
     Conversely, assume that the determinant of the Gram matrix 
\thetag{1.7} is equal to zero. Then the columns of this matrix are 
linearly dependent and, hence, there is a nontrivial linear combination
of them that is equal to zero:
$$
\hskip -2em
\sum^m_{j=1}g_{ij}\,\alpha_j=0.
\tag1.9
$$
Let's denote $\bold v=\alpha_1\cdot\bold v_1+\ldots+\alpha_m\cdot\bold
v_m$. Then consider the following double sum, which is obviously equal 
to zero due to the equality \thetag{1.9}:
$$
\align
0=\sum^m_{i=1}\sum^m_{j=1}\alpha_i\,g_{ij}\,&\alpha_j=
\sum^m_{i=1}\alpha_i\,(\bold v_i\,|\,\alpha_1\cdot\bold v_1+
\ldots+\alpha_m\cdot\bold v_m)=\\
&=(\alpha_1\cdot\bold v_1+\ldots+\alpha_m\cdot\bold v_m\,|
\,\bold v)=(\bold v\,|\,\bold v)=|\bold v|^2.
\endalign
$$
Thus, we get $|\bold v|^2=0$ and, using the positivity of the basic
quadratic form $g$ of the Euclidean space $V$, we derive $\bold v=
\bold 0$. Since $\bold v=\bold 0$, we get the nontrivial linear
combination of the form \thetag{1.8}, which is equal to zero. Hence,
our vector $\bold v_1,\,\ldots,\,\bold v_m$ are linearly dependent.
\qed\enddemo
     Let $\bold e_1,\,\ldots,\,\bold e_n$ be a basis in a 
finite-dimensional Euclidean vector space $(V,g)$. Let's consider
the Gram matrix of this basis. Knowing the components of the Gram
matrix, we can calculate the norm of vectors \thetag{1.1} and the
scalar product of vectors \thetag{1.2} through their coordinates:
$$
\xalignat 2
&\hskip -2em
\qquad |\bold v|^2=\sum^n_{i=1}\sum^n_{j=1} g_{ij}\,v^i\,v^j,
&&(\bold v\,|\,\bold w)=\sum^n_{i=1}\sum^n_{j=1} g_{ij}\,v^i\,w^j.
\tag1.10
\endxalignat
$$
A basis $\bold e_1,\,\ldots,\,\bold e_n$ in a Euclidean space $V$ 
is called an {\it orthonormal basis\/} if the Gram matrix for the basis vectors is the unit matrix:
$$
\hskip -2em
g_{ij}=
\cases
1&\text{\ \ for \ }i=j,\\
0&\text{\ \ for \ }i\neq j.
\endcases
\tag1.11
$$
If the condition \thetag{1.11} is not fulfilled, then the basis 
$\bold e_1,\,\ldots,\,\bold e_n$ is called a {\it skew-angular
basis}. In an orthonormal basis the vectors $\bold e_1,\,\ldots,\,
\bold e_n$ are unit vectors orthogonal to each other. This simplifies
the formulas \thetag{1.10} substantially:
$$
\xalignat 2
&\hskip -2em
|\bold v|^2=\sum^n_{i=1}(v^i)^2,
&&(\bold v\,|\,\bold w)=\sum^n_{i=1}v^i\,w^i.
\tag1.12
\endxalignat
$$\par
     Orthonormal bases do exist. Due to \thetag{1.2} and \thetag{1.7} 
we know that the Gram matrix of the basis vectors $\bold e_1,\,\ldots,
\,\bold e_n$ is the matrix of the quadratic $g$ in this basis.
The theorem~3.3 of
Chapter~\uppercase\expandafter{\romannumeral 4} says that there exists
a basis in which the matrix of $g$ has its canonic form
(see \thetag{3.11} in Chapter~\uppercase\expandafter{\romannumeral 4}).
Since $g$ is a positive quadratic form, its matrix in a canonic form is 
the unit matrix (see theorem~4.1 in 
Chapter~\uppercase\expandafter{\romannumeral 4}).\par
     The theorem~4.8 on completing the basis of a subspace formulated
in Chapter~\uppercase\expandafter{\romannumeral 1} has its analog for
orthonormal bases.
\proclaim{Theorem 1.3} Let $\bold e_1,\,\ldots,\,\bold e_s$ be an
orthonormal basis in a subspace $U$ of a finite-dimensional Euclidean 
space $(V,g)$. Then it can be completed up to an orthonormal basis 
$\bold e_1,\,\ldots,\,\bold e_n$ in $V$.
\endproclaim
\demo{Proof} Let's consider the orthogonal complement $U_{\sssize\perp}$
of the subspace $U$. According to the theorem~4.4 of
Chapter~\uppercase\expandafter{\romannumeral 4}, the subspaces $U$ and 
$U_{\sssize\perp}$ define the expansion of the space $V$ into a direct 
sum:
$$
V=U\oplus U_{\sssize\perp}.
$$
The subspace $U_{\sssize\perp}$ inherits the structure of a Euclidean
space from $V$. Let's choose an orthonormal basis 
$\bold e_{s+1},\,\ldots,\,\bold e_n$ in $U_{\sssize\perp}$ and then
join together two bases of $U$ and $U_{\sssize\perp}$. As a result we
get the basis in $V$ (see theorem~6.3 of 
Chapter~\uppercase\expandafter{\romannumeral 1}). The vectors of
this basis are unit vectors by their length and they are orthogonal
to each other. Hence, this is an orthonormal basis completing the
initial basis $\bold e_1,\,\ldots,\,\bold e_s$ of the subspace $U$.
\qed\enddemo
     Let $\bold e_1,\,\ldots,\,\bold e_s$ and $\tilde\bold  e_1,\,\ldots,
\tilde\bold e_s$ be two orthonormal bases and let $S$ be the transition.
The Gram matrices of these two bases are the unit matrices. Therefore,
applying the formulas \thetag{1.12} of 
Chapter~\uppercase\expandafter{\romannumeral 4}, for the transition matrix
$S$ we derive
$$
\xalignat 2
&\hskip -2em
S^{\tr}\,S=1,
&&S^{-1}=S^{\tr}.
\tag1.13
\endxalignat
$$
Note that a square matrix $S$ satisfying the above
relationships \thetag{1.13} is called an {\it orthogonal\/} matrix.\par
     From the relationships \thetag{1.13} for the determinant of an
orthogonal matrix we get: $(\det S)^2=1$. Therefore, orthogonal matrices
are subdivided into two types: matrices with positive determinant
$\det S=1$ and those with negative determinant $\det S=-1$. This subdivision
is related to the concept of {\it orientation}. All bases in a linear
vector space over the field of real numbers $\Bbb R$ (not necessarily a
Euclidean space) can be subdivided into two sets which can be called
{\tencyr\char '074}left bases{\tencyr\char '076} and 
{\tencyr\char '074}right bases{\tencyr\char '076}. The transition matrix
for passing from a left basis to a left basis or for passing from
a right basis to another right basis is a matrix with positive determinant
--- it does not change the orientation. The transition matrix for passing
from a left basis to a right basis or, conversely, from a right basis to a
left basis is a matrix with negative determinant. Such a transition matrix
changes the orientation of a basis. We say that a linear vector space $V$
over the field of real numbers $\Bbb R$ is equipped with the 
{\it orientation\/} if there is some mechanism to distinguish one of two
types of bases in $V$.\par
\head
\S\,2. Quadratic forms in a Euclidean space. Diagonalization of a 
pair of quadratic forms.
\endhead
\rightheadtext{\S\,2. Quadratic forms in a Euclidean space.}
     Let $(V,g)$ be a Euclidean vector space and let $\varphi$ be
a quadratic form in $V$. For such a form $\varphi$ we define the 
following ratio:
$$
\hskip -2em
\mu(\bold v)=\frac{|\varphi(\bold v)|}{|\bold v|^2}.
\tag2.1
$$
The number $\mu(v)$ in \thetag{2.1} is a real non-negative number.
Note that $\mu(\alpha\cdot\bold v)=\mu(\bold v)$ for any 
nonzero $\alpha\in\Bbb R$. \pagebreak Therefore, we can assume 
$\bold v$ in \thetag{2.1} to be a unit vector.\par
     Let's denote by $\Vert\varphi\Vert$ the least upper bound
of $\mu(v)$ for all unit vectors (such vectors sweep out the {\it
unit sphere\/} in the Euclidean space $V$):
$$
\hskip -2em
\Vert\varphi\Vert=\sup_{|\bold v|=1}\mu(\bold v).
\tag2.2
$$
\definition{Definition 2.1} The quantity $\Vert\varphi\Vert$ determined
by the formulas \thetag{2.1} and \thetag{2.2} is called the {\it norm\/}
of a quadratic form $\varphi$ in a Euclidean vector space $V$. If the
norm $\Vert\varphi\Vert$ is finite, the form $\varphi$ is said to be
a {\it restricted quadratic form}.
\enddefinition
\proclaim{Theorem 2.1} If $\varphi$ is a restricted quadratic form,
then there is the estimate $|\varphi(\bold v,\bold w)|\leqslant
\Vert\varphi\Vert\,|\bold v|\,|\bold w|$ for the values of corresponding
symmetric bilinear form.
\endproclaim
\demo{Proof} In order to calculate $\varphi(\bold v,\bold w)$ we use
the following equality, which, in essential, is a version of the recovery
formula:
$$
\hskip -2em
4\,\alpha\,\varphi(\bold v,\bold w)=\varphi(\bold v+\alpha\cdot\bold w)-
\varphi(\bold v-\alpha\cdot\bold w).
\tag2.3
$$
From \thetag{2.3} we derive the following inequality for the quantity
$4\,\alpha\,\varphi(\bold v,\bold w)$:
$$
\hskip -2em
4\,\alpha\,\varphi(\bold v,\bold w)\leqslant |\varphi(\bold v
+\alpha\cdot\bold w)|+|\varphi(\bold v-\alpha\cdot\bold w)|.
\tag2.4
$$
Now let's apply the inequality 
$|\varphi(\bold u)|\leq\Vert\varphi\Vert\,|\bold u|^2$ derived from
\thetag{2.1} and \thetag{2.2} in order to estimate the right hand
side of \thetag{2.4}. This yields
$$
\hskip -2em
4\,\alpha\,\varphi(\bold v,\bold w)\leqslant \Vert\varphi\Vert\,(|\bold
v+\alpha\cdot\bold w|^2+|(\bold v-\alpha\cdot\bold w)|^2).
\tag2.5
$$
Let's express the squares of moduli through the scalar products:
$$
|\bold v\pm\alpha\cdot\bold w|^2=|\bold v|^2\pm 2\,\alpha\,(\bold v
\,|\,\bold w)+\alpha^2\,|\bold w|^2.
$$
Then we can simplify the inequality \thetag{2.5} bringing it to the
following one:
$$
4\,\alpha\,\varphi(\bold v,\bold w)\leqslant
2\Vert\varphi\Vert\,(|\bold v|^2+\alpha^2\,|\bold w|^2).
$$
Now let's transform the above inequality a little bit more:
$$
f(\alpha)=\alpha^2\,\Vert\varphi\Vert\,|\bold w|^2-2\alpha\,
\varphi(\bold v,\bold w)+\Vert\varphi\Vert\,|\bold v|^2
\geqslant 0.
$$\par
     The numeric function $f(\alpha)$ of a numeric argument $\alpha$
is a polynomial of degree two in $\alpha$. Let's find the minimum
point $\alpha=\alpha_{\min}$ for this function by equating its
derivative to zero: $f'(\alpha)=0$. As a result we obtain
$$
\alpha_{\min}=\frac{\varphi(\bold v,\bold w)}{\Vert\varphi\Vert
\,|\bold w|^2}.
$$
Now let's write the inequality $f(\alpha_{\min})\geqslant 0$ for the
minimal value of this function. This yields the following inequality
for the bilinear form $\varphi$:
$$
\pagebreak
\varphi(\bold v,\bold w)^2\leqslant \Vert\varphi\Vert^2\,|\bold v|^2
\,|\bold w|^2.
$$
Now it is easy to derive the required estimate for $|\varphi(\bold v,
\bold w)|$ by taking the square root of both sides of the above 
inequality. Note that a quite similar method was used when proving 
the Cauchy-Bunyakovsky-Schwarz inequality in the theorem~1.1.
\qed\enddemo
\proclaim{Theorem 2.2} Any quadratic form $\varphi$ in a 
finite-dimensional Euclidean vector space $V$ is a restricted form.
\endproclaim
\demo{Proof} Let's choose an orthonormal basis $\bold e_1,\,\ldots,
\,\bold e_n$ in $V$ and consider the expansion of a unit vector 
$\bold v$ in this basis. For the coordinates of $\bold v$ in this
basis due to the formulas \thetag{1.12} we obtain
$$
(v^1)^2+\ldots+(v^n)^2=1.
$$
Hence, for the components of $\bold v$ we have $|v^i|\leqslant 1$. 
Let's express the quantity $\mu(\bold v)$, which is defined by formula 
\thetag{2.1}, through the coordinates of $\bold v$:
$$
\mu(v)=|\varphi(v)|=\left|\,\sum^n_{i=1}\sum^n_{j=1}
\varphi_{ij}\,v^i\,v^j\right|.
$$
From $|v^i|\leqslant 1$ we derive the following estimate for
the quantity $\mu(\bold v)$:
$$
\hskip -2em
\mu(v)\leqslant \sum^n_{i=1}\sum^n_{j=1} |\varphi_{ij}|
<\infty,
\tag2.6
$$
Right hand site of \thetag{2.6} does not depend $\bold v$. Due to
\thetag{2.2} this sum is an upper bound for the norm $\Vert\varphi\Vert$.
Hence, $\Vert\varphi\Vert<\infty$. The theorem is proved.
\qed\enddemo
\proclaim{Theorem 2.3} For any quadratic form $\varphi$ in a
finite-dimensional Euclidean vector space $V$ the supremum in 
formula \thetag{2.2} is reached, i\.\,e\. there exists a vector
$\bold v\neq\bold 0$ such that $|\varphi(\bold v)|=\Vert\varphi\Vert
\,|\bold v|^2$.
\endproclaim
\demo{Proof} From the course of mathematical analysis we knows that
the supremum of a numeric set is the limit of some converging sequence
of numbers of this set (see~\cite{6}). This means that there is 
a sequence of unit vectors $\bold v(1),\,\ldots,\,\bold v(n),\,\ldots$ 
in $V$ such that the norm $\Vert\varphi\Vert$ is expressed as the 
following limit:
$$
\hskip -2em
\Vert\varphi\Vert=\lim_{s\to\infty} |\varphi(\bold v(s))|.
\tag2.7
$$
Let's choose an orthonormal basis $\bold e_1,\,\ldots,\,\bold e_n$
in $\bold v$ and let's expand each vector $\bold v(s)$ of the sequence
in this basis. The equality 
$$
\hskip -2em
(v^1(s))^2+\ldots+(v^n(s))^2=1
\tag2.8
$$
is derived from $|\bold v(s)|=1$ due to the formulas \thetag{1.12}.
Now the equality \thetag{2.8} means that each specific coordinate
$v^i(s)$ yields a restricted sequence of real numbers:
$$
-1\leqslant v^i(s)\leqslant 1.
$$
From the course of mathematical analysis we know that in each restricted
sequence of real numbers one can choose a converging subsequence. So, in
the sequence of unit vectors $\bold v(s)$ one can choose a subsequence
of unit vectors whose first coordinates form a convergent sequence of
numbers. Let's denote this subsequence again by $\bold v(s)$ and choose
its subsequence with converging second coordinates. Repeating this choice
$n$-times for each specific coordinate, we get a subsequence of unit 
vectors $\bold v(s_k)$ such that their coordinates all are the converging
sequences of numbers. Let's consider the limits of these sequences:
$$
\hskip -2em
v^i=\lim_{k\to\infty} v^i(s_k).
\tag2.9
$$
Denote by $\bold v$ the vector whose coordinates are determined by
the limit values \thetag{2.9}. Passing to the limit $s\to\infty$ in
\thetag{2.8}, we conclude that $\bold v$ is a unit vector: 
$|\bold v|=1$.\par
     Now let's calculate $|\varphi(\bold v)|$ using the matrix of 
the quadratic form $\varphi$ and the coordinates of $\bold v$ in the
basis $\bold e_1,\,\ldots,\,\bold e_n$:
$$
|\varphi(\bold v)|=\left|\sum^n_{i=1}\sum^n_{j=1} \varphi_{ij}\,
v^i\,v^j\right|=\lim_{k\to\infty}\left|\sum^n_{i=1}\sum^n_{j=1}
\varphi_{ij}\,v^i(s_k)\,v^j(s_k)\right|.
$$
On the other hand, taking into account \thetag{2.7}, for
$|\varphi(\bold v)|$ we get
$$
|\varphi(\bold v)|=\lim_{k\to\infty}|\varphi(v(s_k))|=
\lim_{s\to\infty}|\varphi(v(s))|=\Vert\varphi\Vert.
\tag2.10
$$
Thus, for the unit vector $\bold v$ with coordinates \thetag{2.9} 
we get $|\varphi(\bold v)|=\Vert\varphi\Vert$. Multiplying
$\bold v$ by some number $\alpha\in\Bbb R$, we can remove the
restriction $|\bold v|=1$. Then the equality \thetag{2.10} will
be written as $|\varphi(\bold v)|=\Vert\varphi\Vert\,|\bold v|^2$.
The theorem is proved.
\qed\enddemo
\proclaim{Theorem 2.4} For any quadratic form $\varphi$ in a
finite-dimensional Euclidean vector space $(V,g)$ there is an
orthonormal basis $\bold e_1,\,\ldots,\,\bold e_n$ such that
the matrix of the form $\varphi$ in this basis is a diagonal
matrix.
\endproclaim
\demo{Proof} The proof is by induction on the dimension of the space
$V$. In the case $\dim V=1$ the proposition of the theorem is obvious:
any square $1\times 1$ matrix is a diagonal matrix.\par
     Suppose that the proposition of the theorem is valid for all 
quadratic forms in Euclidean spaces of the dimension less than $n$.
Let $\dim V=n$ and let $\varphi$ be a quadratic form in the Euclidean
space $(V,g)$. Applying theorems~2.2 and 2.3, we find a unit vector
$\bold v\in V$ such that $|\varphi(\bold v)|=\Vert\varphi\Vert$. 
For the sake of certainty we assume that $\varphi(v)\geqslant 0$.
Then we can remove the modulus sign: $\varphi(v)=\Vert\varphi\Vert$. 
In the case $\varphi(\bold v)<0$ we replace the form $\varphi$ by the opposite form $\tilde\varphi=-\varphi$ since two opposite forms 
diagonalize simultaneously.\par
     Let's denote $U=\langle\bold v\rangle$ and consider the orthogonal
complement $U_{\sssize\perp}$. The subspaces $U=\langle v\rangle$ and
$U_{\sssize\perp}$ have zero intersection, their sum is a direct 
sum and $U\oplus U_{\sssize\perp}=V$ (see theorem~4.4 in 
Chapter~\uppercase\expandafter{\romannumeral 4}). Let's take an arbitrary
vector $\bold w\in U_{\sssize\perp}$ of the unit length and compose the
vector $\bold u$ as follows:
$$
\bold u=\cos(\alpha)\cdot\bold v+\sin(\alpha)\cdot\bold w.
$$
Here $\alpha$ is a numeric parameter. It is easy to see that $\bold u$ 
is also a unit vector, this follows from the identity 
$\cos^2(\alpha)+\sin^2(\alpha)=1$.\par
     Let's calculate the value of the quadratic form $\varphi$ on the 
vector $\bold u$ and treat it as a function of the numeric parameter
$\alpha$:
$$
f(\alpha)=\varphi(\bold u)=cos^2(\alpha)\,\varphi(\bold v)
+2\sin(\alpha)\cos(\alpha)\,\varphi(\bold v,\bold w)
+\sin^2(\alpha)\,\varphi(\bold w).
$$
According to the choice of the vector $v$, we have the estimate
$\varphi(\bold u)\leqslant\varphi(\bold v)$, and for $\alpha=0$,
i\.\,e\. when $u=v$, we have the equality $\varphi(\bold u)=
\varphi(\bold v)$. Hence, $\alpha=0$ is a maximum point for the
function $f(\alpha)$. Let's calculate its derivative at the point
$\alpha=0$ and equate it to zero. This yields
$$
f'(0)=2\,\varphi(\bold v,\bold w)=0.
\tag2.11
$$
Hence, $\varphi(\bold v,\bold w)=0$ for all vectors $\bold w\in
U_{\sssize\perp}$. Let's apply the inductive hypothesis to the
subspace $U_{\sssize\perp}$ whose dimension is less by $1$ than 
the dimension of the space $V$. Therefore, we can find an
orthonormal basis $\bold e_1,\,\ldots,\,\bold e_{n-1}$ in the
subspace $U_{\sssize\perp}$ such that the matrix of the form
$\varphi$ is diagonal in this basis: $\varphi(\bold e_i,\bold e_j)
=0$ for $i\neq j$. Let's complete the basis $\bold e_1,\,\ldots,\,
\bold e_{n-1}$ with the vector $\bold e_n=\bold v$. The complementary
vector $\bold e_n$ is a vector of unit length. It is orthogonal to
the vectors $\bold e_1,\,\ldots,\,\bold e_{n-1}$. Therefore, the
basis $\bold e_1,\,\ldots,\,\bold e_n$ is an orthonormal basis in
$V$. The matrix of the form $\varphi$ is diagonal in the basis
$\bold e_1,\,\ldots,\,\bold e_n$. This fact is immediate from
\thetag{2.11}. The theorem is proved.
\qed\enddemo
     The theorem~2.4 is known as the theorem on simultaneous 
diagonalization of a pair quadratic form $\varphi$ and $g$. 
For this purpose one of them should be positive. Then the positive
form $g$ defines the structure of a Euclidean space in $V$ and then
one can apply the theorem~2.4. Orthonormality of the basis 
$\bold e_1,\,\ldots,\,\bold e_n$ means that the matrix of 
$g$ is diagonal in this basis (it is the unit matrix). The matrix
of $\varphi$ is also diagonal as stated in the theorem~2.4.\par
\head
\S\,3. Selfadjoint operators. The theorem on the spectrum and the basis
of eigenvectors for a selfadjoint operator.
\endhead
\rightheadtext{\S\,3. Selfadjoint operators.}
\definition{Definition 3.1} A linear operator $f\!:\,V\to V$ in a
Euclidean vector space $V$ is called a {\it symmetric operator\/} or a
{\it selfadjoint operator\/} if for any two vectors $\bold v,\bold w\in V$
the following equality is fulfilled:
$(\bold v\,|\,f(\bold w))=(f(\bold v)\,|\,\bold w)$.
\enddefinition
\definition{Definition 3.2} A linear operator $h\!:\,V\to V$ in a
Euclidean vector space $V$ is called an {\it adjoint operator\/}
to the operator $f\!:\,V\to V$ if for any two vectors $\bold v,\bold w\in V$
the following equality is fulfilled: $(\bold v\,|\,f(\bold w))
=(h(\bold v)\,|\,\bold w)$. The adjoint operator is denoted as follows:
$h=f^{\sssize+}$.
\enddefinition
     In S\,4 of Chapter~\uppercase\expandafter{\romannumeral 3} we have
introduced the concept of conjugate mapping. There we have shown that
any linear mapping $f\!:\,V\to W$ possesses the  conjugate mapping 
$f^*:W^*\to V^*$. For a linear operator $f\!:\,V\to V$ the conjugate mapping $f^*$ is a linear operator in dual space $V^*$. It is related to $f$ by means
of the equality
$$
\hskip -2em
\langle f^*(\bold u)\,|\,\bold v\rangle=\langle\bold u\,|\,f(\bold v)
\rangle,
\tag3.1
$$
which \pagebreak is fulfilled for all $\bold u\in V^*$ and for all 
$\bold v\in V$.
\par
      The structure of a Euclidean vector space in $V$ is determined
by a positive quadratic form $g$. Like every quadratic form, the form
$g$ possesses the associated mapping $a_g\!:\,V\to V^*$ (see \S\,2
in Chapter~\uppercase\expandafter{\romannumeral 3}) such that
$$
\hskip -2em
\langle a_g(\bold v)\,|\,\bold w\rangle=g(\bold v,\bold w)
=(\bold v\,|\,\bold w).
\tag3.2
$$
In the case of finite-dimensional space $V$ and positive form $g$
the associated mapping $a_g$. is bijective. Therefore, for any 
linear operator $f\!:\,V\to V$ we can define the composition 
$h=a_g^{-1}\compos f^*\compos a_g$. Then from \thetag{3.1} and
\thetag{3.2} we derive
$$
\hskip -2em
\aligned
(h(&\bold v)\,|\,\bold w)=\langle a_g\,{\ssize\circ}\,h(\bold v)\,
|\,\bold w\rangle=\\
&=\langle f^*{\ssize\circ}\,a_g(\bold v)\,|\,\bold w\rangle=
\langle a_g(\bold v)\,|\,f(\bold w)\rangle=(\bold v\,|\,f(\bold w)).
\endaligned
\tag3.3
$$
Comparing \thetag{3.3} with the definition~3.2, we can formulate the
following theorem.
\proclaim{Theorem 3.1} For any operator $f$ in a finite-dimensional 
Euclidean space $(V,g)$ there is the unique adjoint operator 
$f^{\sssize+}=a_g^{-1}\compos f^*\compos\,a_g$.
\endproclaim
\demo{Proof} The existence of an adjoint operator is already derived
from the formula $f^{\sssize+}=a_g^{-1}\compos f^*\compos \,a_g$ and
the equality \thetag{3.3}. Let's prove its uniqueness. Assume that $h$
is another operator satisfying the definition~3.2. Then for the
difference $r=h-f^{\sssize+}$ we derive the relationship
$$
\hskip -2em
\aligned
(r(\bold v)\,|\,\bold w)=(h(&\bold v)\,|\,\bold w)-(f^{\sssize+}(\bold v)
\,|\,\bold w)=\\
&=(\bold v\,|\,f(\bold w))-(\bold v\,|\,f(\bold w))=0.
\endaligned
\tag3.4
$$
Since $\bold w$ in \thetag{3.4} is an arbitrary vector, we conclude
that $r(\bold v)\in\Ker g$. However, $\Ker g=\{\bold 0\}$ for a 
positive quadratic form $g$, hence, $h(\bold v)=\bold 0$ for any 
$\bold v\in V$. This means that $h=0$. Thus, we have proved that
the adjoint operator $f^{\sssize+}$ for $f$ is unique. This completes
the proof of the theorem.
\qed\enddemo
\proclaim{Corollary} The passage from $f$ to $f^{\sssize+}$ is an
operator in the space of endomor\-phisms $\End(V)$ of a finite-dimensional
Euclidean vector space $(V,g)$. This operator possesses the following
properties:
$$
\xalignat 2
&(f+h)^{\sssize+}=f^{\sssize+}+h^{\sssize+},
&&(\alpha\cdot f)^{\sssize+}=\alpha\cdot f^{\sssize+},\\
&(f\compos h)^{\sssize+}=h^{\sssize+}\compos f^{\sssize+},
&&(f^{\sssize+})^{\sssize+}=f.
\endxalignat
$$
\endproclaim
Relying upon the existence and the uniqueness of of the adjoint 
operator $f^{\sssize+}$ for any operator $f\in\End(V)$, we can 
derive all the above relationships immediately from the definition~3.2.
The relationship $f^{\sssize+}=a_g^{-1}\compos f^*\compos\,a_g$ can
be expressed in the form of the following commutative diagram:
$$
\CD
V@>f^{\sssize+}>>V\\
@Va_gVV         @VVa_gV\\
V^*@>>f^*> V^*
\endCD
$$
Comparing the definitions~3.1 and 3.2, now we see that a selfadjoint
operator $f$ is an operator which is adjoint to itself: $f^{\sssize+}
=f$.\par
     Let $\bold e_1,\,\ldots,\,\bold e_n$ be a basis in a 
finite-dimensional Euclidean space $(V,g)$ and let $h^1,\,\ldots,
\,h^n$ be the corresponding dual basis composed by coordinate 
functionals. For any vector $\bold v\in V$ we have the following
expansion, which follows from the definition of coordinate functionals
(see \S\,1 in Chapter~\uppercase\expandafter{\romannumeral 3}):
$$
\bold v=h^1(\bold v)\cdot\bold e_1+\ldots+h^n(\bold v)
\cdot\bold e_n.
$$
Let's apply this expansion in order to calculate the matrix of the
associated mapping $a_g$. For this purpose we need to apply $a_g$
one by one to all basis vectors $\bold e_1,\,\ldots,\,\bold e_n$ 
and expand the results in the dual basis in $V^*$. Let's consider
the value of the functional $a_g(\bold e_i)$ on an arbitrary vector
$\bold v$ of the space $V$:
$$
\align
a_g(&\bold e_i)(\bold v)=\langle a_g(\bold e_i)\,|\,\bold v\rangle
=g(\bold e_i,\bold v)=\\
&=g(\bold e_i,h^1(\bold v)\cdot\bold e_1+\ldots+h^n(\bold v)\cdot
\bold e_n)=\sum^n_{j=1} g_{ij}\,h^j(\bold v).
\endalign
$$
Since $\bold v\in V$ is an arbitrary vector, we conclude that the
matrix of the associated mapping $a_g$ in two bases $\bold e_1,\,\ldots,
\,\bold e_n$ and $h^1,\,\ldots,\,h^n$ coincides with the matrix $g_{ij}
=g(\bold e_i,\bold e_j)$ of the quadratic form $g$ in the basis 
$\bold e_1,\,\ldots,\,\bold e_n$ 
    The matrix $g_{ij}$ is non-degenerate (see theorem~1.2 or 
Silvester's criterion in \S\,4 of 
Chapter~\uppercase\expandafter{\romannumeral 4}). Let's denote by 
$g^{ij}$ the components of the matrix inverse to $g_{ij}$. The matrix
$g^{ij}$ is the matrix of the inverse mapping $a_g^{-1}$, i\.\,e\. we
have:
$$
\xalignat 2
&\hskip -2em
a_g(\bold e_i)=\sum^n_{j=1} g_{ij}\,h^j,
&&a_g^{-1}(h^j)=\sum^n_{j=1} g^{ij}\,\bold e_i.
\tag3.5
\endxalignat
$$
The matrix inverse to a symmetric matrix is again a symmetric matrix
(this fact is well-known from general algebra). Therefore $g^{ij}
=g^{ji}$.\par
     Remember that we have already calculated the matrix of the
conjugate mapping $f^*$ (see formula \thetag{4.2} and theorem~4.3 
in Chapter~\uppercase\expandafter{\romannumeral 3}). When applied to
our present case the results of 
Chapter~\uppercase\expandafter{\romannumeral 3} mean that the matrix
of the operator $f^*\!:\,V^*\to V^*$ in the basis of coordinate
functionals $h^1,\ldots,h^n$ coincides with the matrix of the initial
operator $f$ in the basis $\bold e_1,\,\ldots,\,\bold e_n$. Let's
combine this fact with \thetag{3.5} and let's use the formula 
$f^{\sssize+}=a_g^{-1}{\ssize\circ}f^*{\ssize\circ}\,a_g$ from the
theorem~3.1. Then for the matrix of $F^{\sssize+}$ of the adjoint 
operator $f^{\sssize+}$ we obtain:
$$
\hskip -2em
(F^{\sssize+})^i_j=\sum^n_{k=1}\sum^n_{q=1} g^{iq}\,F^k_q\,g_{kj}.
\tag3.6
$$
In matrix form the formula \thetag{3.6} is written as 
$F^{\sssize+}=G^{-1}\,F^{\,\tr}\,G$, where $G$ is the Gram matrix
of that basis in which the matrices of $f$ and $f^{\sssize+}$ are
calculated. The formula \thetag{3.6} simplifies substantially for
orthonormal bases. Here the passage to the adjoint operator means
the transposition of its matrix. The matrix of a selfadjoint 
operator in an orthonormal basis is symmetric. For this reason 
selfadjoint operators are often called symmetric operators.
\par
     Let $f\!:\,V\to V$ be a selfadjoint operator in a Euclidean
space $V$. Each such operator produces the quadratic 
$\varphi\kern-2pt\lower 4pt\hbox{$\ssize f\,$}$ according to the
formula
$$
\hskip -2em
\varphi\kern-2pt\lower 4pt\hbox{$\ssize f\,$}(\bold v)=(\bold v
\,|\,f(\bold v)).
\tag3.7
$$
Conversely, assume that we have a quadratic form $\varphi$
in a finite-dimensional Euclidean space $(V,g)$. The form $\varphi$
determines the associated mapping $a_\varphi$ (see
definition~2.5 in Chapter~\uppercase\expandafter{\romannumeral 4}).
This mapping satisfies the relationship
$$
\hskip -2em
\langle a_\varphi(\bold v)\,|\,\bold w\rangle=\varphi(\bold v,
\bold w)
\tag3.8
$$
for any two vectors $\bold v,\bold w\in V$. The positive quadratic 
form $g$ defining the structure of Euclidean space in $V$ has also
its own associated mapping $a_g$. The mapping $a_g$ is bijective
since $g$ is non-degenerate (see theorem~4.2 in
Chapter~\uppercase\expandafter{\romannumeral 4}). Therefore, we can 
consider the composition of $a_g^{-1}$ and $a_\varphi$:
$$
\hskip -2em
f_\varphi=a_g^{-1}\compos\,a_\varphi.
\tag3.9
$$
This composition \thetag{3.9} is an operator in $V$. It is called
the {\it associated operator\/} of the form $\varphi$ in a Euclidean
space. Since $a_g$ is bijective, we can write \thetag{3.2} as
$\langle\bold u\,|\,\bold w\rangle=(a_g^{-1}(\bold u)\,|\,\bold w)$.
Combining this equality with \thetag{3.8}, we find
$$
\hskip -2em
(f_\varphi(\bold v)\,|\,\bold w)=(a_g^{-1}(a_\varphi(\bold v))
\,|\,\bold w)=\langle a_\varphi(\bold v)\,|\,\bold w\rangle
=\varphi(\bold v,\bold w).
\tag3.10
$$
Now, using the symmetry of the form $\varphi(\bold v,\bold w)$ 
in \thetag{3.10}, we write
$$
\hskip -2em
(f_\varphi(\bold v)\,|\,\bold w)=\varphi(\bold v,\bold w)
=\varphi(\bold w,\bold v)=(f_\varphi(\bold w)\,|\,\bold v)
=(\bold v\,|\,f_\varphi(\bold w)).
\tag3.11
$$
The relationship \thetag{3.11}, which is an identity for all
$\bold v,\bold w\in V$, means that $f_\varphi$ is a selfadjoint
operator (see definition~2.1).\par
     The formula \thetag{3.7} associates each selfadjoint operator
$f$ with the quadratic form $\varphi\kern-2pt\lower 4pt
\hbox{$\ssize f\,$}$, while the formula \thetag{3.9} associates 
each quadratic form $\varphi$ with the selfadjoint operator 
$f_\varphi$. These two associations are one to one and are inverse
to each other. Indeed, let's apply the formula \thetag{3.7} to the
operator \thetag{3.9} and use \thetag{3.10}:
$$
\varphi\kern-2pt\lower 4pt\hbox{$\ssize f\,$}(\bold v)=(\bold v\,|\,
f_\varphi(\bold v))=\varphi(\bold v,\bold v)=\varphi(\bold v).
$$
Now, conversely, let's construct the operator $h=f_\varphi$ for the
quadratic form $\varphi=\varphi\kern-2pt\lower 4pt\hbox{$\ssize f\,$}$.
For the operator $h$ and for two arbitrary vectors $\bold v,\bold w\in 
V$ from \thetag{3.10} we derive 
$$
(h(\bold v)\,|\,\bold w)=\varphi\kern-2pt\lower 4pt\hbox{$\ssize f\,$}
(\bold v,\bold w)=(\bold v\,|\,f(\bold w)=(f(\bold v)\,|\,\bold w).
$$
Since $\bold w\in V$ is an arbitrary vector and since the form $g$
determining the scalar product in $V$ is non-degenerate, from the
above equality we get $h(\bold v)=f(\bold v)$.\par
     Thus, from what was said above we conclude that defining a
selfadjoint operator in a finite-dimensional Euclidean space 
is equivalent to defining a quadratic form in this space.
Therefore, we can apply the theorem~2.4 for describing selfadjoint
operators in a finite-dimensional case.
\proclaim{Theorem 3.2} All eigenvalues of a selfadjoint operator
$f$ in a finite-dimensional Euclidean space $V$ are real numbers 
and there is an orthonormal basis composed by eigenvectors of such 
operator.
\endproclaim
\demo{Proof} For the selfadjoint operator $f$ in $V$ we consider
the symmetric bilinear form $\varphi\kern -4pt\lower 4pt\hbox{
$\ssize f\,$}(\bold v,\bold w)$ determined by the quadratic form
\thetag{3.7}. Let $\bold e_1,\,\ldots,\,\bold e_n$ be an orthonormal
basis in which the matrix of the form $\varphi\kern-4pt\lower 4pt\hbox{
$\ssize f\,$}$ is diagonal. Then from the formula \thetag{3.7} we
derive the following equalities:
$$
\varphi\kern-2pt\lower 4pt\hbox{$\ssize f\,$}(\bold e_i,\bold e_j)=
(\bold e_i\,|\,f(\bold e_j)=\sum^n_{k=1}F^k_j\,g_{ik}=F^i_k.
\tag3.12
$$
As we see in \thetag{3.12}, the matrices of the operator $f$ and of
the form $\varphi\kern-2pt\lower 4pt\hbox{$\ssize f\,$}$ in such basis
do coincide. This proves the proposition of the theorem.
\qed\enddemo
     The theorem~3.2 is known as the theorem on the spectrum
and the basis of eigenvectors of a selfadjoint operator. The 
main result of this theorem is the diaginalizability of 
selfadjoint operators in a finite-dimensional Euclidean space.
The characteristic polynomial of a selfadjoint operator is
factorized into the product of linear terms in $\Bbb R$. Its
eigenspaces coincide with the corresponding root subspaces,
the sum of all its eigenspaces coincides with the space $V$:
$$
\hskip -2em
V=V_{\lambda_1}\oplus\ldots\oplus V_{\lambda_{\ssize s}}.
\tag3.13
$$
\proclaim{Theorem 3.3} Any two eigenvectors of a selfadjoint 
operator corresponding to different eigenvalues are orthogonal
to each other.
\endproclaim
\demo{Proof} Let $f$ be a selfadjoint operator in a Euclidean
space and let $\lambda\neq\mu$ be its eigenvalues. Let's consider
the corresponding eigenvectors $\bold a$ and $\bold b$:
$$
\xalignat 2
&f(\bold a)=\lambda\cdot\bold a,
&&f(\bold b)=\mu\cdot\bold b.
\endxalignat
$$
Then for these two eigenvectors $\bold a$ and $\bold b$ we derive:
$$
\lambda\,(\bold a\,|\,\bold b)=(f(\bold a)\,|\,\bold b)
=(\bold a\,|f(\bold b))=\mu\,(\bold a\,|\,\bold b).
$$
Hence, $(\lambda-\mu)\,(\bold a\,|\,\bold b)=0$. But we know that
$\lambda-\mu\neq 0$. Therefore, $(\bold a\,|\,\bold b)=0$. 
The theorem is proved.
\qed\enddemo
     Assume that the kernel of selfadjoint operator $f$ is nontrivial:
$\Ker f\neq\{\bold 0\}$. Then $\lambda_1=0$ in \thetag{3.13} is one of 
the eigenvalues of the operator $f$ and we have 
$$
\xalignat 2
&\Ker f=V_{\lambda_1},
&&\Img f=V_{\lambda_2}\oplus\ldots\oplus V_{\lambda_{\ssize s}}.
\endxalignat
$$
This means that the kernel and the image of a selfadjoint operator
are orthogonal to each other and their sum coincides with $V$:
$$
\hskip -2em
V=\Ker f\oplus\Img f.
\tag3.14
$$
\head
\S~4. Isometries and orthogonal operators.
\endhead
\definition{Definition 4.1} A linear mapping $f\!:\,V\to W$ from
one Euclidean vector space $(V,g)$ to another Euclidean vector
space $(W,h)$ is called an {\it isometry} if
$$
\hskip -2em
(f(\bold x)\,|\,f(\bold y))=(\bold x\,|\,\bold y)
\tag4.1
$$
for all $\bold x,\bold y\in V$, i\.\,e\. if it preserves the scalar
product of vectors.
\enddefinition
     From \thetag{4.1} we easily derive $|f(\bold x)|=|\bold x|$, 
therefore, $f(\bold x)=\bold 0$ implies $|\bold x|=0$ and $\bold x
=\bold 0$. This means that the kernel of an isometry is always 
trivial $\Ker f=\{\bold 0\}$, i\.\,e\. any isometry is an injective
mapping. Due to the recovery formula for quadratic forms (see formula
\thetag{1.6} in Chapter~\uppercase\expandafter{\romannumeral 4}) 
in order to verify that $f\!:\,V\to W$ is an isometry it is sufficient
to verify that it preserves the norm of vectors, i\.\,e\.
$|f(\bold x)|=|\bold x|$ for all vectors $\bold x\in V$.
\proclaim{Theorem 4.1} The composition of isometries is again an isometry.
\endproclaim
\demo{Proof} Assume that the mappings $h\!:\,U\to V$ and $f\!:\,V\to W$
both are isometries. Hence, $|h(\bold u)|=|\bold u|$ for all $\bold u\in U$
and $|f(\bold v)|=|\bold v|$ for all $\bold v\in V$. Then
$$
|f\compos h(\bold u)|=|f(h(\bold u))|=|h(\bold u)|=|\bold u|
$$
for all $\bold u\in U$. This equality means that the mapping $f\compos h$ 
is an isometry. The theorem is proved.
\qed\enddemo
\definition{Definition 4.2} A bijective isometry $f\!:\,V\to W$ is called
an {\it isomorphism of Euclidean vector spaces}.
\enddefinition
\proclaim{Theorem 4.2} Isomorphisms of Euclidean vector spaces possess
the following three properties:
\roster
\item the identical mapping $\id_V$ is an isomorphism;
\item the composition of isomorphisms is an isomorphism;
\item the mapping inverse to an isomorphism is an isomorphism.
\endroster
\endproclaim
The proof of this theorem is very easy if we use the above theorem~4.1 
and the theorem~8.1 of Chapter~\uppercase\expandafter{\romannumeral 1}.
\definition{Definition 4.3} Two Euclidean vector spaces $V$ and $W$
are called {\it isomorphic} if there is an isomorphism $f:V\to W$ 
relating them.
\enddefinition
     Let's consider the arithmetic vector space $\Bbb R^n$ composed by
column vectors of the height $n$. The addition of such vectors and the
multiplication of them by real numbers are performed as the operations
with their components (see formulas \thetag{2.1} in
Chapter~\uppercase\expandafter{\romannumeral 1}). Let's define a quadratic
form $g(x)$ in $\Bbb R^n$ by setting
$$
\hskip -2em
g(\bold x)=(x^1)^2+\ldots+(x^n)^2=\sum^n_{i=1}(x^i)^2.
\tag4.2
$$
The form \thetag{4.2} yields the {\it standard scalar product\/}
and, hence, defines the {\it standard structure of a Euclidean space\/}
in $\Bbb R^n$.
\proclaim{Theorem 4.3} Any $n$-dimensional Euclidean vector space
$V$ is isomorphic to the space $\Bbb R^n$ with the standard scalar 
product \thetag{4.2}.
\endproclaim
     In order to prove this theorem it is sufficient to choose 
the orthonormal basis in $V$ and consider the mapping $\psi$ that
associates a vector $\bold v\in V$ with column vector of its 
coordinates (see formula \thetag{5.4} in
Chapter~\uppercase\expandafter{\romannumeral 1}).
\definition{Definition 4.4} An operator $f$ in a Euclidean vector
space $V$ is called an {\it orthogonal operator\/} if it is bijective
and defines an isometry $f\!:\,V\to V$.
\enddefinition
     Due to the theorem~4.2 the orthogonal operators form a group
which is called the {\it orthogonal group} of a Euclidean space     
$V$ and is denoted by $\MatGrO(V)$. The group $\MatGrO(V)$ is obviously
a subgroup in the group of automorphisms $\Aut(V)$. In the case $V=\Bbb R^n$
the orthogonal group determined by the standard scalar product in 
$\Bbb R^n$ is denoted by $\MatGrO(n,\Bbb R)$.\par
     Let $\bold e_1,\,\ldots,\,\bold e_n$ be an orthonormal basis
in a Euclidean space $V$ and let $f$ be an orthogonal operator. Then
from \thetag{4.1} we derive
$$
(f(\bold e_i)\,|\,f(\bold e_j))=(\bold e_i\,|\,\bold e_j).
$$
For the matrix of the operator $f$ in the basis $\bold e_1,\,\ldots,
\,\bold e_n$ this relationship yields:
$$
\hskip -2em
\sum^n_{k=1}F^k_i\,F^k_j=
\cases
1&\text{\ \ for \ }i=j,\\
0&\text{\ \ for \ }i\neq j,
\endcases
\tag4.3
$$
When written in the matrix form, the formula \thetag{4.3} means that
$$
\xalignat 2
&\hskip -2em
F^{\,\tr}\,F=1,
&&F^{\,-1}=F^{\,\tr}.
\tag4.4
\endxalignat
$$
The relationships \thetag{4.4} are identical to the relationships
\thetag{1.13}. Matrices that satisfy such relationships, as we already
know, are called {\it orthogonal matrices}. As a corollary of this fact
we can formulate the following theorem.
\proclaim{Theorem 4.4} An orthogonal operator $f$ in an orthonormal 
basis $\bold e_1,\,\ldots,\,\bold e_n$ of a Euclidean vector space 
$V$ is given by an orthogonal matrix.
\endproclaim
     As we have noted in \S\,1, the determinant of an orthogonal matrix
can be equal to $1$ or to $-1$. The orthogonal operators in $V$ with 
determinant $1$ form a group which is called the {\it special orthogonal group\/} of a Euclidean vector space $V$. This group is denoted by
$\MatGrSO(V)$. 
If $V=\Bbb R^n$, this group is denoted by $\MatGrSO(n,\Bbb R)$.\par
     The operators $f\in\MatGrSO(V)$ in two-dimensional case $\dim V=2$ are
most simple ones. If $\bold e_1,\,\bold e_2$ is an orthonormal basis in $V$,
then from \thetag{4.3} and $\det F=1$ we easily find the form of an orthogonal matrix $F$:
$$
\hskip -2em
F=
\Vmatrix
\cos(\varphi) & -\sin(\varphi)\\
\vspace{1.7ex}
\sin(\varphi) & \cos(\varphi)
\endVmatrix
\tag4.5
$$
A matrix $F$ of the form \thetag{4.5} is called a {\it matrix of
two-dimensional rotation}, while the numeric parameter $\varphi$ 
\pagebreak is interpreted as the angle of rotation.\par
     Let's consider orthogonal operators $f\in\MatGrSO(V)$ in the case 
$\dim V=3$. Let $\bold e_1,\,\bold e_2,\,\bold e_3$ be an orthonormal
basis in $V$. A matrix of the form
$$
\hskip -2em
F=
\Vmatrix
\cos(\varphi) & -\sin(\varphi) & 0\\
\vspace{1ex}
\sin(\varphi) & \cos(\varphi) & 0\\
\vspace{1ex}
0             & 0             & 1\\
\endVmatrix
\tag4.6
$$
is an orthogonal matrix with determinant $1$. The operator $f$
associated with the matrix \thetag{4.6} is called the operator
of {\it rotation} about the vector $\bold e_3$ by the angle
$\varphi$.
\proclaim{Theorem 4.5} In a three-dimensional Euclidean vector 
space $V$ any orthogonal operator $f$ with determinant $1$ has 
an eigenvalue $\lambda=1$.
\endproclaim
\demo{Proof} Let's consider the characteristic polynomial of
the operator $f$. This is the polynomial of degree $3$ in 
$\lambda$ with real coefficients: 
$$
P(\lambda)=-\lambda^3+F_1\,\lambda^2-F_2\,\lambda+F_3
\text{, \ where \ }F_3=\det f=1.
$$ 
Remember that the values of a polynomial of odd degree for large 
positive $\lambda$ and for large negative $\lambda$ differ in sign:
$$
\xalignat 2
&\lim_{\lambda\to-\infty}P(\lambda)=+\infty,
&&\lim_{\lambda\to+\infty}P(\lambda)=-\infty.
\endxalignat
$$
Therefore the equation of the odd degree $P(\lambda)=0$ with real
coefficients has at least one real root $\lambda=\lambda_1$. This
root is an eigenvalue of the operator $f$.\par
     Let $\bold e_1\neq\bold 0$ be an eigenvector of $f$ corresponding 
to the eigenvalue $\lambda_1$. Then, applying the isometry condition 
$|\bold v|=|f(\bold v)|$ to the vector $\bold v=\bold e_1$, we get
$$
|\bold e_1|=|f(\bold e_1)|=|\lambda_1\cdot\bold e_1|=|\lambda_1|\,
|\bold e_1|.
$$
Hence, we find that $|\lambda_1|=1$. This means that 
$\lambda_1=1$ or $\lambda_1=-1$. In the case $\lambda_1=1$ the proposition
of the theorem is valid. Therefore, we consider the case $\lambda_1=-1$.
Let's separate the linear factor $(\lambda+1)$ in characteristic 
polynomial:
$$
P(\lambda)=-\lambda^3+F_1\,\lambda^2-F_2\,\lambda+1
=-(\lambda+1)(\lambda^2-\Phi_1\,\lambda-1).
$$
Then $F_1=\Phi_1-1$ and $F_2=-1-\Phi_1$. In order to the remaining
roots of the polynomial $P(\lambda)$ we consider the following quadratic
equation:
$$
\lambda^2-\Phi_1\,\lambda-1=0.
$$
This equation always has two real roots $\lambda_2$ and $\lambda_3$
since its discriminant is positive: $D=(\Phi_1)^2+4>0$. Due to the 
Viet theorem we have $\lambda_2\,\lambda_3=-1$. Due to the same 
reasons as above in the case of $\lambda_1$, for $\lambda_2$ and
$\lambda_3$ we get $|\lambda_2|=|\lambda_3|=1$. Hence, one of these
two real numbers is equal to $1$ and the other is equal to $-1$.
Thus, we have proved that the number $\lambda=1$ is among the 
eigenvalues of the operator $f$. The theorem is proved.
\qed\enddemo
\proclaim{Theorem 4.5} In a three-dimensional Euclidean vector 
space $V$ for any ortho\-gonal operator $f$ with determinant $1$ there
is an orthonormal basis in which the matrix of $f$ has the form
\thetag{4.6}.
\endproclaim
\demo{Proof} Under the assumptions of theorem~4.5 the operator 
$f$ has an eigenvalue $\lambda_1=1$. Let $\bold e_1\neq\bold 0$ be
an eigenvector of this operator associated with the eigenvalue 
$\lambda_1=1$. Let's denote by $U$ the span of the eigenvector
$\bold e_1$ and consider its orthogonal complement $U_{\sssize\perp}$.
This is the two-dimensional subspace in the three-dimensional space
$V$. This subspace is invariant under the action of $f$. Indeed,
from $\bold x\in U_{\sssize\perp}$ we derive $(\bold x\,|\,\bold e_1)
=0$. Let's write the isometry condition \thetag{4.1} for the vectors
$\bold x$ and $\bold y=\bold e_1$:
$$
0=(\bold x\,|\,\bold e_1)=(f(\bold x)\,|\,f(\bold e_1))
=\lambda_1\,(f(\bold x)\,|\,\bold e_1).
$$
Since $\lambda_1=1$, we get $(f(\bold x)\,|\,\bold e_1)=0$. Hence,
$f(\bold x)\in U_{\sssize\perp}$, which proves the invariance of the
subspace $U_{\sssize\perp}$.\par
     Let's consider the restriction of the operator $f$ to the invariant
subspace $U_{\sssize\perp}$. This restriction is an orthogonal operator
in two-dimensional space $U_{\sssize\perp}$, its determinant being equal 
to $1$. Therefore, in some orthogonal basis $\bold e_2,\,\bold e_3$ of $U_{\sssize\perp}$ the matrix of the restricted operator has the 
form \thetag{4.5}.\par
     Remember that $\bold e_1$ is perpendicular to $\bold e_2$ and 
$\bold e_3$. It can be normalized to the unit length. Then three
vectors $\bold e_1,\,\bold e_2,\,\bold e_3$ form an orthonormal basis in
three-dimensional space $V$ and the matrix of $f$ in this basis 
has the form \thetag{4.6}. The theorem is proved.
\qed\enddemo
     The result of this theorem is that any orthogonal operator $f$ with
determinant $1$ in a three-dimensional Euclidean vector space $V$ is an
operator of rotation. The eigenvector $\bold e_1$ associated with the
eigenvalue $\lambda_1=1$ determines the axis of rotation, while the real
parameter $\varphi$ in the matrix \thetag{4.6} determines the
angle of such rotation.\par
\newpage
\topmatter
\title\chapter{6}
Affine spaces.
\endtitle
\endtopmatter
\document
\head
\S\,1. Points and parallel translations. Affine spaces.
\endhead
\leftheadtext{CHAPTER~\uppercase\expandafter{\romannumeral 6}.
AFFINE SPACES.}
\rightheadtext{\S\,1. Points and parallel translations.}
     Let $M$ be an arbitrary set. A {\it transformation\/} of the set $M$ 
is a bijective mapping $p\!:\,M\to M$ of the set $M$ onto itself.
\definition{Definition 1.1} Let $V$ be a linear vector space. We say that
an {\it action\/} of $V$ on a set $M$ is defined if each vector $\bold v
\in V$ is associated with some transformation $p_{\bold v}$ of the set
$M$ and the following conditions are fulfilled:
\roster
\item  $p_{\bold 0}=\id_M$;
\item  $p_{\bold v+\bold w}=p_{\bold v}\compos p_{\bold w}$
       for all $\bold v,\bold w\in V$.
\endroster
\enddefinition
    From the properties \therosteritem{1} and \therosteritem{2} 
of an action of a space $V$ on a set $M$ one can easily derive 
the following two properties of such action:
\roster
\item[3]  $p_{-\bold v}=p_{\bold v}^{-1}$
          for all $\bold v\in V$;
\item     $p_{\bold v}\compos p_{\bold w}=p_{\bold w}\compos
          p_{\bold v}$
          for all $\bold v,\bold w\in V$.
\endroster
\definition{Definition 1.2} An action of a vector space $V$ on a set
$M$ is called a {\it transitive action\/} if for any two elements
$A,B\in M$ there is a vector $\bold v\in V$ such that $p_{\bold v}(A)=B$,
i\.\,e\. the transformation $p_{\bold v}$ takes $A$ to $B$.
\enddefinition
\definition{Definition 1.3} An action of a vector space $V$ on a set
$M$  is called a {\it free action\/} if for any element $A\in M$ the
equality $p_{\bold v}(A)=A$ implies $\bold v=\bold 0$.
\enddefinition
\definition{Definition 1.4} A set $M$ is called an {\it affine space\/}
over the field $\Bbb K$ if there is a free transitive action of
some linear vector space $V$ over the field $\Bbb K$ on $M$.
\enddefinition
     Due to this definition any affine space $M$ is associated with some linear vector space $V$. Therefore an affine space $M$ is often denoted 
as a pair $(M,V)$.\par
     Elements of an affine space are used to be called {\it points}. We
shall denote them by capital letters $A$, $B$, $C$, etc. An affine space
itself is sometimes called a {\it point space}. A transformation 
$p_{\bold v}$ given by a vector $\bold v\in V$ is called a {\it parallel translation\/} in an affine space $M$.\par
     Let $U$ be a subspace in $V$. Let's choose a point $A\in M$ and then
let's define a subset $L\subset M$ in the following way:
$$
\hskip -2em
L=\{B\in M:\exists\,\bold u\ ((\bold u\in U) \and (B=p_{\bold u}(A)))\}.
\tag1.1
$$
A subset $L$ of $M$ determined according to \thetag{1.1} is called a
{\it linear submanifold\/} of an affine space $M$. Thereby the subspace
$U\subset V$ is called the {\it directing subspace\/} of a linear
submanifold $L$. The dimension of the directing subspace in \thetag{1.1}
is taken for the {\it dimension\/} of a linear submanifold $L$.
One-dimensional linear submanifolds are called {\it straight lines};
two-dimensional submanifolds are called {\it planes}. If the dimension
of $U$ is less by one than the dimension of $V$, i\.\,e\. if $\dim(V/U)=1$,
then the corresponding linear submanifold $L$ is called a {\it hyperplane}.
Linear submanifolds of other intermediate dimensions have no special
titles.\par
     Let $U=\langle a\rangle$ be a one-dimensional subspace in $V$. Then
any vector $\bold u\in U$ is presented as $\bold u=t\cdot\bold a$, where
$t\in\Bbb K$. Upon choosing a point $A\in M$ the subspace $U$ determines
the straight line in $M$ passing through the point $A$. An arbitrary point
$A(t)$ of this straight is given by the formula:
$$
A(t)=p_{\,t\cdot\bold a}(A).
\tag1.2
$$
The formula \thetag{1.2} is known as a {\it parametric equation} of a
straight line in an affine space, the vector $\bold a$ is called a {\it directing vector}, while $t\in\Bbb K$ is a {\it parameter}.\par
     If $\Bbb K=\Bbb R$, we can consider the set of points on the straight
line \thetag{1.2} corresponding to the values of $t$ taken from the interval
$[0,\,1]\subset\Bbb R$. Such set is called a {\it segment} of a straight 
line. The points $A=A(0)$ and $B=A(1)$ are ending points of this segment.
One can choose a direction on the segment $AB$ by saying that one of the
ending points is the {\it beginning\/} of the segment and the other is
the {\it end\/} of the segment. A segment $AB$ with a fixed direction
on it is called a {\it directed segment} or an {\it arrowhead segment}.
Two arrowhead segments $\overrightarrow{AB\,}$ and $\overrightarrow{BA\,}$
are assumed to be distinct\footnote{\ If $\Bbb K\neq\Bbb R$, an arrowhead
segment $\overrightarrow{AB\,}$ is assumed to be consisting on two points
$A$ and $B$ only, it has no interior at all.}.\par
\adjustfootnotemark{-1}
     Let $A$ and $B$ be two points of an affine space $M$. Due to the
transitivity of the action of $V$ on $M$ there exists a vector $\bold v
\in V$ that defines the parallel translation $p_{\bold v}$ taking the
point $A$ to the point $B$: $p_{\bold v}(A)=B$. Let's prove that such parallel translation is unique. If $p_{\bold w}$ is another parallel
translation such that $p_{\bold w}(A)=B$, then for the parallel translation
$p_{\bold w-\bold v}$ we have 
$$
p_{\bold w-\bold v}(A)=p_{-\bold v}\,{\ssize\circ}\,p_{\bold w}(A)=
p_{\bold v}^{-1}(p_{\bold w}(A))=p_{\bold v}^{-1}(B)=A.
$$
Since $V$ acts freely on $M$ (see definition~1.3), we have $\bold w
-\bold v=\bold 0$. Hence, $\bold w=\bold v$, this proves the uniqueness
of the vector $\bold v$ determined by the condition $p_{\bold v}(A)=B$.
\par
     The above fact appears to be very useful: if we have an affine
space $(M,V)$, then vectors of $V$ can be represented by arrowhead 
segments in $M$. Each pair of points $A,B\in M$ specifies the unique
vector $\bold a\in V$ such that $p_{\bold a}(A)=B$. This vector can be 
used as a directing vector of the straight line \thetag{1.2} passing through
the points $A$ and $B$. The arrowhead segment with the beginning at
the point $A$ and with the end at the point $B$ is called the {\it geometric
representation} of the vector $a$. It is denoted $\overarrow{AB\,}$\par
     A vector $\bold a$ is uniquely determined by its geometric
representation $\overarrow{AB\,}$. However, a vector $\bold a$ can have
several geometric representations. Indeed, if we choose a point $C\neq A$,
we can determine the point $D=p_{\bold a}(C)$ and then we can construct
the geometric representation $\overarrow{CD\,}$ for the vector $\bold a$.
The points $A$ and $C$ specify a parallel translation $p_{\bold b}$ such
that $p_{\bold b}(A)=C$. Using the property \therosteritem{4} of parallel translations, it is easy to find that the parallel translation $p_{\bold b}$
maps the segment $AB$ to the segment $CD$. So we conclude: various geometric
representations of a vector $\bold a$ are related to each other by means of parallel translations. Note that $p_{\bold a+\bold b}(A)=D$. Therefore,
$\overarrow{AD\,}$ is a geometric realization of the vector 
$\bold a+\bold b$. From this fact we easily derive the well-known rules
for vector addition: the triangle rule $\overarrow{AC\,}+\overarrow{CD
\,}=\overarrow{AD\,}$ and the parallelogram rule $\overarrow{AB\,}+
\overarrow{AC\,}=\overarrow{AD\,}$.\par
     Let $O$ be some fixed point of an affine space $M$. Let's call it
the {\it origin}. Then any point $A\in M$ specifies the arrowhead
$\overarrow{OA\,}$ which is identified with the unique vector $\bold r
\in V$ by means of the equality $p_{\bold r}(O)=A$. This vector 
$\bold r=\bold r_{\kern-2pt\lower 2pt\hbox{$\ssize A$}}$ is called the
{\it radius-vector\/} of the point $A$. If the space $V$ is finite-dimensional, then we can choose a basis $\bold e_1,\,\ldots,\,
\bold e_n$ and then can expand the radius-vectors of all points $A\in M$
in this basis.
\definition{Definition 1.5} A {\it frame\/} or a {\it coordinate system\/}
in an affine space $M$ is a pair consisting of a point $O\in M$ and a basis
$\bold e_1,\,\ldots,\,\bold e_n$ in $V$. The coordinates of the 
radius-vector $r_{\kern-2pt\lower 2pt\hbox{$\ssize A$}}=\overarrow{OA\,}$ 
in the basis $\bold e_1,\,\ldots,\,\bold e_n$ are called the {\it 
coordinates\/} of a point $A$ in the coordinate system $O,\,\bold e_1,\,
\ldots,\,\bold e_n$.
\enddefinition
     Coordinate systems in affine spaces play the same role as bases in
linear vector spaces. Let $O,\,\bold e_1,\,\ldots,\,\bold e_n$ and $O',\,
\tilde\bold e_1,\,\ldots,\,\tilde\bold e_n$ be two coordinate systems
in an affine space $M$. The relation of the bases $\bold e_1,\,\ldots,\,
\bold e_n$ and  $\tilde\bold e_1,\,\ldots,\,\tilde\bold e_n$ is given by
the direct and inverse transition matrices $S$ and $T$. The points
$O$ and $O'$ determine the arrowhead segment $\overarrow{OO'\ }$ and the
opposite arrowhead segment $\overarrow{O'O\,}$. They are associated with 
two vectors $\boldsymbol\rho,\tilde{\boldsymbol\rho}\in V$:
$$
\xalignat 2
&\boldsymbol\rho=\overarrow{OO'\ }
&&\tilde{\boldsymbol\rho}=\overarrow{O'O\,}
\endxalignat
$$
Let's expand $\boldsymbol\rho$ in the basis $\bold e_1,\,\ldots,\,\bold 
e_n$ and $\tilde{\boldsymbol\rho}$ in the basis $\tilde\bold e_1,\,\ldots,
\,\tilde\bold e_n$:
$$
\hskip -2em
\aligned
\boldsymbol\rho&=\rho^1\cdot\bold e_1+\ldots+\rho^n\cdot\bold e_n,\\
\tilde{\boldsymbol\rho}&=\tilde\rho^1\cdot\tilde\bold e_1+\ldots+
\tilde\rho^n\cdot\tilde\bold e_n.
\endaligned
\tag1.3
$$
Then consider a point $X\in M$. The following formulas are obvious:
$$
\xalignat 2
&\overarrow{OX\,}=\overarrow{OO'\ }+\overarrow{O'X\,},
&&\overarrow{O'X\,}=\overarrow{O'O\,}+\overarrow{OX\,},
\endxalignat
$$
By means of them we can find the relation of the coordinates of a point
$X$ in two different coordinate systems
$O,\,\bold e_1,\,\ldots,\,\bold e_n$ and $O',\,
\tilde\bold e_1,\,\ldots,\,\tilde\bold e_n$:
$$
\xalignat 2
&\hskip -2em
x^i=\rho^i+\sum^n_{j=1}S^i_j\,\tilde x^j,
&&\tilde x^i=\tilde\rho^i+\sum^n_{j=1}T^i_j\,x^j.
\tag1.4
\endxalignat
$$
Though the vectors $\boldsymbol\rho$ and $\tilde{\boldsymbol\rho}$
differ only in sign ($\tilde{\boldsymbol\rho}=-\boldsymbol\rho$),
their coordinates in formulas \thetag{1.4} are much more different:
$$
\xalignat 2
&\hskip -2em
\rho^i=-\sum^n_{j=1}S^i_j\,\tilde\rho^j,
&&\tilde\rho^i=-\sum^n_{j=1}T^i_j\,\rho^j.
\endxalignat
$$
This happens because $\boldsymbol\rho$ and $\tilde{\boldsymbol\rho}$
are expanded in two different bases (see the above expansions 
\thetag{1.3}).\par
     The facts from the theory of affine spaces, which we stated above,
show that considering affine spaces is a proper way for geometrization
of the linear algebra. A vector is an algebraic object: we can add 
vectors, we can multiply them by numbers, and we can form linear
combinations of them. In affine space the concept of a point becomes
paramount. Points form straight lines, planes, and their multidimensional
generalizations --- linear submanifolds. In affine spaces we have a quite
natural concept of parallel translations and, hence, we can define the
concept of parallelism for linear submanifolds. The geometry of 
two-dimensional affine spaces is called the {\it planimetry}, the geometry
of three-dimensional affine spaces is called the {\it stereometry}. Affine
spaces of higher dimensions are studied by a geometrical discipline which 
is called the {\it multidimensional geometry}.
\head
\S\,2. Euclidean point spaces.\\Quadrics in a Euclidean space.
\endhead
\rightheadtext{\S\,2. Euclidean point spaces. }
\definition{Definition 2.1} An affine space $(M,V)$ over the field
of real numbers $\Bbb R$ is called a {\it Euclidean point space\/}
if the space $V$ acting on $M$ by parallel translations is equipped 
with a structure of a Euclidean vector space, i\.\,e\. if in $V$ some
positive quadratic form $g$ is fixed.
\enddefinition
     In affine spaces, which we considered in previous section, a very
important feature was lacking: there was no concept of a length and 
there was no concept of an angle. The structure of a Euclidean space
given by a quadratic form $g$ brings this lacking feature in. Let $A$ 
and $B$ be two points of a Euclidean point space $M$. They determine
a vector $\bold v\in V$ specified by the condition $p_{\bold v}(A)=B$ 
(this vector is identified with the arrowhead segment $\overarrow{AB\,}$). 
The norm of the vector $\bold v$ determined by the quadratic form $g$
is called the {\it length\/} of the segment $AB$ or the 
{\it distance\/} between two points $A$ and $B$: $|AB|=|\bold v|
=\sqrt{g(\bold v)}$. Due to the equality $|-\bold v|=|\bold v|$ we
derive $|AB|=|BA|$.\par
     Let $\overarrow{AB\,}$ and $\overarrow{CD\,}$ be two arrowhead
segments in a Euclidean point space. They are geometric representations
of two vectors $\bold v$ and $\bold w$ of $V$. The {\it angle\/}
between $\overarrow{AB\,}$ and $\overarrow{CD\,}$ by definition
is the angle between vectors $\bold v$ and $\bold w$ determined by the
formula \thetag{1.6} of 
Chapter~\uppercase\expandafter{\romannumeral 5}.
\definition{Definition 2.2} A coordinate system $O,\,\bold e_1,\,\ldots,
\,\bold e_n$ in a finite-dimensional Euclidean point space $(M,V,g)$ is
called a {\it rectangular Cartesian coordinate system\/} in $M$ if
$\bold e_1,\,\ldots,\,\bold e_n$ is an orthonormal basis of the Euclidean
vector space $(V,g)$.
\enddefinition
\definition{Definition 2.3} A {\it quadric} in a Euclidean point space
$M$ is a set of points in $M$ whose coordinates $x^1,\,\ldots,\,x^n$
in some rectangular Cartesian coordinate system $O,\,\bold e_1,\,\ldots,
\,\bold e_n$ satisfies some polynomial equation of degree two:
$$
\hskip -2em
\sum^n_{i=1}\sum^n_{j=1}a_{ij}\,x^i\,x^j+
2\,\sum^n_{i=1} b_i\,x^i+c=0.
\tag2.1
$$
\enddefinition
The definition of a quadric is not coordinate-free. It is formulated in 
terms of some rectangular Cartesian coordinate system $O,\,\bold e_1,\,
\ldots,\,\bold e_n$. However, passing to another Cartesian coordinate 
system is equivalent to a linear change of variables in the equation
\thetag{2.1} (see formulas \thetag{1.4}). Such a change of variables
changes the coefficients of the polynomial in \thetag{2.1}, but it does 
not change the structure of this equation in whole. A quadric continues 
to be a quadric in any Cartesian coordinate system.\par
     Let $O',\,\tilde\bold e_1,\,\ldots,\,\tilde\bold e_n$ be some other
rectangular Cartesian coordinate system in $M$. Let's consider the
passage from $O,\,\bold e_1,\,\ldots,\,\bold e_n$ to $O',\,\tilde\bold
e_1,\,\ldots,\,\tilde\bold e_n$. In this case transition matrices 
$S$ and $T$ in \thetag{1.4} appear to be orthogonal matrices 
(see formulas \thetag{1.13} in 
Chapter~\uppercase\expandafter{\romannumeral 5}. We can calculate 
the coefficients of the equation of quadric in the new coordinate 
system. Substituting \thetag{1.4} into \thetag{2.1}, we get
$$
\align
&\hskip -2em
\tilde a_{qp}=\sum^n_{i=1}\sum^n_{j=1}a_{ij}\,S^i_q\,S^j_p,
\tag2.2\\
&\hskip -2em
\tilde b_q=\sum^n_{i=1}b_i\,S^i_q+\sum^n_{i=1}\sum^n_{j=1}
a_{ij}\,\rho^j\,S^i_q,
\tag2.3\\
&\hskip -2em
\tilde c=\sum^n_{i=1}\sum^n_{j=1}a_{ij}\,\rho^i\,\rho^j+
\sum^n_{i=1}b_i\,\rho^i+c.
\tag2.4
\endalign
$$
Now the problem of bringing the equation of a quadric to a
canonic form is formulated as the problem of finding a
proper rectangular Cartesian coordinate system in which 
the equation \thetag{2.1} has the most simple canonic
form.\par
     The formula \thetag{2.2} coincides with the transformation 
formula for the components of a quadratic form under a change of
basis (see \thetag{1.11} in
Chapter~\uppercase\expandafter{\romannumeral 4}). Hence, we conclude
that each quadric in $M$ is associated with some quadratic form
in $V$. The form $a$ determined by the matrix $a_{ij}$ in the basis
$\bold e_1,\,\ldots,\,\bold e_n$ is called the {\it primary quadratic 
form\/} of a quadratic \thetag{2.1}.\par
     Let's consider the associated operator $f_a$ determined by the
primary quadratic form $a$ (see formula \thetag{3.9} in 
Chapter~\uppercase\expandafter{\romannumeral 5}). The operator $f_a$ 
is a selfadjoint operator in $V$; it determines the expansion of the 
space $V$ into the direct sum of two mutually orthogonal subspaces
\ $\Ker f_a$ \ and \ $\Img f_a$:
$$
V=\Ker f_a\oplus\Img f_a
\tag2.5
$$
(see \thetag{3.14} in Chapter~\uppercase\expandafter{\romannumeral 5}).
The matrix of the operator $f_a$ is given by the formula
$$
\hskip -2em
F^i_j=\sum^n_{k=1}g^{ik}\,a_{kj},
\tag2.6
$$
where $g^{ik}$ is the matrix inverse to the Gram matrix of the basis
$\bold e_1,\,\ldots,\,\bold e_n$. Apart from $f_a$, we define a vector
$\bold b$ through its coordinates given by formula
$$
\hskip -2em
b^i=\sum^n_{k=1}g^{ik}\,b_k.
\tag2.7
$$
The definition of $\bold b$ through its coordinates \thetag{2.7} is
essentially bound to the coordinate system $O,\,\bold e_1,\,\ldots,
\,\bold e_n$. This is because the formula \thetag{2.3} differs from
the standard transformation formula for the coordinates of a covector
under a change of basis (see \thetag{2.4} in 
Chapter~\uppercase\expandafter{\romannumeral 3}). Let's rewrite
\thetag{2.3} in the following form:
$$
\tilde b_q=
\sum^n_{i=1}S^i_q\left(b_i+\shave{\sum^n_{j=1}}a_{ij}\,\rho^j
\right).
\tag2.8
$$
Then let's consider the expansion of the vector $\bold b$ into the sum 
of two vectors $\bold b=\bold b^{(1)}+\bold b^{(2)}$ according to the
expansion \thetag{2.5} of the space $V$. This expansion induces the
expansion $b_i=b^{(1)}_i+b^{(2)}_i$, where $b^{(1)}_i$ are transformed
as follows:
$$
\hskip -2em
\tilde b^{(1)}_q=\sum^n_{i=1}S^i_q\,b^{(1)}_i.
\tag2.9
$$
The vector $\bold b^{(2)}$ in the expansion $b=b^{(1)}+b^{(2)}$ can be
annihilated at the expense of proper choice of the coordinate system.
Let's determine the vector $\boldsymbol\rho=\overarrow{OO'\ }$ 
from the equality $\bold b^{(2)}=-f_a(\boldsymbol\rho)$. Though it is 
not unique, the vector $\boldsymbol\rho$ satisfying this equality
does exist since $\bold b^{(2)}\in\Img f_a$. For its components we have
$$
b^{(2)}_i+\shave{\sum^n_{j=1}}a_{ij}\,\rho^j=0,
\tag2.10
$$
this follows from $\bold b^{(2)}=-f_a(\boldsymbol\rho)$ due to 
\thetag{2.6} and \thetag{2.7}. Substituting \thetag{2.10} into
\thetag{2.8}, we get the following equalities in the new coordinate 
system:
$$
\xalignat 2
&\tilde\bold b^{(2)}=\bold 0,
&&\tilde\bold b=\tilde\bold b^{(1)}.
\endxalignat 
$$
The relationships \thetag{2.9} show that the numbers $b^{(1)}_i$
cannot be annihilated (unless they are equal to zero from the very
beginning). These numbers determine the vector $\bold b^{(1)}\in
\Ker f_a$ which does not depend on the choice of a coordinate
system. As a result we have proved the following theorem.
\proclaim{Theorem 2.1} Any quadric in a Euclidean point space 
$(M,V,g)$ is associated with some selfadjoint operator $f$ and 
some vector $\bold b\in\Ker f$ such that in some rectangular 
Cartesian coordinate system the radius vector $\bold r$ of an
arbitrary point of this quadric satisfies the following equation:
$$
\hskip -2em
(f(\bold r)\,|\,\bold r)+2\,(\bold b\,|\,\bold r)+c=0.
\tag2.11
$$
\endproclaim
     The operator $f$ determines the leading part of the equation 
\thetag{2.11}. By means of this operator we subdivide all quadrics 
into two basic types:
\roster
\item {\it non-degenerate quadrics}, when $\Ker f=\{\bold 0\}$;
\item {\it degenerate quadrics}, when $\Ker f\neq\{\bold 0\}$.
\endroster
For non-degenerate quadrics the vector $\bold b$ in \thetag{2.11} 
is equal to zero. Therefore, non-degenerate quadrics are subdivided
into three types:
\roster
\item {\it elliptic type}, when $c\neq 0$ and the quadratic form
      $a(\bold x)=(f(\bold x)\,|\,\bold x)$ is positive or negative,
      i\.\,e\. can be made positive by changing the sign of $f$;
\item {\it hyperbolic type}, when $c\neq 0$ and the quadratic form
      $a(\bold x)=(f(\bold x)\,|\,\bold x)$ is not sign-definite, 
      i\.\,e\. its signature has both pluses and minuses;
\item {\it conic type}, when $c=0$.
\endroster
Degenerate quadrics are subdivided into two types:
\roster
\item {\it parabolic type}, when $\dim\Ker f=1$ and $b\neq 0$;
\item {\it cylindric type}, when $\dim\Ker f>1$ or $b=0$.
\endroster\par
      The equation \thetag{2.1} in the case of a non-degenerate quadric
of {\bf elliptic} type can be brought to the following canonic form:
$$
\frac{(x^1)^2}{(a_1)^2}+\ldots+
\frac{(x^n)^2}{(a_n)^2}=\pm 1.
$$
This is the canonic equation of a non-degenerate quadric of 
{\bf hyperbolic} type:
$$
\frac{(x^1)^2}{(a_1)^2}\pm\ldots\pm
\frac{(x^n)^2}{(a_n)^2}=\pm 1.
$$
The canonic equation of a non-degenerate quadric of {\bf conic}
type is homogeneous:
$$
\frac{(x^1)^2}{(a_1)^2}\pm\ldots\pm
\frac{(x^n)^2}{(a_n)^2}=0.
$$\par
     The equation \thetag{2.1} in the case of a degenerate quadric
of {\bf parabolic} type can be brought to the following canonic form:
$$
\frac{(x^1)^2}{(a_1)^2}\pm\ldots\pm
\frac{(x^{n-1})^2}{(a_{n-1})^2}=2\,x^n.
$$\par
     If $n=\dim M>1$, then in a canonic equation of a quadric 
of {\bf cylindric} type there is no explicit entry of at least 
one variable. Therefore, we can reduce the dimension of the space
$M$. The reduced quadric can belong to any one of the above four
types. If it is again of cylindric type, then we can repeat the 
reduction procedure. This process can terminate in some intermediate 
dimension yielding the reduced quadric of some non-cylindric type.
Otherwise we shall reach the dimension $\dim M=1$. In one-dimensional Euclidean point space there is no quadrics of cylindric type. 
Therefore, the quadrics of {\bf cylindric} type are those which
belong to one of the non-cylindric types in the reduced dimension.
\par
\newpage
\topmatter
\title
References.
\endtitle
\endtopmatter
\document
\setfirstpage
\par\noindent
\myref{1.}{Kurosh~A.~G. {\it Course of general algebra}, 
{\tencyr\char '074}Nauka{\tencyr\char '076} publishers,
Moscow.}
\myref{2.}{Sharipov~R.~A. {\it Course of differential geometry\footnotemark}, 
Publication of Bashkir State University, Ufa, 1996.}
\myref{3.}{Sharipov~R.~A. {\it Classical electrodynamics and the theory
of relativity}, Publication of Bashkir State University, Ufa, 1997;
see online
\blue{physics/0311011} in Electronic Archive 
\blue{http:/\negskp/arXiv.org}.
}
\myref{4.}{Kostrikin~A.~I. {\it Introduction to algebra}, 
{\tencyr\char '074}Nauka{\tencyr\char '076} publishers,
Moscow, 1977.}
\myref{5.}{Beklemishev~D.~V. {\it Course of analytical geometry and linear algebra}, 
{\tencyr\char '074}Nauka{\tencyr\char '076} publishers, Moscow, 1985.}
\myref{6.}{Kudryavtsev~L.~D. {\it Course of mathematical analysis, 
Vol.~\uppercase\expandafter{\romannumeral 1} and \uppercase\expandafter{\romannumeral 2}}, 
{\tencyr\char '074}Visshaya Shkola{\tencyr\char '076} 
publishers, Moscow, 1985.}
\myref{7.}{Sharipov~R.~A. {\it Quick introduction to tensor analysis},
free online publication 
\blue{math.HO/0403252} in Electronic Archive 
\blue{http:/\negskp/arXiv.org}, 2004.}
\footnotetext{\ The references \cite{2} and \cite{3} are added in 1998,
the reference \cite{7} is added in 2004.}
\enddocument
\end